\documentstyle[12pt]{book}

\makeindex

\font\teneufm=eufm10
\font\seveneufm=eufm7
\font\fiveeufm=eufm5
\newfam\eufmfam
\textfont\eufmfam=\teneufm
\scriptfont\eufmfam=\seveneufm
\scriptscriptfont\eufmfam=\fiveeufm
\def\frak#1{{\fam\eufmfam\relax#1}}

\newfam\msbfam
\font\tenmsb=msbm10 scaled \magstep1   \textfont\msbfam=\tenmsb
\font\sevenmsb=msbm7 scaled \magstep1  \scriptfont\msbfam=\sevenmsb
\font\fivemsb=msbm5  scaled \magstep1  \scriptscriptfont\msbfam=\fivemsb
\def\Bbb{\fam\msbfam \tenmsb}

\def\RR{{\Bbb R}}
\def\CC{{\Bbb C}}
\def\QQ{{\Bbb Q}}
\def\NN{{\Bbb N}}
\def\ZZ{{\Bbb Z}}
\def\II{{\Bbb I}}

\def\hexdigit#1{\ifnum#1<10 \number#1\else
 \ifnum#1=10 A\else\ifnum#1=11 B\else\ifnum#1=12 C\else
 \ifnum#1=13 D\else\ifnum#1=14 E\else\ifnum#1=15 F\fi\fi\fi\fi\fi\fi\fi}

\chardef\\="5C                    %% Typesets \ in \tt mode
\chardef\{="7B  \chardef\}="7D    %% also left and right braces

 \def\HollowBoxx #1#2#3{{\dimen0=#1 \advance\dimen0 by -#2       
       \dimen1=#1 \advance\dimen1 by #3                       
        \vrule height 0pt depth #3 width #2                   
       \hskip -#3
       \vrule height #1 depth #3 width #3}}                   
 \def\LeftContraction{\mathord{\kern1.45pt \HollowBoxx{6pt}{3.5pt}{.4pt}}\,}

 \def\HollowBox #1#2#3{{\dimen0=#1 \advance\dimen0 by -#3       
       \dimen1=#1 \advance\dimen1 by #3                       
        \vrule height #1 depth #3 width #3                    
        \vrule height 0pt depth #3 width #2                   
        \hskip -#3}}                                             
 \def\RightContraction{\mathord{\, \HollowBox{6pt}{3.1pt}{.4pt}} \kern1.6pt}

\begin{document}

%% \tableofcontents

\font\bfit=cmbxti10 scaled \magstep1
%% \font\bigbfit=cmbxti10 scaled \magstep3

\begin{center}
{\LARGE \bf A Primer of Mathematical} 
\medskip \\
{\LARGE \bf Writing}
\bigskip \\
\small \sc
Being a Disquisition on \\ 
Having Your Ideas \\
Recorded \\
Typeset \\
Published \\
Read \\
\& \\
Appreciated 
\end{center}
\vspace*{.13in} 

\begin{center}
Second Edition
\end{center}
\vspace*{.13in}

\begin{center}
\large Steven G. Krantz
\bigskip \\
\today
\end{center}

\newpage

\font\msa=msam10

\def\reg{{${}^{\mbox{\msa\char"72}}$}}
\def\bad{{\msa\char"7A}}

%% Dedication

\null \vspace*{2.5in}

\thispagestyle{empty}

\setcounter{page}{5}

\noindent This book is dedicated to Paul Halmos.  \hfill \break
For the example he has set for us all.
\vfill
\eject

\hbox{ \ \ \  }

\thispagestyle{empty}

\newpage

\pagenumbering{roman}

\setcounter{page}{7}

\markboth{TABLE OF CONTENTS}{TABLE OF CONTENTS}

\chapter*{Table of Contents}

\contentsline {chapter}{\numberline {}Preface to the Second Edition}{xi}
\contentsline {chapter}{\numberline {}Preface to the First Edition}{xiii}
\vspace*{.12in}
\contentsline {chapter}{\numberline {1}The Basics}{1}
\contentsline {section}{\numberline {1.1}What It Is All About}{2}
\contentsline {section}{\numberline {1.2}Who Is My Audience?}{3}
\contentsline {section}{\numberline {1.3}Writing and Thought}{4}
\contentsline {section}{\numberline {1.4}Say What You Mean; Mean What You Say}{5}
\contentsline {section}{\numberline {1.5}Proofreading, Reading for Sound, Reading for Sense}{10}
\contentsline {section}{\numberline {1.6}Compound Sentences, Passive Voice}{14}
\contentsline {section}{\numberline {1.7}Technical Aspects of Writing a Paper}{19}
\contentsline {section}{\numberline {1.8}More Specifics of Mathematical Writing}{23}
\contentsline {section}{\numberline {1.9}Pretension and Lack of Pretension}{31}
\contentsline {section}{\numberline {1.10}We vs.\ I vs.\ One}{33}
\contentsline {section}{\numberline {1.11}Essential Rules of Grammar, Syntax, and Usage}{34}
\contentsline {section}{\numberline {1.12}More Rules of Grammar, Syntax, and Usage}{49}
\contentsline {chapter}{\numberline {2}Topics Specific to the Writing of Mathematics}{63}
\contentsline {section}{\numberline {2.1}How to Organize a Paper}{64}
\contentsline {section}{\numberline {2.2}How to State a Theorem}{69}
\contentsline {section}{\numberline {2.3}How to Prove a Theorem}{71}
\contentsline {section}{\numberline {2.4}How to State a Definition}{74}
\contentsline {section}{\numberline {2.5}How to Write an Abstract}{76}
\contentsline {section}{\numberline {2.6}How to Write a Bibliography}{77}
\contentsline {section}{\numberline {2.7}What to Do with the Paper Once It Is Written}{83}
\contentsline {section}{\numberline {2.8}A Coda on Collaborative Work}{92}
\contentsline {chapter}{\numberline {3}Exposition}{93}
\contentsline {section}{\numberline {3.1}What Is Exposition?}{94}
\contentsline {section}{\numberline {3.2}How to Write an Expository Article}{95}
\contentsline {section}{\numberline {3.3}How to Write an Opinion Piece}{98}
\contentsline {section}{\numberline {3.4}The Spirit of the Preface}{100}
\contentsline {section}{\numberline {3.5}How Important Is Exposition?}{102}
\contentsline {chapter}{\numberline {4}Other Types of Writing}{105}
\contentsline {section}{\numberline {4.1}The Letter of Recommendation}{106}
\contentsline {section}{\numberline {4.2}The Book Review}{123}
\contentsline {section}{\numberline {4.3}The Referee's Report}{126}
\contentsline {section}{\numberline {4.4}The Talk}{128}
\contentsline {section}{\numberline {4.5}Your Vita, Your Grant, Your Job, Your Life}{138}
\contentsline {subsection}{\numberline {4.5.1}The Curriculum Vitae}{138}
\contentsline {subsection}{\numberline {4.5.2}The Grant}{141}
\contentsline {subsection}{\numberline {4.5.3}Your Job}{144}
\contentsline {subsection}{\numberline {4.5.4}Your Life}{149}
\contentsline {section}{\numberline {4.6}Electronic Mail}{150}
\contentsline {chapter}{\numberline {5}Books}{157}
\contentsline {section}{\numberline {5.1}What Constitutes a Good Book?}{157}
\contentsline {section}{\numberline {5.2}How to Plan a Book}{160}
\contentsline {section}{\numberline {5.3}The Importance of the Preface}{162}
\contentsline {section}{\numberline {5.4}The Table of Contents}{163}
\contentsline {section}{\numberline {5.5}Technical Aspects: The Bibliography, the Index, Appendices, etc.}{164}
\contentsline {section}{\numberline {5.6}How to Manage Your Time When Writing a Book}{168}
\contentsline {section}{\numberline {5.7}What to Do with the Book Once It Is Written}{171}
\contentsline {chapter}{\numberline {6}Writing with a Computer}{183}
\contentsline {section}{\numberline {6.1}Writing on a Computer}{184}
\contentsline {section}{\numberline {6.2}Word Processors}{189}
\contentsline {section}{\numberline {6.3}Using a Text Editor}{192}
\contentsline {section}{\numberline {6.4}Spell-Checkers, Grammar Checkers, and the Like}{193}
\contentsline {section}{\numberline {6.5}What Is \TeX\ and Why Should You Use It?}{196}
\contentsline {section}{\numberline {6.6}Graphics}{201}
\contentsline {section}{\numberline {6.7}The Internet and {\tt hypertext}}{203}
\contentsline {section}{\numberline {6.8}Collaboration by {\it e}-Mail; Uploading and Downloading}{206}
\contentsline {chapter}{\numberline {7}The World of High-Tech Publishing}{209}
\contentsline {section}{\numberline {7.1}Preprint Servers}{209}
\contentsline {section}{\numberline {7.2}\tt MathSciNet}{211}
\contentsline {section}{\numberline {7.3}Mathematical Blogs and Related Ideas}{212}
\contentsline {section}{\numberline {7.4}{\tt FaceBook}, {\tt Twitter}, {\tt Instagram}, and the Like}{214}
\contentsline {section}{\numberline {7.5}Print-on-Demand Books}{215}
\contentsline {chapter}{\numberline {8}Closing Thoughts}{217}
\contentsline {section}{\numberline {8.1}Why Is Writing Important?}{218}
\vspace*{.12in}

\contentsline {chapter}{\numberline {}Bibliography}{231}
\contentsline {chapter}{\numberline {}Index}{237}
\vfill
\eject

\hbox{ \ \ \ }

\thispagestyle{empty}

\newpage

%% Preface

\setcounter{page}{11}

\markboth{PREFACE}{PREFACE}

\chapter*{Preface to the Second Edition}

The reader response to the first edition of this book has been
gratifying.   Especially because of the Internet, and the cognate
rapid and free dissemination of ideas and writings, people
are now paying more attention to the quality of writing.  And
we are all benefiting from the result.

The essential principles of good writing have not changed for many
years.  In this new edition, I am not going to revise my advice
about grammar and syntax and organization and style.   I will certainly
update and amend and correct certain passages.  But the basic
message will be much as in the first edition.

I will still insist that writing is a yoga, and a healthy
one.  It is a discipline that one must cultivate in oneself,
and it is one that is worth cultivating.  In today's world,
good writers are respected and admired.  They are granted
a considerable measure of respect and prestige.  They are an
important part of our discipline.

What will be truly new in this second edition is an extensive discussion of
new technological developments. Today the Internet is a big part of all of
our lives, and especially of the lives of writers. Correspondingly, we are
all aware of blogs and chat rooms and preprint servers. There are now
electronic-only journals and print-on-demand books and Open Access
journals and joint research projects such as {\tt MathOverflow}.  This is just
a brief sample of the many new developments.  It is really a
whole new world, and it can be overwhelming and confusing. With some
trepidation, I have attempted here to describe and catalog this new landscape
and make it accessible to my readers.

As I lay out this new technological world, I endeavor to be as specific
as possible.  I give lots of concrete examples and plenty of detailed
description.   I leave nothing to the imagination.

The good news is that information and writing is today vastly more accessible
to virtually everyone than in past times.   This basic fact is having
a fundamental impact on communication, on education, and on scholarship.
It affects all of us in profound ways.  It is my hope that this
new edition of my writing book will help the mathematics community
to come to grips with this reality, and to derive the most benefit from
it.

As always, I thank my editor Sergei Gelfand and my readers and
colleagues for their support and their friendly, constructive
criticism. Lynn Apfel, Harold Boas, Robert Burckel, and Jerry
Folland have been particularly helpful. I look forward to
reader feedback on this new edition. 
\bigskip 
\break

\rightline{S.G.K.}
\rightline{St. Louis, Missouri}

%% Preface

\markboth{PREFACE}{PREFACE}

\chapter*{Preface to the First Edition}

The past fifty years have not seen as much emphasis on the quality of
mathematical writing as perhaps one would wish.  Because of
competition for grants and other accolades, we hasten 
our work into print.  An obituary for Hans Heilbronn
(1908-1975) asserted that, after he wrote (by hand) a draft of a
paper, he would put it on the shelf for one year.  Then he would come
back to it with fresh eyes, read it critically, and rewrite it. In
effect, after a year's time, Heilbronn was reading his own work as
though he were unfamiliar with it and had to understand each point
from first principles.  It is perhaps worth dwelling on this exercise
to see what we might learn from it.  

There is no feeling quite like that which comes after you have
proved a good theorem, or solved a problem that you have
worked on for a long time.  Driven by the heat of passion,
the words burst forth from your pen, the definitions get
punched into shape, the proofs are built and bent and patched
and shored up, and out goes that preprint to an appreciative
audience.  The whole paper sparkles---both the correct parts
and the incorrect parts.  A friend of mine, who solved
a problem after working on it to the exclusion of all else
for over fifteen years, used to rise up in the middle of the night
just to caress his manuscript lovingly.  

In circumstances like these, you find it virtually impossible to distance
yourself from the material.  Everything is emblazoned in your own
mind and is crystal clear; you are unable to take the part of the
uninitiated reader.  You are torn between the desire (expeditiously)
to record and validate your ideas, and the desire to communicate.

In today's competitive world, you probably do not feel that you have the
luxury of setting a new paper aside for a year. The paper could be
scooped; the subject could take a different direction and leave your
great advance in the dust; the NSF might cancel your grant; the dean
might not give you a raise; you might not be invited to speak at that
big conference coming up.  

Now let us look through the other end of the telescope.  The harsh
reality is this:  If you prove the Riemann hypothesis, or the\index{Riemann hypothesis!proof of}
three-dimensional Poincar\'{e} conjecture, or Fermat's Last Theorem,
then the world is willing to forgive you a lot.  It would be nice if
your paper were well written, for then more people could benefit from
it more quickly.  But---even if the paper is abysmally written---a
handful of experts will be able to slug their way through it, they
will teach it to others, perhaps more transparent proofs may come
out, textbooks will eventually appear.  Science is a process that
tends to work itself out.  

In fact most of us do not produce work that is
at the high level just described---certainly not consistently so.  If
your work is not written in a clear fashion, so that the reader may
quickly apprehend what the paper is about, what the main results are,
and how the arguments proceed, then there is a considerable
likelihood that the reader will set the paper aside before reading
much of it.  Your work will not have the impact that you had hoped or
intended.

I am certainly not writing this book to advocate that you set aside
each of your papers for a year, in the fashion of Heilbronn,
and then rewrite it.  Rather, I am asking you to consider the value
of learning to write.  Heilbronn had his techniques for sharpening
up his prose.  Each of us must learn his own.

I know many examples of mathematicians $A$ and $B$, of roughly similar
talent, with the property that $A$ has enjoyed much greater success
than $B$, and considerably more recognition for his/her ideas, because
$A$ wrote his/her work in an appealing and readable fashion and $B$ did
not.  The $A$s and $B$s that I am thinking about are not
at the Fields Medal level; Fields Medalists are exceptional in almost
every respect, and probably do not need my advice.  Instead, the examples
of which I speak are several notches down from that august level,
like most of us.
 
Even if you accept my thesis---that it is worthwhile to learn to
write mathematics well---you may feel that fine writing is not an avocation
that you wish to pursue.  Fair enough:  if you had
wanted to become a writer, then probably that is what you would have
done.  But I submit that a reasonable alternative might be to spend
an hour or two with this book, and perhaps another hour or two
considering how its precepts apply to your own writing.  The result,
I hope, will be that you will be a more effective writer and will
derive more enjoyment from the writing process.

As a scholar, or a scientist, you do not make widgets, nor do you
grow wheat, nor do you perform brain surgery.  In fact what you do is
manipulate ideas and report on the results.  Usually this report
is in written form.  What you write is often important, and can have real
impact.  Freshman composition teachers at Penn State like to tell their
students of the engineers at Three Mile Island, who wrote to the
governor of Pennsylvania three times to tell him that a nuclear
disaster was in the making at their power plant.  Their prose was so
garbled that the poor governor could not determine {\it what\/} in the
world they were talking about.  The rest is history.
                                      
The very act of writing has, in the last twenty years, taken on a new
shape and form.  Whereas, years ago, it entailed sharpening a quill
and buying a bottle of ink, nowadays most of us do not even own a
quill knife.  Instead we boot up the computer and create a document
in some version of \TeX.  This being the case, I have decided to
devote a (large) portion of this book to {\it techniques\/} of effective
writing and another (much smaller) portion to the {\it instruments\/}
of modern writing. This book is intended in large part for the novice
mathematician. Fresh from graduate school, such a person must engage
in the struggle of figuring out how to survive in the profession. The
lucky budding mathematician will have gone to a graduate program that
provided experience in technical writing and the use of hardware and
software.  If not, then perhaps the person is presently in a working
environment that makes it easy to learn the technical aspects of
writing. But I think that it is useful to have a reference for these
matters. I intend, with this book, to provide one.

My credentials for writing this book are simple:  I have written 
about one hundred articles and have written or edited about fifteen
books.  I have received a certain amount of praise for my work, and
even a few prizes; and I have received plenty of criticism.   Let me
assure you that one of the most important attributes of a good writer
is an ability to listen to criticism and to learn from it.  Anyone
finds it difficult to read criticism without becoming defensive;
nobody wants to be excoriated. But even the most negative,
uncharitable review can contain useful information.  You profit not
at all by becoming emotional; but if you can use the criticism to
improve your work then you have trumped the critic.

This book is a rather personal tract, containing personal
recommendations that reflect my own tastes.  I have reason to believe
that many others share these tastes, but not all do.  There certainly
are treatments of the art of mathematical writing that are more
scholarly than this one---I note particularly the book [Hig] of
Higham.  He has careful discussions of how to select a dictionary or
a thesaurus, careful catalogings of British usage versus American
usage, a history of mathematical notation, clever exercises for
developing skill with English syntax, tutorials on revision, and
so forth.  Higham's book is a real labor of love, and I recommend it
highly.  But there is no sense for me to duplicate Higham's efforts. 
Here I will discuss how to write, why to write, and when to write. 
However, this is not a scholarly tract, and it is not a text.  The book is
intended, rather, to be some friendly advice from a colleague.  If an
Assistant Professor or Instructor were to come to my office and ask
for suggestions about writing, then I might reply ``Let's go to lunch
and talk about it.''  This book comprises what I would say over the
course of several such meetings.

In this book I shall not give an exhaustive treatment of grammar,
nor of {\it any particular aspect\/} of writing.  When I do go into
some considerable detail, it is usually on a topic not given
extensive coverage elsewhere.  Examples of such topics are {\bf (i)} 
How to organize a paper, {\bf (ii)}  How to organize a book, {\bf
(iii)} How to write a letter of recommendation in a tenure case, {\bf
(iv)} How to write a referee's report, {\bf (v)}  How to write a book
review, {\bf (vi)}  How to write a talk, {\bf (vii)}  How to write
a grant proposal, {\bf (viii)}  How to write your Vita.

I have adopted the practice of labeling {\it incorrect\/} examples of
grammar and usage with the symbol \bad.  I do this so that examples
of what is wrong will not be mistaken for examples of what is right.

I have benefited enormously from many friends and colleagues who were
kind enough to read various drafts of the manuscript for this book.
Their comments were insightful, and in many cases essential.  In some
instances they saved me from myself.  I would like particularly to
mention Lynn S.\ Apfel, Sheldon Axler, Don Babbitt, Harold Boas,
Robert Burckel, Joe Christy, John P.~D'Angelo, John Ewing, Gerald B.\
Folland, Len Gillman, Robert E.\ Greene, Paul Halmos, David Hoffman,
Gary Jensen, Judy Kenney, Donald E.~Knuth, Silvio Levy, Chris Mahan,
John McCarthy, Jeff McNeal, Charles Neville, Richard Rochberg, Steven
Weintraub, and Guido Weiss.  George Piranian generously exercised his
editing skills on my manuscript, and to good effect.  I thank Randi
Ruden for sharing with me her keen sense of language and her sharp
wit; she showed no mercy, and spared no pains, in correcting my
language and my logic. Josephine S.\ Krantz provided valuable moral
support.  I find it a privilege to be part of a community of scholars
that is so generous with its ideas.  Pat Morgan, Antoinette Schleyer,
and Jennifer Sharp of the American Mathematical Society gave freely
of their copy editing skills.  Barbara Luszczynska, our mathematics
librarian, also gave me help in tracking down sources. My work at
MSRI was supported in part by NSF grant DMS-9022140.

It would be impossible for me to enumerate, or to thank properly, all
the excellent mathematical writers from whose work I have learned. 
They have set the example, over and over again, and I am merely
attempting to explain what they have taught us.  Several other
authors have addressed themselves to the task of explaining how to
write mathematics, or how to execute scientific writing, or simply how to
write. Some of their work is listed in the Bibliography.  (See also
[Hig] for a truly extensive enumeration of the literature.)  The
present book interprets some of the same issues from my own point of
view, and filtered through my own sensibilities.  I hope that it is a
useful contribution. 
\bigskip \break

\rightline{S.G.K.}
\rightline{St. Louis, Missouri}
\vfill
\eject

\hbox{ \ \ \ }

\thispagestyle{empty}

\newpage

%% Chapter 1
\chapter{The Basics}

\pagenumbering{arabic}

\begin{quote}
\footnotesize \sl Against the disease of writing one must take special
precautions, since it is a dangerous and contagious disease.
\smallskip \hfill \break
\null \mbox{ \ \ } \hfill \rm Peter Abelard \break
\null \mbox{ \ \ } \hfill \rm Letter 8, Abelard to Heloise
\end{quote}

\begin{quote}
\footnotesize \sl Judge an artist not by the quality of what is framed
and hanging on the walls, but by the quality of what's in the 
wastebasket.
\smallskip \hfill \break
\null \mbox{ \ \ } \hfill \rm Anon., quoted by Leslie Lamport
\end{quote}

\begin{quote}
\footnotesize \sl It matters not how strait the gate, \hfill \break
\sl How charged with punishments the scroll, \hfill \break
\sl I am the master of my fate; \hfill \break
\sl I am the captain of my soul.
\smallskip \hfill \break
\null \mbox{ \ \ } \hfill \rm W.~E.~Henley
\end{quote}

\begin{quote}
\footnotesize \sl Your manuscript is both good and original; but
the part that is good is not original, and the part that is
original is not good.
\smallskip \hfill \break
\null \mbox{ \ \ } \hfill \rm Samuel Johnson
\end{quote}

\begin{quote}
\footnotesize \sl In America only the successful writer
is important, in France all writers are important, in
England no writer is important, and in Australia you have
to explain what a writer is.
\smallskip \hfill \break
\null \mbox{ \ \ } \hfill \rm Geoffrey Cotterel
\end{quote}

\begin{quote}
\footnotesize \sl It may be true that people who are merely mathematicians
have certain specific shortcomings; however, that is not the fault
of mathematics, but is true of every exclusive occupation.
\smallskip \hfill \break
\null \mbox{ \ \ } \hfill \rm Carl Friedrich Gauss \break
\null \mbox{ \ \ } \hfill \rm letter to H.\ C.\ Schumacher [1845]
\end{quote}

\begin{quote}
\footnotesize \sl In fifty years nobody will have tenure but everyone
will have a Ph.D.
\smallskip \hfill \break
\null \mbox{ \ \ } \hfill \rm V.\ Wickerhauser 
\end{quote}

%% Section 1.1

\markboth{CHAPTER 1. THE BASICS}{}
\section{What It Is All About}
\markboth{CHAPTER 1.  THE BASICS}{1.1.  WHAT IT IS ALL ABOUT}

In order to write effectively and well, you must have something to
say\index{say something}. This sounds trite, but it is the single most
important fact about writing. In order to write effectively and well, you
also must have an audience. And you must {\it know consciously\/} who that
audience\index{audience} is. Much of the bad writing that exists is
performed by the author of a research paper who thinks that all his/her
readers are Henri Poincar\'{e}, or by the author of a textbook who does
not seem to realize that his/her readers will be students.

Good writing requires a certain confidence\index{confidence}. You must be
confident that you have something to say, and that that something is worth
saying. But you also must have the confidence to know that ``My audience
is $X$ and I will write for $X$.'' Many a writer of a mathematical paper
seems to be writing primarily to convince himself that his/her theorem is
correct, but not as an effort to communicate. Such an author is
embarrassed to explain anything, and hides behind the details. Many a
textbook author seems to be\index{audience} embarrassed to speak to the
student in language that the student will apprehend. Such an author
instead finds himself making excuses to the instructor (who either will
not read the book, or will flip through it impatiently and entirely miss
the author's efforts).

Imagine penning a poem to your one true love, all the while thinking ``What
would my English teacher think?'' or ``What would my pastor think?'' or
``What would my mother think?'' Have the courage of your
convictions\index{love letter}. Speak to that person or to those people
whom you are genuinely trying to reach. Know what it is that you want to
say and then say it.

For a mathematician, the most important writing is the writing of
a research paper. You have proved a nice theorem, perhaps a great
theorem.  You certainly have something to say.  You also know exactly
who your audience is:  other research mathematicians who are
interested in your field of study.  Thus two of the biggest problems
for a writer are already solved.  The issue that remains is how to
say it.  Remember that, if you pen a love\index{love letter}
letter to yourself, then it will have both the good features and the bad
features of such an exercise: it will exhibit both passion and fervor, but
it will tend to exclude the rest of the world. What do these remarks mean
in practice? In particular, they mean that as you write you must think of
your reader---not yourself. You must consider his/her convenience, and
his/her understanding---not your own.

In the Sputnik era, some years ago, when mathematics departments and
journals were growing explosively and everyone was in a rush to publish,
it was common to begin a paper by writing ``Notation is as in my last
paper.'' Today, by contrast, there are truly gifted mathematicians who
write papers that look like a letter home to Mom: they just start to
write, occasionally starting a new paragraph when the text spills over
onto a new page, never formally stating a theorem or even a definition,
never coming to any particular point. The contents are divine, if only the
reader can screw up the courage to pry them loose.

These last are not the sorts of papers that you would want to read, so
why torment your readers like this? Much of the remainder of this book
will discuss ways to write your work so that people {\it will\/} want to
read it, and will enjoy it when they do.

%% Section 1.2

\section{Who Is My Audience?}

If you are writing a diary, then it may be safe to say that your
audience\index{audience} is just yourself. (Truthfully, even this may not
be accurate, for you may have it in the back of your mind that---like Anne
Frank's diary, or Samuel Pepys's diary---this piece of writing is
something for the ages.) If you are writing a letter home to Mom, then
your audience is Mom and, on a good day, perhaps Pop. If you are writing a
calculus exam, then your audience consists of your students, and perhaps
some of your colleagues (or your chairman, if the chair is in the habit of
reviewing your teaching). If you are writing a tract on handle-body
theory, then your audience is probably a well-defined group of fellow
mathematicians (most likely topologists). Know your audience!

Keep in mind a specific person---somebody actually in your acquaint-
ance---to whom you might be writing. If you are writing to yourself or to
Mom, this is easy. If you are instead writing to your peers in handle-body
theory, then think of someone in particular---someone to whom you could be
explaining your ideas. This technique is more than a facile artifice; it
helps you to picture what questions might be asked,\ or what confusions
might arise,\ or which details you might need to trot out and 
explain.  It enables you to formulate the explanation of an 
idea, or the clarification of a difficult point.

I cannot repeat too often this fundamental dictum:  have something to
say and know what it is; know {\it why\/} you are saying it; finally,
know to whom you are saying it, and keep that audience always in mind.

%% Section 1.3

\section{Writing and Thought}

The ability to think clearly and the ability to write clearly are
inextricably linked. If you cannot articulate a thought, formulate an
argument, marshal data, assimilate ideas, organize a thesis, then you will
not be an effective writer. By the same token, you can use
your\index{writing!and thought} writing as a method of developing and
honing your thoughts (see [Hig] for an insightful discussion of this
concept).

We all know that one way to work out our thoughts is to engage in an
animated discussion with someone whom we respect.  But you can
instead, \`{a} la Descartes, have that discussion with yourself.  And
a useful way to do so is by writing.  When I want to work out my
thoughts on some topic---teaching reform, or the funding of
mathematics, or the directions that future research in several
complex variables ought to take, or my new ideas about domains with
noncompact automorphism group---I often find it useful to
write a little essay on the subject.  For writing forces me to
express my ideas clearly and in the proper order, to fill in logical
gaps, to sort out hypotheses from blind assumptions from conclusions,
and to make my point forcefully and clearly.  Sometimes I show the
resulting essay to friends and colleagues, and sometimes I try to
publish it. But, just as often, I file it away on my hard disc and
forget it until I have future need to refer to it. 

The writing of research level mathematics is a more formal process than
that described in the last paragraph, but it can serve you just as well.
When you write up your latest ideas for dissemination and publication,
then you must finally face the music. That ``obvious lemma'' must now be
treated; the case that you did not really want to consider must be
dispatched. The ideas must be put in logical order and the chain of
reasoning forged and fixed. It can be a real pleasure to craft your latest
burst of creativity into a compelling flood of logic and calculation. In
any event, this skill is one that you are obliged to master if you wish to
see your work in print, and read by other people, and understood and
appreciated.

Once you apprehend the principles just enunciated, writing ceases to be a
dreary chore and instead turns into a constructive activity. It becomes a
new challenge that you can aim to perfect---like your tennis backhand or
your piano playing. If you are the sort of person who sits in front of the
computer screen befuddled, frustrated, or even angry, and thinks ``I know
just what my thoughts are, but I cannot figure out how to say them'' then
something is wrong. Writing should {\it enable\/} you to express your
thoughts, not hinder you. I hope that reading this book will help you to
write, indeed will enable you to write, both effectively and well.

%% Section 1.4

\markright{1.4.  SAY WHAT YOU MEAN; MEAN WHAT YOU SAY}
\section{Say What You Mean; \hfill \break
   Mean What You Say}
\markright{1.4.  SAY WHAT YOU MEAN; MEAN WHAT YOU SAY}

You have likely often heard, or perhaps uttered, a sentence like
\begin{quote}
As a valued customer of XYZ Co., your call is very important
to us. \qquad \bad
\end{quote}
Or\index{obscure expression} perhaps
\begin{quote}
To assist you better, please select one of the following
from our menu:  \qquad \bad
\end{quote}

What is wrong with these sentences?  The first suggests that
``your call'' is a valued customer.  Clearly that is not what
was intended.  A more accurately formulated sentence would
be 
\begin{quote}
You are a valued customer of XYZ Co., and your call is very important to us.
\end{quote}
or perhaps
\begin{quote}
Because you are a valued customer of XYZ Co., your call is very important
to us.
\end{quote}
In the second example, the phrase ``To assist you better'' is
clearly intended to modify ``we'' (that is, it is ``we'' who
wish to better assist you); therefore a stronger construction is
\begin{quote}
So that we may assist you better, please select an
item from our menu \dots .
\end{quote}
or perhaps
\begin{quote}
We can assist you more efficiently if you will make a selection
from the following menu.
\end{quote}

What is the point here? Is this just pompous nit-picking? Assuredly not.
Mathematics cannot tolerate imprecision. The nature of mathematical {\it
notation\/} is that it tends to rule out imprecision. But the {\it
words\/} that connect our formulae are also important. In
the\index{precise use of language} two examples given above, we may easily
discern what the speaker intended; but, in mathematics, if you formulate
your thoughts incorrectly then your point may well be lost. Here are a few
more examples of sentences that do not convey what their authors intended:

\begin{quote}
Having spoken at hundreds of universities, the brontosaurus
was a large green lizard. \qquad \bad
\end{quote}
(Amazingly, this sentence is a slight variant of one 
that was uttered by a distinguished scholar
who is world famous for his careful use of prose.)

\begin{quote}
As in our food, we strive to be creative with keeping the
highest quality in mind, we have in our wine selections also. \hfill \break
\null \quad \qquad \qquad \qquad \qquad \qquad \qquad
\qquad \qquad \qquad \qquad \qquad \bad
\end{quote}
(This sentence was taken from the menu of a rather good St.~Louis 
restaurant.)

\begin{quote}
To serve you better, please form a line.  \qquad \bad
\end{quote}
(How many times have you heard this at your local retailer's, or
at the bank?)
\smallskip \hfill \break
\indent The message here is a simple one:  Make sure that your subject
matches your verb.  Make sure that your referents actually refer to
the person or thing that is intended.  Make sure that your
participles do not dangle.  Make sure that your clauses cohere.  {\it
Read each sentence aloud\/}.  Does each one make sense?  Would you {\it
say this in a conversation?\/}  Would you understand it if someone else
said it?

Use words carefully.  A well-trained mathematician is not likely to
use the word ``continuous'' to mean ``measurable'' nor ``convex'' to
mean ``one-connected''.  However we sometimes lapse into sloppiness
when using ordinary prose.  Treat your dictionary as a close friend: 
consult it frequently.  As a consequence, do not use ``enervate'' to
mean ``invigorate'' nor ``fatuous'' to mean ``overweight'' nor
``provenance'' to denote a geopolitical entity.  When I am
being underhanded, it is not because I am short of help.  

In life, we receive many different streams of ideas simultaneously, and we
parallel-process them in that greatest of all CPUs---the human brain. We
absorb and process information and knowledge in a nonlinear fashion. But
written discourse is linearly ordered. Word $k$ proceeds directly after
word $(k-1)$. The distinction between written language as a medium and the
information flow that we commonly experience is one of the barriers
between you and good writing. As you read this book (which purports to
tell you how to write), you \index{linear ordering of written discourse}
will see passages in which I say ``now I will digress for a moment'' or
``here is an aside.'' (In other places I put sentences in parentheses or
brackets; or I use a footnote.) These are junctures at which I could not
fit\index{writer's block} the material being discussed into strictly
logical order. You will have to learn to wrestle with similar problems in
your own writing. One version of writer's block is a congenital inability
to address this linear vs.\ nonlinear problem. In this situation, nothing
succeeds like success. I recommend that, next time you encounter this
difficulty, address it head on. After you have defeated it a few times
(not without a struggle!), then you will be confident that you can handle
it in the future.
                                                  
I have discoursed on accurate use of language in the technical sense. Now
let me remark on more global issues. As you write, you must think not only
about whether your writing is correct and appropriate, but also about
where your writing will go and what it will do when it gets there.

I have already admonished you to know when to start writing. Namely, you
begin writing when you have something to say and you know clearly to whom
you wish to say it. You also must know when to stop writing.\index{when to
stop writing} Stop when you have said what you have to say. Say it
clearly, say it completely, say it forcefully, say it without leap or
lacuna, but then shut up. To prattle on and on is not to convince further.

And never doubt that language is a weapon. ``Sticks and stones may break my
bones but words will never hurt me'' is perhaps the most foolish sentence
ever uttered. You can inflict more damage, more permanently, with words
than you can with any weapon. You can manipulate\index{language as a
weapon} more minds, and more people, with words than with any other
device. For example, when journalists in the 1960s referred to
``self-styled radical leader Abbie Hoffman'', they downgraded Hoffman in
people's minds. They never referred to Spiro Agnew as a ``self-styled [you
fill in the blank]'' or to Gordon Liddy as a ``self-styled \dots''. This
moniker was reserved for Abbie Hoffman---and sometimes for Jerry Rubin and
Mario Savio---and one cannot help but surmise that it was for a reason. By
the same token, newspapers frequently spoke in the 1960s of ``outside
agitators'' visiting university campuses. They were never described as
``colloquium speakers'' or ``expert political consultants''.

When a policeman addresses you by 
\begin{quote}
Sir, may I see your driver's license?
Did you notice that red light back there?
\end{quote}
then he is sending out one sort of signal.  (Namely, you are clearly
a law-abiding citizen and he is just doing his/her job by pulling you
over and perhaps giving you a ticket.)  When instead a cop
in the station house says

\newpage

\begin{quote}
OK, Billy.  Why don't you spill your guts?  You know that those
other bums aren't going to do a thing to protect you.  All they
care about is saving their own skins.  Jacko already confessed
to the heist and told us that you held the gun, Billy.  Now we
need to hear it from you.   Make it easy on yourself, Billy:
play ball with us and we'll play ball with you.
\end{quote}
then he is sending out a different sort of signal.  (Namely, by using
the first name---and not ``William'', but ``Billy''---he is 
undercutting the addressee's dignity; he is
treating the person like a wayward child.  Further, the policeman
is cutting off the individual from his/her peers, making him feel
as though he is on his/her own.  He is suggesting---albeit
vaguely---that he may be willing to cut a deal.)

When you are a person of some accomplishment, and some clout, then
your\index{manipulative language} writing carries considerable
responsibility. Your words may have great effect. You must weigh the
words, and weigh their impact, carefully.

I am going to conclude this section with a brief homily. (I promise that
there will be no additional homilies in the book; you may even
ignore\index{homily} this one if you wish.) Nikolai Lenin said that the
most effective way to bring down a society is to corrupt its
language.\footnote{A similar statement is attributed to John Locke in {\it
On Human Understanding}.} You need only look around you to perceive the
truth of this statement. When language is corrupted, then people do not
communicate effectively. When they do not communicate effectively, then
they cannot cooperate. When they cannot cooperate, then the fabric of
civilization begins to unravel.

Some of us use the word ``bad'' to mean ``good.'' We use the phrase ``let
us keep our neighborhoods safe and clean'' to mean ``let us segregate our
schools and arm every home''; we use the phrase ``I am for gun control and
freedom of choice'' to mean ``I'm a liberal and you're not.'' We say
``account executive'' when we mean ``sales clerk'' and ``sanitation
engineer'' when we mean ``garbage man.'' We use the words ``interesting''
to mean ``foolish,'' ``imaginative'' to mean ``irresponsible,'' and
``naive'' to mean ``idiotic.'' These observations are not just idle
cocktail party banter. They are in fact indicative of barriers between
certain social groups.

It is just the same in mathematics.  When we use the word
``proof'' to mean ``guesses based on computer printouts'' (see
[Hor]), when we use ``theoretical mathematics'' to mean ``speculative
mathematics'' (see [JQ]), when we use the phrase ``Charles
mathematicians'' to belittle the practitioners of traditional and
hard-won modes of reasoning that have been developed over many
centuries (see [Ati, pp.\ 193--196]), when we use the phrase ``new mathematics'' to
mean ``facile intuition'' (see [PS], [Ati, pp.\ 193--196]), then we are corrupting our
subject.   These are gross examples, but the same type of corruption
occurs in the small when we write our work sloppily or not at all. It
is the responsibility of today's scholars to develop, nurture, and
record our subject for future generations.  
 Good writing is of course not an end in itself; writing is
instead the means for achieving the important goal of communicating
and preserving mathematics.  

%% Section 1.5

\markright{1.5. PROOFREADING}
\section{Proofreading, Reading for Sound, \break
Reading for Sense}
\markright{1.5.  PROOFREADING}

Proofreading is an essential part of the writing process.  And it
is not a trivial one.  You do not simply write the words and
then quickly scan them to be sure that there are no gross
errors.\index{proofreading}  
Paul Halmos\index{Halmos, Paul} [Hig] said that he never published a word
before he had read it six times.  Not all of us are that careful,
but the spirit of his practice is correct:
\begin{itemize}
\item  One proofreading should be to check {\it spelling\/} and simple
{\it syntax\/}\index{spelling} 
errors (software can help with the former, and even
with the latter---see Section 6.4).
\item One\index{accuracy} proofreading should be for {\it accuracy\/}.
\item One proofreading should\index{organization} 
be for {\it organization\/} and for {\it logic\/}.
\item  One proofreading should\index{sense} be for {\it sense\/}, and 
       for the flow of the ideas.
\item  One proofreading should be for\index{sound} {\it sound\/}.
\item  One proofreading should be for overall coherence.
\end{itemize}

The great English stage actor Laurence Olivier used\index{Olivier,
Laurence} to rehearse Shakespeare by striding across the countryside
and delivering his lines to herds of bewildered cattle. 
Understandably, you may be disinclined to emulate this practice when
developing your next paper on $p$-adic $L$ functions---especially if
you live in Brooklyn.  However, note this:  all the best writers whom
I know read their work aloud to themselves.  Reading your words aloud
{\it forces\/} you to make sense of what you have written, and to
deliver it as a coherent whole.  If you have never tried this
technique, then your first experience with it will be a revelation.
You will find that you quickly develop a new sensitivity for sound
and sense in your writing.   You will develop an ``ear.''  You will
learn instinctively what works and what does not.  

Consider these simple examples.  Suppose that the Hemingway novel
{\it For Whom the Bell Tolls\/} were instead entitled {\it Who the Dingdong Rings
For}; or that the Thornton Wilder play {\it Our Town\/} were called
{\it My Turf}.  Even though the sense of the titles has not been
changed\index{sound!and sense} 
appreciably, we see that the alternative titles
eschew all the poetry and imagery that is present in the originals.
{\it For Whom the Bell Tolls} evokes powerful emotions; the
proffered alternative falls flat.   The title {\it Our Town\/} suggests
one value system, while {\it My Turf} brings to mind another.
One fancies that, if {\it The Scarlet Letter\/} had had a less poetic
title (how about {\it Bad Girls Finish Last}), 
then perhaps Hester Prynne would have garnered only
an ``$A$\,--,'' or maybe even an ``Incomplete.''

Mathematicians rarely have to wrestle with these poetic
questions. But we need to choose names for mathematical
objects; we need to formulate definitions. We need to describe
and to explain. My Ph.D. advisor thought very carefully about
his choice of notation and choice of terminology. He figured
that his ideas would have considerable influence and lasting
value, and he wanted them to come out right.

As an instance of these ideas, the word ``continuous'' is a perfect name
for a certain type of function; the alternative terminology
``nonhypererratic'' would be much less useful. The phrase ``the point $x$
lies in a relative neighborhood of $P$'' conveys a world of meaning in an
elegant and memorable fashion. Not by accident has this terminology 
become universal. You should strive for this type of precision and
elegance in your own writing.

William Shakespeare said that ``\dots a rose by any other name would
smell as sweet.''  This statement is true, and an apt observation, in
the context of the dilemma that faced Romeo and Juliet.  But the name
of a person, place, or thing can profoundly affect its future.  There
will never be a great romantic leading man of stage and screen who is
named Eggs Benedict and there will never be a Fields Medalist or
other eminent mathematician named Turkey Tetrazzini.  The name of an
object does not change its properties (consult Saul Kripke's New
Theory of Reference for more on this thought), but it can change the
way that the object is perceived by the world at large.  Bear this notion 
in mind as you create terminology, formulate definitions, and
give titles to your papers and other works.  

Have you ever noticed that, when you are reading a menu or listening
to an advertisement, it never fails that the food being described
contains ``fresh creamery butter'' and ``pure golden honey''?  The
marketing people never say ``this grub contains butter and honey,'' 
for there is nothing appealing about the latter statement.  But the
first two evoke images of delicious food.  As mathematicians, we are
not in the position of hawking victuals.  But we still must make
choices\index{detail} 
to convey most effectively a given message, and the spirit of
that message. We want to inform, and also to inspire.  Consider the
sentence

\begin{quote}
The conjecture of Gauss (1830) is false.   \qquad   \bad
\end{quote}

\noindent Contrast this rather bald statement with

\begin{quote}
The lemmas of Euler (1766) and the example of Abel (1827) led
Gauss to conjecture (1830) that all semistable curves are modular.
The conjecture was widely believed, and more than fifty 
papers were written by Jacobi, Dirichlet, and Galois in support
of it.  To everyone's surprise and dismay, a counterexample was
produced by Frobenius in 1902.  This counterexample opened many doors.
\end{quote}

There is no denying that the second passage puts the entire matter in
context, tells the reader who worked on the conjecture and why, and
also\index{terseness}
how the matter was finally resolved.  There is a tradition in
written mathematics to conform to the terse.  In your own writing,
consider instead the advantages of telling the reader what is going
on.

My advice is not to agonize over each word as you write a first
draft.  Just get the ideas down on the page.  But {\it do\/}
agonize a bit when you are editing and proofreading.
A passage that reads
\begin{quote}
This is a very important operator, that has very specific properties,
culminating in a very significant theorem.  \qquad \bad
\end{quote}
is all right as a first try, but does not work well
in the long run.  It overuses the word ``very.''  It does not
flow smoothly.  It makes the writer sound dull witted.  Consider instead
\begin{quote}
This operator will be significant for our studies.  Its spectral
properties, together with the fact that it is smoothing of order one,
will lead to our first fundamental theorem.
\end{quote}
The second passage differs from the first in that it has {\it
content\/}.  It says something.  It flows nicely, and makes the writer
sound as though he/she has something worthwhile to offer.

An amusing piece of advice, taken from [KnLR, p. 102], is never
to\index{repetitive sounds} 
use ``very'' unless you would be comfortable using ``damn'' in
its place.

A good, though not ironclad, rule of thumb is not to use the same
word,\index{alliteration} 
nor even the same sound, in two consecutive sentences.  Of
course you may reuse the word ``the,'' and the nouns that you are
discussing will certainly be repeated; but, if possible, do not
repeat descriptive words and do not place words that sound similar in
close proximity.

Also be careful of alliteration.  Vice President Spiro Agnew, with
the help of speech writer William Safire, earned for himself a
certain reputation by using phrases like ``pampered prodigies,''
``pusillanimous pussyfooters,'' ``vicars of vacillation,'' and
``nattering nabobs of negativism.''  Whatever \'{e}lan accrued to Agnew by
way of this device is probably not something that you wish to
cultivate\index{Agnew, Spiro}\index{Vietnam war} for 
yourself.  Lyndon Johnson led us 
into an escalated Vietnam war by    
deriding ``nervous nellies.''  The alliterative device is often
suitable for polemicism or poetry, but rarely so for mathematics.  For
example 
\begin{quote} 
This semisimple, sesquilinear operator serves to show sometimes that
subgroups of $S$ are sequenced.  \qquad \bad 
\end{quote} 
does not sound like mathematics.  The typical reader probably will pause,
reread the sentence several times, and wonder whether the writer is
putting him/her on.  Better is 
\begin{quote}  
Observe that this operator
is both semisimple and sesquilinear.  These properties can lead to
the conclusion that if $G$ is a subgroup of $S$ then $G$ is sequenced.  
\end{quote} 
Notice\index{alliteration} how simple syntactical tricks are used
to break up the alliteration, and to good effect.

The last two points---not to repeat words or sounds, and to avoid intrusive
alliteration---illustrate the principle of ``sound and
sense.''  If you read your work aloud as you edit and revise, then
you will pick out offending\index{sound!and sense}
passages quickly and easily. With practice, you also will learn how to
repair them. The result will be clearer, more effective writing.

%% Section 1.6

\section{Compound Sentences, Passive Voice}

It would be splendid if we could all write with the artistry of Flaubert,
the elegance of Shakespeare, and the wisdom of Goethe. In
mathematical\index{Flaubert, Gustave} writing, however, such an abundance
of talent is neither necessary\index{Shakespeare, William}\index{Goethe,
Johann} nor called for. In developing an intuitionistic ethics ([Moo]),
for example, one\index{intuitionistic ethics} presents the ideas as part
of a ritualistic dance: there is a certain intellectual pageantry that
comes with the territory. In mathematics, what is needed is a clear and
orderly presentation of the ideas.

Mathematics is already, by its nature, logically complex and subtle. 
The sentences that link the mathematics are usually most effective
when they are simple, declarative sentences.  Compound sentences 
should\index{run-on sentences}
be broken up into simple sentences.  Avoid run-on sentences
at all cost.  Here are some examples:

Rather than saying

\begin{quote}
As we let $x$ become closer and closer to 0, then $y$ tends
ever closer to $t_0$.  \qquad \bad
\end{quote}

\noindent instead say

\begin{quote}
When $x$ is close to 0 then $y$ is close to $t_0$.
\end{quote}

\noindent Of course mathematical notation allows us to write $\lim_{x
\rightarrow 0} y = t_0$ instead of either of these; this\index{notation}
abbreviated presentation will, in many contexts, be more desirable.

Rather than saying 

\begin{quote}
If $g$ is positive, $f$ is continuous, the domain of $f$ is
open, and we further invoke Lemma 2.3.6, then the set of points
at which $f\cdot g$ is differentiable is a set of the second category,
provided that the space of definition of $f$ is metrizable and separable.
\qquad \bad
\end{quote}

\noindent instead say 

\begin{quote}
Let $X$ be a separable metric space.
Let $f$ be a continuous function that is defined on an open
subset of $X$.  Suppose that $g$ is any positive function.
Using Lemma 2.3.6, we see that the set of points at which
$f\cdot g$ is differentiable is of second category.
\end{quote}

\noindent An alternative formulation, even clearer, is this:

\begin{quote}
Let $X$ be a separable metric space.
Let $f$ be a continuous function that is defined on an open
subset of $X$.  Suppose that $g$ is any positive function.
Define $S$ to be the set of points $x$ such that $f \cdot g$
is differentiable at $x$.  Then, by Lemma 2.3.6, $S$ is of second category.
\end{quote}

\noindent Note the use of the words ``suppose'' and ``define''
to break up the monotony of ``let.''   Observe how the formal
definition of the set $S$ clarifies the slightly awkward 
construction in the penultimate version of our statement.

It is tempting, indeed it is a trap that we all fall into,
to overuse a single word that means ``hence'' or ``therefore.''
An\index{hence}\index{therefore}\index{so}
experienced mathematical writer will have a clutch of
words (such as ``thus,'' ``so,'' ``it follows that,'' ``as
a result,'' and so on) to use instead.  A paragraph in which
every sentence begins with ``therefore,'' or with ``let,''
or ``so'' can be uncomfortable to read.  Have
alternatives at your fingertips.

In general, you should avoid introducing unnecessary notation.
Mary\index{notation!unnecessary} Ellen Rudin's famous statement 

\begin{quote}
Let $X$ be a set.  Call it $Y$.
\end{quote}

\noindent is funny because it is so ludicrous.  But this example is not
far from the way we write when we are seduced by notation.
A statement like 

\begin{quote}
Let $X$ be a compact metric subspace of the space $Y$.  
If $f$ is a continuous,
$\RR$-valued function on that space then it assumes both 
a maximum and a minimum value.  \qquad \bad
\end{quote}

\noindent suffers from giving\index{notation!overuse of} 
names to the metric space, its superspace, the
function, and the target space,
and then never using any of them.  Slightly better is

\begin{quote}
Let $X$ be a compact metric space.  If $f$ is a continuous,
real-valued function on $X$ then $f$ assumes both a maximum 
and a minimum value.
\end{quote}
  
\noindent Better still is

\begin{quote}
A continuous, real-valued function on a compact metric space assumes both
a maximum value and a minimum value.
\end{quote}

\noindent  The last version of the statement uses no notation,
yet conveys the message both succinctly and clearly.

Paul Halmos\index{Halmos, Paul} [Ste] asserts that mathematics should be
written so that it\index{walk in the woods, Halmos style} reads like a
conversation between two mathematicians who are on a walk in the woods.
The implementation of this advice may require some effort. If what you
have in mind is a huge commutative diagram, or the determinant of a big
matrix whose entries are all functions, then you will likely be
unsuccessful in conveying your thoughts orally. You must think in terms of
how you, or another reasonable person, would {\it understand\/} such a
complicated object. Of course such understanding is achieved in bits and
pieces, and it is achieved conceptually. That is how you will communicate
your ideas during a walk in the woods.

One corollary of the ``walk in the woods'' approach to writing is that you
should write for a reader who is not necessarily sitting in a library,
with all the necessary references at his/her fingertips. To be sure, most
any reader will have to look up a few things. But if the reader must race
to the stacks, or boot up the computer and do
a {\tt Google} search, at every other sentence, then you are making the job too
hard.  Your paper is far too difficult to follow. Supply the necessary
detail, and the proper heuristic, so that even if the reader is not sure
of a notion he/she will be able temporarily to suspend his/her disbelief
and move on.

Most authorities believe that writing in the passive voice is\index{passive
voice} less effective than writing in the active voice. To write in the
active voice is to identify the agent of the action, and to emphasize that
agent (see [Dup] for a powerful discussion of active voice vs.\ passive
voice). For example,

\begin{quote}
The manifold $M$ is acted upon by the Lie group $G$ as follows:  \hfill \break
\null \quad \qquad \qquad \qquad \qquad \qquad \qquad
\qquad \qquad \qquad \qquad \qquad \bad
\end{quote}

\noindent is less direct, and requires more words, than

\begin{quote}
The Lie group $G$ acts on the manifold $M$ as follows:
\end{quote}

\noindent Likewise, the statement

\begin{quote}
It follows that the set $Z$ will have no element of the 
set $Y$ lying in it. \qquad \bad
\end{quote}

\noindent can be more clearly expressed as

\begin{quote}
Therefore no element of $Y$ lies in $Z$.
\end{quote} 
             
\noindent Even better is

\begin{quote}
The sets $Y$ and $Z$ are disjoint.
\end{quote}

\noindent or

\begin{quote}
Therefore $Y \cap Z = \emptyset$.
\end{quote}

Notice that the last version of the statement used one word,
while the first version used fifteen.  Also, a mathematician
much more readily apprehends $Y \cap Z = \emptyset$ than
he/she does a string of verbiage.  Finally, coming up with the
succinct fourth formulation required not only restating
the proposition, but also thinking about its meaning.
The result was plainly worth the effort.

In spite of these examples, and my warnings against passive voice, I must
admit that passive voice\index{passive voice} gives us certain latitude
that we do not want to forfeit. If, in the first example, you have reason
to stress the role of the manifold $M$ over the Lie group $G$, then you
may wish to use passive voice. In the second example, it is unclear how
the use of passive voice could add a useful nuance to your thoughts. As
usual, you must let sound work with sense to convey your message.

As I have already noted, no rule of writing is unbreakable.  
The active voice is usually more effective than the passive
voice.  It is easy to criticize Lincoln's Gettysburg address
for over-use of the passive voice.
But\index{Gettysburg address}\index{Lincoln, Abraham}
Lincoln had a good ear.  If he had begun the speech with

\begin{quote}
Our ancestors founded this country 87 years ago.
\end{quote}

\noindent then he would have certainly followed the dictum of using
the active voice and using simple declarative sentences.  However he
would not have set the beautiful pace and tone that ``Fourscore and
seven years ago our fathers brought forth on this continent, a new
nation, \dots'' invokes.  He would have jumped too quickly into the
rather difficult subject matter of his speech.  (See [SW] for
the provenance of these last ideas.)

As mathematicians, we rarely will be faced with a choice analogous
to Lincoln's.  But the principle illustrated here is one worth
appreciating.

%% Section 1.7

\section{Technical Aspects of Writing a Paper}

Even when your paper is in draft form, your name should be on it.
A\index{name on work}\index{numbering!pages}
date is helpful as well.  Number the pages.  Write on one
side of the paper only.  Give the paper a working title.

Is all this just too compulsive?  No.

First, you must always put your name on your work to identify it as
your own.  If it contains a good idea, then you do not want someone
else to walk off with it.  Because you tend to generate so many
different drafts and versions of the things that you write, 
you should date your work.  I have even known mathematicians who
put\index{date on work} a time of day on each draft.  
(Of course a computer puts a date
and time stamp on each {\it computer file\/} automatically; here I am
discussing hard copy or paper drafts.)

You should write your affiliation---even on the draft. If you are usually
at\index{affiliation!on work} Harvard, then write that down. If instead
you are spending the year in Princeton, write that down. The draft could,
at some point, be circulated. People need to know where to find you. With
this notion in mind, include your current {\it e}-mail address.

If your writing is highly technical, and you are deeply involved in working
out a complicated idea, then you do not want to burden yourself with not
knowing in which order the pages go. Be sure to number them. The numbering
system need not be ``$1 \ \ 2 \ \ 3 \ \ 4 \ \ 5 \dots$.'' It could be
``$1A \ \ 1B \ \ 1C \dots$'' or ``$1_{\rm cov} \ \ 2_{\rm cov} \ \ 3_{\rm
cov}$ \dots'' (to denote your subsection on the all-important covering
lemma). In a rough draft, self-serving numbering systems like these can be
marvelously useful. On the preprint that you intend to circulate, use a
traditional sequential method for numbering the pages.

Take a few moments to think about the numbering of
theorems, definitions, and so forth.  This task is important both in
writing a paper and in writing a book.  Some authors number their
theorems from $1$ to $n,$ their definitions from $1$ to $k$, their lemmas
from $1$ to $p,$ their corollaries from $1$ to $r$---each item having its
own numbering system.  Do not laugh: this describes the default in
\LaTeX.  As a reader, I find this method maddening; for the upshot is
that I can never find anything.  For instance, if I am on the page
that contains Lemma 1.6, then that gives me no clue about where
to find Theorem 1.5.  If, instead, all displayed items are numbered in
sequence---Theorem 1.2 followed by Corollary 1.3 followed by
Definition 1.4, etc.---then I always know where I am.

Having decided on the logic of your numbering system, you also need
to decide how much information you want each number to contain. 
What\index{numbering!schemes} 
does this mean?  My favorite numbering system (in writing a
book) is to let ``$\langle\langle\hbox{\it Item\/}\rangle \rangle$ 3.6.4'' 
denote the fourth displayed item
(theorem or corollary or lemma or definition) in the sixth section of
Chapter 3.   If there is a labeled, displayed equation in the
statement of the $\langle\langle\hbox{\it Item\/}\rangle \rangle$ 
then I label it $(3.6.4.1)$.  The good
feature of this system is that the reader always knows precisely where
he/she is, and can find anything easily.  The bad feature is that the
numbering system is a bit cumbersome.  Other authors prefer to number
displayed items within each section.  Thus, in Section 6 of Chapter 3
the displayed items are numbered simply 1, 2, 3, \dots.  When
reference is later made to a theorem, the reference is phrased as
``by Theorem 4 in Section 6 of Chapter 3'' or ``by Theorem 4 of
Section 3.6.''  As you can see, this ostensibly simpler numbering
system is cumbersome in its own fashion.

The main point is that you want to choose a numbering system that
suits your purposes, and to use it consistently.  You want to make
your book or paper as easy as possible for your reader to study.  
Achieving this end requires that you attend to many small details. 
Your numbering system is one of the most important of these.

A final point is this: do not number every single thing in your
manuscript.\index{numbering!systems} This dictum applies whether you are
writing a paper or a book. I have seen mathematical writing in which every
single paragraph is numbered. Such a device certainly makes navigation
easy. But it is cumbersome beyond belief. Likewise do not number all
formulas. You will only be referring to some of them, and the reader knows
that. To number all formulas will create confusion in the reader's mind;
he/she will no longer be able to discern what is important and what is
less so.

When writing your draft (by hand), write on one side of the paper only. If
you do not, and if you are writing something fairly technical and
complicated (like mathematics), then you can become hopelessly confused
when trying to find your place. In addition, you must frequently\index{one
side of paper, writing on} set two pages side by side---for the sake of
comparing formulas, for instance. This move is easy with a manuscript
written on one side, and nearly impossible with one that is not.

If you are scrupulous about not wasting paper, and insist on using
both sides, then my advice is this:  write drafts of your
mathematical papers on one side of fresh paper.  When that work is
typed up and out the door, boldly $X$- out the writing on the front
side of each page of your old drafts.  Turn the paper over, and use
it as scratch paper, or for your laundry list.

I suggest writing in ink.  Pencil can smear, erasing can tear the
page,\index{pen vs.\ pencil}
and it is difficult to read a palimpsest.  
Also pencil-written material does not
photocopy well.  Blue pens do not photocopy well either.  I always
write with a black pen on either white or yellow paper.  I write
either with a fountain pen or a rolling writer or a fiber-tip pen so
that the pen strokes are {\it dense\/} and {\it sharp\/} and {\it dark\/}. 
I write with a pen that does not skip or blot.  If it begins to do
either, I immediately discard it and grab a new one.

Of course you cannot erase words that are written with a pen; but you
can cross them out, and that is much cleaner.  It is easier to
read a page written in bold black ink, and which includes some
crossed out passages, than to decipher a page of chicken scratch layered
over erased smears written with a pencil or written with a pen that
is not working properly.

Be sure that your desk is well stocked with paper, pens, Wite-Out\reg,
Post-it\reg\ notes, a stapler, staples, a staple remover, 
cellophane tape, paper clips,
manila\index{desk supplies} 
folders, manila envelopes, scissors, a dictionary, and
anything else you may need for writing.  Have them all at your
fingertips.  You do not want to interrupt the precious writing
process by running around and looking for something trivial.

Do not write much on each page.  I advise writing {\it large\/},
and double or triple spaced.  The reason?  First, you want
to\index{spacing on page} be able to insert passages, make editorial remarks, make
corrections, and\index{writing!large} 
so forth.  Second, a page full of cramped writing
on every line is hard to read.  Third, you can more easily 
rearrange material if there is just a little on each page.  
For example, if one
page contains the statement of the main theorem and nothing else,
another contains key definitions and nothing else, and so forth,
then you can easily change the location of the main theorem in
the body of the paper.  If the main theorem is buried in a page
with a great deal of other material, then moving it would involve
either copying, or photocopying, or cutting with scissors.

Do not hesitate to use colored pens.  For instance, you could be
writing\index{colored pens} 
text in black ink, making remarks and notes to yourself
(like ``find this reference'' or ``fill in this gap'') in red ink,
and marking unusual characters in green ink.  This may 
sound compulsive, but it makes the editing process much easier.
                                                     
A good bibliography is an important component of scholarly work (more
on\index{bibliography}\index{bibliographic!references} bibliographies can
be found in Sections 2.6, 5.5). Suppose that you are writing a paper with
a modest number of references (about 25, say), and you are assigning an
acronym to each one. For instance, [GH] could refer to the famous book by
Griffiths and Harris. When you refer to this work while you are writing,
use the acronym. Keep a sheet of notes to remind yourself what each
acronym denotes. Do not worry about looking up the detailed bibliographic
reference while you are engaged in writing; instead, compartmentalize the
procedure. When you are finished writing the paper, you will have a
complete, {\it informal\/} list of all your references. You can go to
{\tt MathSciNet} (Section 7.2) OnLine and find most of your references in an instant.
You can also go to your library's catalog OnLine to find locally
obtainable references. \LaTeX\ can be a great help in eliminating much of the
tedium of assembling bibliographies. See the discussion in Sections 2.6
and 5.5.

Let me make a general remark about the writing 
process.  As you are writing a paper, there will be several
junctures at which you feel that you need to look something up:  either
you cannot remember a theorem, or you have lost a formula, or you
need to imitate someone else's proof.  My advice is {\it
not\/} to interrupt yourself while you are writing. Take your red pen
and make a note to yourself about what is needed.  But {\it keep
writing\/}.  When you are in the mood to write, you should take
advantage of the moment and do just that.  Interrupting yourself to
run to the library, or for any other reason, is a mistake.

Write on a desk that is free of clutter.  It is romantic, to
be\index{clutter on desk}
sure, to watch a film in which the writer labors furiously on a
desk that is awash with papers, books, hamburger bags, ice cream
containers, old coffee cups, last week's underwear, and who knows
what else.  Leave that stuff to the movies.   Instead imagine tearing
into page 33 of your manuscript and accidentally spilling a week-old
cup of coffee and a piece of pepperoni pizza all over your project. 
Think of the time lost in mopping up the mess, separating the pages,
trying to read what you wrote, copying your pages, and so forth.  
Enough said.

If you are going to drink coffee or a soda or eat a sandwich while
you work, I suggest having the food on a small separate side table. 
This little convenience will force you to be careful, and if you do
have an accident then it will not make a mess of your work.

Write in a place where you can concentrate without interruption.
Whether\index{concentration} 
you have music going, or a white noise machine playing,
or a strobe light flashing is your decision.  But if you are going
to concentrate on your mathematics, it may take up to an hour
to get the wheels turning, to fill your head with all
the ideas you need, and to start formulating the necessary assertions.
After you have invested the necessary time to tool up, you want
to use it effectively.  Therefore you do not want to be interrupted.
Close the door and unplug the telephone if you must.  Victor Hugo
used to remove all his clothes and have his servant lock him in a room
with\index{Hugo, Victor}
nothing but paper and a pen.  Moreover, the servant guarded
the door so that the great man would not be interrupted by so much
as a knock.  This method is not very practical, and is perhaps not
well suited to modern living, but it is definitely in the right spirit.

%% Section 1.8

\markright{1.8.  MORE SPECIFICS OF MATHEMATICAL WRITING}
\section{More Specifics of \hfill \break 
Mathematical Writing}
\markright{1.8.  MORE SPECIFICS OF MATHEMATICAL WRITING}

For the most part, the writing of mathematics is like the writing of
English prose.  Indeed, it {\it is\/} a part of the writing of English. 
({\it Caveat:}  I hope that my remarks have some universality, and apply 
even if you are writing mathematics in Tagalog or Coptic or
Tlingit.)  If you read your work aloud (I advocate this practice in
Section 1.5), then you should be reading complete sentences that flow
from one to the next, just as they do in good prose.

It\index{writing!mathematics vs.\ writing English}
is all too easy to write a passage like

\begin{quote}
Look at this here equation:
$$
x^n + y^n = z^n .  \qquad \mbox{\bad}
$$
\end{quote}

\noindent Much smoother is the passage

\begin{quote}
The equation
$$
x^n + y^n = z^n
$$
tells us that Fermat's Last Theorem is still alive.
\end{quote}

\noindent Another example of good sentence structure is

\begin{quote}
Since
$$ 
A < B
$$
we know that \dots.
\end{quote}
Notice that the the sentence reads well aloud:  ``Since $A$ is less
than $B$ we know that \dots.''
 
An aspect of writing that is peculiar to mathematics is the use
of\index{notation!use of}\index{notation!abuse of} notation. Without
good notation, many\index{notation!good} mathematical ideas would be
difficult to express. Indeed, the development of mathematics in the middle
ages and the early renaissance was hobbled by a lack of notation. With
good notation, our writing has the potential to be elegant and compelling.

A common misuse of notation is to put it at the beginning of a sentence or a clause.
For example,
\begin{quote}
Let $f$ be\index{notation!at the beginning of a sentence}
a function.  $f$ is said to be {\it semicontinuous\/} if \dots  \hfill \break
\null \quad \qquad \qquad \qquad \qquad \qquad \qquad
\qquad \qquad \qquad \qquad \qquad \bad
\end{quote}

\noindent and

\begin{quote}
For most points $x$, $x \in S$.  \qquad \bad
\end{quote}

Even in these two simple examples you can begin to apprehend the 
problem:  the eye balks at a sentence or clause that is begun
with a symbol.  You find yourself rereading the passage a couple
of times in order to discern the correct sense.  Much better is:

\begin{quote}
A function $f$ is said to be {\it semicontinuous\/} if \dots
\end{quote}

and

\begin{quote}
We see that $x \in S$ for most points $x$.
\end{quote}

\noindent Observe that both of these revisions are easily
comprehended the first time through.  That is one of the goals of
good writing.

Mathematical notation is often so elegant and compelling that we are
tempted to overuse, or misuse, it.  For example, the notation in the
sentence ``If $x > 0$, then $x^2 > 0$'' is no hindrance, is easy to
read,\index{notation!overuse of} 
and tends to make the sentence short and sweet (nonetheless,
there are those who would tender cogent arguments for ``If a number
is positive then so is its square.'').  By contrast, the phrase
\begin{quote}
Every real, nonsquare $x < 0$ \dots  \qquad \bad
\end{quote} 
is objectionable.  The
reason is that it is not clear, on a first reading, what is meant.
Are you saying that ``Every real, nonsquare $x$ is negative'' or are
you saying ``Every real, nonsquare $x$ that is less than zero has
the additional property \dots.''  By strictest rules, the notation \ $<$
\ is a {\it binary connective\/}.  The notation is designed for expressing the
thought $A < B$.  If that is not the exact phrase that fits into your
sentence, then you had best not use this notation.

When you are planning a paper, or a book, you should try to plan your
notation in advance. You want to be consistent throughout\index{notation!planning} the work in question. To be sure, we have all seen works that,
in Section 9, say ``For convenience we now change notation.'' All of a
sudden, the author stops using the letter $H$ to denote a subgroup and
instead begins to use $H$ to denote a biholomorphic mapping. Amazingly,
this\index{notation!consistency of} abrupt device actually works much of
the time---at least with professional mathematicians. But you should avoid
it. If you can, use the same notation for a domain in Section 10 (or
Chapter 10) of your work that you used in Section 1 (or Chapter 1). Try to
avoid local contradictions---like suddenly shifting your free variable
from $x$ to $y$. Try not to use the same character for two different
purposes.\footnote{When Andr\'{e} Weil was writing his book {\it Basic
Number Theory} [Wei], he strove mightily to follow this advice. He used up
all the roman letters, all the Greek letters, all the fraktur letters, all
the script letters, all the Hebrew letters, and all the other commonly
used characters that are seen in mathematics. He ended up resorting to
Japanese characters.}

The last stipulation is not always easy to follow.  Many of us
commonly use $i$ for the index of a series or sequence:
$\sum_{i=1}^\infty a_i$ and $a_i$.  No problem
so far, but suppose that you are a complex analyst, and use $i$
to denote a square root of $-1$.  And now suppose that this
last $i$ occurs in some of your sequences and series.  You can
see the difficulties that would arise.  It is probably best to use
$j$ or $k$ as the index of your sequence or series.  A little planning can
help with this problem, though in the end it may involve a great
deal of tedious work to weed out all notational ambiguities.

Many a budding mathematician is seduced by mathematical notation.
There was a stage in my education when I thought that all of
mathematics\index{notation!overuse of} 
should be written without words.  I wrote
long, convoluted streams of $\forall\, ,\, \exists\, ,\, \ni:\, , \,
\Rightarrow\, , \, \equiv$\, ,\, and so forth.  
This style would have served me well had
I been invited to coauthor a new edition of
{\it Principia Mathematica\/} (see [WR]).  In modern
mathematics, however, you should endeavor to use English---and to
minimize the use of cumbersome notation.  Why burden the reader with
$$
\forall x \exists y , x \geq 0 \Rightarrow y^2 = x  \qquad \mbox{\bad}
$$
when you can instead say
\begin{quote}
Every nonnegative real number has a square root.
\end{quote}

The most important logical syllogism for the mathematician is 
{\it modus ponendo ponens\/}, or ``if \dots then.''  If you
begin\index{modusponendo@{\it modus ponendo ponens}} 
a sentence\index{if-then}  
with the word ``if,'' then do not forget to include
the word ``then.''  Consider this example:

\begin{quote}
If $x > 4$, $y < 2$, the circle has radius at least 6, 
the sky is blue, the circle can be squared. \qquad \bad
\end{quote}

\noindent Which part of this sentence is the hypothesis and which
the conclusion?  After a few readings you may be able to figure
it out. If it were sensible mathematics then the mathematical meaning
would probably give you some clues.  But it is clearer to write

\begin{quote}
If $x > 4$, $y < 2$, the circle has radius at least 6, and the sky is
blue, then the circle can be squared.
\end{quote}

Following the dictum that shorter sentences are frequently preferable
to longer\index{short sentences vs.\ long sentences} 
ones, you can express the preceding thought even more
succinctly as

\begin{quote}
Suppose that $x > 4$, $y < 2$, the circle has radius at least 6,
and the sky is blue.  Then the circle can be squared.
\end{quote}

\noindent The word ``then'' is pivotal to the logical structure here.
It acts both as a connective and as a sign post.  The reader can 
(usually) figure out what is meant if the word ``then'' is omitted.
But the reader should not {\it have\/} to do so.  Your job as the
writer is to perform this task {\it for\/} the reader.  Mathematicians
have a tendency to want to jam everything into one sentence. 
However, as the last example illustrates, greater clarity can often
be achieved by breaking things up; this device also forces you to
think more clearly and to organize your thoughts more effectively.

Mathematicians commonly write ``If $f$ is a continuous
function, then prove $X$.''  A moment's thought shows that
this\index{thenprovethat@``then prove that''} is not the intended meaning:  the
desire to prove $X$ is not contingent on the continuity of $f$.  What
is intended is ``Prove that, if $f$ is a continuous function, then
$X$.''  In other words, the hypothesis about $f$ is part of what
needs to be proved.

Sometimes you need to write a sentence that treats
a word\index{word!treated as an object} 
as an object.  Here is an example:

\begin{quote}
We call $\Gamma$ the {\it fundamental solution\/} for the partial
differential operator $L$.  We use the definite article ``the'' because,
suitably normalized, there is only one fundamental solution.
\end{quote}

I have oversimplified the mathematics here to make a typographical
point.  First, when you define a term (for the first time), 
you should italicize the word or phrase that is being defined.  Second, when
you refer to a word (in this case ``the'') as the object of discussion,
then put that word in quotation marks.  For a variety of psychological
reasons,\index{Quine, W.\ V.\ O.} 
writers often do not follow this rule.  It is helpful
to recall W.\ V.\ O.\ Quine's admonition:  `` `Boston' has six letters.
However Boston has 6 million people and no letters.'' 

The phrase ``if and only if'' is a useful mathematical device.  It
indicates logical equivalence of the two phrases that it connects.
While\index{if and only if}
the phrase is surely used in some other disciplines, it plays a 
special role in mathematical writing; we should take some care to
treat it with deference.  Some people choose to write it as ``if, and
only if,''---with two commas.  That is perfectly grammatical, if a
little stilted.  One habit that is unacceptable (because it sounds
artificial and is difficult to read) is to begin a sentence with this
phrase.  For instance,

\begin{quote}
If and only if $x$ is nonnegative, can we be sure that the real
number $x$ has a real square root.  \qquad \bad
\end{quote}

\noindent That is a painful sentence to read, whether the reading is done
aloud or {\it sotto voce\/}.  Better is

\begin{quote}
A real number $x$ has a real square root if and only if $x \geq 0$.
\end{quote}

An alternative form, not with universal appeal (but better than beginning
a sentence with ``if and only if''), is

\begin{quote}
Nonnegative real numbers, and only those, have real square roots.
\end{quote}

The neologism ``iff,'' reputed to have been popularized by Paul
Halmos,\index{iff}\index{Halmos, Paul}
is a generally accepted abbreviation for ``if and only if.''
This is a useful bridge between the formality of ``if and only if''
and the convenience of ``if.''	It is also common to use the
symbol $\Longleftrightarrow$ for ``if and only if.''

Word order can have a serious, if subtle, effect on the meaning
(or at least the nuance)\index{word!order} of a sentence.  The examples

\begin{quote}
Yellow is the color of my true love's hair.
\end{quote}
\begin{quote}
My true love's hair has the color yellow.
\end{quote}
\begin{quote}
The hair, which is yellow, of my true love \dots
\end{quote}
each say something different, as they emphasize a different
aspect---either the color, or the person, or the hair---that
is being considered.  (As an exercise, insert the word ``only'' 
into all possible positions in the sentence
\begin{quote}
I helped Carl prove quadratic reciprocity last week.
\end{quote}
and watch the meaning change.)

In mathematics, word order can seriously alter the meaning of a
sentence, with the result that the sentence is not immediately
understood---if at all.  When you proofread your own work, you tend
to supply\index{word!order} 
meaning that is not actually present in the writing; the 
result is that you can easily miss obscurity imposed by word order.
Reading your work aloud can help cut through the problem.

Do not overuse commas.  I become physically ill when I see
a sentence like
\begin{quote}
We went\index{commas, overuse of} 
to the store, to buy some potatoes.  \qquad \bad
\end{quote}
Slightly more subtle, but still irksome, is
\begin{quote}
Now that we have our hypotheses in place, we state our theorem,
with the point in mind, that we wish to understand the
continuity, of functions in the class ${\cal S}$. \qquad \bad
\end{quote}
We certainly use a comma to indicate a pause.
But the comma indicates a {\it logical pause\/}, not a lack
of air or lack of good sense.
Read the last displayed sentence out loud, with suitable pauses where
the commas occur.  It sounds like someone huffing and puffing;
the pauses have no reason to them.
This sentence is not 
a representative example of the way that we speak, 
hence it is not indicative of the way that
we should write.  Much more attractive is
\begin{quote}
Our hypotheses are now in place, and we next state our theorem.
The point is to understand the continuity properties of functions
belonging to the class ${\cal S}$.
\end{quote}

Mathematicians like the word ``given.''  We tend to overuse
and misuse it---especially in instances where the word can
be\index{given} 
discarded entirely.  Consider the example ``Given a metric
space $X$, and a point $p \in X$, we see that \dots.''  More
direct is ``If $X$ is a metric
space and $p \in X$, then \dots.''  We are often tempted
to transcribe spoken language and call that written language;
such laziness should be defeated.
Our misuse of ``given'' is an example of such sloth.

Whenever possible, use singular constructions rather than plural.
Consider\index{singular constructions vs.\ plural} the sentence
\begin{quote}
Domains with noncompact automorphism groups
have orbit accumulation points in their boundaries.  \qquad \bad
\end{quote}

\noindent First, such a construction is quite awkward:  should it be
``groups'' or ``group''?  More importantly, do all the domains
share the same automorphism group, or does each have its own?
Does each domain have several orbit accumulation points, or just one?
Clearer is the sentence
\begin{quote}
A domain with noncompact automorphism group has
an orbit accumulation point in its boundary. 
\end{quote}

When you are putting the final polish on a manuscript, look it over
for general appearance.  In mathematical writing, several consecutive
pages of dense prose are not\index{prose vs.\  mathematics}
appealing, nor are several
consecutive pages of tedious calculation.  For ease of reading, the
two types of mathematical writing should be interwoven.  It requires
only a small extra effort to produce a paper or book with comfortable
stopping places\index{stopping places for reader}
on every page.  The reader needs to take frequent
breathers, to survey what he/she has read, to pause and look back.  Make
it easy for him/her to do so.

While you are thinking about the counterpoint between prose and
formulas, think also about the use of displayed math versus in-text
math\index{displayed math vs.\ in-text math}
[in \TeX\ (see Section 6.5), the former is set off by double dollar signs \verb@$$@
while the latter is set off with single dollar signs \verb@$@].  
Long formulas are usually better displayed, for they are difficult to
read when put in text.  Of course {\it important\/} formulas should be
displayed no matter what their length---and provided with numbers or
labels if they will be mentioned later.  Do not display every single
formula, for that will make your paper a cumbersome read.  Also do not
put every formula in text, as that will make your writing tedious.  A
little thought will help you to strike a balance, and to use the two
formats to good effect.

And now a coda on the role of English in mathematical writing. More
and\index{English, role of in mathematics} more, English is becoming the
language of choice in mathematics. Therefore those of us who are native
speakers set the standard for those who are not. We should exercise a bit
of care. I have a good friend, also an excellent mathematician, who is
widely admired; his fans like to emulate him. He is fond of saying
(informally) ``What you need here is to cook up a function $f$ such that
\dots.'' Mathematicians of foreign extraction, who have been hearing him
make this statement for years, have now developed the habit of saying
``Take a function $f$. Now cook it for a while \dots.'' It is a bit like
having your children emulate (poorly) all your bad habits. A word to the
wise should suffice. 

%% Section 1.9

\section{Pretension and Lack of Pretension}

Avoid the use of big words when small ones will do.  Do not say
``peregrinate'' for ``walk,'' nor ``omphaloskepsis'' for ``thought,''
nor\index{big words vs.\ small words} 
``floccinaucinihilipilificate'' for ``trivialize'' unless the
longer word conveys some important nuance that the shorter word does
not.  The urge to so bloviate should be resisted.  To indulge in
hippopotomonstrosesquipedalian tergiversation is not to show your
erudition; rather, it is to be superficial.  Also remember that many of
your readers will be foreign born, not native English speakers.  Make
some effort to write simple, straightforward English that they will
easily apprehend.  Save your high-flown rodomontade for ceremonial
occasions.  

Likewise---and I have said this elsewhere in the book---stick to
simple sentence structures.  Even the subjunctive mood can lead to
confusion when it is\index{simple sentence!structures}
used in mathematical writing.   Let the
mathematics speak for itself; do not try to dress it up with fancy
language.

You can have some fun peppering your prose with {\it bon vivant\/} and {\it
Gem\"{u}tlichkeit\/} and {\it ad hominem\/} and {\it samizdat\/}, but
the\index{foreign words and phrases}
careless use of foreign words and phrases does not add anything
to most writing.  And it will confuse many readers.  Use
foreign phrases sparingly.  If you do use them, typeset them in
italics.  (An exception should be made for foreign words like
``etc.'' (short for {\it et cetera\/}), which have become standard
parts of the English language and should be set in roman.)  The books
[Hig], [Por], [SG], [Swa] give more detailed treatments of this topic.

Good mathematics is difficult.  Do not let your writing be a device
for making it more so.  Use simple, declarative sentences---short
ones.\index{simple sentence!vs.\ complex sentence}  
Use short\index{short paragraphs vs.\ long paragraphs} 
paragraphs, each with a simple point.  To understand
my meaning, put yourself in the position of the reader.  You are
slugging your way through a tough paper.  You come to the proof of
the main theorem.  After killing yourself for a couple of hours, you
finally come to the crux of the argument.  And it is a single, dense
paragraph that is two pages long.  Such a daunting prospect is truly
depressing.  You do not want to abuse your readers in this fashion. 
Break up the ideas into palatable bites.

And now a note on flippancy.  A friend of mine once wrote a truly
elegant---and important---book that included the phrase ``the reader
should\index{flippancy}
review enough functional analysis so that he does not barf
[{\it sic\/}] at the sight of a Banach or Frechet space.''  At the
reviewer's insistence, the phrase was toned down before publication.
Another friend published a book with the phrase ``we leave the details
of this proof for the mentally infirmed.''  I would advise against
this sort of sarcasm.  This suggestion is not simply a nod to
propriety.  You want to be proud of your work.  Remember that your
thesis advisor and the authorities in the field are likely to look at
it.  Such puerile prose is not what you want them to see.  Most
likely, ten years hence, you will wish fervently that you
had not included such phrases.  Anyone who continues to grow
intellectually will look on his/her work of ten years ago with some
disdain.  But there is no percentage in adding embarrassment to the
mix.

Suit your tone, and your choice of words, to the subject at hand.
It\index{tone}\index{choice of words}
might be suitable to use the phrase ``He had all the efficiency
and dexterity of a ruptured snail'' to describe a clumsy waiter;
this is probably not appropriate language for describing the pope.

Finally, stay away from faddish prose.  If you say ``fraternally
affiliated,\index{faddish prose} 
ethically challenged young male'' to mean ``gang member'' or
``peregrinating, fashion-challenged, pulchritudinally advanced hostess''
to mean ``prostitute,'' then you may be politically correct
today but you will be strictly out to lunch tomorrow.
Today, many a writer or speaker wants to work the word ``dis''
(gang talk for ``disrespect'') or ``flame'' (yuppie talk
for ``disrespect'') into his/her prose.  This practice is a mistake,
because in ten years the words will have no meaning.

By the same token, avoid old-fashioned modes of expression.  In 1827
it\index{old-fashioned prose}
was appropriate for a physician to diagnose a patient with
``falling crud and palpitation of the pluck''; in 1930 it was
fashionable for a woman to complain of ``the fantods.''  Today these
phrases are meaningless.  It might exhibit devotion to Fermat to use
``ad\ae quibantur'' instead of \ ``$=$'' \ (as did he), but such a 
practice would lead to boundless confusion today.

Some American writers think that it is tony to pepper their writing
with\index{British English} 
British English.  They use ``humour'' for ``humor,''
``lorry'' for ``truck,'' and ``spanner'' for ``wrench.''
Such language is out of place, and can only lead to obfuscation.  It would
be just as foolish for an American cookbook to give recipes
for spotted dick, bubble-and-squeak, and stodge.
Nobody would know what the author was talking about.

For the same reasons I advise against using contractions,
using abbreviations, or using slang---at least in formal writing. 
Even acronyms are dangerous (see Section 1.12); use them with
caution.\index{contractions, use of}  
We write because we want our thoughts to last, and to be
comprehensible both now and in years hence.  Do not let language stand
in the way of that goal.

%% Section 1.10

\section{We vs.\ I vs.\ One}

When I was a child, I once asked a mathematician why mathematics was
usually written in the first person plural:  ``We now prove this'';
``Our next task is thus''; ``We conclude our story as follows.''
The rejoinder that I received was ``This is so that the reader will think
that\index{ivswevsone@`I' vs.\ `we' vs.\ `one'} there are a lot of you.''

More seriously, when you are writing up mathematics, then you must
make a choice.   You can say ``I will now prove Lemma 5'' or
``We will now prove Lemma 5'' or ``One may now turn one's 
attention to Lemma 5.''  Which is correct?

As with many choices in writing, this one involves a degree of
subjectivity.  I shall now tell you what I think about the matter.
The first option is rarely chosen.  Most people consider it pompous
and inappropriate.  The only instance where I find the first person
singular to be a comfortable choice is the following:  sometimes at
the end of a paper one says ``At this time I do not know how to prove
Conjecture A.'' The choice is appropriate for this particular
statement because in fact you are imparting to the reader some
specific information about what you yourself know.  It would be
misleading, and a trifle affected, to say ``At this time {\it one\/}
does not know \dots.''  Likewise for ``At this time {\it we\/} do not
know \dots.''  However, you could say, ``At this time it is not
known whether \dots.''

The custom in modern mathematics is to use the first person plural,
or ``we.''  It stresses the\index{weuseof@`we', use of} 
participatory nature of the enterprise,
and encourages the reader to push on.  Moreover, since ``we'' is what
people are accustomed to hearing, it is less likely to jar their
ears, or to distract them, than one of the other choices.   The use
of third person singular, or ``one,'' often leaves the writer
struggling with awkward sentence structures.  If you endeavor to
write in that mode, then you will likely find yourself soon breathing
a sigh of relief as you abandon it.  If you read with sensitivity,
you also will likely learn that first person singular, or ``I,'' is
irritating; therefore you will not use it.

With a little craftsmanship, you can avoid entirely the use
of the first person in your writing.  Rather than say ``We now
turn to the proof of Lemma 4,'' instead say ``Next is the proof
of Lemma 4'' or perhaps ``The next task is the proof of Lemma 4.''
Rather than say ``We see that the proof is complete,'' say
``The proof is now complete'' or ``This completes the proof.''
The book [Dup, Ch.\ 2] has a sensible and compelling discussion of
the question of ``We'' vs.\ ``I'' vs.\ ``one.''

Sound and sense will dictate which of the words ``I,'' ``we,'' or 
``one''---or\index{sound!and sense} 
perhaps none of these---you wish to use.  I am offering
``we'' as the default.  But the sense of what you are writing may
dictate another choice.

%% Section 1.11

\markboth{CHAPTER 1.  THE BASICS}{1.11.  ESSENTIAL RULES OF GRAMMAR}
\section{Essential Rules of Grammar, Syntax, and Usage}
\markboth{CHAPTER 1.  THE BASICS}{1.11.  ESSENTIAL RULES OF GRAMMAR}

I have intentionally put this discussion of the rules of grammar and
syntax and usage at the end of Chapter 1.  The reasons are
several. First I want, in a gentle way, to de-emphasize them.  I am
not one of those who says ``the battle against `hopefully' is lost,''
``the battle for `which' vs.~`that' is lost,'' ``the battle for `lay'
vs.~`lie' is lost,'' and so forth.  I find such statements facile, and
they miss the point that careful writing requires some precision. The
argument ``You know what I mean; whether I use `that' or `which' is
incidental'' abrogates the fact that accurate writing, and accurate
expression\index{rules of grammar} 
of your thoughts, requires accurate use of language. But
you do not develop skill as a writer by concentrating on the rules of
language; they are merely a set of tools that are used in the
process.

The intent of this book is that you should learn to write
logically and cogently; to say precisely what you mean, using just
the right number of words; to eschew obfuscation.  You want to
develop an ear, so that clear writing becomes natural.  Exact use of
the language is a part of the process.  But it is not the main
point.

Most of the rules of English usage are succinct and logical.  A
particularly concise enunciation of the basic rules appears in [SW]. 
Since I cannot improve on that presentation, I certainly shall not
repeat the rules of grammar here.  It is a revelation for any adult
writer to review the rules of usage and to see what eminent sense
they make.  Here I shall mention just a few sticky points that come
up frequently in mathematical and other writing.  I hope that you
will find this section, and the next, to be a useful
``quick-and-dirty'' reference. With that goal in mind, I have presented
the topics in alphabetical order.  See also [Chi], [Dup], [Fow],
[Fra], [Hig], and [MW] for a more thorough treatment of issues of
grammar, syntax, and usage.

Bear in mind, as you read these precepts, that no rule of English 
grammar is etched in stone.  There will certainly be times that a
sentence or phrase formed according to the strictest rules will sound
just awful.  In such an instance, you must override the rules and
use your good sense and taste.  More will be said about this
technique as the book develops.

Now for some rules:

\begin{itemize}
\item{\bfit All, Any, Each, Every}\ \ In
mathematics\index{allanyeach@`All', `Any', `Each', `Every'}
we commonly formulate statements such as ``Show that any
continuous function $f$ on the interval $[0,1]$ has a point $M$ in
its domain such that $f(M) \geq f(x)$ for $x \in [0,1]$.''  For
cognoscenti it is clear that, when we say ``any'' here, we mean
``all.''  But for others---for students, or for nonnative
speakers---this slight abuse of language could cause confusion.  For
example, a student reading this sentence could (perfectly correctly)
construe it to mean ``Demonstrate that for {\it some\/} function $f$
\dots.''  Thus, if this sentence were part of an exercise, the student
might answer 
\begin{quote}
The function $f(x) = - (x - 1/2)^2$ is continuous on
$[0,1]$ and the point $M = 1/2$ satisfies the conclusion.
\end{quote}
The lesson is to avoid using ``any'' when ``all'' or ``each'' or
``every'' is intended.

\null \quad \ Conversely, even when you are writing for experts you can cause
confusion by misusing quantifiers.  Sentences like
\begin{quote}
All continuous functions have a maximum. \qquad \bad
\end{quote}
are far too common in mathematical writing.  Notice that the sentence
suggests that all continuous functions share the {\it same\/}
maximum.  Of course what was intended was
\begin{quote}
Every continuous function has a maximum.
\end{quote}
or, more precisely,
\begin{quote}
Each continuous function has a maximum.
\end{quote}
(Once again we see the advantage, from the point of view of clarity,
of the singular over the plural.)
As you proofread your work, you must learn to take the part of
the reader (who is not {\it a priori\/} sure of what is being said) 
in order to weed out misused quantifiers.

\item{\bfit Brevity}\ \  Endeavor\index{brevity} to formulate your thoughts 
briefly and succinctly.  For
example, you {\it could\/} say
\begin{quote}
In point of fact, we devolved upon the decision to solicit
opinions, form an enumeration, and produce a tally. \hfill \break
\null \mbox{} \qquad \qquad \qquad \qquad \qquad \qquad \qquad
\qquad \qquad \qquad \qquad \qquad \bad
\end{quote}
Such a sentence sounds mellifluous, sanguine, and high toned.
But why not instead say
\begin{quote}
We decided to take a vote.
\end{quote}
The second sentence says in 6 words what the first said in 19; and it
presents the message more clearly and forcefully.  Strunk and White
[SW] give a thorough and engaging treatment of the topic of brevity,
and they speak particularly cogently of eliminating extra or extraneous
words.  Mathematics is difficult to read under the best of
circumstances.  Do not make the reader's job even more difficult by
weighing down your prose with excess baggage.

I once saw a sign in the elevator of a Washington, D.C. hotel that said
\begin{quote}
Do not carry lighted tobacco products in the elevator.
\end{quote}
I can only suppose that some politician created this sign.  Why
not just say
\begin{quote}
No smoking.
\end{quote}

\item{{\bfit cf., e.g., i.e., n.b., q.v.,} \bf and the like}\ \
These\index{cf.}\index{e.g.}\index{i.e.}\index{n.b.} 
are\index{q.v.} abbreviations for specific Latin expressions:  {\it confer\/},
{\it exempli gratia\/}, {\it id est\/}, {\it nota bene\/}, {\it quod
vide\/}.  They have particular meanings, and you should
strive to use them accurately. In particular, ``cf.'' is often
misused to mean ``see.''  It actually means ``compare.''  Sometimes ``e.g.''\ 
and ``i.e.''\ are interchanged in error; the first of these
means ``for example,'' and the second means (literally) ``the favor
of an example'' or (more familiarly) ``for the sake of example.'' It is
difficult to use ``n.b.'' with grace.  If you are unsure, then use
the English equivalent of which you {\it are\/} sure.

\null \quad \  In fact it is difficult to make a compelling case for
``i.e.'' in favor of ``that is,'' or for any of the other Latin
substitutes in favor of their English equivalents.  The punctuation
and font selection questions connected with these Latin abbreviations
are tricky (see [Hig] or [Fow] or [Chi] or [SK]).  To repeat, use
them with care.
  
\item{\bfit Comprise {\bf vs.} Compose}  \ \ People use the word ``comprise''
because they think it makes them sound tony.  Unfortunately, because
most everyone misuses the word, they instead sound uneducated.
The correct use of the word ``comprise'' is
\begin{quote}
The standing committee comprises two women, three men, and a donkey.
\end{quote}
The formula is ``$A$ comprises $B$.''  What people often say, or write,
instead is
\begin{quote}
The committee is comprised of two women, three men, and a donkey.
\end{quote}
What {\it should} have been used in this last instance is ``composed,''
{\it not} ``comprised.''
			
\item{\bfit Contractions}\ \  Do {\it not\/} use contractions in 
formal\index{contractions} writing.  Thus the
words ``don't,'' ``can't,'' ``shouldn't,'' ``I'm,''
``you're,'' etc., are taboo.  Of course you should never
write ``ain't.'' You also should avoid abbreviations.
Particularly avoid using informal abbreviations
like ``cuz'' for ``because,'' ``tho'' for ``though,'' and
so forth.  You will probably never be tempted to
work ``bar-b-q'' into your next paper on para-differential
operators; but you might be tempted to use ``rite inverse.''
Please resist.

\null \quad \  Occasionally you will find it suitable to use contractions in
various kinds of {\it informal\/} writing.  It can be a way of drawing
in your audience, or of warming yourself up to your subject.  For
example, in the book [Kr2] I intentionally used an occasional
contraction in an effort to create a friendly air about the book.
By contrast, the present book is a book about writing, and I wish to
set a more formal example---so there are no contractions.

\item{\bfit Denote} \ \ Use the word ``denote'' carefully.  It has
a\index{denote@`denote'} special purpose in mathematics (and in 
logical positivist philosophy and modal logic) 
and we should take care to preserve it as a useful tool. 
Suppose that a certain mathematical symbol $A$ stands for, or
represents, the item or set of ideas $B$ (ideally, you should be able
to excise any occurrence of $A$ and replace it with $B$ and preserve
exactly the intended meaning).  Under these circumstances, and {\it
only under these circumstances\/}, do we say that ``$A$ denotes $B$.'' 
For example,

\begin{quote}
Let $X$ denote the set of all semisimple homonoids
with stable quonset hut.
\end{quote}

\null \qquad   \ There is a shade, but an important shade, of difference between
the statements
\begin{quote}
{\bf (1)}\ Let $f$ be a continuous function.
\end{quote}
and
\begin{quote}
{\bf (2)}\ Let $f$ denote a continuous function.
\end{quote}
The intended meaning of the first sentence 
here is ``let $f$ be {\it any old\/} continuous
function.''  Thus the first statement is both customary and correct.
The second is neither customary {\it nor\/} correct.  For we use ``denote''
when we want to say that a certain specific item stands for some
other specific item.  That is not what we are trying to say here.

\null \quad \  Lack of familiarity with English, or lack of familiarity with
the precise meaning of ``denote,'' sometimes leads to dreadful abuses of
the word. A common one is ``Denote $X$ the set of all left-handed
polyglots.'' I leave it to you to decide whether failing English or
failing intellect might be the correct provenance of such a sentence; the
lesson for you is not to use ``denote'' in such a fashion.

\null \quad \  The word ``connote,'' rarely used in mathematical writing, can be
(but should not be) confused with ``denote.''  The dictionary teaches
us that ``$A$ connotes $B$'' means that $A$ {\it suggests\/} $B$, but not
in a logically direct fashion.  For example,
\begin{quote}
To a young man, ``love'' connotes flowers, beautiful
music, and happiness.
\end{quote}
is an appropriate use of the word ``connote.''

\item{\bfit enervate} \ \ Often we are lazy, and we use a word
according to how it sounds, rather than according to what it actually
means.  This text offers ``enervate'' as an instance of this phenomenon.
What the word actually {\it means} is ``to lessen the vitality or
strength of.''  But, intuitively, we confuse ``enervate'' with
``energize'' and give it essentially the opposite meaning.  The
lesson is to be careful with words that are unfamiliar.

\item{\bfit He {\bf and} she} \ \ It used to be the custom that,
if one referred to an abstract person in one's writing, then
one used the pronoun ``he'' or ``him'' or ``his.''  Now this is
considered to be politically incorrect.  One must treat women the
same as men.  

\null \quad \ One way to address the problem is to replace
``he'' with ``he/she,``, ``him'' with ``him/her,`` and
``his'' his ``his/her.''  But this is a bit awkward.  Another
popular means is to replace ``he'' with ``she,'' replace
``him'' with ``her,'' and replace ``his'' with ``her.''  This
does not really seem to solve the problem; instead it replaces one
conundrum with another.  A third possibility, commonly taught
at colleges and universities, is to replace the gender-specific
pronouns with ``they,'' ``them,'' and ``their,''  This unfortunately 
results in some rather awkward constructions.  A fourth possibility is to use the
words dreamed up by Michael Spivak to replace the offending pronouns.
Spivak replaces ``he'' and ``she'' with ``e,'' he replaces ``him'' and ``her''
with ``em,`` and he replaces ``his'' and ``her'' with ``eir.''

\null \quad The really best way to solve this problem, though it requires some
extra time and effort, is to phrase your sentences so that they
omit pronouns altogether.  As an example, instead of saying
\begin{quote}
One might replace his notebook with a computer.
\end{quote}
one could instead say
\begin{quote}
The notebook can be replace with a computer.
\end{quote}

\item{\bfit Hyphen {\bf vs.}\ en dash} \ \  It is common
in mathematics, if two mathematicians have proved a theorem,
to call it something like ``the Riemann-Lebesgue lemma.''
Nowadays this is considered to be inappropriate.  The use
of the hyphen here may suggest that Riemann and Lebesgue have
more than a professional relationship.  So the politically correct
thing to do is to write ``the Riemann--Lebesgue lemma.''
What is the difference?  In the second example I used
the so-called {\it en dash} rather than the hyphen.  The en dash
is a dash that is about the width of the letter ``n`` in the current
font, and it is typically used to denote a range of numbers
(as in ``pages 324--386'').  So it is considered to carry less
emotional baggage and is therefore a better choice to denote
a mathematical collaboration.

\null \quad \ You may think this discussion ludicrous, but I can
tell you that if you do not conform to the prescription described
in the last paragraph then your copy editor will change all your
hyphens to en dashes.

\item{\bfit If\/ {\bf vs.}\ Whether}\ \ The words ``if'' and ``whether''
have\index{ifvswhether@`if' vs.\ `whether'} 
different meanings, and are suitable for different contexts.
Follow the example of master editor George Piranian:
\begin{quote}
Go to the window and see {\it whether\/} it is raining; {\it if\/} 
it is raining, then let Fido inside.
\end{quote}

\item{\bfit Infer {\bf and} Imply}\ \ The words ``infer'' and ``imply''
are\index{inferandimply@`infer' and `imply'} often confused in everyday usage. It should
not be difficult for a mathematician to keep these straight. A set of
assumptions can {\it imply\/} a conclusion. But one {\it infers\/} the
conclusion from the assumptions. It is that simple.

\item{\bfit Its {\bf and} It's}\ \  Use ``it's'' only 
to denote\index{itsandits@`its' and `it's'} the contraction for ``it is.''
Otherwise use ``its.''  For example ``Give the class its exam'' and
``A place for everything and everything in its place.''  Compare with
``It's a great day for singing the blues.''

\null \quad \   More generally, the apostrophe is never used to denote the possessive
of a pronoun:  what is correct is ``its,'' ``hers,'' ``his,'' ``theirs.''

\item{\bfit Lay {\bf and} Lie}\ \ ``Lay'' is a transitive 
verb\index{layandlie@`lay' and `lie'} and ``lie'' is intransitive.  This
means that ``lay'' is an action that you perform on some object,
while ``lie'' is not.  For instance, ``Lay down your weary head,''
``Now I will lay down the law,'' and ``I shall lay responsibility
for this transgression at your feet''; compare with
``I am tired and I shall lie down'' or ``Let sleeping dogs lie.''
Note, however, that the past tense of ``lie'' is ``lay.''  Therefore
you may say ``Yesterday I was so tired that I 
laid down my books and then I lay down.''

\item{\bfit Less {\bf and} Fewer}\ \ How many times have you been in 
the\index{lessandfewer@`less' and `fewer'} grocery store and gravitated
toward the line that is labeled {\it Ten Items or Less\/}?  Of course
what is intended here is {\it Ten Items or Fewer\/}, and I have
a special place in my heart for those few grocery stores that
get it right.  The word ``fewer'' is for comparing two numbers
while ``less'' is for comparing quantity.\footnote{Another way to
think about the matter is that ``fewer'' is used to compare 
discrete sets while ``less'' is used to compare continua.}  Mathematics deviates a bit
from these rules, because we certainly say ``3 is less than 5''
(of course the meaning here is that ``the number 3 represents
a quantity that is less than the quantity represented by
the number 5'').
Avoid saying ``3 is smaller than 5,'' because ``smaller'' is a
word about {\it size\/}:  perhaps the digit 3 is smaller than the
digit 5.  It also {\it could\/} be correct to say
``5 is smaller than 3'' if comparison of digit size is what is
intended:  {\small 5} versus {\Large 3}.

\item{\bfit Lists Separated with Commas {\bf (the Serial Comma)}}\ \  
When\index{serial comma} you are presenting a list, separated with commas, then you
should put a comma after every item in the list except the last.  For
example, say ``the good, the bad, and the ugly'' rather than ``the
good, the bad and the ugly.''  A moment's thought reveals that the
former conveys the intended meaning; the latter may not, for the
reader could infer that ``bad'' and ``ugly'' are somehow linked.

\item{\bfit Numbers}\ \ Some sources will tell you that (whole)
numbers\index{numbers!and numerals} 
less than 101 should be written out in words; larger numbers
should be expressed in numerals (other sources will put the cutoff at
twenty or some other arbitrary juncture).  A discursive discussion
appears in [SG].  Such considerations are, for a mathematician, next
to ludicrous.  The main thing, and this advice applies to spelling
and to many other {\it choices\/}, is to select a standard and to be
consistent.

\item{\bfit Obviously, Clearly, Trivially}\ \ These words 
have\index{obviously@`obviously'}\index{clearly@`clearly'}\index{trivially@`trivially'} 
become part of standard mathematical jargon.  This is too bad. In
the best of circumstances, when you use these phrases you are
endeavoring to push the reader around.  In the worst of circumstances
you are throwing up a smoke screen for something that you yourself
have not thought through.  It would be embarrassing to count the
number of major published mathematical errors that have been prefaced
with ``Obviously'' or ``Clearly'' or ``Trivially'' (no doubt the supreme
deity's way of reminding us that ``pride goeth before the fall''). 
The use of these words is one of the ways that we have of kidding
ourselves.\footnote{In a moment of exasperation, a friend
of mine said of her soon-to-be-ex-husband mathematician, ``You
look at anything and you either say that it is `very interesting'
or `trivial.' ''}

\null \quad \  As you proofread your manuscript relentlessly, and
endeavor to weed out superfluous words, pay particular attention to
the use, abuse, and overuse of these trite words.  They add nothing
to what you are saying, and are frequently a cover-up.

\item{\bfit Overused Words}\ \ Certain words in the English language
are\index{overused words} grossly overused.  Among these are ``very'' and ``most'' and
``nice'' and ``interesting.''  It is certainly very pleasant and most
insightful to express great appreciation for a very nice and
supremely interesting theorem; but I encourage you not to do so---at
least not with these banal words.  If such language
represents how you wish to express
yourself, then perhaps you have nothing to say.  
Instead think
carefully about what you really mean to say, and then say it.

\null \qquad   \ In fact the language is littered with overused words
that come into and out of fashion.  The words ``awesome,''
``totally,'' ``dude,'' and ``righteous'' are current examples.  The
phrase ``today I'm not 100\%,'' foisted upon us by some semiliterate
sports announcer, is currently the bane of our collective existence. 
Each field of mathematics has its own set of stock phrases and
tiresome clich\'{e}s.  Endeavor not to propagate them.  

\null \quad \  A good general principle is to put every word in every
sentence under the microscope:  What does it add to the sentence? 
Will the sentence lose its meaning if the word is omitted?  Can the
thought be expressed with fewer words?  Strunk and White [SW] have a
splendid discussion of the concept of weighing each word.

\item{\bfit Plural Forms of Foreign Nouns}\ \ We all grind our teeth
when\index{plural forms of foreign nouns}
we hear our freshmen say ``And this point is the {\it maxima\/} of the
function.''  To no avail we explain that ``maxima'' is plural, and
``maximum'' is singular.  Yet we make a similar error when we do not
differentiate ``data'' (plural) from ``datum'' (singular) and
``criteria'' (plural) from ``criterion'' (singular).  As usual, exercise
special care when dealing with foreign words.

\item{\bfit The Possessive}\ \  When you express the possessive of a
singular\index{possessives} noun, always use \ 's.  Thus you should say 
``Pythagoras's society,'' ``the dog's day,'' ``Stokes's theorem,'' 
``Bliss's book,'' ``baby's bliss,'' and
``van der Corput's lemma.''  The terminal ``s'' is omitted when you
are denoting the possessive of a plural noun:  ``the boys' trunk,''
``the dogs' food,'' ``the students' confusion.''
\vspace*{-.15in}

\null \quad \  ``Collective nouns'' are treated in a special manner. 
For instance, we write ``the people's choice'' and ``the children's
folly'':  even though the nouns are plural, we denote the possessive
{\it with\/} a terminal ``s.''

\null \quad \ Just because we frequently see such misuse in advertising
and other informal writing, we sometimes get sucked into using extraneous
apostrophes.  As an example, one often sees expressions like ``This
sentence contains a lot of TLA's.''  Here a TLA is a ``three-letter acronym.''
What is wrong with this sentence?  The last ``word'' in the sentence is
supposed to be a plural---{\it not} a possessive.  So the apostrophe is
out of place.  It should be ``TLAs,'' not ``TLA's.''  Likewise, do not
write, ``I surely miss the 1960's.''  It should be ``1960s,'', {\it not} ``1960's.''

\item{\bfit Precision and Custom}\ \ At times, the goal of
precision\index{precision and custom}
in writing flies in the face of custom.  Antoni Zygmund
once observed that the World Series of American baseball might more
properly be called the ``World Sequence.''  I am inclined to agree (in
no small part out of fealty to my mathematical grandfather), but I
must be over-ruled by custom:  if you use the phrase ``World Sequence''
then nobody will know what you are talking about.  Bear this thought
in mind when you are tempted to invent new terminology or new
notation (see also the remarks in Section 2.4 on terminology and
notation).

\item{\bfit Subject and Verb, Agreement of}\ \ Make sure that subject
and verb match in your sentences.  A mismatch not only grates on the
sensitive\index{subject!and verb, agreement of}
ear, but can seriously distort meaning.  Consider the
example ``The set of all morphisms are compact.'' This syntax is
incorrect.  The {\it subject\/} (that is, the person or thing
performing the action) in this sentence is {\it set\/}.  We should
conjugate  the verb ``to be'' so that it
agrees with this subject.   As a result,
the grammatically correct statement is ``The set of all
morphisms is compact.''  (Note, in passing, that the original form of
the sentence might have misled the reader into thinking that the
writer was---rather clumsily---discussing a collection of compact
morphisms.)
\vspace*{-.25in}

\null \quad \ Of course the test is easy:  omit the prepositional
phrase ``of all morphisms'' and analyze the root sentence.  Clearly
``The set is compact'' is correct while ``The set are compact'' is
not.  You will find the device of focusing on the root statement, or
breaking into pieces (see our analysis of Subject and Object below),
to be a valuable tool in analyzing many grammatical questions.  
\vspace*{-.14in}

\null \quad \ As a parting exercise, consider the phrases ``the sequence
$z_n$ {\it converges\/} to $p$'' while ``the numbers $z_n$ {\it
converge\/} to $p$.'' Think carefully about why both statements are
correct.

\item{\bfit This {\bf and} That}\ \ 
We often hear\index{thisandthat@`this' and `that'}, 
especially in conversation, phrases
like ``Because of this, we decided that.'' 
If we exercise the full force of logic then
we must ask `` `Because' of {\it what\/}?''
and `` `we decided' {\it what\/}?''  And this niggling query
raises an entire body of common errors that I would like to point
out.  This corpus is not composed so much of errors in English
usage, but rather errors in logic and precision.  Consider the
following examples:

\begin{quote}
Shakespeare was an important writer.  This tells us a lot
about English literature.  \qquad \bad
\end{quote}

\begin{quote}
A triangle is a three-sided polygon.  This means that \dots  \qquad \bad
\end{quote}

\begin{quote}
The day was bright and beautiful.  Because of this, Mary smiled.  \qquad \bad
\end{quote}

\null \quad \  In each of these sample sentences, my objection is 
```this' what?''
(Notice that I did {\it not\/} say ``In each of these, my objection 
is \dots.''
I was careful to say {\it precisely\/} what I meant.)   The following
passages convey the same spirit as the preceding three, but they
actually {\it say\/} something:

\begin{quote}
Shakespeare was an important writer.  The forms of
his plays and poems as well as his use of language have had a strong
influence on English literature.
\end{quote}

\begin{quote}
A triangle is a three-sided polygon.  The trio
of sides satisfies the important {\it triangle inequality\/}.
\end{quote}

\begin{quote}
The day was bright and beautiful.  Observing the
weather caused Mary to smile.
\end{quote}
\vspace*{-.1in}

\null \quad \   Here is a delightful example that was contributed 
by G. B. Folland:
\begin{quote}
Saddam Hussein was determined to resist attempts to force
Iraqi troops out of Kuwait, although George Bush made it clear
that he did not want to be seen as a wimp.  This caused the Gulf
War.
\end{quote}
If you were to ask someone to which clause ``this'' refers,
then the answer you received would probably depend on that person's
politics.
\vspace*{-.1in}

\null \quad \   The message here is fundamental:  as a default, do not use
``this'' or ``that'' or ``these'' or ``those'' without a
clear point of reference.  When the occurrence  of
``this'' or ``that'' is fairly close to the referent, then
the intended meaning is often clear from context.  When
instead the distance is greater (as in Folland's example),
then confusion can result.
\vspace*{-.1in}

\null \quad \   Repetition is a good thing, so repeat your nouns rather than
refer to them with a potentially vague pronoun or other word.  There
{\it will\/} be cases where the casual use of ``this'' or ``that'' is
both natural and appropriate, but such instances will be exceptions.
\vspace*{-.1in}

\null \quad \   Copy editor Rosalie Stemer says that a hallmark of good writing is
that it answers more questions than it raises.  Applying this
philosophy will lead naturally to many of the points raised in this
book, including the present one.
\vspace*{-.1in}

\item{\bfit Where}\ \ One of the most common types of run-on
sentence\index{where@`where'} in mathematics is a statement with a dangling
concluding phrase such as ``where $A$ is defined to be \dots.''  
An example is

\begin{quote}
Every convex polynomial function is of even degree, where we
define a function to be convex if \dots  \qquad \bad
\end{quote}

\noindent We see this abuse so often that we are rather accustomed to it.
This is also an easy crutch for the writer:  he/she did not bother to 
plant the definition before this statement, so he/she just
tacked the definition onto the end.

\null \quad \    That is sloppy writing and there is no excuse for it:  
before you
use a term, define it.  You need not use a formal, displayed definition.
But you must put matters in logical order.  The example I have
given is quite trivial; but in serious mathematical writing it is
taxing on the reader to have to pick up definitions on the fly.
Especially if you are writing with a computer, it is very easy
for you to scroll up and put the needed definition where it
belongs.

\item{\bfit Who {\bf and} Whom; {\bf Subject and Object}}\ \  
Be conscious\index{whoandwhom@`who' and `whom'}\index{subject!and object}
of the difference between ``who'' and ``whom.''
The word ``whom'' is an object; used properly, it denotes a person
that is {\it acted upon\/}.  An example of the common misuse 
of the word ``whom'' is 
\begin{quote}
The pastor, whom expected 
a large donation, smiled warmly. \qquad \bad
\end{quote}
Here the issue
is what is the correct subject to put in front of the
verb ``expected.''  The word ``whom''
cannot act as a subject.  The correct word is ``who'':  ``The pastor,
who expected a large donation, smiled warmly.''   In the same
vein, it is correct to say ``To whom am I speaking?''  and
``Is he the man to whom the Nobel Prize was awarded?''

\null \quad \   Also do not confuse ``I'' and ``me.''  The latter is
an object, the former not.  For example, ``The teacher was addressing
Bobby and I'' is plainly wrong, since here ``I'' is used 
incorrectly as the object of the verb ``addressing.''  President
Clinton's famous mis-statement ``Give Al Gore and I a chance to
bring America back'' is a dreadful error; nobody would say
``Give I a chance \dots.''  That sort of sentence analysis---breaking
a sentence down to its component parts---is the 
method you should use to detect the error.
The sentence
\begin{quote}
Him and me proved the isotopy isomorphism theorem in
1967.  \qquad \bad
\end{quote}
is an abomination.  Unfortunately, even smart people make 
mistakes like this.  Anyone can see that ``Him proved the isotopy
isomorphism \dots'' and ``Me proved the isotopy isomorphism \dots''
are incorrect.  But, somehow, the ganglia are more prone to misfire
when we put the two sentences together.  Conclusion:  test the
correctness of a sentence with compound subject (or any compound
element) by breaking it into its component sentences.
\end{itemize}

%% Section 1.12

\markright{1.12.  MORE RULES OF GRAMMAR}
\section{More Rules of Grammar, Syntax, and Usage}
\markright{1.12.  MORE RULES OF GRAMMAR}

\null \quad \ Here I include additional rules of grammar and syntax that
are dear to my heart.  They come up frequently in general
writing, less so in specifically mathematical writing.
They should prove useful in your expository work, and
sometimes in your research work as well.

\begin{itemize}
\item{\bfit Adjectives vs.~Adverbs}\ \ An adjective is designed to
describe, or\index{adjectives}\index{adverbs} to modify, a noun.  An adverb is designed to describe, or
to modify, a verb.  Correct is to say 

\null \ \ \ \  This is a good book. 

\noindent and
\vspace*{-.13in}

\null \ \ \ \  This is an expensive car. 

\noindent and 
\vspace*{-.13in}

\null \ \ \ \  The quick, brown fox jumped over the stupid, lazy dog.
\vspace*{.1in}

\noindent  because ``good,'' ``expensive,'' ``quick,''
``brown,'' ``stupid,'' and ``lazy'' are adjectives.  They modify the
nouns ``book,'' ``car,'' ``fox'' (twice), and ``dog'' (twice),
respectively.   You may also say 
\smallskip 

\null \ \ \ \  She shouts loudly. 

\noindent and
\vspace*{-.13in}

\null \ \ \ \  He sings beautifully.

\noindent and 
\vspace*{-.13in}

\null \ \ \ \  She strove sporadically to master her homework thoroughly.
\smallskip 

\noindent  because ``loudly,'' ``beautifully,''
``sporadically,'' and ``thoroughly'' are adverbs.  They modify the
verbs ``shouts,'' ``sings,'' ``strove,'' and ``master.''  
Learn to distinguish between adjectives and adverbs, and learn to 
use both correctly.  

\null \quad \ After Paul Halmos had seen an early draft of the first edition 
of this book, he sent me the message ``You write good.''  One can
guess effortlessly that he was joking mischievously about this silly,
little book.

\item{\bfit Alternate {\bf vs.}\ Alternative}\ \ The words ``alternate''
and ``alternative'' (used as adjectives) have different meanings, though
they\index{alternatevsalternative@`alternate' vs.\ `alternative'} are often, and erroneously,
used interchangeably. The word ``alternate'' (most commonly used in the
form ``alternately'') refers to some pair of events that occur repeatedly
in successive turns; the word ``alternative'' refers to a choice between
two mutually exclusive possibilities. For example:

\begin{quote}
Pierre alternately dated Mimi and Fifi.  He had
considered monogamy, but had instead chosen the alternative lifestyle
of a concupiscent lothario.
\end{quote}

\item{{\bf The Verb} \bfit To Be}\ \ ``The verb `to be' can never
take\index{tobe@`to be'} 
an object.''  Probably you have been hearing this statement
all your life.  What does it mean?

\null \quad \   When you formulate the sentence 
\begin{quote}
I hit the ball.
\end{quote}
then ``I'' is the subject 
(of the verb ``hit'') and ``ball'' is the object (of the verb ``hit'').
But when you formulate the sentence 
\begin{quote}
I am the walrus.
\end{quote}
then ``I'' is the subject (of the verb ``to be,'' conjugated
as ``am'') but ``walrus'' is the {\it predicate nominative\/}
(also sometimes called the {\it predicate noun\/} or {\it subjective
complement\/}).
The word ``walrus'' is {\it not\/} the object of an action.
It has a different grammatical role in this sample sentence.
\vspace*{-.1in}

\null \quad \   By the same token, it is incorrect to answer the query (over the
telephone) ``Is this Napoleon Bonaparte?'' with the answer ``This is
me.'' The word ``me'' is supposed to be used as the object of an
action.  The verb ``to be,'' however, does not take an object.  Thus
the correct rejoinder is ``This is I''  or ``This is he.''
\vspace*{-.1in}

\null \quad \   To make a long story short, your writings should not include the
statement ``The person who proved Fermat's Last Theorem is me.'' 
Grammatically correct is ``The person who proved Fermat's Last
Theorem is I'' or ``It is I who proved Fermat's Last Theorem'' or ``I
am the one who proved Fermat's Last Theorem.''  You should not, however, 
pen any of these statements unless you are Andrew Wiles.

\item{\bfit Compare {\bf and} Contrast} \ \  The words ``compare'' and
``contrast'' have\index{compareandcontrast@`compare' and `contrast'} 
different meanings.  One compares two items in
order to bring out their similarities; one contrasts two items in
order to emphasize their differences.  For instance, 
we can compare groups and semigroups because they are both
associative.  We can contrast them because one contains
all inverse elements and the other does not.

\item{\bfit Different from {\bf and} Different than}  The phrase 
``different from'' is correct,
while ``different than'' is not.  Examples are ``His view of
grammar\index{differentfromvs@`different from' vs.\ `different than'} 
is different from mine'' and ``His syntax is different
from what I expected.''  The classical rationale here is that
the word ``different'' demands a preposition and a noun.  Modern
treatments (see [Fra, p. 266]) suggest that this classical dictum
is too restrictive and that ``different than'' (without the noun) is
more comfortable.  You will have to decide which usage you prefer,
but do be consistent.

\item{\bfit Due to, because of, {\bf and} through}\ \  Mathematicians 
commonly\index{duetovsbecause@`due to' vs.\  `because' vs.\ `through'} 
use the phrase ``due to,'' and we often use it incorrectly. 
We sometimes say ``due to the fact that'' when instead ``because''
will serve nicely.  The phrase ``due to'' tempts us to wordiness
that is best resisted.

\null \quad \   A good rule of thumb (thanks to G.\ Piranian) is to use
``due to'' only to introduce an {\it adjectival clause\/}---never an
adverbial clause.  In fact the grammatical issues at play here---including
the use of ``through''---are
rather complex, and not suitable for discussion in this book.  See
[MW] for a detailed treatment.

\item{\bfit Farther {\bf and} further}\ \ It is common to interchange 
the\index{farthervsfurther@`farther' vs.\ `further'}
words ``farther'' and ``further,''
but there is a loss of precision when you do so.  The word ``farther''
denotes distance, while ``further'' suggests time or quantity.
For example, one might say ``I wish to study {\it further\/} the question
of whether Lou Gehrig could throw the baseball {\it farther\/} than
Ty Cobb.''

\item{\bfit Good Taste and good sense} \ \ Suit your prose to the
occasion.\index{good!taste}\index{good!sense}  
The writer of a Harlequin romance novel might write
\begin{quote}
Clutched in the gnarled digits of the syphilitic Zoroastrian homunculus
was a dazzling Faberg\'{e} egg.
\end{quote}
while Raymond Chandler would have written something more like
\begin{quote}
The dwarf held a gewgaw.
\end{quote}
In mathematics, simpler is usually better.  Flamboyant writing 
is out of place.

\item{\bfit Hopefully {\bf and} I hope}\ \  With due homage to 
Edwin R. Newman [New], I note that it
is incorrect to use ``hopefully'' (at the beginning of a sentence)
when\index{hopefully@`hopefully'} you mean to say ``It is hoped that'' or (more sloppily)
``I hope.''  The word ``hopefully'' is an
adverb.  It is intended to modify a verb.  For example, consider the sentence

\begin{quote}
She wanted so badly to marry him, and she looked at him hopefully
while she waited for a proposal.
\end{quote}

\noindent Note that the word ``hopefully'' modifies ``looked.''  It is
incorrect to say
\begin{quote}
Hopefully the weather will be better today.  \qquad \bad
\end{quote}
\noindent when what you mean to say is
\begin{quote}
I hope that the weather is better today.
\end{quote}
By the same token, do not say ``This situation looks hopeful.''
People can be hopeful, objects or things never.

\null \quad \ Monty Python tells us that ``Mitzi was out in the garden,
hopefully kissing frogs.''  If you are comfortable with the common
misuse of ``hopefully'' then you will probably misunderstand this
sentence.

\null \quad \ The reference [KnLR, p. 57] offers a detailed analysis of the
history of the word ``hopefully,'' and another, more liberal,
point of view about its use.  See also [MW].

\item{\bfit Infinitives, Splitting of}\ \  As a general rule, do not
split infinitives.  For example, do not say ``He was determined to
immensely enjoy his\index{infinitives, splitting of}
food, so he smothered it in ketchup.''  The correct
version (though one may argue with the sentiment) is ``He was
determined to enjoy his food immensely, so he smothered it in ketchup.'' 
Here the infinitive is ``to enjoy'' and the two words should not be
split up.  Curiously, the reason for this rule is an atavism:  some
of the languages that contributed to the formation of modern English,
such as Latin and French, combine these two words into one.  Our
rule not to split the infinitive carries on that tradition.

\null \quad \   There are a number of opinions on this matter.  The ``modern'' point
of view is that it is acceptable to split an infinitive when it sounds
right; otherwise it is not.  For example, sometimes a mathematical
sentence will resist the suggested rule.  G.~B.~Folland supplies the
example ``Hence we are forced to severely restrict the allowable
range of values of the variable $x$.'' Strictly speaking, the word
``severely'' splits the infinitive ``to restrict.''  But where else
could you put ``severely'' while maintaining the precise meaning of the
sentence? 

\null \quad \   On a more personal level, the sentence 
\begin{quote}
To really love someone requires a lot of effort.
\end{quote}
evinces a particular sentiment while
\begin{quote}
To love someone really requires a lot of effort.
\end{quote}
conveys a different meaning.  

\null \quad \ Arguably, it would
be better to formulate the sentence differently (how about
``To love a person with passion and intensity requires
a lot of effort.''?).   But if one were wedded to the ``really''
construction then one would have a problem.  Use your
ear, and use sound and sense, to convey your message clearly
and forcefully.

\item{\bfit In terms of}\ \  Sentences of the form
\begin{quote}
Who\index{intermsof@`in terms of'} is he, in terms of surname?  \qquad \bad
\end{quote}
and
\begin{quote}
How is she doing, in terms of her math classes? \qquad \bad
\end{quote}
are simply dreadful.  Usually the phrase ``in terms of'' is gratuitous,
and can be omitted entirely.  Consider instead
\begin{quote}
What is his surname?
\end{quote}
and
\begin{quote}
How is she doing in her math classes?
\end{quote}

\item{\bfit Need Only; Suffices to} \ \ 
In written mathematics, we often find it 
convenient to say ``We need only show that \dots'' or ``It
suffices to show that \dots.''  These are lovely turns of phrase.
Strive not to overuse\index{needonly@`need only'}\index{sufficesto@`suffices to'}
them, or to misuse them.  Too often we see
instead ``We only need to show that \dots'' or ``Suffices it to show
that \dots.''  With these misuses, the message still comes across---but in a
more halting and less compelling manner.

\item{\bfit Parallel Structure} \ \ 
The\index{parallel structure}
principle of parallel structure is that proximate clauses which have
similar content and purpose are (often) more effective if
they have similar form.
The use of parallel structure is an advanced writing
skill:  good writing can be made better, more forceful,
and more memorable with the use of parallel structure.  Consider
the dictum
\begin{quote}
Candy is dandy but liquor is quicker.
\end{quote}
Whether you approve of the sentiment or not, the thought is
memorably expressed---using a quintessential example of parallel
structure.  As an exercise, try expressing the thought with
more desultory prose, and see for yourself what is lost in
the process.

\null \quad \   The first inspirational quotations (from Sir Francis Bacon) in
Chapters 3 and 5 provide less frivolous examples of parallel
structure.

\item{\bfit Participial phrases}\ \ Participial phrases are a frequent 
cause for discomfort.
For\index{participial phrases} example, 
\begin{quote}
Shining like the sun, the man gazed happily upon the 
heap of gold coins.  \qquad \bad
\end{quote}
The participial phrase ``shining like the sun'' modifies ``man,''
whereas it was clearly intended to modify ``the heap of gold
coins.''  Better would be

\begin{quote}
The man gazed happily upon the heap of gold coins, which
shone like the sun.
\end{quote}

\null \quad \   Harold Boas contributes the following useful maxim:  ``When 
dangling, don't use participles.''

\item{\bfit Prepositions, Ending a Sentence with}\ \ As a general rule,
do not\index{prepositions, ending a sentence with}
end a sentence with a preposition.   Do not say ``Where do we
stop playing at?'' Instead say ``At what point do we stop
playing?''  Better still is ``When do we stop playing?''  Do not say
``What book are you speaking of?'' Instead say ``Of which book do
you speak?'' or ``Which book is that?''

\null \quad \   Often, when you are tempted to end a sentence with a preposition,
what is in fact occurring is that the errant preposition is a spare 
word---not needed at all.  The preceding examples,  and  the
suggested alternatives, illustrate the point.

\null \quad \   An old joke has a yokel trying to find his way across the
Harvard campus.  A Brahmin student corrects him sternly
for posing the question ``Excuse me.  Where's the library at?''
After the Harvardian explains at length that one does not
end a sentence with a preposition, the yokel tries
again:  ``Excuse me.  Where's the library at---{\it jerk\/}?''  This
is perhaps a bizarre example of sound working with sense.

Harold Boas cautions:  ``Watch out for prepositions that sentences end with.''

\item{\bfit Quotations}\ \  We do not often include quotations in
mathematics\index{quotations} papers.  If you decide to include a quotation, then be
aware of the following technicality.  Logically, it makes sense to
write a sentence of the following sort:

\begin{quote}
As Methuselah used to say, ``When the going gets tough, 
the tough get going''.
\end{quote}

\noindent What is logical here is that the quotation itself is a
proper subset of the entire sentence; therefore it stands to reason
that the terminal double quotation mark should be {\it before\/} the
period that terminates the sentence.  Unfortunately, 
logic fails us here.
Admittedly typesetters are still
debating this point, but the current custom in the United States is to put the closing
double quotation mark {\it after\/} the period.  Open any novel and
see for yourself.  Thus the sentence {\it should\/} be written

\begin{quote}
As Methuselah used to say, ``When the going gets tough, 
the tough get going.''
\end{quote}

\null \quad \   The fact is that the complete rule is even a bit more complicated
than has already been indicated.  By the rules of {\it American\/}
usage, commas and periods should be placed inside quotation marks,
and colons and semicolons outside quotation marks (see [SG, p.~222]
and [Dup, p.~192]).  Placing exclamation points and question marks
inside or outside of quotation marks depends on context.  British
usage is even more ambiguous.  This is all a bit like the infield fly
rule in baseball.  But do be consistent, and be prepared to
arm-wrestle with your publisher or with your copy editor if you
have strong opinions in the matter.

\null \quad \ If your quotation is $n$ paragraphs in length, then there is an
opening double quotation mark on every paragraph.  There is no
closing double quote on paragraphs $1$ through $(n-1)$; but there certainly
{\it is\/} a closing double quote on paragraph $n$.  Again, check any
published novel to see that this is the case.

\item{\bfit Redundancy}\ \  Logical redundancy, used with discretion,
can be\index{redundancy} a powerful teaching device.  By contrast,
avoid (local) verbal redundancy. 
The phrases ``old adage,'' ``funeral obsequies,'' ``refer back,''
``advance planning,'' ``strangled to death,'' 
``invited guest,'' ``body of the late,''
and ``past history'' display an ignorant and superfluous use of
adjectives.  Avoid constructions of this sort.

\item{\bfit Shall {\bf and} Will}\ \ In common speech, the words
``shall'' and ``will'' are often used interchangeably, or according
to what appeals to the\index{shallandwill@``shall'' and ``will''} 
speaker.  In formal writing there is a simple
distinction:  when expressing belief regarding a future action or
state, ``shall'' is used for the first person (``I'' or ``we'') and
``will'' is used for the second person (``you'') or third person (``he,''
``she,'' ``it,'' or ``they'').  When the first person is expressing
determination, then it is appropriate to use ``will.''  These rules,
taken from [SW], are illustrated whimsically in that source by

\begin{quote}
{\bf Bather in Distress:}  ``I shall drown and no one will save me.''
\end{quote}

\noindent but

\begin{quote}
{\bf Suicide:}  ``I will drown and no one shall save me.''
\end{quote}

\item{\bfit That {\bf and} Which}\ \ The word ``that'' is used to denote
{\it restriction\/} while\index{thatandwhich@``that'' and ``which''} 
the word ``which'' denotes {\it amplification\/}. 
For example, ``I am speaking of the vase that sits on the table'' and
``The book that is by Gibbons is in the study.'' Compare with ``The
vase, which is red, sits on the table'' and ``The book, which is by
Gibbons, is fascinating.''

\null \quad \   Here is a mathematical example:  ``A holomorphic function that vanishes
on $S$ must be identically zero.'' Compare with ``A holomorphic function
which vanishes on $S$ must be identically zero.'' Which is correct? 
Think about the logic.  What we are saying is that a holomorphic
function $f$ such that $f(z) = 0$ for $z \in S$ must be identically
zero.  (For the mathematics, note that, in one complex variable, a set
$S$ with an interior accumulation point will suffice for the truth of
the statement.)  Phrased in this way, the statement is
restrictive:  a holomorphic function with a certain additional
property must be zero.  Thus the correct choice is ``that'' rather
than ``which.''

\null \quad \   Modern grammarians approve of the use of  ``which''  for ``that''
in suitable contexts.  Consult a grammar book, such as [SG],
for the details.

\null \quad \   I have already noted that it is sometimes 
useful to let ``sound and
sense'' overrule the strict code of grammar.  In particular, there
are times when ``which'' sounds more weighty, or more formal, than
``that.'' Thus some writers will make the technically incorrect
choice, just to achieve a certain effect.

\end{itemize}

As already noted, the rules of grammar and syntax are not absolute. 
English usage is constantly evolving.  While some current aspects of
usage are fads and nothing more, others become common and are finally
adopted by the best writers and speakers.  Those tend to stay with
us.  But there is a more subtle point.  Sometimes a sentence formed
according\index{rules of grammar!flexibility of} to 
the strict rules of usage {\it sounds awkward\/}.  A
classic\index{rules of grammar!strictness of}
example (usually attributed to Winston Churchill) is

\begin{quote}
That is the sort of behavior up with which I will not put.
\end{quote}

\noindent Notice that the speaker is going into verbal contortions to
avoid ending the sentence with a preposition.  The result
is a sentence that is so ludicrous that it defeats the main
purpose of a sentence---to {\it communicate\/}.  Better is
to say 

\begin{quote}
That is the sort of behavior that I will not put up with.
\end{quote}

\noindent While technically incorrect---because the preposition
is at the end of the sentence---this statement nevertheless will
not grate on the ears of the listener, will convey the sentiment
clearly, and will get the job done.  Of course it would be even
better to say

\begin{quote}
I will not tolerate that sort of behavior.
\end{quote}

\noindent This sentence conveys exactly the same meaning as the first
two.  But it has the advantage that it is direct and forceful.  
In most contexts, the last sentence
would be preferable to the first two.  This is again a matter of 
sound working with sense.  And here is a point that I will make
several times in this book:  often it is a good idea {\it not\/}
to wrestle with a sentence that is not working; instead, reformulate
it.  That is what we did with the last example.

As an exercise, find a better way to express the following sentence
(which ends with five prepositions, and which I learned
from Paul Halmos\index{Halmos, Paul} by way of [KnLR]):
\begin{quote}
What did you want to bring that book I didn't want to be read
to out of up for?
\end{quote}

Do not use acronyms, abbreviations, or jargon unless you are dead
certain\index{acronym}\index{abbreviations}\index{jargon} 
that your audience knows these shortcuts.  Speaking of an
ICBM, the NAFTA treaty, ARVN, and MIRV is fine for those well read in
the current events of the past twenty-five years---and who have an excellent
memory to boot.  But most of us need to be reminded of the meanings of
these acronyms.  The best custom is to define the acronym
parenthetically the first time it is used in a piece of writing.  For
example, 
\begin{quote}
The SALT (Strategic Arms Limitation Talks) were progressing
poorly, so we broke for lunch.  A few hours later, we resumed
our efforts with SALT.  
\end{quote}

I have served on many AMS (American Mathematical Society)
committees, and am somewhat horrified at the extent to which I have
become inured to certain acronyms.  How many of these do you know:  CPUB,
COPROF, JSTOR, LRPC, ECBT, COPE?  I am conversant with them all, and
none has done me a bit of good.  In practice, you may not even safely
assume that your reader knows what the AMS is---what if he/she is
Turkish?

I was once at a meeting to discuss the writing of a new grant
proposal---to apply for renewal of funds from a generous source which,
we hoped, would be inclined to give again. One of the PIs (``PI''
denotes ``Principal Investigator'') said, in all seriousness, ``I think
that we are going to need more blue sky in this proposal if we want
to generate more bottom line.''  Of course his meaning was ``We must
endeavor to paint an enlarged picture of long-term goals and
anticipated achievements if we want to increase the size of this
grant.''  The first mode of expression might be appropriate among
venture capitalists, who are inured to such language.  It is probably
inappropriate among academics.

%%%%%%%%%%%%%%%%%%%%%%%%%%%%%%%%%

%% Chapter 2
\chapter{Topics Specific to the Writing of Mathematics}

\begin{quote}
\footnotesize \sl What I really want, doctor, is this.  On the day when
the manuscript reaches the publisher, I want him to stand up---after
he's read it through, of course---and say to his staff: ``Gentlemen,
hats off!'' 
\smallskip \hfill \break
\null \hfill \rm  Albert Camus \break
\null \hfill {\it The Plague} \break
\end{quote}
\vspace*{-.13in}
\begin{quote}
\footnotesize \sl You don't write because you want to say something;
you write because you've got something to say.
\smallskip \hfill \break
\null \hfill \rm F.~Scott Fitzgerald \break
\end{quote}
\vspace*{-.13in}
\begin{quote}
\footnotesize \sl So I'm, like, ``We need to get some food.''
And he's, like, ``I don' wanna go th' store.  How `bout some
`za?''  And I'm, like, ``Well, we gotta eat, dude.  I could get 
like totally into a pizza.''  And he's, like, ``No biggie.''  
And I'm, like, ``This guy is grody to the max.  Gag me with a spoon.'' 
\smallskip \hfill \break
\null \hfill \rm A Valley Girl \break
\end{quote}
\vspace*{-.13in}
\begin{quote}
\footnotesize \sl 
We have read your manuscript with boundless delight. \hfill \break
If we were to publish your paper, \hfill \break
it would be impossible for us to publish any work of lower standard. \hfill \break
And as it is unthinkable that in the next thousand years \hfill \break
we shall see its equal, we are, to our regret, \hfill \break
compelled to return your divine composition and to beg \hfill \break
you a thousand times to overlook our short sight and timidity. 
\smallskip \hfill \break
\null \hfill \rm  Memo from a Chinese Economics Journal \break
\null \hfill \rm  From {\it Rotten Rejections\/} (1990)
\end{quote}
\vspace*{.05in}

\begin{quote}
\footnotesize \sl 
Having imagination, it takes you an hour to write a paragraph that,
if you were unimaginative, would take  you only a minute.  Or  you
might not write the paragraph at all.
\smallskip \hfill \break
\null \hfill \rm Franklin P. Adams \break
\null \hfill \rm {\it Half a Loaf\/} (1927) \break
\end{quote}

%% Section 2.1

\section{How to Organize a Paper}
\markboth{CHAPTER 2.  TOPICS SPECIFIC TO MATHEMATICS}{2.1.  HOW TO ORGANIZE A PAPER}

To begin, a mathematics paper has certain technical components.
It requires a title, and that title should convey some information
to the reader.  If it does not, then the reader is likely to
move\index{math paper, components of}
on to more stimulating reading matter, without looking
any further at your work.  Thus a title like {\sl On a theorem
of Hoofnagel} says almost nothing.  The title
{\sl On differentiation of the integral} is only slightly better;
but\index{title, importance of}
at least now the reader knows that the paper is about analysis, and
he/she has a rough idea what sort of analysis.  The title
{\sl Quadratic convergence of Lax/Wendroff schemes with optimal 
estimates on the error term} is ideal.  In a nutshell, this
title tells the reader exactly what the paper is about and, further,
what point it makes.

Of course an equally important component of your paper is
the identification of the author or authors.  At the beginning of
your career, pick a name for yourself and stick to it.  And I do not
mean a name like ``Stud'' or ``Juicymouth.'' I might have called myself
Steven George Krantz or S.\ G.\ Krantz or S.\ Georgie Krantz or any
number\index{author!name} of other variants.  I chose Steven G.\ Krantz, just as it
appears on the title page of this book.  When an abstracting,
indexing, or reviewing service endeavors to include your works, you
want it to be a zero-one game:  it should retrieve all your works or none
of them.  You do not want any to be left out, and you should
leave no doubt as to your identity.

Here is a quick run-down of other technical components that
belong in most papers:  
\begin{enumerate}
\item[{\bf (1)}]  affiliations of authors,
\item[{\bf (2)}]  postal\index{affiliation!of author}
addresses\index{postal address of author} 
and {\it e}-mail addresses of authors,
\item[{\bf (3)}] date,\index{date on paper} 
\item[{\bf (4)}] abstract, 
\item[{\bf (5)}] key words,
\item[{\bf (6)}] AMS subject\index{abstract!of paper} 
classification numbers, 
\item[{\bf (7)}]\index{key words} 
thanks\index{thanks to granting agencies} to 
granting agencies and\index{subject!classification numbers}
others.  
\end{enumerate}

Items {\bf (1)}, {\bf(2)}, 
and {\bf (3)} require no discussion; topic {\bf (4)} is discussed in Section 2.5.
Let us say a few words about {\bf (5)}--{\bf (7)}.  

The key words are provided so that {\it Math.\ Reviews\/} and other
archiving services can place your paper correctly into a database. 
Endeavor to choose words that reveal what your paper is about;
that is, you want words that will definitely lead a potential reader
to your paper.  Thus ``new,'' ``interesting,'' and ``optimal'' are
not good choices for key words.  Instead, ``pseudoconvex,'' ``Cauchy
problem,'' and ``exotic cohomologies'' {\it are\/} good choices.  

Similar comments apply to the AMS subject classification numbers.
The American Mathematical Society has divided all of mathematics into
97 primary classification areas (rather like {\it phyla} in the
classification of animals) and these in turn into subareas.
Assigning\index{MR subject classification numbers}
the\index{Mathematical Reviews} correct 
classification numbers to your paper is a
reliable way to put your paper before the proper audience.  It also
will help to ensure that your paper is classified correctly.  The
American Mathematical Society publishes an elegant little book,
which can be found in most mathematics libraries, that lays out the
AMS subject classification scheme.  This information is also readily
available on the Web.  The key words and classification numbers usually
appear in footnotes on the first page of your paper; some journals
instead specify that title, abstract, key words, and classification
numbers appear on a separate ``Title Page.''

Finally item {\bf (7)}:  often it will be appropriate to thank other
mathematicians for helpful conversations or specific hints.
Sometimes you will thank someone for reading an earlier draft of the
paper, or for catching errors.  If you want to do things strictly by
the book, you should ask a person before you thank him/her in public
(because, for example, most people would not want to be
thanked heartily in a paper that turned out to be hopelessly incorrect). 
However, as a matter of fact, most people do not engage in this formality;
and most of those who are thanked do not object.  Occasionally you may
wish to thank the referee for helpful comments or suggestions
(best is to do this {\it after\/} the paper has been
refereed---not in advance, or in anticipation of a friendly referee);
sometimes you will need the editor's help in handling this particular
``thank you'' correctly.  Sometimes one thanks one's spouse for
forbearance, or one's typist for a splendid job with the
manuscript.  Sometimes one thanks one's department for time off to
complete the work, or for the opportunity to teach an advanced
seminar in which the work was developed.  The one particular form of
thanks that you are honor bound to include is thanks to any
agency---government, university, or private---that has provided you
financial (or other) support.  In some cases, this thanks is
mandatory; in all cases it is an eminently appropriate courtesy.

Now let us turn to the contents of the paper.
A mathematical paper is not a love letter to yourself
(in\index{love letter to your self}
content it might be, but in form it definitely should
not be).  You are writing about a topic on which you
have become expert.  You have made an advance, and you
want to share it with the mathematical community.  This
should be your point of view when writing up your results.

The simplest way to write a paper is to introduce some notation,
state your theorem, and begin the proof (for simplicity I am
supposing\index{paper!contents of} 
that this is a ``one-theorem paper'').  Such a procedure
probably involves the least effort on your part, it gets the theorem
recorded for posterity, and it might even get the paper published. 
But this methodology is the least effective if you genuinely want
your work to be read and understood, and if you want the ideas
disseminated to the broadest possible audience.

In point of fact a good mathematics paper is {\it not\/} necessarily
written in strict logical order. The reason lies in theories of learning
due to Piaget and others. The point is\index{paper!ordering material in}
simply this: While it can be useful---when recording mathematics for
archiving in the literature---to develop ideas \`{a} la Bourbaki/Hilbert
in strict logical order, {\it this is not the way that we learn.\/} It is
not the way that a typical human being---even a mathematician---apprehends
ideas. This is the case even if the reader is an expert in the subject,
just like yourself.

Reading a mathematics paper is work, and a typical reader approaches the
task with caution. Most people will not read more than a couple of math
papers per month---I mean really read them, verifying all the details.
However, those same people will {\it look\/} at several dozen papers each
month. We all receive a great many preprints in the mail and over the
Internet or by way of preprint servers. We must make choices about which
ones to {\it read\/}.

Having established this premise, let us think about what sort of
paper will encourage the potential reader to plunge in, and what sort
will not.  If the first couple of pages of the paper consist of
technical definitions and technical statements of theorems, then I
would wager that most potential readers will be discouraged. Imagine
instead\index{paper!first paragraph of} 
a paper written as follows.  The first paragraph or two
summarizes the main results of the paper, in nontechnical language. 
The next several paragraphs provide the history of the problem,
describe earlier results, and state exactly what progress the current
paper\index{paper!introduction to} 
represents.  This introduction concludes, perhaps, with
acknowledgements and an outline of the organization of the paper
(either in Table of Contents form or paragraph form).

A reader faced with the latter organizational form has many
advantages.  This person knows {\bf (i)} what the paper is about, 
{\bf (ii)} why the result of the paper is new, {\bf (iii)} what is the
context into which the paper fits, and {\bf (iv)} whether he/she wants
to read on.

One person whom you must keep in sharp focus as you craft your paper
is the referee.  You cannot, indeed you must not, assume that the
referee will compensate for your shortcomings.  If {\it you\/} do not
explain what the paper is about, why you wrote it, why your
theorems\index{paper!referee of} 
are new, why this paper makes an interesting contribution,
why its techniques are original---then nobody else is going to do it
for you.  And the referee, who really does not want to do the full job
of reading the {\it entire paper\/}, will (if the introductory portion
of your paper is not up to snuff) conclude quickly that this is just
another piece of second rate drivel and will reject it.  

Back to the chase:  Imagine that, having concluded the 
introductory section of the paper, you (the writer)
turn to the necessary technical definitions and a formal
statement\index{definitions!placement of}
of results.\index{theorem!statement of}  This central material
would be the substance of Section 2
of the paper (assuming that the introductory material, discussed
in the last paragraphs, was Section 1).  Now the reader---the expert
who has desired to slog this far---knows precisely what he/she is getting
himself/herself into.  Turning to Section 3, you (the writer) can now
dive into all the pornographic details of the proof.  Right?
Wrong.

Reading a difficult mathematical proof in strict logical order is an
onerous task.  If the first five pages of Section 3 consist of a
great many technical lemmas, with nary an indication of where things
are going, of what is\index{proof!organizing}
important, and of what is not, then many 
readers will be discouraged.  Let me now describe a better way.

It is more work for the writer, but definitely a great favor to the
reader, to organize the paper as follows.  Section 3
should consist of the ``big steps'' of the proof.   Here you should
formulate the technical lemmas (provided that the reader can
understand\index{paper!organization of} 
them at this point) and then you should describe how they
fit together to yield the theorem.  You should push the nasty details
of the proofs of the lemmas to the end of the paper.

It should be clear by now why the proposed organizational scheme
makes sense.  First, the reader can decide at each of these signposts
how far he/she wants to get into the paper.  Each new epsilon of effort
on his/her part will yield additional and predictable benefit.  And the
hardest and most technical parts are left to the end for the real
die-hard types.   This writing style is of course beneficial for the
reader; it will also aid the writer.  It disciplines the writer,
forces him/her to evaluate and predigest what he/she has to say, and will
tend to reduce errors.

The principles of writing a math paper that have been described here
do not apply to every paper that is, has been, or ever will be written.
They probably do not literally apply to the Feit/Thompson paper
on the classification of finite, simple groups (an entire issue
of the {\it Pacific Journal\/}) or to Andrew Wiles's proof of
Fermat's Last Theorem (an entire issue of the {\it Annals\/}).
They certainly apply to a twenty-page, ``one-theorem'' paper.
And the general principles described above probably apply in some form
to virtually any mathematics paper.

And now a word about redundancy.  In general, redundancy is a good
thing.\index{redundancy}  One fault that all mathematicians have is this:
we think that when we have said something once clearly then that is the
end of it; nothing further need be said.  This observation explains why
mathematicians so often lose arguments.  You must repeat. 
Help the reader by recalling definitions---especially if the definition
was given 50 pages ago.  If you need to use the definition {\it now\/},
and\index{repetition} 
if you have not used it for quite a while, then give the reader
some help.  Give a quick recap or at least a cross-reference; 
do likewise for a theorem or a lemma that you need to recall.  Think of
how much you would appreciate this assist if you were the reader.

%% Section 2.2

\section{How to State a Theorem}
\markboth{CHAPTER 2.  TOPICS SPECIFIC TO 
MATHEMATICS}{2.2.  HOW TO STATE A THEOREM}

There are some mathematical subjects---geometric measure
theory\index{theorem!how to state}
is\index{geometric measure theory}
one of them---in which the custom is for the
statement of a theorem to occupy one or more pages, and for
the enumerated hypotheses to number 25 or more.  This practice is
too bad.  It makes the subject seem impenetrable
to all but the most devoted experts.  People who present
their theorems in the fashion just described are wont to claim
that their subject prevents any other formulation of the theorems,
that this is just the nature of the beast.  I would like to take
this opportunity humbly to disagree.  It may require extra
effort on the part of the writer, but I claim that you never
need to state a theorem in this tedious manner.

You should strive to hold the statement of a theorem to fewer than
ten lines, and preferably to five lines.  (Some books on writing
assert that a theorem should consist of only a single sentence!)  How
can you do this if there are twenty-five hypotheses?  First of all,
the assertion that there are twenty-five hypotheses is only a
manifestation of what is going on in the writer's mind.  Mathematical
facts are immutable and stand free from any particular human mind, but
the way that we describe them, verify them, and understand them is
quite personal.  In particular, the way that a theorem is presented
is subject to considerable flexibility and massaging.  Let us
consider a quick and rather artificial example:  
\medskip 

\begin{quote}
\noindent {\bf Theorem:}  Let $f$ be a function satisfying the
following hypotheses:

\begin{enumerate}
\item The function $f$ has domain the real number line;
\item The function $f$ is positive;
\item The function $f$ is uniformly continuous;
\item The function $f$ is monotone;
\item The function $f$ is convex;
\item The function $f$ is differentiable except possibly
       on a set of the first category;
\item The function $f$ has range that is dense in the positive real numbers;
\item The function $f$ has no repeated values;
\item The function $f$ is a weak solution of the differential equation
$$
 L f = 0 
$$
(where the operator $L$ has been defined earlier in the paper);
\item The function $f^2$ is a subsolution of $L f = 0$.
\end{enumerate}

Then $f$ operates, in the sense of the functional calculus, on
all bounded linear operators on a separable, real Hilbert space $H$.
\end{quote}
\vspace*{.12in}

This sample ``theorem'' has only ten hypotheses, and these
assumptions are not
all that difficult to absorb; but it
serves to illustrate our
point.  Here is a more efficient, and more user-friendly,
manner in which to state the theorem.

Suppose that, prior to the statement of the theorem, we formally
define a function to be {\it regular\/} if it is defined on the
real line, uniformly continuous, convex, monotone, and positive.
Further, we define a function to be {\it amenable\/} if
it has range dense in the positive reals and has no repeated values.
Finally, let us say that a function $f$ is {\it smooth\/} if it
is differentiable except possibly on a set of the first category,
it is a weak solution of $L$ and, in addition, $f^2$ is
a subsolution of $L$.  Each of these should be stated
as a formal definition, prior to the formulation of the theorem.
Moreover, we should state that, until further notice,
$H$ will designate a separable, real Hilbert space
and ${\cal L}(H)$ the bounded linear operators on $H$.  With this
groundwork in place, we can now state the theorem as follows.
\medskip \hfill \break

\noindent {\bf Theorem:}  If $f$ is a regular, amenable, smooth
function, then it operates on ${\cal L}(H)$ in the sense of the 
functional calculus.
\medskip \hfill \break

Notice that, by planning ahead and introducing the terms ``regular,''
``amenable,'' and ``smooth,'' we have grouped together cognate
ideas.  We are not just engaging in sleight of hand; in fact we are
providing organization and context.  We are also 
helping the reader by keeping the statement of the theorem short and
sweet.  The reader will come away from reading the theorem
remembering that {\bf (i)} there is a hypothesis about $f$ involving
continuity, convexity, and so forth, {\bf (ii)} there is a hypothesis
about the value distribution of $f$, and {\bf (iii)} there is a
hypothesis about the way that $L$ acts on $f$.  The conclusion
is that $f$ operates on ${\cal L}(H)$.  You, the writer, have done
some of the work for the reader, and given him/her something to take
away. The reader can always refer to the text for details as
they are needed.  But if the theorem is recorded in the first form
rather than the second then, most likely, the reader will
not quite know what he/she has read, nor when and where he/she 
can use it.\index{theorem!stating in one sentence}

Also note that we managed to state the theorem in one sentence,
and in just two lines.

%% Section 2.3

\section{How to Prove a Theorem}
\markboth{CHAPTER 2.  TOPICS SPECIFIC TO MATHEMATICS}{2.3.  HOW TO 
       PROVE A THEOREM}

What I mean here, of course, is ``how to {\it write the proof
of a theorem\/}.''  You are not doing your job---unless the proof is
short and fairly simple---to begin at the beginning and
charge through to the end.  A proof of more than a few
pages should be broken\index{theorem!how to prove}
into lemmas and corollaries
and organized in such a fashion that the reader can
always tell where he/she has been and where he/she is going.

A useful device in writing up a proof is the ``Claim.''
This tool is often used in the following manner.  You have
set up the basic pieces of your proof; that is, you have
defined the sets and functions and other objects 
that\index{claim!use of} 
you need.  You are poised to strike.  Then you write ``We claim
that the following is true.''  State the claim.  Then you say
``Assuming this claim for the moment, we complete the proof.''

Used correctly, this technique is a terrific psychological device.
It allows you to say to the reader ``This is the crux of the proof,
but its verification involves some nasty details.  Trust
me on this for the moment, and let me show you how the crux
leads to a happy ending.''  The reader, having arrived at the
end of the proof (modulo the claim), will feel that progress
has been made and he/she will be in a suitable mood either to
study the details of the claim or to skip them and come back to
them later.

Another useful device---nearly logically equivalent to the ``claim''---is
to enunciate a technical lemma right at the point where you need
to use\index{lemma!use of} it (sometimes a good idea because to enunciate it well in
advance would make almost no sense to the reader), but then
to say ``Proof Deferred to Section 8.''  If you indulge in this trick,
be sure\index{deferred proof} that your paper is well organized and that the
different parts of the paper are well labeled.  Do not leave
your poor reader with a head full of dangling claims and unproved
lemmas to sort out.  A good rule of thumb ([Gil, p.\ 8]) is to be sure
that your reader always knows the {\it status\/} of every statement that 
you make.\index{status of your statements}

In Section 2.1 I have advocated that a paper should be
organized so that the technical stuff is at the back and the
explanatory stuff at the front.  The paper should proceed, by
gradations, from the latter to the former.  The proof of a theorem
should proceed in roughly the same way.  You, as the author and
creator of the theorem, have the whole thing jammed into your head;
it has no beginning and no end---it just resides there.  Part of the
writing process is to transfer\index{proof!organizing} 
this organic mass from your head to
someone else's.  Thus, as you write, try to provide signposts so that
the reader always knows where he/she has been and where
he/she is going.  This writing goal is best achieved by pushing the
technicalities to the end.

Many books, and some papers, are written as follows:  the author
rattles on for several pages---chatting about this and that---and
abruptly says 
\begin{quote}
Note that we have proved the following theorem:
\medskip \hfill \break
{\bf Theorem [The Riemann Mapping Theorem]:}\ \ Let $\Omega$
be a simply connected, proper subset of the complex plane \dots. \qquad \bad
\end{quote}

\noindent Good heavens!  What a disservice to the reader.  The 
Riemann mapping theorem is a milestone in mathematical thought,
perhaps even in human thought.  Each of the steps in its proof---the
extremal problem, the normal families argument, etc.---is a subject
in itself.  The writer must lay these milestones
out for the reader and must pay due homage to each.  The offhand ``Note
by the way that we have proved the Riemann mapping theorem'' is a real
travesty, and ignores the author's duty to {\it explain}.  Rise above
the idea that it suffices for the writer to
somehow record the thoughts on the page; if you, the author,
have not crafted them and worked them and, indeed, handed them to the
reader, then you have not done your job.

And here is a small note about proofs by contradiction.  Some
mathematicians\index{proof!by contradiction} 
begin a proof by contradiction with

\begin{quote}
Not.  Then there is a continuous function $f$ \dots  \qquad \bad
\end{quote}

\noindent Others begin with 

\begin{quote}
Deny.  Then there is a continuous function $f$ \dots  \qquad \bad
\end{quote}

\noindent This is all rather cute; the first of these is perhaps a
tribute\index{Saturday Night Live}
to John Belushi and the {\sl Saturday Night Live} gang. But
both examples (and these are {\it not\/} made up---people actually
write this way) hinder the task of {\it communicating\/}.  A
preferred method for beginning a proof by contradiction is 

\begin{quote} 
Seeking a contradiction, suppose that $f$ is a
continuous, real-valued function on a compact set $K$ that does not
assume a maximum.  Then \dots 
\end{quote}

%% Section 2.4

\section{How to State a Definition}
\markboth{CHAPTER 2.  TOPICS SPECIFIC TO MATHEMATICS}{2.4.  HOW TO STATE A DEFINITION}

Definitions are part of the bedrock of mathematical writing and thinking.
Mathematics is almost unique among the sciences---not to mention
other disciplines---in insisting on strictly rigorous
definitions\index{definitions!importance of} 
of terminology\index{definition!how to state} and concepts.
Thus we must state our definitions as
succinctly and comprehensibly as possible.  Definitions
should not hang the reader up, but should instead provide
a helping hand as well as encouragement for the reader to
push on.

As much as possible, state definitions briefly and
cogently.  Use short, simple
sentences rather than long ones.  To avoid excessively complex
and introspective definitions, endeavor to {\it build\/}
ideas in steps.  For instance, suppose that you are writing an
advanced calculus book.  At some\index{definition!statement of}
point you define what a function
is.  Later you say what a continuous function is.  Still later you
state what the intermediate value property for continuous functions
is.   Further on, you use the latter property to establish the
existence of $\sqrt{2}$.  You do not, all at once, attempt
to spit out all these ideas in a single sentence or a single
paragraph.  In fact you build stepping stones leading to the key idea,
so that the reader is given a chance to internalize idea $n$ before going
on to idea $(n+1)$.

Just how many definitions should you supply?  If you are writing a
paper on von Neumann algebras (algebras of bounded operators on
Hilbert space), then you certainly need not say what a Hilbert space
is, nor what a bounded linear operator is.  Every graduate student
who has passed through the qualifying exams is familiar with these
ideas, and\index{definitions!how many to supply} 
you may take these for granted (that is why we have
qualifying exams).  Define ${\cal L}(H)$ (see Section 2.2) only if
you think that readers likely will misinterpret this
(rather standard) notation.  Of course you would have to define
``regular,'' ``amenable,'' and ``smooth'' (the terminology that we
introduced in Section 2.2). Those terms are not standard, and have
been given other specialized meanings elsewhere.

What I am describing here is another of many subjective matters
that pertain to writing.  If your paper supplies too few,
or poorly written, definitions then both the referee and the readers
will lose their patience.  If your paper supplies too many
definitions, then you also will irritate your audience.  For standard
terminology, you could give a well-known reference like Dunford
and Schwartz [DS] or Griffiths and Harris [GH] or Birkhoff and
MacLane [BM] or Kuratowski [Kur].  This habit is preferable to taking up
valuable journal space with a rehash of well-known ideas.\footnote{But do
your reader the kindness of telling him/her to what {\it part} of [DS] or [GH]
to refer.  Best is to give a specific definition number.  But giving
the page number is fine.   Or even the section number is OK.}  
Less kind is to refer to a\index{definitions!how many to include}
semi-obscure journal article for terminology.  If
that is the best reference for definitions, then you should probably
repeat them.

There is some terminology that you simply cannot take the space to
repeat or define, even though it is rather advanced.  For example,
you cannot rehash---for the convenience of your readers---the
standard theory of elliptic partial differential equations, nor the
basics\index{definitions!which not to include} 
of $K$ theory, nor the guts of the Atiyah-Singer Index
Theorem.  (In writing a book you in fact {\it can\/} indulge in such a
review; I treat book writing elsewhere.)  Try to refer the reader to
a good source for the important ideas on which you are building.

I have advocated (Section 2.2) the tasteful use of terminology to clump ideas
together, thus making them more palatable for the reader.  However,
try to avoid introducing any more new terminology than is necessary.
If your paper contains a plethora of unfamiliar language, then it may
cause your reader to suspect that you actually have nothing to say.
And\index{terminology as organizational tool} 
if there is a standard bit of notation or terminology for what
you are saying, then by all means use it.  I once saw a paper in a
standard mathematics journal of good repute that defined the space
$Q^{17}_{\rm reg}$ to be the set of all bounded holomorphic functions
on the unit disc in the complex plane.  Of course the well-known
notation $H^\infty(D)$ describes this space of functions, and it is
virtually mandatory to use {\it that\/} notation.  The proposed alternative
notation is just crazy, {\it unless\/} the author is introducing a
whole new scale of function spaces in which $H^\infty$ arises in a
natural way.  If that is the case, then the author should certainly
mention this relationship explicitly.  (For example, all the standard
function spaces---$L^p$, Lipschitz, Sobolev, Hardy, Besov,
Nikol'skii, etc.---are special cases of the Triebel-Lizorkin spaces
${\stackrel{\cdot}{F}}{}^{\alpha,q}_p$.  Thus, in certain contexts,
it would be appropriate to refer to the Lebesgue spaces, or the
Sobolev spaces, using the Triebel-Lizorkin notation.)

Good notation is extremely important, sometimes as important as a
theorem.  As an example, the notation of differential forms is a
small miracle.  Large parts of geometric analysis would be completely
obscure without it.  Of course you cannot perform at the level of Elie
Cartan\index{notation!importance of} 
every time you dream up a piece of notation, but you can
consider following these precepts: {\bf (1)} Do not create new
notation if there already exists well-known notation that is suitable
for the job at hand; {\bf (2)}  If you must introduce new notation,
then think about it carefully; {\bf (3)}  Strive for simplicity and
clarity at all times.  

Fiddle with several different notations before
you make a final decision.  Consult the standard references in
the field to see whether they give you any ideas.  If possible,
try\index{notation!choosing} 
your new notation out on a colleague, or on one of your
graduate students.

Technically speaking, a definition should almost always be 
formulated in ``if and only if'' form.  For example

\begin{quote}
A function $f$ on an open interval $I$ is said to be
{\it continuous\/} at $c \in I$ if and only if, for every
$\epsilon > 0$, there is a $\delta > 0$ such that \dots
\end{quote}

\noindent In practice, we generally replace the phrase ``if and only
if'' in\index{definition!use of ``if and only if'' in}
 this definition with ``if.'' We do so partly out of laziness,
and partly because the ``if'' phraseology is less cumbersome than
``if and only if.''  The price that we pay for this convention is
that we must teach our students to read definitions; the fact is that
we {\it do not\/} write what we mean.  

Although nobody will punish you for writing ``if and only if'' in
your definitions, and some will appreciate it, it is usually best to
follow mathematical custom and simply to write ``if.''  A useful,
and modern, compromise is to use Paul Halmos's\index{Halmos, Paul}
invention ``iff'' (see Section 1.8).  The word `iff'' captures the
brevity of ``if'' but carries the precision of ``if and only if.''
\vspace*{-.23in}

%% Section 2.5
\section{How to Write an Abstract}
\markboth{CHAPTER 2.  TOPICS SPECIFIC TO MATHEMATICS}{2.5.  HOW TO WRITE AN ABSTRACT}

Many journals now require that, when you submit a paper,
you include an abstract of the paper.  The abstract, usually
not\index{abstract!how to write}  
more than ten lines, is supposed to convey on a quick
reading what the paper is about.  According to the strictest
standards, the abstract should be self-contained, should not
make any bibliographic references, and should contain a minimum
of notation and jargon.

A rough rule of thumb is that any reader who looks at your paper will
read the abstract, only 20\% of those will read the introduction, and
perhaps one fourth of that 20\% will dip into the body of the paper. 
This being the case, your abstract is obviously of preeminent
importance.\index{abstract!notation in} 
Many\index{abstract!references in}
indexing and reviewing services will rely on your
abstract.  So it had better give a clear picture of what is in the
paper.

As usual, endeavor to employ simple, short, declarative sentences in
your abstract.  Eschew nasty details.  Do not say, with a plethora of
$\epsilon$s and $\delta$s, exactly what interior elliptic estimate
you are proving; instead state that you are proving a new interior
elliptic estimate in the Nikol'skii space topology and that it
improves upon classical results of Nirenberg.  State that it has
applications to certain free boundary problems.  The interested
reader\index{abstract!simplicity of} 
can then move on to the introduction, where further details
are provided.

If your abstract is too long or too short, then the editor of the
journal will likely make you rewrite it.  The ``Instructions to
Authors'' section in the journal should give you an idea of what is
required for an abstract in {\it that\/} journal.  Study several
abstracts\index{Instructions!to Authors} 
in the journal to which you plan to submit to get an idea
of what is suitable.

%% Section 2.6

\markboth{CHAPTER 2.  TOPICS SPECIFIC TO 
MATHEMATICS}{2.6.  HOW TO WRITE A BIBLIOGRAPHY}
\section{How to Write a Bibliography}
\markboth{CHAPTER 2.  TOPICS SPECIFIC TO 
MATHEMATICS}{2.6.  HOW TO WRITE A BIBLIOGRAPHY}

The bibliography, or list of references, is one of the most important
components of a mathematical work. This assertion is true for research
articles,\index{bibliography!how to write} for books, and for
expository articles as well. The bibliography tells the reader where you
are coming from and where you are going, it keeps you honest, and it
provides critical assistance for those readers not already familiar with
the subject.

Real sticklers---mavens of good scholarly form---will tell you that a
bibliography should {\it only\/} be assembled from primary sources. The
book\index{primary sources, use of} 
[Hig, pp. 87--8] has several examples of bibliographic
inaccuracies in the literature that have been propagated for dozens
of years because reference $(n+1)$ was always copied from reference $n$. 
The book [Hig] also advises you never to retrieve information about a
paper either from the cover of the journal issue or from the Table of
Contents since information is frequently misrepresented in both
places.  Even if it is not, you could easily get the first or last
page of an article wrong.  Accuracy and scholarship are best served
when you gaze upon the actual paper; and you will also be
able to say truthfully that you have ``looked'' at the paper.  

Purists also will tell you
that each reference should include an {\bf MR} (or {\it
Math.~Reviews\/}) number.  Such an addition is often quite convenient
for\index{Math Reviews number} 
the reader, and a lot of extra work for the writer (though with
the advent of MathSciNet (Section 7.2),\index{MathSciNet}
the web service available from the American
Mathematical Society, the task has become much easier).

You should only list references in your bibliography that you 
also cite in the text.  We are frequently tempted to include
extra references either for sentimental reasons or because
we think that these references might be handy for the reader.
The\index{references, which to cite}
former motivation is spurious, and the latter misguided.
If you give the reader no advice on the value of a reference,
then you are offering nothing by listing it.

What some writers (at least in writing a book) do is, in addition to
writing a Bibliography, they also provide a list called ``Additional
Reading.'' This list consists of books, and perhaps papers, that are {\it
not} cited in the text.

There are many possible formats for bibliographic entries.  If you use
AMS-\TeX, then your bibliographic entries are formatted for you
automatically\index{bibliography!in \TeX} 
in the approved AMS style.  (Similar comments apply to
\LaTeX's\index{latex@\LaTeX} treatment of bibliographic entries.)  
But if you do not then
you must make some choices.  At the beginning of my career
I picked a favorite journal and adopted its bibliographic style.
I chose a format that is commonly used, and it has served me well.
Here\index{bibliography!style of} 
are two bibliographic references formatted in that style:

\begin{quote}
[Bat] \ \ Gill Bates, {\it How I Made My First Billion}, 2nd
ed., \hfill \break
\quad \quad \  Acquisitive Press, New York, 1986.
\end{quote}

\begin{quote}
[Beh] \ \ Viscount Hugh Behave, Some theories on the gentle \hfill \break
\quad \quad \   art of belching in public, {\it 
The Journal of Eminently} \hfill \break
\quad \quad \  {\it Forgettable Theories\/} 
42(1976), 35--53.
\end{quote}

\noindent The first of these is a book, and the second a paper in a
journal.  Notice that the information provided for a book is
different\index{bibliographic!style, for a book} 
from\index{bibliographic!style, for a paper} 
that provided for a journal article or paper.  For a book, the author,
title, edition number (if this is not the first edition),
publisher, city of publication, and date of publication are
usually considered complete bibliographical data.  These are shown,
in order, in the example given.  For a paper, the author, title,
journal, volume number of the journal, year, and pages are usually
considered complete bibliographical data.  Some people include
the issue number of the journal.  Of course the protocol for
a preprint, for a conference proceedings, for an unpublished
manuscript, for a translated paper, and for a Ph.D.\  thesis are all a
bit different.  I shall not go into all the details here.  The software
Bib\TeX\ (part of \LaTeX) provides particular formats for all
these special types of references.  See
[SG, pp.~407-410] or [Hig] or [VanL] for particulars.

Not everyone likes the use of acronyms for citing elements of the
bibliography.   Some people prefer to number the elements of the
bibliography from [1] to [$n$].  The method of enumeration has the
disadvantage\index{acronym!use of in bibliographic references} 
that, if you add or delete a reference late in the game,
it throws off all your numbering.  However using good software can
circumvent that problem (see below).  Both the numbering scheme and
the acronym scheme have the disadvantage that even a one-character
typographical\index{numbers!use of in bibliographic references}
error can make it virtually impossible for the reader
to tell which reference was intended.

One excellent scheme for bibliographic references, and one that is
virtually essential when the bibliography is long, is illustrated 
in the following example.  It lists three works by John Q.\ Public,
just as they might appear in a bibliography.

\begin{quotation}
\noindent John Q.\ Public  \hfill \break
\vspace*{-.2in}

\begin{tabbing}
\noindent [1987] \ \ \ \= {\it Why I Never Vote}, Ignoramus Press, Brooklyn.   \\
\noindent [1992a]  \> {\it The Less I Know, the Better}, Rosicrucian Press, \\
\null              \> Poughkeepsie.   \\
\noindent [1992b]  \> On Doctoring Polls, {\it The Smart Pollster\/} 31, 59-71.
\end{tabbing}
\end{quotation}

\noindent If you use this system (known as the
{\it Harvard system}), 
then when\index{Harvard system for bibliographic references} you refer to a
bibliographical item in the text you say ``By J. Q. Public [1992b],
we know that \dots.''

Note that in mathematics we do not usually put bibliographical
references in footnotes (however it {\it is\/} customary in
certain statistical work, and it was fairly common in
mathematics one hundred years ago). This habit came about in
part because typesetters objected to the expense and trouble
of typesetting copious footnotes. With the advent of \TeX,
that particular objection is moot. However, the rule persists.
In fact, if you were to submit to most mathematics journals a
paper with all the references in footnotes, then you would
most likely be asked to reformat it. The trouble with using
footnotes in a mathematics paper is that the footnote tags can
be mistaken for exponents.

If you are writing your paper in \LaTeX, then you have the option
of using \LaTeX's bibliographic utilities.   One of \LaTeX's
tools allows (see, for instance, Bib\TeX) you to assign a nickname to each of your bibliographical
references.\index{bibliographic!references in \LaTeX}  
Then, in the text, you can cite any reference by
its nickname.  When you compile your {\tt *.TEX} file, each nickname
citation is replaced by the appropriate preassigned acronym or
number; the full bibliographic citation occurs at the end of the
document as usual (see Section 6.5 for more on \TeX\ and \LaTeX).

Slightly more sophisticated is \LaTeX's bibliographic database
system.  With this device, you never write another bibliography as
such.  You simply have an ever-growing database of bibliographic
references.  Whenever a new reference comes to hand, you add it
directly to the database.  Each reference has a preassigned acronym
and a preassigned nickname.  Then, when you are writing a new
document, you make a reference by referring to the appropriate
nickname in the database (if you cannot remember all the
nicknames---perhaps your database has thousands of items in it!---you
can just pull the database into a window with your text editor and
check it).  When you compile the document, a beautiful bibliography
is created for you, with the requisite information pulled in from the
database.

If the last systems do not appeal to you, then you also can keep the \TeX\
files of all your papers in a single directory.  Most of us tend to
use many of the same references repeatedly.  Thus, when you are
writing a new paper and need a reference, you can open a window
with your text editor, pull in an earlier paper that has the
reference, and cut and paste the reference into your new document.

Incidentally, \LaTeX\ also allows you to assign nicknames to your
equations and theorems.  You can refer to them, during the writing
process, by nickname.  Then, when the paper compiles, the correct
line numbers and theorem numbers are inserted for you automatically.

The advantage of these \LaTeX\ devices is that you no longer have
to worry about numbering of theorems or of bibliographic references.
If you insert a new theorem, you no longer have to renumber all the
old theorems.  \LaTeX\ does it for you.

I know mathematicians---excellent ones---whose bibliographies look
like this:

\begin{quote}
{\bf 1.}  Knuth, 1992.

{\bf 2.}  Lister, 1991.

{\bf 3.}  Machedon, 1988.
\end{quote}

\noindent This is it!  No titles, no journal names, no volume numbers, no
page references.  This scheme in effect takes the \LaTeX\ device to the
limit:\index{bibliographic!reference, unacceptable}  
you just supply the nicknames but none of the details.

The practice of listing abbreviations in lieu of correct
bibliographic references is irresponsible.  In truth, such
sloppiness should have been caught by the editor, who should
have demanded that the author rectify the matter.  As indicated at the
beginning of this section, the bibliography is part of your paper
trail.  You hold the responsibility for providing complete bibliographic
information.  It should be complete in the sense that you have cited
everyone who merits citation, but it also should be complete in the
sense that all the information is there.  The bibliographic sample
just provided might mean something to a few experts for a few years.
In fifty years it would not mean anything to anyone.

And speaking of ``meaning nothing to anyone,'' do {\it not\/} give
in-text bibliographic references that have the form ``see Dunford
and Schwartz'' (for those not in the know, [DS] is a three volume
work totaling more than 2500 pages).  The only  
correct and thorough way to give a reference is to cite the specific theorem
or the specific page.  Sometimes, to conserve space and to prevent
repetition, we say ``by a variant of theorem thus and such'' or
``by a variant of the argument in this paper''  (the subject of
analysis, in particular, seems to be littered with references
of this nature).  If you find such references necessary in your
own work, be as specific as you can so that the reader may follow
your path. 

Modern technology enables a marvelous writing environment---at
least in principle.  If, for example, I am a {\tt Windows}\reg\
 user, then I\index{windows@{\tt Windows}}
can have my text editor going in one window (this is where I actually
do my writing), a thesaurus and dictionary on CD-ROM in another, the
library's\index{technology and the bibliography} 
on-line catalog in another, and MathSciNet\index{MathSciNet} on line in a
fourth.  Passing from one environment to the next requires only a
keystroke or a mouse click or two.  Clearly such an environment makes
tedious trips to the library a thing of the past, and makes
assembling a bibliography relatively quick and easy.

Now let us treat styles for citation.  In this section, I have
spoken of bibliographic references with the assumption that they
will\index{citation!styles} 
occur on the fly, right in the text.  For example:
\begin{quote}
By a theorem of Steenrod [Ste], we know that every instance
of generalized nonsense is a generalization of specific nonsense.
\end{quote}
The good feature of this methodology is that it tells you right
away what the source is.  The bad feature is that it clutters up
the text a bit.  In most mathematics {\it papers}, the on-the-fly
style is used.  You make a reference either by acronym, or by number,
or\index{citations!on-the-fly} 
by author surname, but the reference occurs at the moment
of impact.

In [Ste], Steenrod fulminates against this 
bibliographic style for the
writing of a book.  His preference is to have a paragraph or more at
the end of each chapter detailing the genesis, development, and
sources\index{citations!Princeton style}
for the theorems in that chapter.  This methodology is commonly
known as the ``Princeton style.'' 

Many books in the
Princeton book series {\it Annals of Mathematics Studies} handle
bibliographic references in this fashion.  These little 
end-of-chapter essays can be
quite informative and, if well written, can give the reader a sense
of the historical flow of thought that in-context references (as
indicated above) do not.  I would say that the down side of this
end-of-chapter approach is the following. It serves the big shots
well.  If you are annotating a chapter on singular integrals, then
you will certainly not overlook Calder\'{o}n, Zygmund, Stein, and the
other major figures.  But you might overlook the smaller
contributors.  The advantage of the in-text, on the fly reference
method is that it systematically holds you accountable:  you state a
theorem, and you give the reference; you recall an idea, and you give
the reference.  You are much less likely to give someone short shrift
if you adhere to this more pedestrian methodology.  Of course the
final choice is up to you.

%% Section 2.7

%\markboth{CHAPTER 2.  TOPICS SPECIFIC TO MATHEMATICS}{2.7.  WHAT TO DO ONCE
%    THE PAPER IS WRITTEN}
\section{What to Do with the Paper Once It Is Written}
\markboth{CHAPTER 2.  TOPICS SPECIFIC TO MATHEMATICS}{2.7.  WHAT TO DO ONCE
    THE PAPER IS WRITTEN}

Ours is a profession where, by and large, we are left on our
own to figure out how to function.  Nobody shows us how to teach,
nobody tells us how to write a paper, and nobody tells us how
to get published.  This section addresses the last issue.

So imagine that you have written a paper that you think is good.  How
do\index{paper!how to write} you know it is good?  Being a
mathematician is a bit like being a manic depressive:  you spend your
life alternating between giddy elation and black despair.  You will
have difficulty being objective about your own work:  before a
problem is solved, it seems to be mightily important; after it is
solved, the whole matter seems trivial and you wonder how you could
have spent so much time on it.  How do you cut through this
imbroglio?

If you are smart, you have told some colleagues about your results. 
Perhaps you have given some seminars about it. You have sent 
preprints around (either
by {\it e}-mail or\index{paper!publishable}
 by snail mail) to colleagues.  If
you have kept your ears open, you have some sense of how receptive
the world is to your ideas.  Are your listeners surprised, impressed,
confused, bored?  Sometimes they will suggest changes.  Consider all
criticisms and suggestions carefully, and make appropriate changes to
your paper.  Now you must decide where to submit it.

Before you make that momentous decision, let me back-pedal a minute
and\index{paper!where to submit} 
address the question of how to decide when you have something
that is worth writing up.  This is a confusing issue, and one
that every mathematician must learn to face.

We all know that the keys to success in this profession of ours
include intelligence, perseverance, drive, and hard work (not
necessarily in that order).  Some may deny it, but there is also an
art to the business.  Let me explain.  Ideally, the working
mathematician sets a problem for himself:  solve the (restricted) Burnside 
problem, or calculate the dual of the Hardy space $H^1$, or prove the
corona theorem in several variables.  We all know that there are
extraordinary mathematicians who can actually do just this:  E.\ Zelmanov
did the first and C.\ Fefferman did the second.  Nobody has done the
third, although many of us have tried.  In practice,
this point-and-shoot technique is rarely the way that mathematics is
successfully practiced.  

A somewhat more modest way to get one's feet wet is this:  become
completely immersed in a subject, and then formulate a program.
Determine to assume hypotheses $A, B, C$ and endeavor to prove
conclusion $X$.  Sadly, this {\it modus operandi\/} is also only
occasionally successful.

In fact what happens in practice is that we try a great many things.
Some succeed and some do not.  Along the way, hypotheses are 
constantly being altered and substituted and strengthened; conclusions
are redirected or transmogrified or reversed.  The theorem that you
end up proving is rarely the theorem that you set out to prove.
This is a perfectly reasonable way to proceed.  Columbus sought
a new passage to India and instead found America.  Jonas Salk
discovered the polio vaccine by accident.  Milnor 
discovered multiple differentiable
structures on the 7-sphere because calculations on another problem
were not working out as planned.

One of the chief differences between a successful mathematician and an
also-ran is that the former can take his/her partial results and his
tries---and yes, even his/her failures\footnote{A twentieth century
Hungarian philosopher once said that a mathematician is nothing but a
collection of statements that he/she cannot prove.}---and turn them into
an attractive tapestry of theorems and corollaries and partial results and
conjectures; the latter instead takes two years of hard work and dumps it
in the trash.

As you read these words, do not suppose that I am advocating any degree of
chicanery, or self-promotion, or hype. I am instead encouraging you to
have the confidence and fortitude to make something of your work. Part of
doing mathematics successfully is to get in there and calculate and reason
and think and ponder. But another part is to evaluate and organize and
deduce. What I am describing is a bit easier to imagine for a laboratory
scientist. He/she performs a huge experiment that may take a year or two
and may cost a few million dollars. No matter how things turn out, he/she
must make a show of it. He/she must report to his/her granting agency and
write papers about how his/her laboratory has been spending its time and
effort. The message here is that a mathematician must do something
similar, but his/her wherewithal is somewhat more tenuous; indeed it is all in
his/her head. Part of training yourself to\index{assessing your work}
survive in this profession is coming to terms with the reality that I have
described.

I cannot conclude this digression without also noting that another key to
success is actually making some progress. It just will not do
to\index{progress in mathematics} tell yourself (and the world) that for
the next twenty years you will work on the Riemann hypothesis, {\it unless\index{Riemann hypothesis}
you can arrange to have something to show along the way}. You do not get
tenure, or a promotion, or an invitation to the International Congress by
advertising that you are working on a great problem and telling people
that they should contact you a generation later to see how things worked
out. I have a friend who has a twenty-five step program for proving the
Riemann hypothesis: ``Count to twenty-four and then prove the Riemann
hypothesis.'' There is wisdom in this little joke. The successful
mathematician knows how to manage his/her research program so that it will
proceed incrementally, so that he/she can report progress along the
way---including writing up papers and giving talks and showing the world
just what he/she is up to. By the same token, the good mathematician knows
how to determine when he/she is {\it not\/} making progress, when his/her
program is {\it not\/} paying off, when it is time to move on to something
else.

Now let us return to more pedestrian matters.  Let us suppose
that you have organized some of your material and turned it
into a paper.  You believe that this is a worthy piece of work.
You want to get it published.  The next move is yours.

Keep in mind that the one hard and fast rule in this business is that
you can\index{paper!where to submit}
submit a paper to just one journal at a time.  {\it Never
consider\index{paper!how to submit} 
deviating from this policy.}  In the words of Clint Eastwood,
``Don't even think it.''  If you do send the same paper to two
different journals simultaneously, then that paper is liable to be
sent to the same referee by both journals; thus you will be caught
red handed!   Agonizing though it may be, you must wait for a
decision from journal $n$ before you submit to journal $(n+1)$.  As a
result, there is considerable motivation to exercise wisdom when
choosing a journal.

There is a distinguished mathematician, now retired, who in his heyday
wrote about a dozen papers per year.  He submitted them all to {\it
The Annals of Mathematics}.  Several of his papers were accepted by
the {\it Annals}.  Others were either rejected or else the author was
asked to perform various revisions.  Now, writing twelve papers per
annum as he did, this mathematician had no time for revisions.  So,
in cases two and three, he sent his papers to a well-known journal
that was reputed to have minimal standards (what the famous computer
scientist Dijkstra would call a ``write-only'' journal).  Thus this
esteemed man has a publication list, emblazoned in {\tt MathSciNet\/} 
for all to see, consisting of several citations in the {\it
Annals\/} alternating with citations in this other ``catch-all''
journal.

Another famous mathematician was in the habit of bringing his latest
preprint to the departmental secretary, together with a list of
journals to which it might be submitted.  Her job was to cycle
through the journals on the list, one by one, and to inform Herr
Doktor Professor when his paper was finally accepted.  In this way
the good Doktor was spared the grief of dealing with surly referees
and uncooperative editors.

The preceding two strategies are amusing but probably unwise for most
mathematicians.  The working mathematician should have a
sense of which are the very best journals, which are at the next
level, and which are of average quality.  How can one gauge which
journals are which?  They all look rather elegant, and all profess to
have high standards.  They all have distinguished people on their
editorial boards.  What is the trick?

Begin by considering where cognate results have appeared.  The {\it
Journal of Algebra\/} will probably not consider papers on singular
integrals.  The {\it Journal of Symbolic Logic\/} probably does not
publish papers on Gelfand-Fuks cohomology.  Certain journals have
become\index{paper!where to submit} 
the default forum for work on operator theory or several
complex variables or potential theory.  Consider
those if your work fits.  You will naturally consider which
editors will understand what your paper is about and will know how to
select a referee.  You need not actually {\it know\/} the editor, but
it is comforting to know where the editor is coming from.

If you submit your work to a journal of the highest rank, then you might
pay in several ways:  {\bf 1)}  the refereeing process may take an
extra long time, {\bf 2)}  the journal may have a huge backlog,  {\bf
3)} the paper may be rejected for almost any reason.  Thus the entire
process of getting your work published could drag on for two years or
more.  If you are fighting the tenure clock, this could be a
problem.  In some ways it is better to err on the low side.  Usually
mathematical work is judged on its own merits.  Nobody will downgrade
your work, or you, if your theorems are not published in the optimal
journal.  But do not publish in an obscure journal that nobody ever
reads.  

Part of the secret to success in this profession is to talk to
people.  Doing so, you will quickly learn that {\it Acta
Mathematica}, the {\it Annals}, {\it Inventiones}, and the
{\it Journal of the American Mathematical Society} are four of the
pre-eminent mathematics journals. This choice of four reflects my
prejudices as an analyst.  Others might name {\it The Journal of
Differential Geometry\/} or {\it The Journal of Algebra\/} or the {\it
The Journal of Symbolic Logic\/} as being at the top.  Opinions
will vary. Perhaps {\it Duke}, the {\it Transactions of the AMS},
the {\it Journal of Geometric Analysis}, and
several\index{journal!selecting one for your paper}
others are at the next level.  And on it goes.   There are
prestigious journals and there are excellent journals.  Many journals
fit into both categories, and many fit into neither.
There are nearly 2000 mathematics journals in the world,\footnote{By contrast,
geophysics has only about five journals.} so you have many choices
as to where to publish.

You can form your own opinion of journals by seeing what papers they
publish and by which authors; you can look at how many truly eminent
people (and from which universities) are on their editorial boards,
and\index{journal!ranking of} 
you can learn something just by submitting your papers to various
journals and seeing what happens.

Since the latter strategy is costly---in terms of time, and
perhaps your bruised feelings as well---you should
develop a sense of what is a typical {\it Annals\/}
paper, what is a typical {\it Transactions\/} paper, and 
what is a typical {\it Rocky Mountain Journal\/} paper.
If you are in doubt, ask someone with more experience.
If someone whom you respect and trust has read your preprint,
then he/she would be an ideal person to ask for suggestions
as to where to submit.
 
Also of interest in considering journals is the backlog of the
publication, the turnaround time from submission to acceptance (or
rejection),\index{journal!backlog of}
and similar data.  Fortunately, the {\it Notices of the
AMS\/} publishes a detailed
analysis---containing just this sort of information---of all the major
journals at least once per year.  In the end, you have the
responsibility to pick a suitable journal for your work; and the
choice is not a trivial matter, since a year may pass while you are
waiting for an acceptance or rejection.

Of course you will learn from experience.  You also will have to
decide for yourself whether to shoot high and take your chances,
or to shoot low and optimize your likelihood of a quick acceptance.
If your tenure case is a few years down the road, then this
choice should not be taken lightly.  Deans tend to know which
are the good journals and which are not.  (In fact I know of
several universities where the dean has circulated a ranked list
of mathematics journals.  The implication is that ``If you want
to\index{journal!selecting one for your paper} 
get promoted then you had better publish in these journals
but not in those journals.'')
They are not impressed by a young assistant professor whose work is all
submitted to ``gimme'' journals.  They are also not impressed
by a dossier with most papers ``submitted'' but not yet accepted.

\def\email{{\it e}-mail}

Most journals have a section called ``Instructions to Authors'' or
``Instructions for Submission.''  Before you submit paper $X$ to
journal $A$, you should read those instructions.  They will tell you
how\index{Instructions!to Authors} 
many copies are needed, whether the title and abstract and other
data should be on a separate page, whether the journal requires key
words\index{Instructions!for Submission}
 and AMS subject classification numbers, what languages the
journal will accept (English, French, and German are the most
common---though there {\it are\/} mathematics journals that
will take papers in Latin or Esperanto or Japanese), any
formatting requirements, length restrictions, where to send
the paper (to the Editorial Office, or to an Associate Editor
of your choosing, or perhaps another option), whether the
journal prefers submissions in \TeX, whether the journal has a
\TeX\ style file that you should use, whether the journal
accepts electronic submissions, and so forth. You will annoy
the editors, and cause unnecessary delays and confusion, if
you do not follow these readily available instructions. Of
course many journals now take electronic submissions (i.e.,
uploading at an Internet site), and this system has the
advantage that it {\it forces} you to make the right choices.
Some journals just ask you to send the paper as an \email\
attachment to one of the editors.

It must be noted that, in today's world, it is common to submit
a paper to a journal by way of a Web site.  You simply fill
out some Web forms and upload your paper and that is it.  With
this system, you certainly do not need to worry about how many
copies to submit.  And you will be prompted for all the ingredients
that are needed.  It is a reliable system, and it works.

The journal will assume that the ``communicating author'' is the
person who submitted the paper---unless you explicitly tell
it otherwise.  All further correspondence will be conducted with that
person\index{communicating author} 
at that {\it e}-mail address.  

Even if you submit your paper at a Web site, you will often be asked
to supply a cover letter.
Some authors think that the cover letter is an opportunity to make a
pitch for the paper.  Such an author will fill the cover letter with
fulsome\index{paper!cover letter for}
praise of what is in the preprint, why it improves on the
existing literature, and who might be a suitable referee.  Most
editors will not find such remarks helpful, and many will find them
annoying.  By naming potential referees, you may in fact be ruling
them out in the mind of the editor (since he/she may think that they are
your pals).  Best is to keep the cover letter simple and
dispassionate.

Here is what you can expect after you have sent your paper to a journal.
After a short delay, the journal will notify you that it has received
the paper; this is usually done by {\it
e}-mail. The journal will often assign a manuscript number to your paper,
and will advise you to use this number in all future correspondence. I run
a journal, and I can tell you that this number is valuable. The journal
office can easily misfile a paper with multiple authors;
also,\index{paper!acknowledgement of receipt}
 since the paper is passed from Managing Editor to Associate Editor to one
or more referees, the paper can be misplaced. It helps significantly when
authors and editors use the manuscript number. Such a number might be
``JGEA-D-16-243,'' indicating that this is the $243$rd paper received in
2016 by the {\it Journal of Geometric Analysis}. The {\it e}-mail will
conclude by saying something like ``Don't call us; we'll call you.'' In
other words, you may have to wait a while for the referee's report; so sit
tight.

Expect to wait four to six months for a report. After that wait, you are
well within your rights to send a polite note to the editor to whom you
submitted the paper;\footnote{This is the stage where it is important
that you have the editor's name and the manuscript number at hand.} 
simply state that you submitted the paper on thus and
such a date, received an acknowledgement on another date, and you are
wondering if there has been any progress in the matter. Most editors
appreciate a gentle reminder, and will in turn nudge the appropriate
Associate Editor or referee.

Eventually you will receive a referee's report.  It may be a
paragraph or it may be five pages or more.  It may say ``This paper
is terrific.  Publish it as quickly as you can.''  Or it may say
``This paper is dreadful.  Stay as far away from it as you can.'' 
Most\index{paper!dealing with referee's report}
often it will say something in between these extremes.

If the paper is rejected, then you will have to ply your wares
elsewhere.  A rejection does not necessarily mean that the paper is
bad, or that its results have no value.  Many journals suffer from a
serious backlog, and send most papers back unread (this is, properly
speaking, not a rejection---for the paper has not even been examined
or evaluated); sometimes the editor picks the wrong referee, or a
referee with an ax to grind, or a referee who did not understand the
paper; sometimes the editor misunderstands the referee's report;
sometimes\index{referee!errors by} the referee is just plain wrong. 
Some of my own most influential papers have\index{referee!dealing
with} been treated rather shabbily.  I know even Fields Medalists who
tell horror stories of papers rejected.  One of the secrets to
success in the academic game is perseverance.  If your paper is
accepted the first time around, then congratulations.  If not, you
should try to be objective and figure out why.  Then act
intelligently on that new information.

If your paper is accepted, then the referee will most likely have offered
comments and suggestions. Some referees go so far as to suggest
alternative proofs, different references, or entirely different
approaches. Some editors will instruct you to read the referee's remarks,
make those changes that you wish, and then to send the final version of
the manuscript,\index{editor, dealing with}
labeled ``revised'' and with a new submission date, to the journal; other
editors will explicitly make final acceptance conditional on your
responding in detail to everything that the referee has said. In this last
case, if you want to continue doing business with the journal (you always
have the fallback option of withdrawing the paper), then you are honor
bound to respond to {\it each of the referee's remarks}. The best way to
respond is to treat the referee's remarks one by one, and to record in a
cover letter to the editor a brief description of just what you did in
each instance. In some cases, you may say ``the referee is mistaken and
here is why.'' Or you could even say ``this is a matter of taste and I
respectfully disagree.'' In most instances you can expect the referee's
comments to be accurate and useful and you will probably want to implement
them in some form.

If you feel that the referee has been particularly helpful, then you may
wish to add\index{referee!thanking} a sentence to the paper---alongside
your other acknow- \hfill \break ledgements---saying that you thank the
referee for useful suggestions. You will find it awkward to endeavor to
determine the identity of the referee, so do not plan to mention the
referee by name.

When you are finished with your revision, assuming that a revision is
what has been requested, then make the usual number of copies of the
revised manuscript, mark each of them ``Revised'' and put the date,
and then return these to the editor along with a new cover letter. 
Your new cover letter should state plainly that this is a revision of
a previously submitted paper, that you have responded to the
referee's remarks, and that you consider this to be the final copy.\footnote{Of
course if you are dealing with the journal by way of a Web site then
the procedure is different.  You will submit your revision by uploading
it to a suitable node.  You will also fill out various forms telling
the journal just what you are doing.} 
Please note, however, that the editor {\it might\/} send the revised
paper to the original referee---or to some other referee---again, and
you may be asked to make even further changes. You can expect to
receive an acknowledgement of your new submission, together with a
clear statement of whether this is the end of the road or whether you
will be hearing again from a referee.

Alternatively, you may be able to submit your revision using the Internet.
Doing so, you will receive an immediate acknowledgement of receipt.  You
will probably still have to wait a bit for the editor to render a
final decision.

And now, as the Managing Editor of a journal, I would like to ask you
a favor.  It is quite common for a young mathematician to learn
\TeX\ on the local computer system at his/her university.  That
system will no doubt have in-house macros and fonts that
everyone at University $X$ uses.  The trouble with this
is that the \TeX\ file for the paper will compile on the
computer at {\it that} university, but not on other computers.

So you must learn to make your papers self-contained, so that
all the macros and font calls are in the source code file
for the paper.  Usually, when you first submit a paper electronically,
you will only send in a {\tt *.pdf} file.  This is because
the referee and/or the editor will likely ask for revisions so
there is no sense to send the \TeX\ source code file at that
time.  But, after the paper is accepted, you {\it will} submit
your \TeX\ source code file.  And you want the journal to be
able to compile your file.  So pay attention to the point
made in this paragraph.

After a suitable number of iterations of the procedures just described, you
and the journal will reach some closure. Then you must wait---this wait
could be from six months to two years or more---for the galley proofs or
page proofs of your paper.\footnote{It should be noted, however, that
some journals now instantly publish your paper on the Web as
soon as it is accepted.  Most authors are pleased at this
eventuality.  But it is understood that the archival copy of the paper is
the hard copy version, the one that has undergone copy editing
and other vetting.} These you must\index{paper!galley proofs for}
 proofread meticulously, both for mathematical accuracy and for
typesetting accuracy.  There also will be ``Author Queries,'' noted
by hand, on the proof sheets.  You must respond to each of these.
Often you will handle these matters by email.  Although sometimes
you will be sent a {\tt *.pdf} file and you will annotate that
file with ``electronic sticky notes.''

Again I must make a nod to the Internet.  These days there are some
journals that, once your paper is accepted, will publish the
paper electronically {\it immediately}.  Since the hard copy
version of the journal is published in a queue, there may be
a considerable delay before the printed version of your paper appears.
It is also the case these days that some journals do no copy editing
whatsoever.  Whatever the liabilties of such a policy, it certainly 
streamlines the process.

You will always be asked to turn your proof sheets around rather
quickly---often within 48 hours. You will sometimes be asked to sign a
statement saying that you approve this version of the manuscript going
into print.\index{galley proofs, return of} {\it And you will be asked at
this time how many reprints you want.\/} The truth is that, these days,
reprints are something of an atavism. Typically a journal will give you a
{\tt *.pdf} file of the final version of your article. And then you can
print all the reprints that you may want. After you send back your
response to the galley proofs, then your job is done. Just wait for your
paper to show up in the library, and you will know that your paper
is now part of the permanent archive.

%% Section 2.8

\section{A Coda on Collaborative Work}
\markboth{CHAPTER 2.  TOPICS SPECIFIC TO MATHEMATICS}{2.8.  A CODA
ON COLLABORATIVE WORK}

I have written a great many collaborative papers, and some
collaborative\index{collaborative work}
books as well.  I know others who have never
collaborated.  And there are others still who have collaborated a few
times and would never do so again.  Which characteristics lead to a
successful and happy collaboration and which do not?

First, if you agree to collaborate on a project (and both parties had
better agree at the outset; do not leave this question until the
project is finished!)\ then set aside all questions of priorities.  At
the end of the collaborative process, it is both painful and
inappropriate\index{priority disputes} 
for one author to say ``Well, you didn't contribute
very much.  My name should go first''  or, worse, ``Your name should
not appear at all.''  If, at the end of the first paper, either or
both participants deems the collaboration unsatisfactory, then the
authors\index{paper!order of author names}
can go their separate ways.  But, in my view, an agreement to
collaborate is an {\it a priori\/} contract to see things through to
the end.

Some of my collaborations involve multiple papers; in one case the joint
work amounts to thirteen papers and a two books. Another collaboration of
mine involves six papers and four books. In these monumental
collaborations, both my collaborator and I know that on some papers he/she
contributed more and on others I contributed more. Same for the books. I
can honestly say that neither of us dwells on the matter. Taken as a
whole, we are both quite pleased and proud of the {\it oeuvre}. As it
happens in the first of these collaborations, one of us has lost interest
in this subject area and the other one has pushed on a bit further, either
writing papers alone or in collaboration with others. This has worked out well,
because each of us respects the other.

And\index{author!disputes}
this last point is the real key.  I know of collaborations in
which one author purposely introduced errors into the joint paper in
order to see whether the other author was truly reading the paper or
not.  I know of a collaboration that got to the stage of the paper
being\index{author!mutual respect of} 
submitted to a journal; after a time the authors had a dispute,
and one author unilaterally withdrew the paper and resubmitted it
elsewhere under his name alone.  I know of a collaboration
between two lifelong friends who were developing their twentieth
paper together; they could not agree on whether to call the first
result Theorem 1A or Theorem A1; the matter ended with lawyers, death
threats, and guns brandished in the air (this {\it really\/}
happened).  In all these cases the base of the difficulties was that
the authors did not respect each other.
  
None of the scenarios described in the last paragraph should take
place when you engage in a collaborative effort.  You should enter a
collaboration with the view that this is an adventure and you will each
see what you can derive from it.  Each of you should respect the
other(s),\index{collaboration!as adventure} 
and will take great pains to be courteous and helpful to
the other(s).  The goal is to produce a nice piece of work---{\it not\/}
to squabble over credit, {\it not\/} to argue over whose name should go
first (alphabetical is almost always best), {\it not\/} to argue over
whether future papers will be joint or will be written separately.

This last point can lead to sticky wickets, even between the most
well meaning of participants.  Imagine this scenario:  
Mathematicians $A$ and $B$ write one or more joint papers.  The
collaboration then seems to go into remission.  That is, each author
goes\index{author!disputes} 
his/her own way, and they have little contact for a couple of
years.  Then one of the authors (say $A$) cooks up another idea and
writes a new paper, by himself, which in some sense builds on the
ideas in the old series of joint papers.  Mathematician $B$ gets wind
of this new paper, feels that his contribution to the earlier work
justifies his name appearing on this new paper, and relates this 
feeling to $A$ in no uncertain terms.  Mathematician $A$ feels that
the joint work was long ago and far away.  The main reason for the
existence of the new paper is {\it his\/} new idea---which is due to
him alone.  Mathematician $A$ feels that $B$ has already received
adequate credit for the joint work; no further credit is due $B$.  As
you can imagine, a major fight ensues.

This situation is most unfortunate.  Nobody is right and nobody is
wrong.\index{collaboration!etiquette}  Here is what 
{\it should\/} have happened---in the best of all
possible worlds.  Realizing that this new paper builds on old joint
work with $B$, mathematician $A$ should have phoned $B$ and told him
about it and then said ``I think it would be appropriate for this new
paper to be joint between us.  What do you think?''  Mathematician
$B$, ever the gentleman, should then have said ``Oh no, this is
your idea.  Write the paper by yourself.  You can thank me in the
introduction if you like.''  Having participated in transactions of
this nature, I can tell you that this is a most satisfactory way to
handle the matter.  Typically, mathematician $B$ is not hungry for
another paper; he/she just wants his/her due.  Typically, mathematician $A$
is not anxious to offend $B$; he/she just wants credit for his/her new idea. 
(Of course the human condition is such that there are
always more complex forces at play.  Perhaps $A$ feels that, in the
world at large, $B$ is generally given more credit for the
collaborative work than $A$.  Perhaps $B$ feels that $A$ never pulled
his weight in the first place and therefore $A$ owes $B$. 
Fortunately, this is not a tract on psychology, so I shall not
comment further on these complexities.)   

By touching base with your collaborators 
in a courteous fashion, you can usually avoid
friction.  And the effort is worth it.  To go
to a conference and run into a former collaborator is a pleasure.  If the
relationship is healthy and friendly, then there is plenty to discuss
and the potential for future joint work always lies in the offing. If
instead there is friction and resentment between you and your former
collaborator, then meeting again could be perfectly dreadful.  This is
another skeleton in your closet.  Bend over backwards to avoid such a
liability.

\section{Professional Societies}

Sad to say, young mathematicians these days do not show a great
interest in belonging to our professional societies.  This is a
great mistake.  Our professional societies are the lifeblood of
our profession.  They look out for our best interests, they
promote and nurture the profession, and they provide both
financial and spiritual support for everything that we do.

Just to take an instance, the American Mathematical Society engages
in the following activities (and these are only some of the many):
\begin{itemize}
\item Organizes important national and regional (and sometimes
international) conferences.
\item Publishes a number of frontline research journals.
\item Has several significant book series.
\item Produces {\it Math Reviews} and {\tt MathSciNet}.
\item Has a lobbyist in Washington to promote mathematical
interests.
\item Sponsors a number of important mathematics prizes.
\item Sponsors several fellowships and postdocs.
\item Keeps careful account of mathematical infrastructure and
our role in society.
\end{itemize}
The list goes on and on.  The mathematics profession---not just
in the United States but worldwide---benefits immensely from
the activities of the American Mathematical Society.

Among the professional mathematical organizations in this country
there are:
\begin{itemize}
\item the American Mathematical Society
\item the Mathematical Association of America
\item the Association for Women in Mathematics
\item the Society for Industrial and Applied Mathematics
\item the Society for Advancement of Chicanos/Hispanics 
and Native Americans in Science
\item the American Statistical Association
\item the American Mathematical Association of Two-Year Colleges
\item the Association for Symbolic Logic
\item the National Council of Teachers of Mathematics
\item the Society of Actuaries
\end{itemize}
and there are many more.  

You have probably heard of most of societies adumbrated in the
last paragraph.  But there are smaller and less prominent mathematical
societies that have played significant roles in people's lives.
An example is the Resource Modeling Association.  This organization
only has 200 members.  But it is a real touchstone for the people
involved with it.  Its statement of purpose is this:
\begin{quotation}
The RMA is an international association of scientists working at the
intersection of mathematical modeling, environmental science, and natural
resource management. We formulate and analyze models to understand and
inform the management of renewable and exhaustible resources. We are
particularly concerned with the sustainable utilization of renewable
resources and their vulnerability to anthropogenic and other disturbances.

We hold an annual conference and we publish the quarterly journal Natural
Resource Modeling in association with Wiley-Blackwell.
\end{quotation}

Another small but incisive organization is the River Management Society.
There are a great many of these, and they play a significant role
in people's lives.

I would urge you to join one or more professional societies.  The potential
benefits are immense, and this will certainly enrich your life in
a number of ways.
\vfill
\eject

\hbox{ \ \ \ }

\thispagestyle{empty}

\newpage

%%%%%%%%%%%%%%%%%%%%%%%%%%%%%%%%%%%%%%%%%%%%%%%%%%%%%%%%%%%%%%

%% Chapter 3
\chapter{Exposition}

\begin{quote}
\footnotesize \sl Reading maketh a full man, conference a ready man,
and writing an exact man.
\smallskip \hfill \break
\null \mbox{ \ \ } \hfill \rm Francis Bacon \break
\null \mbox{ \ \ } \hfill \rm {\it Essays\/} [1625], Of Studies
\end{quote}

\begin{quote}
\footnotesize \sl You can fool all of the people all of the
time if the advertising is right and the budget is
big enough.
\smallskip \hfill \break
\null \mbox{ \ \ } \hfill \rm Joseph E. Levine 
\end{quote}

%\begin{quote}
%\footnotesize \sl Never spend more than a year on anything.
%\smallskip \\
%\null \hfill \rm Jeff Ullman 
%\end{quote}

\begin{quote}
\footnotesize \sl
When Kissinger can get the Nobel Prize, what is there left for satire?
\smallskip \hfill \break
\null \mbox{ \ \ } \hfill \rm Tom Lehrer
\end{quote}

%\begin{quote}
%\footnotesize \sl The commonest thing is delightful if only one hides it.
%\smallskip \\
%\null \hfill \rm Oscar Wilde
%\end{quote}

\begin{quote}
\footnotesize \sl
Life should be as simple as possible, but not one bit simpler.
\smallskip \hfill \break
\null \mbox{ \ \ } \hfill \rm ascribed to Albert Einstein
\end{quote}

\begin{quote}
\footnotesize \sl If you have one strong idea, you can't help 
repeating it and embroidering it.  Sometimes
I think that authors should write one novel and then be put
in a gas chamber.
\smallskip \hfill \break
\null \mbox{ \ \ } \hfill \rm John P.~Marquard 
\end{quote}

\begin{quote}
\footnotesize \sl Writing comes more easily if you have something to say.
\smallskip \hfill \break
\null \mbox{ \ \ } \hfill \rm Scholem Asch 
\end{quote}

\begin{quote}
\footnotesize \sl Considering the multitude of mortals that handle the
pen in these days, and can mostly spell, and write without glaring
violations of grammar, the question naturally arises:  How is it,
then, that no work proceeds from them, bearing any stamp of
authenticity and permanence; of worth for more than one day?
\smallskip  \hfill \break
\null \mbox{ \ \ } \hfill \rm Thomas Carlyle  \break
\null \mbox{ \ \ } \hfill \rm {\it Biography\/} (1832)
\end{quote}

%% Section 3.1

\markboth{CHAPTER 3.  EXPOSITION}{3.1.  WHAT IS EXPOSITION?}
\section{What Is Exposition?}
\markboth{CHAPTER 3.  EXPOSITION}{3.1.  WHAT IS EXPOSITION?}

Perhaps the highest and purest form of mathematical writing is the
research paper.  A research paper, in its best incarnation,
contributes something useful and insightful to our collective
mathematical\index{exposition!importance of}
knowledge.  If it is very good, then the contribution may
live for a long time.  The creation and publication
of research is what mathematics is all about.

But our profession involves us in other types of writing.  We must
write letters of recommendation.  We must write referee's reports.  
We must review cases for tenure and promotion.  We write
surveys.  We sometimes write book reviews.  We may be called on to
write opinion pieces.  The present chapter concentrates on such
{\it expository writing}.\footnote{We save a discussion of bookwriting
for Chapter 5.}

In its simplest form, mathematical exposition could take the form
of a survey of a field on which you are an expert.  Or it could
be\index{exposition!what is} 
a text or monograph on some specific area of mathematics.
The new challenges present in such a writing task are these:
{\bf (i)}  you are attempting to reach a broader audience than that
which would read one of your research papers; {\bf (ii)}  you must
strike a balance between how much mathematical detail to give and 
how much explanation and/or handwaving to provide; {\bf (iii)}  you
must be open to the idea that this is a new type of writing with
new goals and new audiences.

The reader of an expository article does not want to work as hard as
the reader of a research article.  Envision your reader sitting on a
park bench reading your expository article, or putting his/her feet up
and drinking a cup of coffee while reading.  {\it Do not\/} imagine
your reader with a pencil gripped in his/her fist, slaving away over each
detail of your paper.  Thus, if you are writing an article about the
influence of the Atiyah-Singer Index Theorem on modern mathematics,
you certainly will not prove the theorem.  To be sure, you will refer
to some of the excellent books on the subject.  You will explain how
the result is a far-reaching generalization of de Rham's theorem and the
Riemann-Roch theorem.  You will describe the ingredients of the
proof, and will give a rough sketch of its structure.  
But you will certainly not {\it prove\/} the theorem.

You also will not assume that your reader already knows all the jargon
in\index{exposition!level of difficulty in}
the subject.  You will not assume the reader to be
expert in $K$-theory or pseudodifferential operators.
Nor will you assume that your reader is familiar
with the motivation for, and the applications of, 
the subject.  {\it You should not assume that your reader
has the perfect background to read what you are writing.}

So you have your work cut out for you.  Expository writing is a lot
like teaching.  You frequently must anticipate your audience's
shortcomings, and make suitable adjustments in your presentation. But
in expository writing you must be smarter than you are when you are a
calculus teacher.  In the latter situation, your audience is before
you and is sending you signals.  When you are writing, your audience
is (if you are lucky) only in your head.

%% Section 3.2

\section{How to Write an Expository Article}

For the purposes of this section, the phrase ``expository article''
means a survey article.  Such an article might be a survey of some field on
which\index{exposition!how to write} 
you are expert.  Perhaps you are one of the pre-eminent
experts, and therefore the canonical person to be writing
this survey.  Or your colleagues have called upon you to plant a
flag for the subject.  One scenario is that you have been invited to
give a one-hour talk at a national meeting of the American
Mathematical Society.  It is natural, and commonplace, for you
to turn such a talk into a survey that will appear in the {\it
Bulletin} of the AMS.  But such an august occasion is not a necessary
condition for the writing of a survey.  You may simply feel that a
survey is needed, that certain wrongs need to be righted, or that the
time is ripe.  More and more journals are soliciting expository
articles; there is a market for high-level exposition done well.

In order to write a good survey article you will need a
detailed outline before you.  You will be covering a lot of ground,
and you do not have the luxury of hiding behind the details of the
proofs.  In fact, you will most likely not be presenting any proofs
in their entirety.  When you do choose to present a proof, it will
probably\index{outline, importance of} 
be a sketch, or a pseudo-proof.  If you are clever, you
will present a well-chosen example, work it through, and then
say ``this reasoning also shows you how the proof works.'' Or you
might\index{pseudo-proofs} say 
\begin{quote} This example is in fact the enemy.
The proof shows that this example represents the only 
thing that could possibly go
wrong, and then systematically\index{exposition!use of examples in}
shows that the hypotheses rule the
example out.
\end{quote}

A good survey should build to a climax:  Poincar\'{e} looked at
this example, Lefschetz looked at that example, eventually people realized
what was going on, and the Eilenberg-Steenrod axioms were formulated
(I am thinking here of the genesis of algebraic topology).  
Alternatively:  first there was the Laplacian, then there was
general elliptic theory, then there was the $\overline{\partial}$-Neumann
problem, then pseudodifferential operators evolved.  Simply
to begin citing technical results in chronological sequence is not
to\index{survey!writing of} write an effective survey.  You are telling a story, and you
must create a tapestry.

A good survey should have a stirring conclusion.  By this
I do not mean ``That's all, folks!!'' or a hearty cry for
more and better research on Moufang loops.  Instead, your survey
should conclude by taking a look back at what ground has been
covered and where the subject\index{survey!conclusion of} 
might go in the future.  It should
note the historical turning points (which you have, I hope, described in
the body of your survey), and make speculations about what
the future milestones might be.  If appropriate, it should
sum up what this subject has taught us so far, and what
it might show us in the years ahead.

A good survey should have an extensive bibliography.  You are
not doing your job if you merely say ``The three standard books in the
subject\index{survey!bibliography for} 
are these, and you should look in their bibliographies for
all the technical references.''  By all means mention the three
standard treatises, and extol the virtues of their bibliographies. 
But you {\it must\/} create your own bibliography.  Your list of
references is your detailed definition of what the subject is, what
are the most important papers, and what is the latest hot stuff. 
Compiling a good bibliography is a lot of work (though, with the aid
of modern technology, not nearly so arduous as in years past---see
Sections 2.6, 5.5).  But this effort is a necessary part of the process, and the
result will be a valuable tool for both you and your colleagues in
the\index{survey!how to write} years to come.

Writing a good survey---one to which people will refer for many
years---is one of the hardest writing tasks there is.  Getting all
the basic ideas on paper, and in the proper order, is just the
first step.  Once that task is completed, then you must craft the
piece into a compelling tale with introduction, entanglements,
climax, denouement, and finale (much as in a Shakespearean tragedy). 

{\it Be absolutely certain that you have not slighted any of the
players in the subject, nor inadvertently misrepresented their
contributions.}  Almost certainly, in the course of writing
your survey, you will be saying (perhaps {\it sotto voce}) ``Here is
the right way to see things (implying, perhaps, that
some others are not\index{survey!giving credit in a} 
the right way) and here are the important
contributions.''  Whether you do this consciously or not, 
you certainly will do it.  Take extra care that you are
diplomatic, and that you let everyone's voice be heard.

It is generally a good idea, when you write a survey, so send it around to
all the key players in the field {\it before publication}. Give them a
chance to respond to it and offer ideas. Most scholars are pleased to be
given such an opportunity.

There was something of a mathematical mini-crisis a few years
ago when an important mathematician wrote a survey of a topic
in harmonic analysis.  He made the mistake of {\it not} showing
it to people before it appeared in the {\it Notices of the AMS}.
In those days the {\it Notices} had a policy that no article
could have more than ten references.  Well, you can certainly
imagine that a survey will have many more than ten references.
So, when people saw the published article, they were alarmed
and offended.

What the author {\it should} have done is this:  First, he should
have circulated the article {\it before} publication.  Second,
he should have put a line in the article saying that he
was restricted to just ten references, but that a complete
list of references appears on thus and such a Web site.  That
Web site could also contain ancillary examples, figures, and other
material that would enhance the article.  This {\it modus operandi}
would have solved the problem and left everyone happy.

The reader of your survey should come away from it feeling that he/she
has been given an entr\'{e}e to some new mathematics.  He/she should
have {\bf (i)} learned some new facts and {\bf (ii)} seen some
new techniques, {\bf (iii)} learned some interesting history.  
Forty pages of descriptive prose, without any
substance, will not wash with a mathematical audience.  You must
sketch how the ideas unfold, and endeavor to give some
indications of the proofs.  When deciding what to prove, you must
balance\index{survey!purpose of} 
what is instructive against what is feasible in a short
space.   Often you can prove only a special case, or you might say 
``to simplify matters, we add some hypotheses.''  As an instance,
proving the inverse function theorem for a $C^1$ function is hard. 
But if you assume that the function is $C^2$, then you can use the
remainder in the Taylor expansion to good effect and the proof
suddenly becomes easy (see the details in [Kr1]).  The key ideas are
still present, and they come out much more clearly.  

Just as when giving a talk, you can fudge a bit in your survey.
State theorems precisely and correctly but, when presenting the
proof, say ``For simplicity, we consider only a special case''
or ``For a quick proof, let us assume that the function is actually
real\index{survey!imprecision in} 
analytic.''  Readers will appreciate being given a nugget
of knowledge, without the gory details.  A research paper should
contain complete proofs---proofs that are categorical, and leave
no doubt of their correctness.  A survey acts as a pointer to the
research literature; it is not usually the final word on the proofs.

%% Section 3.3

\section{How to Write an Opinion Piece}

Most mathematicians agree that writing good exposition is
considerably more difficult than writing good mathematics.  As has
already\index{opinion piece!how to write an}
been described in this book, the latter activity makes few
demands on your abilities as a creative writer.  You need only
exercise some taste in organizing and presenting the ideas.  However,
when you are expositing, then you are less engaged in statement and
proof and more engaged in description, explanation, and opinion
formation.  There is much more latitude and therefore you, as a
writer, must exercise more control.

Let us turn now to the writing of an opinion piece.  The writing of
an effective opinion piece will involve all the skills already noted
as mandatory for good expository writing.  But the opinion piece must
also, if it is to be effective, have fire and life and drive.  It
must capture the reader's attention, and it must {\it convince\/} him/her
of something.  How does one go about pulling this off?

A parody of midwestern political oratory has the would-be congressman
declaiming
\begin{quote}
Agriculture is important. \hfill \break
Our rivers are full of fish.  \hfill \break
The future lies ahead.
\end{quote}
I hope that, when you write your opinion piece, your thoughts
have more focus than this politician's, and your message is
more incisive and more substantive.  

First, to repeat one of the main themes of this book, 
if you are going to write an essay expressing an opinion then you
must have something to say.  And you must know clearly and
consciously\index{opinion piece!having something to say in an} 
what that something is, and how you propose to formulate
it and to defend it.  It does not wash to say, in your mind's voice,
``I am going to write an essay in opposition to the teaching of
calculus in large lectures because it is a bad idea and I hate it.'' 
In point of fact most of us agree with this thesis, but it turns out
to be extremely difficult to harness facts and arguments to support
the thesis. Couple this lack with the fact that there are articulate and
serious people---who spend their professional lives studying such
matters---who disagree vehemently with the thesis (see [Dub]) and can
marshal forceful arguments against it, and your obvious essay in
support of an obvious contention suddenly becomes quite painful.  

In fact the statistics that bear on the ``large lecture'' question
are a bit unsettling.  They tend to suggest that students taught in
small classes feel better about themselves and about the subject
matter (than do students taught in large classes); they do {\it not\/}
tend to suggest that such students will turn in a better
performance.   That is the trouble with facts:  they
sometimes force you to conclusions that differ with your intuition.

What I am suggesting in the preceding paragraphs is that the writing
of a position paper or an opinion piece often involves considerable
research.  This is not the same sort of research that one performs in
order to prove the Riemann hypothesis.  But it needs to be done, and\index{Riemann hypothesis!proof of}
done\index{position paper, research for} 
thoroughly.\index{opinion piece!research for}
There is no substitute for knowing what you are
talking about.

You must have a deliberate and explicit formulation of your thesis
and your contentions.  Best is to enunciate that thesis in your
first paragraph.  The thesis\index{opinion piece!thesis of an} 
could constitute your first
sentence, or it could be the culmination of suitable background
palaver that lays the history and orients the reader's mind toward
the main point of your essay.  But the thesis you are defending
should be put forth---so that it cannot be mistaken---at the
outset of your opinion piece.

The next (major) portion of your essay should consist of cogent
presentation of material gathered in support of the previously enunciated
thesis. This prose could include facts, reports of studies, anecdotes,
logical arguments,\index{opinion piece!arguments in} and other materials
as well. Note that a defense consisting entirely of anecdotes is at first
entertaining but, in the end, not convincing. On the other hand, a defense
consisting only of dry facts and logical arguments does not generally hold
the reader's attention and is not forceful. (The validity of this last
statement depends, of course, on your audience, on your subject matter,
and on the context. Obviously, most any mathematical research paper
contains just facts and logical arguments. But such a recondite exercise
is directed toward specialized researchers with an {\it a priori\/}
interest in what the writer has to say. The audience for an expository
paper or opinion piece is more diffuse, less well prepared, and less
patient.) 

The last portion of your position paper should sum up the major points you
have made, repeat the most important ideas, and force the desired
conclusion. These are the final thoughts with which you will leave the
reader. They are analogous to the closing arguments in a
jury\index{opinion piece!summation in} trial. Weigh each word carefully.
Remind the reader what he/she has read, and why.\index{opinion piece!conclusion of}

I do not mean to suggest that persuasive writing is formulaic. This
activity is not like laying bricks. Some of the best position papers
conform only loosely, if at all, to the rubric just laid out. But the
points I have made, and the issues I have raised, are salient to any
polemic, no matter what its exact form.

%% Section 3.4

\section{The Spirit of the Preface}

Many writers spend little time in writing a preface; some forget to
write it at all.  This is a mistake.  Your prefatory remarks are
often the most important part of your writing.\index{preface!what is a}
They tell the reader why you write what you write, what your goals are,
and what you intend to accomplish. They state what you assume, and what
you conclude. These principles apply whether you are writing a book (which
has a separate, formal preface) or an article (which may have a prefatory
section, or collection of paragraphs) or a letter (which may have just one
prefatory paragraph). The preface is your statement of purpose. It is
vital to your mission.

When\index{preface!purpose of the} an editor at a publishing house
receives a mathematical manu- script (for a book, say) for consideration,
he/she usually seeks advice from one or more experts in the field. When I
am asked by a publisher to review such a manuscript, the first thing I
look for is a Preface or Prospectus (a marketing version of the Preface),
and a Table of Contents. These two items, if present, will give me a quick
overview of the project: What material is covered? At what level? Who
constitutes the intended audience? What are the prerequisites?
What\index{preface!importance of the} need does this book fill? What are
the book's competitors? (In modified form, these queries also apply to an
expository article. If I receive a 100-page expository article to review,
then I hope that it---or at least the cover letter 
that accompanies it---\index{preface!questions that it answers} contains the
information that a book Preface and Table of Contents usually provide.)

Without this information, I have no idea what I am reading.  The
manuscript could start out with sophomore-level differential
equations, and before long be doing canonical transformations
for Fourier integral operators.  As a result, I have no
idea what the author is trying to accomplish.

Whether you are writing a research announcement, a research paper, an
expository paper, a book, or virtually anything of a scholarly nature, you
should always ask yourself the questions in the second paragraph in this
section. Most importantly, you must decide {\it in advance\/} the book's
intended audience and you must, at all times, keep that group clearly in
focus. If you are writing a calculus book, then presumably the audience is
freshmen, and therefore you must resist the temptation to indulge in
asides to the professor. If you are writing a research article, then
presumably the audience is fellow researchers, like yourself---not Gauss
and God. If I may be permitted a little hyperbole, I will say
now\index{sense of audience} that having a strong sense of your audience
is the single most important attribute of an effective writer.

I have sung the praises of the Preface. But the Table of Contents (known as
the TOC in the publishing industry) is nearly as important. Here is why.
When you are writing a book, which is a big project, you should have the
entire scope of your endeavor firmly planted before\index{Table of Contents!importance of} your\index{TOC} mind's eye. In this way you can measure
your task, you can see what progress you are making, and you
can\index{Table of Contents!role in the writing process} keep the affair
in perspective. Sometimes you can have fun just sitting down and starting
to write, or just seeing where your thoughts will lead you, or modifying
your project every time an interesting new preprint comes across your
desk. But let me assure you that these methods are a sure way to guarantee
that your book will never be completed. Writing the TOC addresses this
impasse.

Now you may not be interested in writing something like a book. Such a
writing project is awesome and onerous; the task is\index{book!writing of
a} not for everyone. But the principles in the last paragraph apply even
to writing a twenty-five-page research paper. You need to have the full
scope of the paper in your mind so that you can endow your working methods
with a pace and give yourself a sense of incremental accomplishment. This
sort of organization is also just a simple device for keeping yourself
from becoming completely confused. 

I often begin a book by writing the Preface, because it helps me to
\noindent organize my thoughts and to orient myself toward the project. I
refer to it frequently as my work on the project progresses. At the same
time, it also makes sense to write the Preface last: for when the book is
complete, then you know in detail what you have written and\index{preface!as writing compass} you can describe it lovingly to the reader. My
recommendation is to do both. Write a version of the Preface before you
begin the book. When you are finished, write it again.

This section on the Preface may seem like a digression, but it is not. Even
if you are writing an opinion piece (Section 3.3), or a letter of
recommendation (Section 4.1), or a book review (Section 4.2), your piece
should contain prefatory remarks. Such remarks are good for your reader,
so he/she knows what this piece of writing is supposed to be about. But,
most importantly, they are good for you: they keep you honest, and keep
you on your course.

%% Section 3.5

\section{How Important Is Exposition?}

There are those who will argue that mathematical exposition is not
important at all. The one true pursuit, whose fruits are recognizable and
of lasting value, is mathematical research. Some would\index{exposition!importance of} go so far as to claim that even the writing and publishing
of research results is an activity suitable only for hacks. Instead, for a
good mathematician, it suffices to prove the theorems, tell at least one
other person about them, and then let the word spread.

I happen to think that the attitudes described in the preceding paragraph
are counterproductive. Scholars are not monks. They are an active and
engaged component of society. Along with the universities and research
institutes, they are the vessels in which mankind stores its accumulated
knowledge and civilization. Thus scholars must communicate. They can do so
by giving lectures; the importance of lectures cannot be overemphasized.
But scholars also must write. The written word---unlike the spoken
word---lasts for ages, and can influence generations.

Exposition is important because it reaches a broader audience than do
specialized research articles. Thus good expository articles disseminate
information quickly, and they are much more likely to spawn collaboration
between different fields than are specialized articles.  An outstanding
expository article will cause even the experts to reorganize the subject
in their own minds.

In my own work, I have found that expository writing is a device for
teaching myself.  It forces me to organize my thoughts, and to be
sure that I understand how a subject is constructed---from the ground
up.  This is also a device for teaching my students: after I have
explained the same idea several times to several different students
(perhaps over a period of years), I find it useful to write something
down.  Then, when the next student comes along, I can give him/her
something to keep.

I can still remember, many years ago, reading an article by Freeman Dyson
called ``Missed Opportunities'' [Dys]. I have never seen anything like it
before or since. In this article, the author makes statements like the
following: ``In 1956, $X$ proved this and in 1957, $Y$ discovered that. If
only I had been alert, I could have combined these results with ideas of
my own and with Heisenberg's uncertainty principle and\index{Dyson,
Freeman} I could have done thus and such. Instead, $Z$ combined the ideas
in a different way; for this work he later won the Nobel Prize.''

I found this article to be an inspiration in several respects. First, I was
amazed that anyone could understand his subject so well that he could
recombine its parts in ways that never actually occurred. Second, I was
amazed by Dyson's candor. Third, Dyson helped me see what creativity is.
Fourth, he gave me a sense of the scope of knowledge.

Dyson's article is a sterling example of what good exposition can do. 
In a lifetime, you probably will not read more than half a dozen
articles that are this good.  But half a dozen is enough.  If you are
truly fortunate, and extremely talented, then perhaps you will write
one.

%%%%%%%%%%%%%%%%%%%%%%%%%%%%%%%%%%%%%%%%%%%%%%%%%%%%%%%%%%%%%%%

%% Chapter 4
\chapter{Other Types of Writing}

%\begin{quote}
%\footnotesize \sl Not of the letter, but of the spirit:  for the
%letter killeth, but the spirit giveth life.
%\smallskip \\
%\null \hfill \rm The Holy Bible, the New Testament \\
%\null \hfill \rm The Second Epistle of Paul the Apostle \\
%\null \hfill \rm to the Corinthians.  Chapter 3, Verse 6.
%\end{quote}

\begin{quote}
\footnotesize \sl Stand firm in 
your refusal to remain conscious during algebra.  In real
life, I assure you, there is no such thing as algebra.
\smallskip \hfill \break
\null \mbox{ \ \ } \hfill \rm Fran Lebowitz  
\end{quote}

\begin{quote}
\footnotesize \sl Neither can his Mind be thought to be in Tune, whose
words do jarre; nor his reason in frame, whose sentence is
preposterous.
\smallskip \hfill \break
\null \mbox{ \ \ } \hfill \rm Ben Jonson  \break
\null \mbox{ \ \ } \hfill \rm {\it Explorata---Timber},   \break
\null \mbox{ \ \ } \hfill \rm or Discoveries Made upon Men and Matters 
\end{quote}
\begin{quote}
\footnotesize \sl The flabby wine-skin of his brain \hfill \break
Yields to some pathological strain, \hfill \break
And voids from its unstored abysm \hfill \break
The driblet of an aphorism
\smallskip \hfill \break
\null \mbox{ \ \ } \hfill \rm The Mad Philosopher, 1697 \break
\null \mbox{ \ \ } \hfill \rm in {\it The Devil's Dictionary}  \break
\null \mbox{ \ \ } \hfill \rm \quad by Ambrose Bierce
\end{quote}

\begin{quote}
\footnotesize \sl What is written without effort is in general
read without pleasure.
\smallskip \hfill \break
\null \mbox{ \ \ } \hfill \rm Samuel Johnson 
\end{quote}

\begin{quote}
\footnotesize \sl Sometimes a cigar is only a cigar.
\smallskip \hfill \break
\null \mbox{ \ \ } \hfill \rm Sigmund Freud
\end{quote}

\begin{quote}
\footnotesize \sl The good writing of any age has always been the
product of {\it someone's\/} neurosis, and we'd have a mighty dull
literature if all the writers that came along were a bunch of happy
chuckleheads.
\smallskip \hfill \break
\null \mbox{ \ \ } \hfill \rm William Styron   \break
\null \mbox{ \ \ } \hfill \rm interview, Writers at Work (1958)
\end{quote}

\begin{quote}
\footnotesize \sl Close your eyes and think of England.
\smallskip \hfill \break
\null \mbox{ \ \ } \hfill \rm ---a Victorian mother, giving advice to her daughter \break
\null \mbox{ \ \ } \hfill concerning behavior on the wedding night.  
\end{quote}

%% Section 4.1

%Cherries Jubilee

\markboth{CHAPTER 4.  OTHER TYPES OF WRITING}{4.1.  THE LETTER OF RECOMMENDATION}
\section{The Letter of Recommendation}
\markboth{CHAPTER 4.  OTHER TYPES OF WRITING}{4.1.  THE LETTER OF RECOMMENDATION}

Once you have become an established mathematician, you are likely to be
asked for letters of recommendation. Such a document could be a letter of
recommendation for a tenure case, or for a promotion, or
for\index{recommendation, letter of} both. It could be a letter
recommending a young person for a first or second job. It could be a
letter recommending a senior person for an endowed Chair Professorship, or
for the Chairmanship of a department. (For the sake of this discussion, I
will call these ``professional letters.'') It also could be a
letter\index{recommendation letter!professional} of recommendation for a
student (such letters are treated a bit differently from professional
letters---see below). There are many variants; here I would like to
distill out some unifying principles on writing letters of recommendation.

When you are asked to write a letter as described in the first paragraph,
you are in effect set a task. You have become a one person ``taskforce.''
What makes a taskforce different from a committee is that a taskforce is
not supposed to debate the task at hand; instead, the taskforce is
supposed to perform the designated task. In the present instance, you are
supposed to offer your professional opinion on a certain matter.

In my view, it is both unprofessional and irresponsible to dodge the
assigned task. Let me be more precise. There certainly will arise
circumstances where you either cannot write or should not write. Perhaps
you have had a fight with the candidate in question and feel that you
cannot offer an objective opinion; perhaps you have a conflict of
interest; perhaps you are unfamiliar with the general area in which the
candidate works; or perhaps you do not know the candidate well at all. In
any of these cases, or similar ones, you should quickly and plainly write
to the person (the dean or chairman) who requested the letter and say that
you cannot write. Best is if you can give the reason, but it is acceptable
if you cannot. Do not agonize over the task for six months and {\it
then\/} decline to write; take care of the matter right away.

The circumstances described in the last paragraph should be considered to
be extreme exceptions. They will come up less than one percent of the
time. In most instances, you will be asked to write about some particular
person for some particular circumstance and you should say ``yes'' and
then you should do it.

I know mathematicians who will agree to write an important letter and then
not do it. This paradox usually occurs for one of two reasons: {\bf (1)}
the putative letter\index{recommendation letter!reasons not to write}
writer is pathologically disorganized and forgets, or {\bf (2)} the
putative letter writer has nothing nice to say about the issue or person
at hand and does not want to say it. I have already addressed the second
of these {\it conundra}. The first of these situations is not likely to
arise if the request to write was submitted to you as a formal
letter---from a dean, for instance. For then the piece of paper is sitting
somewhere on your desk and you will probably get to it eventually. The
paradox {\it can} occur if instead a student pokes his/her head in your
door and asks for a recommendation to graduate school. You give a cheery
``yes'' and then the entire matter vacates your head. To avoid this error,
ask the student to put the following on a slip of paper: his/her name, any
classes he/she took from you or other pertinent data, and the address to
which the letter is to be sent. (This ruse also helps you to avoid the
embarrassment of having to ask the student's name.) Now you have it in
writing. Also ask the student to come back in a week or so and remind you.
I usually find it convenient to write the letter right away
(if\index{recommendation letter!not writing} the student has poked
his/her head in the door then it is likely my office hour and I might as
well be doing {\it something}). For once the request has been tendered, I
am probably already thinking about what I am going to say; I might as well
write it down and be done with it.

Some clever people create a Web page for students who want a letter of
recommendation.  Then, when a student comes to you and asks for a letter,
you tell him/her to fill out the Web page.  This Web page asks
the student for all sorts of pertinent information---about courses taken,
grades received, personal interactions, and so on.  If you construct
this Web page carefully, then the student in effect writes the letter
for you.

Make a point of writing requested letters in a timely fashion.  It
is the professional thing to do, and you would appreciate
such consideration if the letter were about you.

Having decided to perform the task---that is, to write the requested
letter---you must do what you have been asked to do. That is, you must
formulate an opinion, state it clearly, and defend it. The standard format
will be explained below.

In the first few sentences, state plainly the question that you
are addressing.  For example:

\begin{quote}
\indent The purpose of this letter is to support the tenure and promotion
of Zoltan A.~Beelzebub.   I have known the candidate and his work
for\index{recommendation letter!introductory paragraph}
a period of six years, and have been
impressed with his originality and his productivity. I indeed
think that tenure and promotion are 
appropriate.  My detailed remarks follow.  
\end{quote}

\noindent Alternatively:

\begin{quotation} 
\indent You have asked for my opinion on the tenure, and promotion to
Associate Professor, of Dr.\ Aloysius K.~Foofnar. Dr.~Foofnar is now six
years from the Ph.D., and in that time has produced nothing but some
rotten teaching evaluations and a letter to the editor of the {\it
Two-Year College Math Journal}.  Based on that track record, my
opinion is that he is worthy of neither tenure nor of promotion.  
\end{quotation}
The bulk, or body, of the letter follows, and it should support in
detail the thesis enunciated in the first paragraph.  I shall comment
below on what might constitute that support.  First, let me conclude
these beginning thoughts.

Once the body of the letter is written---and this could comprise one
or two (or even more) pages---then you must write a concluding paragraph.
You {\it must\/} write it.  You must sum up the point you have made, and
restate\index{recommendation letter!summation of} 
your thesis.  A sample of this practice is

\begin{quotation}
\indent In view of the stature of Laszlo Toth
%he busted up the statue of the Pieta with a sledge hammer
in the field of computational algebraic geometry, and considering
his accomplishments as a teacher and as a scholar, I can recommend
him without reservation for promotion and tenure in your department.
\end{quotation}

\noindent (I am assuming that you have in fact described
Toth's status and accomplishments, in a favorable manner, in
the preceding paragraphs.)  Another possibility is:

\begin{quotation}
\indent In sum, I feel strongly that Seymour Schlobodkin should not be 
promoted or tenured.  Indeed, I cannot imagine the circumstances
in which such a move could be considered appropriate.
\end{quotation}

There are those who, although experienced letter writers, do
not adhere to the general scheme just described.  One of the standard
rationales for this behavior is that, in many states and at many
institutions, it is (theoretically) possible for the candidate to
have access to the complete text of his/her letters of 
recommendation---including the identities of the writers.  If such is
the case, then the soliciting school will inform the writer at the
time the letter is solicited.  Of course the letter writer is
offered the option up front of declining to write if he/she is
uncomfortable with this ``freedom of information'' situation.  

There are those who, still uncomfortable, agree to write but are afraid to
say anything. The most negative thing that they are willing to do is to
``damn with faint praise.'' Not only does this artifice undercut the
responsibility\index{damning with faint praise} of the letter writer, but
it puts on those evaluating the case the onus of trying to figure out what
the writer was trying to (but did not) say. In the best of all possible
circumstances, someone at the soliciting institution will phone the letter
writer and just {\it ask\/} what the letter was meant to say. In the worst
of circumstances, the evaluators are left to guess what was meant. Given
that someone's life and career are in the balance here, it is a genuine
shame for such a circumstance to come to pass.

Enough preaching. I will now give some advice about the body of the letter.
If you want your (professional) letter to have some impact, and to be
taken seriously, then you must do two things:
\begin{enumerate}
\item[{\bf (i)}]  make some specific
comments about specific work or
specific\index{recommendation letter!specificity in}
\index{recommendation letter!describing scholarly work 
in} papers\index{recommendation letter!binary comparisons in} of the candidate, 
\item[{\bf (ii)}]  make binary
comparisons. 
\end{enumerate}

You may also wish to discuss other qualifications of the
candidate. No matter what these may be, you should heed these principles:
be {\it precise}, speak of {\it particular\/} attributes, and speak only
of those topics of which you have {\it direct knowledge}. Now let me
explain.

Your letter had better say more than ``Ahmenhotep Smith
is a hail fellow, well met.  Give him whatever he wants.''  First,
such a letter does not say anything.  Second, given the circumstances
described above, in which some letter writers attempt to avoid
litigation by ``damning with faint praise,'' such a vague letter could
be construed as {\it sotto voce\/} damnation.  If your comments are
instead detailed and specific then it is difficult for people to
misconstrue them.  

Thus you should dwell, for a page or
more, on specific virtues of the candidate's scientific work.
Make detailed remarks about specific papers:  Why is this result 
important?  How does it improve on earlier work?  How does the work
advance the field?  Who else has worked on this problem?  This
material should not be a self-serving introspection.  Remember that
most of the readers of the letter will be nonspecialists. Many,
including the dean and members of his/her committee, will not even be
mathematicians.  Thus attempt, briefly, to give background and
motivation.  Drop some names.  For example, say that Ignatz of MIT
worked\index{recommendation letter!scientific work described in} 
on this problem for years and obtained only feeble partial
results.  The candidate under review murdered the problem.  If
appropriate, point out that the candidate publishes in the {\it
Annals} and {\it Inventiones}---and that these are eminent,
carefully refereed journals.

It is astonishing, but true, that even highly placed people,
who write dozens of influential letters every year, seem to be
unaware of the need for binary comparisons.  To put it bluntly, an
important letter that is to have a strong effect {\it must\/} compare
the candidate being discussed to other people, of a similar age and
career level, at other institutions.  The comparison should be with
people---preferably other academic mathematicians---whose names the
informed reader will recognize.  Thus, if the candidate is an
algebraic geometer and you say in your letter that ``this candidate
is comparable to Mumford when Mumford was the same age'' then most
algebraic geometers will know exactly what you mean and will be
extremely impressed; they will in turn explain to their colleagues
the significance of your remarks.  If instead you say ``this candidate
is comparable to Prince Charles when Charlie was a student
at Gordonstoun'' then nobody will know what you are talking about---and
you can be sure that they will not be impressed.

To come to the point, if you are writing an important letter that you
want people to notice, then you must say something like

\begin{quote}
The five best people under the age of 35 in this area are
$A$, $B$, $C$, $D$, $E$.
\end{quote}

\noindent In the best of all circumstances, the candidate under
consideration in your letter is one of $A$ through $E$---and
you should point out that fact.  Alternatively,
you could say

\begin{quote}
Two of the best people in this field, at the beginning tenure stage,
are Jones and Schmones.  Candidate Bones fits comfortably between
them.  Bones is surely more original than Schmones and more 
powerful than Jones.
\end{quote}
Or you could say that the candidate falls into the next group.  Or
that the candidate is so good that it would be silly to compare him
to\index{recommendation letter!binary comparisons in}
the usual five best.  Say what you think is appropriate.  But {\it
say something}.  If you do not, then the readers will notice the
omission and infer that, between the lines, you are saying that this
guy is not any good.   Better to say that he/she is number 15 than to
say nothing at all.\footnote{A {\it caveat\/} is in order if the
letter that you are asked to write is {\it not\/} solicited from a
research institution.  If the candidate is in fact at a four-year
college, where the primary faculty activity is teaching, then the
school probably demands a lot of classroom activity---and not so 
much scholarship.  These days, almost every school wants
its permanent faculty to have some sort of academic profile; but a
teaching college can hardly expect its instructors to stand up to 
hard-nosed binary comparisons.  The lesson is this:  read the
soliciting letter carefully; speak to people in {\it their
language}, and tell them what {\it they\/} want to know.  If the
soliciting letter is from a teaching institution, then it is probably
most appropriate for you to write about teaching, curriculum,
publications in the {\it Monthly}, and letters to the editor of {\it
UME Trends}.  A disquisition on Gelfand-Fuks cohomology is probably
less apropos.}

Tailor your binary comparisons to the circumstances.  It would be
inappropriate to compare a candidate two years from the Ph.D.\ with a
sixty-year-old member of the National Academy of Sciences (unless the
comparison is favorable, and you are trying to knock the reader's
eyes out).  It would probably be inappropriate to compare {\it
anyone} with Gauss (although I {\it have\/} seen favorable comparisons
with Gauss!).\index{recommendation letter!binary comparisons in}   
Note also that, if you are recommending a senior
person for (just as an instance) an honorary degree, then binary
comparisons might be entirely out of place, and uncomfortable as
well.  If the person is already a Chair Professor at Harvard, then to
whom will you compare him/her?  And to what end?

Your letter of recommendation can contain other specifics and details
that might grab the reader's attention.  You could say that the candidate
gives excellent talks at conferences.  You could say that he/she is a
wonderful collaborator.  You could say that the candidate has
beautiful insights, and that talking mathematics with this person is
a\index{recommendation letter!things to say in} 
pleasure.\footnote{I saw one letter of recommendation, by a very
famous mathematician about another famous mathematician, that said,
``Talking mathematics with $X$ is like talking to Enriques.''  This
written by someone who was too young to have ever met Enriques.}
You could describe in glowing and heartfelt terms the
process of proving a theorem, or of writing a paper, with the
candidate.

These days, credible evidence that the candidate
is a good teacher will certainly help the case.  Of course
you are probably not in the same department as the candidate, so you
very well may not be able to discuss his/her teaching.  If the candidate is a
truly outstanding teacher, then perhaps you have heard his/her colleagues
mention his/her talents, or perhaps you know that he/she has won a teaching
award.  It makes quite an impression on letter readers if Professor
$A$, from University $X$, can comment knowledgeably and in detail 
on the teaching of Professor $B$ from College $Y$.

Here are some travesties that I have seen (all too frequently) in 
letters\index{recommendation letter!mistakes in} 
of recommendation.  You should certainly not emulate any
of these mistakes:

\begin{enumerate}
\item The writer begins in one of the fashions indicated above.
Then he/she says       

\begin{quote} Nefertiti Prim has proved the following 
theorem about pseudographs (state the theorem).  This is a nice
result.  The theorem is
based on some old ideas of mine.  [{\it And the rest of
the letter consists of a description of the letter writer's own work!}]
\end{quote}

\noindent Such a letter violates all the precepts laid out above, and
marks the writer as a thoroughly self-absorbed fool.  Of course
this letter does nothing to help, nor to hurt, the candidate; but it
gives a rather poor impression.

\item  The writer discusses the candidate, discusses the candidate's
work, makes binary comparisons, and mentions specific papers.  In
short, the writer makes all the right moves.  In the concluding
paragraph, he/she says
\begin{quote}
I am going to make no specific recommendation as to whether you
should promote Mergetroyd Plotz or not.  After all, you know what the needs of your 
department and your school are.  You
can use the information that I have provided to come to an
appropriate decision.  
\end{quote}

\noindent Rubbish!  Imagine taking your car to a mechanic and hearing
him say ``Your transmission runs at half speed and your rear
wheels turn forward.  Your stroke is short and your valves 
rattle.  I am not going to make any specific recommendation for
repairs because, after all, it is your car and you know what your
needs are.''  Or imagine your physician saying ``Your heart will give
out any day now, and you are also a prime candidate for a stroke or
for total paralysis.  However, I will make no specific
recommendations.  It is your body, and you know best \dots.'' {\it You
are a professional; you are expected to render an opinion.}

\item The writer neglects to address explicitly the question at hand.
This omission is sometimes committed inadvertently, but this omission
is a dreadful error.  If you are asked whether Sara Glockenspiel should be
tenured, or promoted, or given a certain post, or a grant, then you
must say point blank what your opinion is {\it about that question
as it applies to that candidate}.
If you neglect to say, then your letter (taken as a whole)
is likely to be read as the worst sort of ``damning with faint
praise.''  Whether you intended it or not, you may have buried the
candidate.  

\item The writer faces the following request (and blows
it):  In a school that fancies that it wants to make hard decisions,
and elicit the {\it bona fide\/} truth from the letter writers, it is
common for the dean to include in his/her solicitation letter a query
like ``Would you tenure Marilyn Montpelier in {\it your\/} department?'' 
If the person being asked for the letter works at Harvard, and if the
institution soliciting the letter is a four-year teaching college,
then such a dean is just looking for trouble.              

\null \quad \ Even if the letter of solicitation does not explicitly ask this
question, we letter writers are often tempted to answer it.  Unfortunately,
the answer sometimes comes out like this:

\begin{quote}  \begin{tabbing}
$(*)$ \ \ \ \= Dr. Morris Fischbein is not good enough for us,  \\
\null       \> but he is certainly good enough for you.  \qquad \bad
\end{tabbing}
\end{quote}
Rarely is a letter writer clumsy enough to phrase things quite this
bluntly, but I have seen many a letter in which this sentiment comes
through loud and clear.   

\null \quad \ This is just too bad.  The person writing such a statement (or a
euphemistic paraphrase of it) probably thinks that he/she is being
frank and helpful.  He/she is being neither.  Instead, he/she is insulting
the maximum number of people in the least constructive possible
fashion.  A word to the wise should be sufficient:  proofread your
letter of recommendation to be sure that you have not inadvertently
(or intentionally) made statement $(*)$.  The inclusion of such
an assertion in your letter will vex the readers, and render
your letter ineffectual, so that it will not count.  I presume that
this effect is not the one that you want.

\null \quad \ If in fact you are at a place like Harvard, and if your letter is
solicited\index{recommendation letter!for a non-research school}
from a much more humble institution, and if you {\it must\/}
address this difficult question, then you should endeavor to tell the truth.
Say that Harvard's math department is usually ranked in the top three;
you only tenure people who are world leaders, indeed great historical
figures; such standards would be inappropriate to apply at an institution
like the one which has solicited the present letter.  However, you
certainly would recommend this candidate for tenure at Bryn Mawr,
or Swarthmore, or some other institution that you choose
for comparison.
\end{enumerate}

That concludes my enumeration of woeful errors. Now let me cut to the
chase. When you are writing a letter for a candidate, then a heavy
responsibility rests on your shoulders. The dean or chairman who solicits
the letters of recommendation is not simply casting his/her net and taking
a vote: this person\index{recommendation letter!significance of a
negative} wants a {\it mandate}. He/she will {\it not\/} weigh good
letters against bad: he/she wants to be socked between the eyes. A tough
dean once told me ``If a case is not overwhelming then I turn it down. If
the candidate is any good, he'll land on his/her feet. If not, then we are
better off without him/her.'' Thus if your letter says

\begin{quote}
Herkimer Nixon is no good.  Don't do it.
\end{quote}

\noindent then you may as well face the music and realize that {\it your
letter alone\/} will have killed the case---at least for now. I cannot
repeat this point too strongly: it is dead wrong to say to yourself ``This
is a negative letter that I am writing, but it will not count unless all
the other letters are negative too.'' Baloney! One negative letter will
usually stop the case cold. That is all there is to it.

A letter with inadvertent errors (of the sort mentioned above) will
not\index{recommendation letter!errors in} necessarily bury the
candidate, but it certainly will not help him/her.

In the closing paragraph of your letter, you will typically indicate a
degree of enthusiasm for the case at hand. Here is a
graded\index{recommendation letter!enthusiasm in} list of
examples---taken from letters that I have actually seen:

\begin{quote}
Igor Stravinsky has done a workmanlike job with his research program. 
\smallskip \hfill \break

Ayatollah Bono is a reasonable case for tenure.  You would
not go wrong to tenure him.
\smallskip \hfill \break

I recommend Zigamar Pistachini warmly for tenure and promotion.
\smallskip  \hfill \break

I recommend Rufus P. Quackenbush enthusiastically for ten- 
ure and promotion.  
\smallskip \hfill \break

The case for Guy de Maupassant Rabinowitz is overwhelming.  
I recommend him without reservation.  
\smallskip \hfill \break

I give Chicken \`{a} la King my strongest possible recommendation.  Phone
me if you require further details on the case.
\end{quote}

In case my admonitions have not sunk in, let me beat you over the
head with them.  The first two of these statements are in fact
negative.  Whether they were {\it intended\/} to be negative, or are
simply an articulation of the writer's loss for words, that is how
they will be read.  You might as well take the candidate out and
shoot him. The third passage is a little better (many evaluators will
read ``warmly'' as ``lukewarmly''), but does not convey passionate
affirmation.

By contrast, the fourth example will definitely be construed
favorably.  The adverb ``enthusiastically'' conveys the positive
nature of the assertion.  The last two sample sentences represent the
sort of forcefulness that is virtually mandatory if you want to argue
for the tenuring of a candidate at any of the best institutions.

The writing of letters of recommendation is not formulaic.  Indeed,
if all letters of recommendation fit a pattern and sounded the same,
or if all {\it your\/} letters look the same, then they will 
eventually be ignored.  Mathematicians keep a mental database on
letter writers in the same way that good baseball pitchers keep a
database on batters.  After several years, we know who ``tells it like
it is'' in his/her letters, who spins tales, and who simply cannot be
trusted.  We know who always writes the same letter for everyone. And
we act accordingly.\footnote{In fact there is an eminent mathematician
who has had many students and writes a great many letters of
recommendation.  They are so similar that you could hold any two of
them up to the light, one behind the other, and most of the words
would line up.  But then he scribbles his real opinion in the margin
by hand.}

You will develop your own style of writing letters.  Mathematics is
a sufficiently small world that, after several years, people will
recognize your letters of recommendation at a glance.  But,
no matter how you write your letters, you will want to take into
account the issues raised in this section.

During times when jobs are hard to come by, letters of recommendation
tend to become more and more inflated.  Everyone feels that he/she must
try harder if he/she is going to land a job for that special someone. 
Here\index{recommendation letter!inflation in} 
are examples of lines that have been used to describe specific,
rather famous, job candidates.  I do not necessarily recommend that
you use any of them; if you do, the readers might think that you are
eating with only one chopstick.  But these examples will give you an
idea of what some people have done to draw attention to what they are
saying, or to remove their particular letter from the ranks of the
humdrum.  (Of course names have been changed to protect the
innocent.)

\begin{quote}
Beef E. Wellington has a good idea every other day and writes a
brilliant paper every week. 
\mbox{ \ \ } \medskip 

Potatoes au Gratin knows both classical analysis and modern analysis.
He is the natural successor to Hardy and Littlewood.
\mbox{ \ \ } \medskip 

Talking to Leon Czogolsz is like talking to Enriques. 
%this is the guy who assassinated William McKinley 
% Federigo Enriques, 1871-1946.
(An inspiring thought, written by one too young to have ever
spoken with Enriques.)
\mbox{ \ \ } \medskip 

Cherries Jubilee is the most mathematically intelligent
person that I have ever met.
\mbox{ \ \ \ }  \medskip 

Rootie N.\ Kazootie is the greatest mathematician since Gauss.
\mbox{ \ \ }
\end{quote}

Although there is an art to writing a ``professional letter,'' it is also
the case that at least you are dealing with familiar territory, and
speaking of matters on which you are expert. Any professional
mathematician for whom you might write has a publication list, and a track
record in teaching, and a reputation as a lecturer, and some {\it gestalt}
as a collaborator. When you are writing for a student, by contrast,
matters are more nebulous. The student has none of the professional
attributes that you are comfortable discussing. Yet, if you want your
letter to be memorable, and to be perceived positively, then you still
want to say something noteworthy about the student.\index{recommendation
letter for a student}

While the precepts of organization that I have stated above still apply
in a letter for a student, some of the other particulars do not.  For
example, you most likely cannot remark on the student's scientific
work, and you most likely cannot make binary comparisons.  In fact
any attempt at binary comparison is likely to be ludicrous.  Imagine
saying ``I am delighted to recommend Sacajawea Smith.  She is every
bit as good as Euthanasia Jones, whom I recommended five years ago to
a different institution.'' If in fact you previously recommended a
student who turned out to be a well-known star---or at least a 
well-known star at the institution to which your letter is
addressed---then by all means make a binary comparison involving that
person if such a comparison is appropriate. Usually it is not
appropriate, so no such comparison should be included.

Thus in practice you must try a bit harder to say something specific about
the student for whom you are writing a letter. After you have been
teaching for several years, it may be the case that you have actually
taught a few thousand students (this would be true, for example, if you
have taught calculus in large lectures several times). It becomes
difficult\index{recommendation letter!specificity in} to distinguish
students---even good ones---in your memory, much less to say something of
interest about any of them. If you apply yourself to the task, then you
can nevertheless come up with some noticeable things to say. Here are some
examples, taken from genuine letters:

\begin{quote} 
Fig Newton is one of the five most talented
undergraduates that I have encountered in twenty years of teaching.
\mbox{ \ \ } \medskip

Iphiginea Mandelbrookski is hard working and perseverant.  She
can think on her feet---at the blackboard---just like a
mathematician.  She is original and imaginative.
\mbox{ \ \ }  \medskip

In order to test her creative abilities, I have given Cleopatra Jones
extra work outside of class.  She discovered a new proof of
Gronwall's inequality, discovered Euler's equation in the calculus of
variations on her own, and has also posed numerous interesting
problems of her own creation.  Needless to say, she breezed through
all the standard class work.
\end{quote}

\noindent As usual, the point is to say {\it something}---and that
something should be quite specific.  The view of letter {\it readers\/} is
that if the letter writers cannot say anything unambiguous and
remarkable about a student, then there is probably nothing remarkable
about that student.  So what if the student can earn mostly $A$s in his
classes?---that is no big deal, and in any event can be gleaned from
the transcript.

Sometimes a student, or someone else, will ask you for a letter about
himself/herself and you do not feel that you can write a good one.  Either
you have nothing to say, or you have nothing good to say, or you have
some other valid reason for not writing.  (Note that this case is
different\index{recommendation letter!declining to write}
from the one in which a dean is asking
you for a letter about one of his/her faculty.  Now the candidate
himself/herself is standing before you and asking for a letter
{\it about himself/herself}.)
You always have the option
of agreeing to write, and then writing a negative letter.  Often,
however, you bear the candidate no malice and think that he/she
deserves a chance.   In that case, the honorable thing to do is to
say to the candidate ``I'm sorry.  I frankly do not feel that I could
write a good (or supportive) letter for you. Perhaps you should ask
someone else.''  The rotten thing to do---and this happens far too
often---is to say ``Oh yes, fine'' but with no intention of ever
writing {\it anything.\/}  Note that the lack
of your letter in the dossier will make that dossier incomplete; in
many cases the candidate will not, as a result, be considered at
all.  If that is the effect you want, then you should have the
courage to say something in a letter.  If it is not the effect you
want, then you should have the courtesy to take a ``pass.''

One of the most critical, and delicate, types of letter that you will
have to write is a letter seeking a job for a student completing his
M.A.\ or Ph.D.\ under your direction.  Your statements are {\it a priori\/}
suspect\index{recommendation letter!for your own Ph.D.\ student}
because you obviously have a vested interest in finding this
student a job, and in seeing him/her succeed.  Thus you must strive to
put into practice the precepts described above:  {\bf (i)}  say why
this student is good, {\bf (ii)}  say what this student has
accomplished, {\bf (iii)}  if possible, compare the student favorably
with other recent degree holders, {\bf (iv)}  say something about the
student's ability to teach, {\bf (v)}  say something about the
student's ability to interact with other mathematicians.   

A meat-and-potatoes job application from a fresh Ph.D.\ has a
detailed letter from the thesis advisor that conforms, at least in
spirit, to the suggestions just adumbrated.  This detailed letter is
accompanied by two or three additional letters from other instructors
at the same institution, each of which is rather vague and says in
effect ``Doo dah, doo dah; see the letter by the thesis advisor.'' If
you want your student's dossier to stand out, and to really garner
attention, then you should strive to help the student make his/her
dossier rise above this rather dreary norm.  Endeavor to ensure that
the other writers know something about what is in the thesis. If
possible, convince someone from another institution to write a letter
for the student. Make sure that the dossier includes detailed letters
about the student's teaching abilities.

When you write a letter of recommendation, tell the truth. If all your
letters\index{recommendation letter!truth in} read ``This candidate is
peachy, and a dandy teacher too. Give him/her $X$'' (where $X$ is the plum
that the candidate is applying for), then after a while nobody will pay
any attention to what you say. I presume that if you take the trouble to
write letters, then that is not the result that you wish. The
infrastructure has a memory. It will remember whether you are a person who
can make tough decisions, or whether you are wishy-washy. If you want your
letters to count, then you must call it as you see it. It is hard to be
hard, but that is what the situation demands.

One issue that we, as letter writers, often must address is whether or not
a job candidate can speak English, and how well (this question could even
apply to an undergraduate student---especially if that
student\index{English ability of candidate} is applying to graduate school
and might be considered for a Teaching Assistantship). In this matter we
are, in the United States, cursed by our group dishonesty over the past
forty years. Too often have we said in a letter that ``this candidate
speaks excellent English, can teach well, and is a charming
conversationalist to boot.'' In a more frank mode, we might have said
``This candidate speaks better than average English'' (recalling Garrison
Keillor's statement about the town of Lake Woebegon, in which ``all the
children are above average''). When the candidate arrived to assume
his/her position, the hiring institution often found that he/she could not
understand even simple instructions and had no idea how to teach.

It is difficult, but you must endeavor to be honest about the candidate's
fluency in English (again, your credibility---which will follow you
around all your life---is at stake).  You could say, for example,

\begin{quote}
This candidate speaks English like a home-grown American, with no
trace of an accent.  Listening to him/her is like listening to Walter
Cronkite.
\end{quote}

\noindent  This would be the ideal thing to write, and would dispel all
trepidations about the candidate's fluency.  Unfortunately, if the
question needs to be addressed at all, then this statement probably
is not true.  You could instead say 

\begin{quote} Luisa Longshoremanska has been taking ``English as a second
language'' and has taught several lower-division courses successfully.
Her English is accurately formulated and clearly enunciated.  Students
have no trouble understanding her.  
\end{quote}

\noindent Unfortunately, you cannot always be so enthusiastic.  Sometimes
you must say something like

\begin{quote} 
Mr.\ Anthrax Xlpltqlpl has been working hard on his
English, and has made substantial progress.  One still needs to
concentrate in order to understand him.  
\end{quote}

\noindent Or you might say

\begin{quote}
It takes students three or four days to become accustomed \hfill \break
to Ms.~Imelda Rasputin's English, but her charming personality helps
them along.  As a result she is a most successful teacher.
\end{quote}

\noindent The thought that I am trying to formulate here is that Mr.\
Xlpltqlpl's English or Ms.~Rasputin's English is not perfect.  But
Mr.\ X and Ms.\ R are real troupers.  They try hard, and the students
(at least in Ms.~R's case) forgive them a lot. 

Of course you can plainly see that I am trying to suggest ways to
avoid saying  ``This person cannot speak English and refuses to
learn.  He/she is only suitable for a nonteaching position.''  But
sometimes---presumably not in the case of Mr.~Xlpltqlpl or 
Ms.~Rasputin---it must be said.

At the risk of repeating myself, let me say that when you address the
candidate's ability with English then you should not be
formulaic.  If all your letters about foreign candidates say 
\begin{quote}
$X$'s English is just fine.  He/she is a good teacher.
\end{quote} 
then, after a while, the world will mod out by that portion of your
letters.  Try to say something original, apt, and true about each
candidate.  I once wrote the following about a fresh Ph.D., from
a foreign country, who was applying for a job:

\begin{quote} I consider myself to be rather a good teacher, but I
really learned something when watching Mr.\ Frangi Pani with his
class.  He moves skillfully among his students, looks at their work,
makes insightful remarks, and does a marvelous job of eliciting class
participation. It is clear that the students like and respect him.  
\end{quote}

\noindent This passage addresses the language issue implicitly, for
it confirms that the candidate can {\it teach.\/}  Moreover, it is
not just a bunch of pap.  It says something particular and notable 
about the candidate's abilities. 

Occasionally, you will have to address a truly thorny matter in one
of your letters of recommendation.  As an instance, I was once
writing on behalf of a young mathematician who was applying to
several dozen first class universities for a position.  I thought
that I knew this person quite well.  But, a few days before I was
going to draft my letter, I learned that the candidate was undergoing
a sex change.  I had to decide whether I should mention this fact in
my letter.  I reasoned as follows:  if he/she were changing from
Catholicism to Judaism, or from Democrat to Republican, or from
carnivore to vegetarian, I certainly would not consider discussing
the matter in my letter of recommendation for a mathematical post; so
why should I treat trans-sexuality?  And I did not.  Some time later, I
discussed the matter with one of my mentors.  He told me that I had
erred.  In stern terms, he informed me that a matter like this could
affect the candidate's ability to teach, and his ability to function
as a colleague; therefore I was morally obligated to mention the
matter.  I still do not know what the correct course of action should
have been.  I only hope that I will not be faced with another choice
like this one any time soon.   

Just for fun, let me conclude this long section by quoting from a letter
for tenure that was written (truly!) about twenty years ago for a candidate
in a French department.  Call the candidate Mr.~de Gaulle.

\begin{quote}
Surely Mr.~de Gaulle is now wiser than he once was.
\end{quote}

\noindent {\it That was the entire text of the letter!---No
introduction, no conclusion, no binary comparison, no
exegesis of the candidate's scholarly work.  
Just\index{recommendation letter!brevity in} 
the one sentence.}  Although the
letter does not follow the precepts described in this section, it
definitely gets its point across.

%% Section 4.2

\section{The Book Review}

As with most topics in the subject of writing, there is some
disagreement over what constitutes a good book review.  When Paul
Halmos\index{book!review, how to write a}
was the book reviews editor of the {\it Bulletin of the AMS}, he
sent every reviewer a set of instructions.  The gist of these
instructions was that a book review is not a book report.  It should
{\it not\/} say ``Chapter 1 says this, while Chapter 2 says that. 
Chapter 3 is a bore, and Chapter 4 is too hard.''

Instead, according to Halmos, a book review on a book about $X$ is
an excuse to write an essay about the subject $X$.  Look at the
book reviews in the {\it New York Review of Books}.  
On the whole they are a delight to read, and they conform to
Halmos's view of what a book review is and does.  These reviews
tell you about the book, but they paint the picture on a large
canvas.

To reiterate: If you are reviewing a book on harmonic analysis, then you
should write about the history of the subject, what the\index{book!review,
what should be in a} milestone books and theorems have been, who the major
players are and were, and what the big problems are. Drop some names. Make
some assertions and conjectures. Having laid considerable groundwork, then
finally focus on the book under review. Describe where it fits into the
infrastructure you have outlined. Indicate its strengths and weaknesses.
Suggest who would profit from reading it, and why. Touch on areas where
there is room for improvement. Do not, however, use my suggestions here as
an excuse to write an opinionated essay and virtually ignore the book. The
book review is supposed to be {\it about the book;\/} but it should be
about the book in the context of the subject matter, not the other way
around.

Here\index{book!review, check list for} are some other issues that your
book review might address:

\begin{quote}
$\bullet$ \  Will students benefit from reading the book?
\mbox{ \ } \smallskip

$\bullet$ \  Are there exercises?
\mbox{ \ } \smallskip  

$\bullet$ Are there lists of open problems? 
\mbox{ \ } \smallskip  

$\bullet$ Is there an accurate and complete bibliography?
\mbox{ \ } \smallskip  

$\bullet$ Is there an index?
\mbox{ \ } \smallskip  

$\bullet$ Is there a list of notation?
\mbox{ \ } \smallskip  

$\bullet$ Is there sufficient review material?  Does the book begin
at a reasonable level?
\mbox{ \ } \smallskip  

$\bullet$ Does the author provide an adequate amount of 
detail in the book?  Does the book make too many demands on the reader?
\mbox{ \ } \smallskip  

$\bullet$ Are the proofs complete, clear, and accurate?
\mbox{ \ } \smallskip  

$\bullet$ Is the book organized in an intelligent fashion that
is useful to the reader?  Can the beginner navigate his
way through the book?
\mbox{ \ } \smallskip  

$\bullet$ Is the history correct?  Are attributions complete and 
accurate?
\mbox{ \ } \smallskip  

$\bullet$ Does the book bring the reader up to the cutting edge of
research?
\end{quote}
If you think about the issues that I have raised here, then you will
realize that I have described what a potential reader of the book
will want to know when he/she is making a decision as to whether to
buy the book and whether to read the book.  One of the main purposes
of your review is to inform such decisions.

Many mathematical book reviewers---writers for the {\it Bulletin of
the AMS}, for instance---feel obligated to write a {\it positive\/} or
upbeat book review, no matter what they really think of the book. 
They are afraid to be critical.  In my opinion, this attitude is an
error.  Not all books are good, and not all good books are entirely
good.\index{book!review, positive}\index{book!review, negative}
You will help the audience, and the author as well, if
your review points out inadequate features of the book, or omissions,
or errors, or items that can be improved. You should tender your
criticisms in a constructive fashion:  in this manner you will
increase the likelihood that people will attend to what you have to
say, and your thoughts will perhaps make friends and influence people
(rather than the opposite).  

On the other hand, there is the occasional reviewer who lets it all hang
out. Books seem to have a sort of permanence that papers do not. An
incorrect or wrong-headed paper is, after all, ultimately buried in a
bound journal volume and hidden away. But a book is always right there on
the shelf, staring us all in the face. And, as previously noted, a book
reaches out to a larger audience than does a paper. As a result, emotions
can run high over a book. I have seen a book review that (literally) began
by questioning the editorial decision to publish the book and asserted
further that the book completely misrepresented the subject matter; the
reviewer spent the rest of the review describing what the subject was {\it
really} (in his opinion) about, with nothing further said about the book
itself. And I have also seen a book review [Blo] that compared the subject
matter of the book to rather delicate portions of the female anatomy. A
recent (and rather controversial) book review [Kli] asserts that the book
under review is obviously about a weak subject, as one can see by
examining the Bibliography and noting the substandard journals in which
the relevant papers are published; the reviewer neglects to point out that
he/she is or has been on the editorial board of most of the relevant
journals (see [NoS] for an incisive reply). While these essays are briefly
diverting they are, in retrospect, embarrassing for us all. As you write
your review, pretend that you are reviewing the book of a friend: you want
to be honest, and you want to be helpful, but you also want to be
scholarly and dignified. Brutality is almost never the order of the day.

In 1978 there appeared a marvelous book on algebraic geometry that
is almost universally admired, but that is famous for having a large number
of errors: either slight misstatements, or omissions, or incomplete
proofs.  The fact remains that everyone loves this book, and there
is no other like it.  (Heck, I may as well tell you:  it is [GH].)
One reviewer [Lip] praised the book to the heavens, but
felt that he had to say something about the hasty writing and the
density of errors. So he wrote in part 
\begin{quote} If it makes you
feel better, think of this book as a set of lecture notes, or even as a
fantastic collection of exercises, with copious hints.  
\end{quote}
Thus the reviewer did his duty:  he certainly said something
critical, but he said it in a charitable manner, and with good
humor.  Even the authors must have chuckled over these remarks, and
everyone learned something.

The harshest book review that I have ever seen appears in [Mor]. Mordell in
fact uses this review to trot out his frustration with the French school's
abstraction of his beloved number theory. He attacks not just the book,
but he attacks its author on a rather personal and visceral level. A now
famous letter was subsequently written by C.\ L.\ Siegel to Mordell,
praising the review and heaping even more calumny upon the book's author.
A discussion of that interchange, and its significance, appears in [Lan].
The trouble with such a review is that any flow of scholarly thought or
criticism is lost in the morass of venom and vituperation. No constructive
purpose is served by such a review. It is also virtually impossible to
have any useful dialogue following upon such a review.

If you are called upon to review a book, and you are tempted to trash it,
then I suggest that you set the draft of your review aside for a month (a
year would be too long!). Let the ideas gel, and let the words mellow.
Show it to a few trusted friends. After a month, you will probably be
inclined to take the long view, and to express your ideas in a more
temperate fashion. The result will be a better review, and one that you
will still be proud of ten years after it appears.

%% Section 4.3

\section{The Referee's Report}

When you are asked to write a referee's report, then you are 
being requested to
offer\index{referee's report!how to write a}
your opinion as an expert.  If you agree to write
the report (and you {\it should\/}---refereeing is an important part of your
professional duties), then you should adhere to the following
precepts: 
\begin{itemize} 
\item Write the report in a timely
manner---if possible within the time frame suggested by the editor.
\item Tell it as you
see it.  Just as in a letter of recommendation, enunciate your
opinion clearly and succinctly, defend it, and then summarize your
findings.  
\item Defend your opinion in detail.  You need not find a
new proof\index{referee's report!what to include in a}
of each lemma, nor read every bibliographic reference.  But
you must read enough of the paper so that you can comment on it
knowledgeably.  While you may not have checked every detail
in the paper, you should at least be confident of your opinion
as to the paper's correctness and importance.

\null \quad \ If somebody asked you whether you liked your car, and
whether you would recommend that they buy one, you would not (in all
likelihood) tell how each bolt 
was\index{referee's report!level 
of detail in a} installed in the chassis, nor how
the finish was applied to the body.  You would instead summarize the
overall performance and features of the automobile.  Just so, when
you evaluate a paper you should address Littlewood's three precepts:
{\bf (1)}  Is it new? {\bf (2)} Is it correct? {\bf (3)}  Is it
surprising?\index{Littlewood, J.\ E.!precepts}
You should speak to its contribution to our knowledge,
and to the literature.  
\item Provide constructive criticisms of the writing, 
or of the paper's organization.  You may
enumerate spelling and grammatical errors (if you wish
to do so).  You should certainly
point out mathematical errors, or places where the reasoning is
unclear.  But you should not be captious.  ({\it Exercise:}  Look up this
word in your Funk \& Wagnall's and think about its relevance to the
present discussion.)  
\item Place the paper in context:  How does it
compare to other recent papers in the field?  Where does it fit?  
Does it represent progress?  If you were not the referee, then is it
a paper that you would want to read?  
\end{itemize}

Of course your report should be tailored to the journal to which
the paper was submitted.  For instance, the {\it Annals of Mathematics\/}
professes to publish papers of great moment, written for the ages.
Other journals have the more modest goal of publishing papers that
are correct and of some current interest.  Still other journals have
no standards\index{referee's report!keying to a given journal} 
at all.  You must speak to people in their own
language---language that they will understand.  Likewise, when you
evaluate a paper for a journal, base your assessment on {\it that 
journal's value system}.

A typical referee's report is anywhere from one to five pages (or, in
rare instances, even more).  Its most important attribute should be
that it makes\index{referee's report!key points in}
a specific recommendation.  Everything else that you
say is for the record: it is important, but it is secondary.

%% Section 4.4

\font\teneufm=eufm10
\font\seveneufm=eufm7
\font\fiveeufm=eufm5
\newfam\eufmfam
\textfont\eufmfam=\teneufm
\scriptfont\eufmfam=\seveneufm
\scriptscriptfont\eufmfam=\fiveeufm
\def\frak#1{{\fam\eufmfam\relax#1}}

\newfam\msbfam
\font\tenmsb=msbm10  scaled \magstep1  \textfont\msbfam=\tenmsb
\font\sevenmsb=msbm7  scaled \magstep1   \scriptfont\msbfam=\sevenmsb
\font\fivemsb=msbm5   scaled \magstep1   \scriptscriptfont\msbfam=\fivemsb
\def\Bbb{\fam\msbfam \tenmsb}

\def\RR{{\Bbb R}}
\def\CC{{\Bbb C}}
\def\QQ{{\Bbb Q}}
\def\NN{{\Bbb N}}
\def\ZZ{{\Bbb Z}}
\def\II{{\Bbb I}}

%%%%%%%%%%%%%%%%%%%%%%%%%%%%%%%%%%%%%%%%%%%%%%%%%%%%%%%%%%%%%%
%%%%%%%%%%%%%%%%%%%%%%%%%%%%%%%%%%%%%%%%%%%%%%%%%%%%%%%%%%%%%%

\def\hexdigit#1{\ifnum#1<10 \number#1\else
 \ifnum#1=10 A\else\ifnum#1=11 B\else\ifnum#1=12 C\else
 \ifnum#1=13 D\else\ifnum#1=14 E\else\ifnum#1=15 F\fi\fi\fi\fi\fi\fi\fi}

\chardef\\="5C                    %% Typesets \ in \tt mode
\chardef\{="7B  \chardef\}="7D    %% also left and right braces

 \def\HollowBoxx #1#2#3{{\dimen0=#1 \advance\dimen0 by -#2       
       \dimen1=#1 \advance\dimen1 by #3                       
        \vrule height 0pt depth #3 width #2                   
       \hskip -#3
       \vrule height #1 depth #3 width #3}}                   
 \def\LeftContraction{\mathord{\kern1.45pt \HollowBoxx{6pt}{3.5pt}{.4pt}}\,}

 \def\HollowBox #1#2#3{{\dimen0=#1 \advance\dimen0 by -#3       
       \dimen1=#1 \advance\dimen1 by #3                       
        \vrule height #1 depth #3 width #3                    
        \vrule height 0pt depth #3 width #2                   
        \hskip -#3}}                                             
 \def\RightContraction{\mathord{\, \HollowBox{6pt}{3.1pt}{.4pt}} \kern1.6pt}             

\section{The Talk}

Giving a talk is different from writing. But it is relevant to the writing
process. We ordinarily do some writing to prepare for a talk. And what we
write\index{talk!how to give a} will strongly influence the talk itself.
So this topic is fair game for the present book.

A talk is more flexible than a paper. In a talk, you may indulge in
informalities, whimsicalities, and a little imprecision; it helps the
audience a lot to tell of things tried, and things that failed.
You\index{talk!flexibility of a} may work trivial examples, and use them
as a foundation on which to build ideas.

A talk is also less flexible than a paper. Because the audience receives
the talk in linear order---it\index{talk!inflexibility of a} cannot
rewind or speed ahead to check on things---it is therefore at your mercy.
You are at a great advantage, when preparing a talk, if you are aware of
the limitations of the medium. Endeavor to be gentle.

John Wermer [Wer] makes an excellent case that many mathematics talks are
not as effective as they might be because the lecturer is\index{Wermer,
John} speaking to an imaginary audience located inside his/her head. This
audience is one that knows all the necessary motivation, can pick up on
fifty new technical definitions quickly and easily, can follow a technical
proof (without explanation) in a jiffy, and can fill in the logical gaps
and potholes left by the speaker. Of course\index{talk!apocryphal
audience for a} such an audience is apocryphal, and thus we are often left
with a communication gap between speaker and audience. This section will
give you some advice on how not to be like Wermer's {\it idiot-savant}.

What are the ingredients of a good mathematics talk? First, you must know
your subject cold. This does not mean simply that you know it well enough
to communicate it to another expert like yourself, but rather that you
know it well enough to teach it: that you know the background,\index{talk!ingredients of a} the biases, the reasons for the questions, the good and
the bad attacks on the problems, and the current state of the art.
However, just because you know all these things does not mean that you
need to say all of them. A good mathematics lecture is an exercise in
self-restraint. Never mind impressing the audience with your profound
erudition, your spectacular vocabulary, your extensive professional
connections, or your readiness to cite last week's hot results. Instead
showcase a nugget of knowledge and insight, and shore it up with crisp
comments and incisive examples.

If your talk is scheduled for fifty minutes, then the first twenty should
be accessible\index{talk!time management in} to a graduate student who
has passed the quals. My statement is a strong one. Such a student is not
expert at anything. He/she knows the basics of real and complex analysis,
algebra, and perhaps a little geometry. This student has (we hope) an open
mind and wants to learn. But your talk in those first twenty minutes
should presuppose no specialized knowledge beyond what has just been
mentioned. This explains why a nice example or two can be so useful. With
an example, God is in the details. The playing field is level, and
everyone can benefit. The example(s), of course, should lead to some
definitions and the formulation of the questions that you wish
to\index{talk!use of examples in} treat in the body of the talk.

The next twenty minutes of the talk should be pitched at a mathematically
literate person who is not a specialist. By this I mean\index{talk!organization of} that, if your talk is about some part of analysis, then
the second twenty minutes should be comprehensible {\it not just to a
specialist in another part of analysis}, but to an algebraist. So you can
mention more sophisticated ideas---sheaf theory, or elliptic regularity,
or wave front sets, or singular integral operators---but you should not
beat them to death.

The last ten minutes can be for the experts, for God, and for you (not
necessarily in that order). Every speaker should have a chance to strut
his/her stuff, and the\index{talk!time segments in} end of the talk is
when you should do so. Mention some gory details. Make speculations,
formulate technical corollaries, sketch the key ideas in the proof. Forget
the neophytes and address yourself to the people who might read your
papers. In fact if you do not use the last ten minutes of your talk in
this way, then you might leave the impression that you are a lightweight,
or that you have nothing to say.

Attempt to finish with a bang. Too many math talks begin with ``Well, what
I want to talk about today is \dots'' and then a definition\index{talk!conclusion of a} goes onto the blackboard. Too many math talks end with
``Well, I guess that's all I have to say'' or ``I see that I'm out of time
so that's it'' or ``I guess I'll stop here; thank you.'' Surely you can
devise a more creative and informative way to conclude your discourse. You
would never end a paper in this fashion. Of course when you write a paper
you have time to sit and think of a nice turn of phrase for your
conclusion. You should do the same when composing a talk: prepare the
introductory sentence or two in advance; likewise prepare the concluding
sentence or two. You could finish with a few courteous words of thanks for
the opportunity to visit your hosts and to enjoy the hospitality of their
department; or you could end with a few mathematical sentences---of real
substance---that summarize your enthusiasm for your subject matter. But do
end by {\it saying something}.

The preceding discussion may make it seem that giving a good talk is a
piece of cake; that it requires only an acclimatization to certain simple
proprieties. Not so. Many other parameters figure into the process.

In fact there are many types\index{talk!types of} of talks: the
colloquium, the seminar, and the ``job'' talk (in which you are showcasing
yourself before a department that is considering offering you a job) are
three of these. The colloquium is supposed to be for the entire department
and perhaps for the graduate students as well; the seminar is for a group
of specialists, probably your friends; and the job talk is a set
piece---something like Kabuki theater---in which you show yourself. An
entire separate book could be written on the art of giving talks. In the
interest of brevity, my remarks below will center around colloquium talks.
Seminar talks are less demanding and job talks more so. The remarks below
apply in some form to {\it any\/} talk; the trick in interpreting my
advice for a particular circumstance is to {\it know your audience.} As
you read my detailed\index{talk!checklist for} remarks below, keep this
unifying principle in mind.

\begin{enumerate} 
\item  {\it Showcase one theorem}, or
perhaps a single cluster\index{talk!focus of} of theorems.  
There is no point to giving
a talk on five truly different theorems, because the audience cannot
absorb so much material in one sitting.  On the other hand, if you
cannot build your talk around one theorem then perhaps you have
nothing to say.  Here is what Gian-Carlo Rota 
thought about the matter:

\begin{quote}
Every lecture should state one main point and repeat it over
and over, like a theme with variations. An audience is like a
herd of cows, moving slowly in the direction they are being
driven towards. If we make one point, we have a good chance
that the audience will take the right direction; if we make
several points, then the cows will scatter all over the field.
The audience will lose interest and everyone will go back to
the thoughts they interrupted in order to come to our lecture.
\end{quote}

\null \quad \ If the talk is a survey, then you should temper this
last advice to suit the occasion.  Better to give a survey of a
particular aspect of semi-Riemannian geometry than to endeavor to
survey the entire subject of geometry.  And do suit the talk to the
audience.  You can survey non-Euclidean geometry for junior/senior
mathematics majors, or you can do it for seasoned mathematicians. But
you would do it differently for each of these audiences.
\vspace*{-.08in}

\item {\it Have an attractive title.}  A casual
observer, seeing the title ``Subelliptic estimates for
a quasi-degenerate, semilinear partial differential operator
satisfying a weak symplectic condition with applications
to the hodograph technique of H\"{o}rmander,'' 
will probably be 
more tempted to head for a late afternoon\index{talk!title of} beer 
than to attend your talk.  The title ``A new attack on a class
of nonelliptic equations'' conveys the same spirit and is
likely to suggest to a broader class of people that the talk
may contain something for them.
\vspace*{-.08in}

\item {\it Prove something}. It leaves a bad taste\index{talk!proofs in}
in everyone's mouth if you talk about a subject but do not get in there
and do it. One good strategy is to prove a special case, or work out an
example, in some detail; then use this prolegomena to sketch the key ideas
in the proof of the main result.

\item {\it Structure your talk so that everyone will take something away from
it.} Ideally, a member of the audience who is questioned\index{talk!structuring of} that evening about that day's colloquium should be able to
say ``The talk was about this'' or ``The main theorem was that'' or ``The
speaker was relating geometry to combinatorial theory in a new way.'' If
you bear this thought in mind while composing your talk, then it will have
a strong, and salubrious, influence on your entire approach to the
process.

\item {\it Be specific.} Heed this advice, both when you are writing and
when you are speaking.\index{talk!specificity in} Nobody wants to listen
to an hour of vague fluff. Nobody wants to perceive that you are dodging
the main point of the discussion. If you appear to be evasive then, at
best, you could make people think that you are sloppy and imprecise; at
worst, you could leave people with the impression that you are faking
it---indeed that you do not know how to prove these theorems.

\null \quad \ I once heard a mathematician dedicating a lecture to an
eminent person, on the occasion of that man's sixtieth birthday.  In brief,
the dedicator said ``In my country the tradition in lectures,
especially with my thesis advisor (whom we are too polite here to name) has been
to deal in vague generalities.  This man (to whom I am dedicating
my remarks) has taught us to present concrete
examples, and to work through them completely.''  The value of
showing your audience the inner workings of the material you are
presenting cannot be overemphasized.  This process helps to draw in
students (both young and old), and shows them how the subject works. 
It also helps to involve those who are not already expert.

\item {\it Do not be afraid to dream.}  I say this cautiously, for I have already
warned you not to prevaricate or mislead. But a talk is a different
vehicle from a formal piece of\index{talk!dreaming in} writing. Standing
before a group and speaking is an opportunity for you to tell the audience
what you tried, what did not work, and what might work in the future. It
is absolutely impossible in mathematics to publish a paper that says
``Today I woke up and tried to prove the Riemann hypothesis and I\index{Riemann hypothesis}
failed.'' In a mathematical {\it talk}, you can dandle such thoughts
before your audience and not only survive, but in fact heighten the
audience's appreciation for you and your insights.

\item {\it Do not be afraid to be informal.} One of the most effective devices
that I have seen is for the speaker to say ``If we assume these three
explicitly stated hypotheses, together with some other\index{talk!informality in} technical things that I shall not enunciate, then the
following conclusion holds.'' Often the technical items that are left
unspoken are of great interest to the deep-down experts; but to everyone
else they would be meaningless, indeed confusing. It takes real insight,
and a dash of courage, to be able to say to the audience that you are
sloughing over some difficult points. Of course you should never lie; but
you may certainly downplay some of the technical points in your subject.

\null \quad \ These comments also apply when you are presenting a
proof.  In a specialized seminar, it might be appropriate to slug
your way through every technical detail of your argument.  In a
colloquium, such arcana are virtually never appropriate.  If your
theorem has any substance at all, then its proof may consist of ten
or twenty or more pages of dense argument.  It could take a serious
reader a week to digest thoroughly the inner workings of your
reasoning.  Thus it could never work to present the entire theorem,
with its proof, in a colloquium talk.  Hit the high points, say a
word about what you are omitting and oversimplifying.  Proud as you
are of the cute argument you cooked up for the proof of Sublemma
3.1.5, do not trot it out during your colloquium.

\item {\it Prepare your talk with multiple entry points and multiple exit
points.}\index{talk!entry points}  What does this mean?  
Rare is the listener who can pay rapt
attention for the full space\index{talk!exit points} 
of 50--60 minutes.  Many members of your
audience will drift in and out.  If you say something interesting,
then certain people will begin to think their own thoughts, or try to
produce a necessary example or lemma.  Make it easy for such people
to ``re-enter'' the lecture.  Provide several doorways.

\null \quad \  Likewise, there is no way to predict how a given talk
will go before a given audience.  If you are lucky, there will be
fortuitous interruptions and serendipitous comments.  Time will not
be used in just the way that you had planned; you could easily be
caught short.  Therefore you should create several
propitious\index{talk!time management in} points at which you can
make a gracious exit from the talk.  As already noted, a hasty
``Egad, I'm out of time'' is not a savvy way to finish your
colloquium.  In any event, do not run overtime---at least not by more
than a couple of minutes.  First, to do so is rude; second, colloquia
are at the end of the day and people have other things to do (such as
going home to dinner); third, people simply have no patience for a
talk that exceeds the allotted time.  At my university, we had a
leading job candidate who, in his ceremonial talk, ran out of time. 
He gave us a big smile, went to the clock, and pushed the big hand
back twenty-five minutes.  And then he used them!  Suffice it to say
that there was no further discussion of his candidacy.

\item {\it Prepare, prepare, and prepare some more.}  You should have
thorough notes before you, but you should rarely refer to them.  
Your\index{talk!preparation of} talk should have an edge:  you
want to be thinking through the ideas with your audience, and you want to
be {\it talking\/} to the people in the room.  You are not giving a
recitation to your buddy in the front row.  You are not lecturing to
the fictitious audience that is engraved in your frontal lobes.  You
are talking to the individuals who are breathing the same
air as you.  Pick them out as you speak; look at them; change your
focus and your depth perception as the talk develops.  Pace around.
Step backward and forward.  Use your body and your voice to lure
the audience into the talk.  Do {\it not\/} be a slave to your notes.

\null \quad \ Several years ago I watched an eminent mathematician prepare
to give a colloquium on a topic that I personally had seen him lecture on
at least four times previously. He had probably lectured on it fifty times
in total. I had attended his course on the subject. He {\it owned\/} this
material; he had created it from whole cloth. He could have given this
talk in his sleep. Nonetheless for this, his fifty-first performance, he
insisted on sitting quietly in a room for an hour and writing out
everything that he was going to say. During the talk, he cast not a glance
at his notes. At the end of the talk, he summarily dumped the notes into
the trash.

\null \quad \ This process made a tremendous impression on me, and I have reminisced
frequently about what I observed. {\it Writing out his talk was his
mantra.} He used this process to prepare himself psychologically for the
talk.\index{talk!preparation of} Some\index{talk!mantra for} people will
prepare by just staring off into space and walking, mentally, through the
talk. Others will stroll to the student center and buy a cup of coffee.
Still others will spend the entire afternoon in the library sweating over
the literature, and worrying about questions that someone might ask but in
fact never will. It does not matter what you do to psyche yourself up for
your talk, to guarantee that you are prepared. The main point is to {\it
do something}: find a technique that works for you and use it.

\item {\it Be careful in your talk to give credit where it is due.} Do not give
attributions only when your name is involved. In fact most\index{talk!giving credit in} speakers tend to do the opposite. When presenting
someone else's theorem, the speaker is careful to write out all the
relevant (sur)names in full. When it comes to his own theorem, the speaker
just writes something like

\begin{quote}
{\bf Theorem:} \ [Fu-Isaev-K]  
\vspace*{.03in}
 
Let $\Omega \subseteq \CC^n$ be a smoothly bounded, pseudoconvex \hfill \break
domain with noncompact automorphism group \dots
\end{quote}

\noindent This citation is an example taken from my own life, where Siqi Fu
and Alexander Isaev are my collaborators and ``K'' is yours truly.
\end{enumerate} 

Now that I have listed the ten commandments, let me turn to a discussion of
general principles. Many technical skills are necessary for giving a good
talk. I have already mentioned eye contact and organization. Let me also
discuss blackboard technique. Even if you are a great expert in your
subject, and have a charming and erudite delivery, you will be putting a
substantial barrier between yourself and your audience if your writing is
an incoherent mess, or if you fill the blackboard with a chaotic barrage
of longhand. Learn to write in straight lines, horizontally, from left to
right. Write large, and write neatly.\index{talk!blackboard use in} Do
not put much on each blackboard. Give the audience a chance to copy what
you have written before you erase it.

Do not stand in front of what you
have written.  As you write, read the sentences aloud. 
Learn to draw your figures accurately and skillfully.  Isolate
material that will require later reference and {\it do not
erase it.\/}  Plan in advance how you will use the blackboard, so that
you can be sure that you will always have room for what you
want to write next.  Just as the director of a play knows
in advance where each actor will be at each moment, and plans
every movement on the stage in considerable detail, so you
should plan the moves of your talk in advance.  The
audience will grow phenomenally frustrated watching a forlorn
speaker pace back and forth in front of his/her blackboards---for several
minutes!---trying to decide what to erase, or what to save.

Some people solve the blackboard problem by using overhead slides
(transparencies) instead.  The very act of creating slides in advance
addresses virtually all the issues that I have raised about
blackboards in the last paragraphs.  Slides, in the hands of a
skilled\index{talk!overhead slides in} 
user, can be a powerful tool.  (The blackboard is sometimes
inescapable, however, so you should learn to use it.)

If you do use slides, then learn to use them wisely.  Each slide
should contain one thought, or one idea.  Each slide should contain
about six to eight lines, and should have wide margins.
The bottom 2 inches of each slide
should be left blank---because this portion of the slide is often
blocked from the vision of those in the back of the room. 

Do not write out complete, long sentences on your slides.  Abbreviate
wherever possible.  If you are going to \TeX\ your slides, then
consider using Sli\TeX\ (which is a version of \TeX\ that contains
extra large fonts and other artifacts that are useful for preparing
overhead slides). Note, however, that a neatly prepared
handwritten slide is often as effective as a \TeX ed slide---and
handwritten slides give you the additional flexibility of colors,
arrows, and other graphic tricks. You should have only about one
slide per two to three minutes of speaking.  

Clearly the use of slides---or of software like {\tt PowerPoint} 
or {\tt Beamer}\footnote{{\tt Beamer} is a German product that is a
\LaTeX\ macro for giving talks.  It is very much like {\tt PowerPoint},
and has many of the same capabilities; but, since it is a \TeX\ product,
you can formulate all the mathematics that you like.}
solves several of the problems indicated above.  If you are projecting
your material on a raised screen, then you cannot be standing in front
of it.  And you will probably be using nice fonts, so handwriting
is not an issue.  In addition, the material will likely be formatted
in a nice way.  So ``blackboard technique'' does not come into the picture.
On the other hand, presenting your stuff on a screen presents new problems:
you might have too much on each screen, you might change screens too quickly,
you might have problems going back and forth to refer to earlier material.

I have seen talks in which the speaker simply printed out the text of
a fifty-page \TeX\ document onto transparencies---in 10-point 
or 12-point type.  Moreover, the speaker showed every single slide to the
audience.  What a disaster! First, this
is far too much material per slide---and none of it can be read. Second,
this is too many slides for a fifty- or a sixty-minute talk.

Edward Tufte is a notable advocate of good speaking, and of the skillful
presentation of graphics to illustrate quantitative information (see [Tuf1], [Tuf2]).
He also, since his retirement from the Statistics Department at 
Yale University, gives day-long presentations at hotels on how
to give a talk.  In these presentations he rails against {\tt PowerPoint}.
In particular, he makes these points about the software:		  
\begin{itemize}
\item Its use to guide and reassure a presenter, rather than to enlighten the audience;
\item Its unhelpfully simplistic tables and charts, resulting from the low resolution of early computer displays;
\item The outliner's causing ideas to be arranged in an unnecessarily deep hierarchy, itself subverted by the need to restate the hierarchy on each slide;
\item Enforcement of the audience's lockstep linear progression through that 
hierarchy (whereas with handouts, readers could browse and relate 
items at their leisure);
\item Poor typography and chart layout, from presenters who are poor 
designers or who use poorly designed templates and default 
settings (in particular, difficulty in using scientific notation);
\item Simplistic thinking---from ideas being squashed into bulleted lists; 
and stories with beginning, middle, and end being turned into a collection 
of disparate, loosely disguised points---presenting a misleading facade 
of the objectivity and neutrality that people associate with science, 
technology, and ``bullet points.''
\end{itemize}
See \verb@\https://en.wikipedia.org/wiki/Edward_Tufte@ for more on this matter. 
									  
One of the most important skills that you need to develop, both as a
teacher and as a colloquium or seminar speaker, is time management. 
You need to fit what you have to say into the time allotted.  People
will be monumentally irritated to watch a mature mathematician spend the
last thirty minutes of his/her fifty-minute talk pacing back and forth, scowling
at\index{talk!time management in}
the clock, and declaiming that he/she does not have sufficient time to
present his/her thoughts.  I have seen many such a speaker act as though
it were the audience's fault, or the university's fault, or perhaps
his/her host's fault, that he/she did not have more time.  What nonsense. The
speaker knew when he/she was invited---probably many months before---what
the parameters were.  Giving a fifty- or a sixty-minute 
colloquium talk is part of the academic game.  Learn to play by the rules.                                
       
Perhaps you are saying to yourself---or have said to yourself
in the past---``all good and well, but this speech-making stuff is
for joke-tellers and hams and showoffs.  I am a scholar.  I am not an
actor.''  Such a statement is a cop-out (if you will pardon the vernacular).
Nobody expects you to be Bob Hope.  Part of a scholar's existence
is to communicate---both in writing and in speaking.  The thoughts
in this section are intended to help you to enhance your abilities
with the latter.  Giving a talk is a personal affair; you should
do it in the fashion that best suits you.  But I hope that the ideas
presented here will help you to sharpen your wits and your technique.

%% Section 4.5

\section{Your Vita, Your Grant, Your Job, Your Life}

\subsection{The Curriculum Vitae}

\indent A businessman has his/her resum\'{e} and an academic has his/her Curriculum
Vitae (or {\it Vita} for short).  The Vita is your professional
history---it should give a quick sketch of who you are, where you
were\index{Vita!importance of the} 
educated, your professional experience, any honors that you have
earned, your scholarly accomplishments, and related materials. 
Usually you will include your publication list with your Curriculum
Vitae. 

Your Vita should {\it not\/} read like this:
\vspace*{-.1in}

\begin{quote}
Born on a mountain top in Tennessee. \hfill \break
Greenest state in the land of the free. \hfill \break
Raised in the woods so he knew every tree.  \hfill \break
Killed him a b'ar when he was only three.\footnote{From
{\it The Ballad of Davey Crockett,\/} Walt Disney Productions.}  \hfill \break
\end{quote}

\noindent All quite charming, but a Vita should {\it never\/} be in
paragraph (or stanza) form.  The material should be laid out in a
tableau so that the reader can quickly pick out the information he/she
needs.  Your name should be in boldface at the center top.  (I
recommend that you use your official name---the one on your birth
certificate.  Your friends may call you ``Goober,'' but you should save
that information for another occasion.)  Quickly following should be
your date of birth, your educational information, your address and
phone numbers and {\it e}-mail address, your employment record, key
honors earned, and so on.  An example of the first page of a Vita
appears later in this section.

Your publication list should be a separate section of the Vita.
Those who are especially careful separate published works from 
unpublished (or to-be-published) works and separate items in refereed
journals from items in nonrefereed publications (such as conference
proceedings); books are often listed separately; some people even
list class notes they have\index{Vita!parts of the}
prepared, or software that they have
written (if you are a numerical analyst or a specialist in algorithms
then this last would be essential).  Usually the items in a
publication list are given in approximate chronological order,
although some people use reverse chronological order.

Another section of the Vita lists grants or funding that the person
has received over the years.  For each grant, you usually list the
funding agency, the title and number of the grant, the amount of
money in the grant, and the year(s).

Often a Vita will include a section of invited talks or, if you
are quite senior, of major invited talks (that is, an hour speaker
at a national AMS meeting, or principal speaker at a CBMS 
conference, or a speaker at the International Congress).

Yet another section will list graduate students (Masters and Ph.D.) that
you have directed.  Another could list material describing your
teaching experience (courses taught, curricula developed, and so
forth).  Indicate your expertise
with computers---either software developed, or courses taught,
or other accomplishments.

You may wish to say something about your skill with foreign languages.
Have you done any translation work?  Are you well traveled?
Have you taught in another country?

Finally, some Vitae have a catch-all section with editorial activities,
service to professional societies, or anything else that the person
writing the Vita thinks may be of interest.

Your Vita is no place to be humble.  This document is the {\it
gestalt} that you present to the world.  Certainly do not
prevaricate---or even exaggerate---but be sure to tell the reader
everything that you want him/her to know about yourself.  

At the risk of sounding preachy, let me expand a bit on one of the points
in the last paragraph. When preparing the Vita, we all want to present
ourselves in the best possible light. There is a tendency to\index{Vita!honesty in} dress things up---beyond what is strictly kosher. Perhaps you
did not complete that French course---but you ate quiche Lorraine
once---so you write that you speak French. The letter from the journal to
which you submitted your latest paper says ``if you make the following
twelve changes then the referee will have another look at it,'' and you
list the paper as accepted. The NSF tells you that you are on the
``maybe'' list for a grant, and you put on the Vita that you have a grant.
People who have made these slips are not liars; they are just trying too
hard. Strive to avoid such exaggerations. Most departments check facts
carefully. Many schools only believe in publications that have appeared,
and for which there is a {\it bona fide\/} reprint (many schools have been
burned once too often in the past). If the Funded Projects
office at your school does not have the letter from the NSF, then your
grant does not exist. Worse, if you make such claims in your Vita and the
claims do not wash with the people evaluating your case, then the
situation will weigh against you. My advice is to be extra careful. 

\newpage

\begin{center}
{\tt SAMPLE VITA}  \hfill \break
-----------------------------------------------------------  \hfill \break
{\large \bf CURRICULUM VITAE}
\smallskip \break
for  \break  
\bf Clemson Ataturk Kadiddlehopper\index{Vita!sample of}
\end{center}
\noindent {\bf Date of Birth:}  March 15, 1947
\smallskip  \hfill \break
{\bf Home Address:}  17 Poverty Row, Faculty Ghetto, Iowa 50011  
\smallskip \hfill \break
{\bf Current Academic Affiliation:}  Department of
Mathematics \hfill \break
Walmart A\&M, Sam's Clubville, Iowa 50012  
\vspace*{-.03in}
\begin{tabbing}
{\bf Telephone:} \ \ \ \ \ \= \ \ \ \ \ \= (515) 294-6021\ \  \= (office) \\
               \>  \> (515) 373-3286     \> (home)    \\
               \>  \> (515) 294-6047     \> (FAX)  \\
{\bf e-mail address:}  \> \> {\tt cak@math.sam.edu} \\
{\bf Graduate Education:}\ \ \ \ \= {\sl Ph.D., Mathematics} \\ 
             \>  Montana Institute for the Tall, 1974 \\
	            \>  Thesis directed by Charles Ulmont Farley  \\
                           \vspace*{.06in} 
             \> {\sl M.S., Mathematics} \\
             \> Frisbee State University, 1971  \\
{\bf Undergraduate Work:} \> {\sl B.A., Mathematics}  \\
                     \> Joe's Bar and University, 1969        \\
{\bf Academic Positions:}  \= Assistant Professor, College of the Yodeling \\
                           \> \qquad Yuppie, 1974--1979  \\
                           \> Associate Professor, Steland Lanford  \\
                           \> \qquad University, 1979--1988 \\
                           \> Professor, Walmart A\&M, 1988--present  \\
{\bf Honors:} \=Neural Sediment Fibration Graduate Fellow, 1971--1974  \\
           \> Visiting Professor, Callipygean Institute of Tectonics, 1977 \\
           \>	Shinola Fellow, College of Good Hair, 1979  \\
           \>	Visiting Professor, Upper College of \\
           \> \qquad  Lower Academics, 1980  \\
           \>	Visiting Professor, University of Basic Bourgeoisie, 1986    \\
           \>	Visiting Professor, Hahvahd University, 1986  \\
           \> Honorary Lecturer, Crab Louie College, 1987   \\
\null 
------------------------------------------------------------  \hfill
\end{tabbing}

Now let us return to matters prosaic.
You must tailor your Vita to the circumstances.  I have been
teaching for 42 years.  Thus it would be crazy for me to list every
course that I have ever\index{Vita!what not to include in a} 
taught.  It would make more sense for my Vita to
list courses that I have created, or textbooks that I have written. 
On the other hand, if you are just starting out in the profession,
then you should indicate the depth and range of your
teaching experience and certainly indicate your facility with
computers, both in the classroom and outside it.  A beginner will
probably have directed few if any graduate students.  That is not a
problem, since such activity is not expected.  Do, however, be
complete in describing your other activities.
\medskip 

\subsection{The Grant}

\indent Funding is available for many different types of
activities that a mathematician might undertake. These range
from quite specific, goal- \break oriented
projects\index{grants!availability of} that are funded by
industry all the way to grants from the NSF (National Science
Foundation) to encourage pure research in abstract
mathematics. There is also funding from the Department of
Defense, from DARPA (an arm of the CIA or Central Intelligence
Agency), from NASA (the National Aeronautics and Space
Administration), from NIH (the National Institutes of Health),
from DOE (the Department of Energy), from NSA (National
Security Agency), from the Simons Foundation, and from many
other sources as well. Granting agencies such as the NSF have
considerable funds to encourage work on the mathematics
curriculum---from developing new ways to teach calculus to
revising substantial blocks of undergraduate mathematics
education.

Given the range of activities that granting agencies are
willing to fund, and given the variety of different potential
sources of grants, I could discuss grantsmanship at length. I
shall content myself here with a few general precepts that
should apply to virtually any grant application that you may
write.

Read the prospectus for the program to which you are applying.
Doing so, you will learn what the program is looking for, what
particulars should\index{grant!prospectus for} be itemized in
the proposal, what page limits will be enforced, and when the
proposal is due. Learn about what type fonts are acceptable,
what margins the pages should have, how long the Curriculum
Vitae portion of the proposal should be, how long the
references section should be, how many pages should address
previous work, how many pages should address new work \dots
{\it and so forth.\/} Grant proposal writing is not a free
form activity. Get the rules straight before you begin.

The main issue that is in the air when your grant proposal is being
evaluated is your credibility:  {\it can\/} you do the work that is
being proposed, and {\it will\/} you do the work that is being
proposed?\index{grant!credibility of your proposal for}
Given your stature, your abilities, and your track record,
is it clear that you can work on these problems (be they research or
education)?  Can you solve them or make progress on them?  Are
you capable of evaluating your own progress?  Finally, can a case
be made that you are {\it the right person\/} to work on this project?
Or will the work be done as a matter of course by others (if indeed
it is worth doing at all)?

You must walk a delicate line here. On the one hand, you want to make it
clear to the potential granting agency that you know this subject inside
and out, that you know the existing literature, and that you have a good
program for proceeding. You want to demonstrate that you are already
engaged in some version of the proposed activities. On the other hand, you
do not want to make it sound as though you have already solved the
proposed problem. You also want to give the strong impression that you are
working on substantial problems of real significance; these should be
problems for which even partial results will be of interest. But it should
be plausible that you are up to the task. In particular, if you propose to
prove the Riemann hypothesis, then you will have a difficult time making\index{Riemann hypothesis}
your case. After all, many of the bigshots are working on this problem; if
they cannot make inroads then how will you?

Generally speaking, granting agencies will not provide funds to help
you to learn something new, or to retool.  Thus 
you should make a case that you are already engaged in
the proposed project, that you have a grip on it, and that you have a
viable program.  If those considerations entail your learning
something about nonlinear elliptic PDEs, then by all means you
should say so.  But a grant proposal that reads (in effect) ``I'm
tired of studying finite groups so now I'm going to do symplectic
geometry'' just will not fly.

As already noted, you must prove that you are up on the relevant
field---not just what is in the published literature but what is
available in preprint and other tentative form.  For the most part,
grants are refereed by your peers.  These will be peers who are on
the cutting edge.  They will judge you by their own standard---the
standard by which they themselves would expect to be judged.

When writing a grant proposal, you must walk another delicate line. You
will usually have a page limit. You simply cannot go on at length, or have
extensive digressions, or have verbose introductions\index{grant!page
limitation for application} or chatty conclusions. But you must make the
proposal as easy to read, and as self-contained, as possible. If your
proposal engages the referee's interest, and teaches him/her something,
and does not force him/her to keep running to the library to figure out
what you are talking about, then you will be at a real advantage. If,
instead, your prose is a bore and the referee has to slug his/her way
through it, then your proposal is likely to be penalized.

Do not be afraid to telephone the granting agency to which you are applying
and to talk to the program officer. Many grant programs, and many program
officers, encourage this activity. By talking to a program officer, you
can better focus and tailor your\index{grant!program officer of} proposal
to the goals of the intended program, and you will not waste the program's
time with a proposal that is completely off the wall.

Proposals in mathematics education and curriculum often require a section
on ``self-evaluation'' and a section on ``dissemination.'' You should
not\index{grant!self-evaluation in} (in\index{grant!dissemination in}
the self-evaluation portion of your proposal) say ``We'll see how happy
the students are at the end by distributing teaching evaluations.'' You
also should not (in the dissemination portion of your proposal) say ``I'll
go to conferences and talk about this stuff with my buddies.'' I have seen
both of these in serious proposals, and they do not work. Both approaches
are too facile, and show no imagination and no effort. Good
self-evaluation programs often involve motivational psychology experts
from your institution's School of Education, tracking of students after
they leave the experimental program, exit interviews, and many other
devices. Good dissemination programs often involve writing a textbook for
publication, creating a newsletter, setting up a Web page, organizing
workshops, and so forth. I am not necessarily advocating any of these
devices. I am merely explaining how the world works. 

Self-evaluation and dissemination play an implicit role in a research
proposal as well. Your report on previous work will give an indication of
your ability to evaluate your own progress. The scientist who says ``In
the last five years I tried a lot of things but nothing panned out'' shows
both bad judgment and an inability to learn from his/her work. The
dissemination aspect of a research proposal is reflected in your
publication record, your invitations to speak at conferences or colloquia,
and your collaborative activities.  NSF proposals now require
a section on dissemination and impact.

Before you submit your proposal, run it through a spell-checker. Check and
recheck the grammar. Show it to a senior colleague. Proofread\index{grant!proofreading of the proposal} it more times than you think could possibly
be necessary, and then proofread it once more. The reviewer will be
phenomenally irritated to read a proposal that appears to have been
prepared hastily or sloppily. Your proposal should be as slick as glass.
It should be a pleasure to read, and it should get the reviewer excited
about and interested in what you are proposing to do. 
\smallskip 
\hfill \break

\subsection{Your Job}

\indent At one time or another, most of us will have to apply for a job.
Let me first say a few words about applying for an academic job.

When\index{job!applying for a} 
you apply for a job at a college or university, you send in your
Vita (discussed above) and a cover letter.  The cover letter should
be brief (well under a page); it should identify you, your 
present\index{job!cover letter for application} 
position, the type of position you seek, and your areas of
interest.  No more.  A cover letter with a multitude of exclamation
points, shamelessly extolling your virtues as a teacher and your
bonhomie as a colleague, is highly inappropriate.  A sample cover
letter appears later in the section.

I would be remiss not to note that the AMS (American Mathematical
Society) has created an OnLine utility called {\tt mathjobs}.  This 
is a device for applying for a job.  If you are a candidate, then
you arrange for all your materials (your letters of recommendation,
your Teaching Statement, your Research Statement, your Vita, your cover letter, and
so forth) to be uploaded to the {\tt mathjobs} Web site.  Then
any school that is interested in you has complete access to all
your materials.  There are many advantages to this system:  {\bf (i)}  all
your materials are in one place, {\bf (ii)}  the materials cannot
be lost or misplaced, {\bf (iii)}  several people can view your
materials at once.  It is safe to say that {\tt mathjobs} has become
a cornerstone of the job marketplace.  Not every school participates in
{\tt mathjobs}, but a great many do.  If you are a senior person
applying for a position, then using {\tt mathjobs} may be inappropriate.
The school to which you are applying will let you know.

If you are a beginner in the profession, with a short publication
record, then you might include some of your preprints with the
job\index{Vita!references in}
 application. Your Vita should also list your ``references'' or
``recommenders.'' This point is vital, and many job applicants
overlook it. Before applying for any job you should approach
three or four people (for a senior job it could be six or
eight, or even more) and ask whether they are willing to write
in support of your application. Ideally, these should be
prominent people in your field, whose names will be recognized
by those evaluating your dossier.\footnote{If you are a very
senior candidate for a position, such as I was for my most
recent job, then you do not contact the letter writers. The
hiring institution will do so.} Once you induce a group of
such people to agree to write,\footnote{The notion of the
candidate {\it asking\/} people to write for him/her seems to
be a peculiarly American custom. In many
countries---especially in Europe---the hiring institution does
all the solicitation of letters.} then include their names,
business addresses, business phone numbers, {\it e}-mail
addresses, and fax numbers in a section of your Vita called
``References.'' These days---especially if you are a job
candidate just a few years past the Ph.D.---you should have
one or two letters from people who can say something specific
and positive about your teaching. Their names should be listed
in your References section as well. It is now commonplace for
job candidates in the United States to include the ``AMS
Standard Cover Sheet'' in the dossier; this form may be
obtained\index{AMS Standard Cover Sheet} from most issues of
the {\it Notices of the AMS}.\footnote{Many people,
especially young job candidates, include in their job
application a one or two page statement describing their
research\index{teaching philosophy, statement of} program;
many young people also include a ``statement of teaching
philosophy.'' The first of these can be quite helpful to a
nonexpert who is endeavoring to evaluate the dossier: a good
research statement can at least guide such a reviewer to an
appropriate expert colleague who can comment on the case in
detail. I dare not say whether a ``teaching philosophy''
statement has any real value; there is little grass growing in
this subject area, and you will have trouble finding anything
interesting or original to say. On the other hand, some
schools {\it require\/} a statement of teaching philosophy; in
such a circumstance you must do your best to write something
thoughtful and thought provoking.}

Some job candidates arrange to have a sample of their teaching evaluations,
or\index{job!teaching credentials for}
passages from their present institution's Teacher Assessment
Book (such a book is often published by the campus student
organization) to be included in their dossier.  Such an inclusion
can help to lift the dossier out of the ordinary, and will
add substance to the letters that praise the candidate's teaching
abilities.

If you use your imagination, you can probably think of all sorts
of things that might be included in your dossier in candidacy for
a job.  I recommend that you consider each one cautiously.  
People on hiring committees these days often must wrestle
with 500 job applications in a season.  A big, fat dossier 
will\index{job!what to include in the application for a}
just turn off a weary committee member.\footnote{It is quite
common these days for a high school student applying to college
to include a video, or a videotape, in his/her application.
Such an item would be quite inappropriate in a professional
job application.}  So do not leave
anything important out of your dossier, but think carefully about
what you do include. 

If you make your application in the manner described in the preceding
paragraphs, then a school to which you apply will know just
how to process your paperwork.  Once it has your cover letter and
Vita, it can start a file on you.  Then it has a place to put the
letters of recommendation as they come in.  And, because you have
included a list of references, the school will know when your dossier
is complete.

If your application is for any position beyond a beginning lectureship, and
if you make the ``short list,'' then you will likely be invited to give a
talk and to meet your potential future\index{job!short list for}
colleagues. Let me not mince words: this is a make-or-break situation.
Dress well (not as though you were entering a ballroom-dancing contest,
but rather as though you are taking the situation seriously). Give a
polished, well prepared talk (see Section 4.4 on how to give a talk).
Think in advance about some of the topics of conversation that may come up
when you meet your new colleagues. Be prepared to describe your research
conversationally to a small group of nonexperts; be able to say in five or
ten minutes what you do, how it fits into the firmament, who are some of
the experts, what are some of the big questions.

{\it Be prepared to say who in this new department has interests in
common with yourself; with\index{job!talk given in candidacy for}
whom you think you might talk mathematics; who might become your
collaborator.\/} Do not underestimate the significance of this circle of
questions. Do {\it not\/} say ``Oh, I talk to everybody; I'm the Leonardo
da Vinci of modern mathematics.'' Such a statement is not credible; utter
it and you will surely send yourself plummeting to the cellar of the short
list.

Think over your ideas about teaching, about the teaching reform
movement, about teaching with calculators or computers, about
teaching students in interactive groups, and about any other topics
that may arise.  Some schools have special problems connected with
the teaching of large lectures; be prepared to share your views on
that topic.  Other schools have special tutorials for calculus
students; be prepared to chat about that topic as well.  Some schools
like to conduct a formal interview, with a few of the
senior\index{job!interview for} faculty asking you direct questions
about your research, your teaching, your attitudes about curriculum
and reform, about teacher/student rapport, or anything else that is
in the air at the time.  It makes a dreadful impression 
if you are inarticulate, do not seem to know your own mind, or simply 
have not given any thought to these matters.  I am not advocating 
that you go to your ceremonial job interview
with a sheaf of notes {\it in your hand;\/} I am instead
advocating that you go with a few note cards {\it in your head}.

When a school decides to offer you a job, the chairman will usually
telephone you, or send you an {\it e}-mail message
followed\index{job!offer of a}
by a phone call.  At that time he/she may discuss salary, teaching load,
computer equipment, startup funds, health insurance, the retirement annuity, or
other perquisites.  You may wish to take the opportunity to ask about these or
about other concerns.  

Important information about you can be lost in the Vita---especially
if your Vita is long.  If you are a graduate student applying
for a first job, and if you have won a teaching award for ``Best TA,''
then certainly mention that encomium in your cover letter.  If you are 
a few years from the Ph.D.\ and the
holder of a Sloan Fellowship, or an NSF Postdoc, then you should
mention these honors in your cover letter.  You do not want your cover
letter\index{cover letter, what to include in a}
to look like a flyer for your local supermarket, 
but you want it quickly to lead the reader to your strong points.
			       
\begin{center}
{\tt SAMPLE COVER LETTER}  \hfill \break
------------------------------------------------------------ \hfill \break
\end{center}
\rightline{November 22, 1996}

\leftline{Noodles Romanoff, Chairman} 
\leftline{Department of Mathematics}
\leftline{Little Sisters of the Swamp College}
\leftline{Sanctuary, Oklahoma 23094}
\vspace*{.1in}

\leftline{Dear Professor Romanoff:}
\leftline{ \ }

I wish to\index{cover letter, sample of} 
apply for a faculty position, at or near tenure,
in your department.  I received the Ph.D.\ in Mathematics in 1974. 
I am a geometric analyst, with specializations in complex function
theory, harmonic analysis, and partial differential equations.  My
Vita is enclosed.  It includes my list of references.  I currently
hold the position of Instructor of Mathematics at Brouhaha Subnormal
School in Wichita Falls.

Please note that my research is supported by a grant from the
Normative Sodality Agency.  I am also a recipient of the
Mudville Distinguished Teaching Award.  I have strengths in
research, teaching, and curriculum. 

I look forward to hearing from you.  
\hfill
\bigskip \break

\leftline{Sincerely,}
\vspace*{.6in}

\leftline{Shrimp Chop Suey}
\leftline{Instructor of Mathematics}
\begin{center}
------------------------------------------------------------  \hfill \break
\end{center}

If you are applying for a job in the private sector---say at Texas
Instruments, or AT\&T, or Aerospace Corporation, or Microsoft---then the
application process may be a bit different from the process in an
academic setting.  An industrial organization is probably not
interested in letters of recommendation from the Universit\'{e} de
Paris, nor in binary comparisons with famous young algebraic
geometers.  A business resum\'{e} is different from an academic
Vita.  Go to your local bookstore and purchase a book on how to write
an effective resum\'{e} (see, for example, [Ad3]), how to write a cover
letter (see [Ad1]), and how to apply for a job (see [Ad2]).   
Although working in industry certainly will involve communication
skills, it probably will not involve much classroom teaching.  
The interview\index{communication skills}  
for an industrial job will likely be even
more crucial than the interview for an academic position.  Consult
acquaintances who have been through the process so that you can be
well prepared. 
\smallskip \hfill \break

\subsection{Your Life}

\indent In the abstract, the rewards for good writing may seem far off and
vague; instead you can see clearly how a well-written Vita or grant
proposal could lead to just deserts. Your Vita is a tool for helping you
to find employment, or a promotion, or to achieve some other goal. Your
grant proposal is a way to seek funding. I have also discussed how to find
a job. The principles of good writing described in other parts of this
book apply just as decisively to these practical matters: express yourself
directly, cogently, and briefly; do not show off; know what you are
talking about; and (paraphrasing Jimmy Cagney) plant both feet on the
ground and tell the truth.

I have seen many a Vita in my time. One of these contained a page entitled
``Cities Beginning with the Letter `Q' in which I Have Spoken Fewer than
Five Times.'' Another listed an uncompleted mystery novel. Yet another
listed\index{Vita!things not to include in a} forty (count `em)
collaborative papers that were incomplete and in progress. One Vita by
Mathematician $X$ listed poetry, both published and unpublished, that was
written to $X$, by $X$, about $X$. Another Vita listed the subject's (not
very happy) marital history. Your Vita is a business document. This piece
of paper is a pr\'{e}cis of your professional life. Think carefully about
what you put into it and how you organize it.

Your grant proposal is a manifestation of your professional values,
what you are all about, and what you are trying to do.  As you
develop it, read it with the eyes of your potentially most critical
reviewer.   

Here we are discussing writing with immediate impact, and with a direct
effect on your life. This is writing that you wish to succeed because it
must. Even more than in your other writing, you will want to strive to
make each word count, and to force each sentence to say precisely what is
intended. The critical skills discussed in this book should help you in
these tasks.

%% Section 4.6

\def\email{{\it e}-mail}

\section{Electronic Mail}

For many of us, electronic mail (or \email\ for short) has become an
important part of life. The technology of \email\ has enabled us to carry
on extended conversations with people all over the globe. We
can\index{electronic mail!significance of} engage in topic-specific
discussion groups, conduct business, develop friendships, and even have
fights via the Internet. Perhaps more significantly, we can conduct
mathematical collaborations with people 10,000 miles away, in some cases
with people whom we have never even met. You may actually (though I
encourage you to exercise this option with discretion) send an \email\
blind to a professor\index{electronic mail!uses of} at MIT and say
``Hello, I'm so and so. Do you know the answer to the following
question?'' I have occasionally engaged in this speculative activity and,
more often than not, I have received a useful answer.

Several years ago I was writing a series of papers with two collaborators,
one of whom is usually in Los Angeles and the other usually in Canberra,
Australia. During this last year, one of us spent a leave in Berkeley,
another took a leave in Wuppertal, Germany, and the third changed jobs. We
did not miss a beat, because \email\ is oblivious to\index{electronic
mail!collaboration via} these moves. Marshall McLuhan [McL2] died too
soon: the global village is finally here in spades.

G.~H.~Hardy and J.~E.~Littlewood carried on what is by now the most
famous, and certainly the most prolific, mathematical collaboration
in history.  Usually in two different locales (one in Cambridge, the
other Oxford), they conducted their collaboration by regular post
(now\index{Hardy, G.\ H.}\index{Littlewood, J.\ E.} known as ``snail
mail'').  Their hard and fast rule was that if one of them received a
letter from the other, he was under no obligation to open it---right
away or at any time.  Many a letter was thrown into a pile, not to be
read then or perhaps ever; this to guarantee that the recipient could
think his/her own thoughts, and not be interrupted.  To my mind, \email\
is a bit different:  it gets right in your face, once or several
times a day.  Once you have determined (perhaps by looking at 
the Subject line) what a particular \email\
message is and whence it came, then you are looking at it.  The
temptation is to read it.  As McLuhan taught us [McL1], ``the medium
is the massage.''

For many purposes, communicating via \email\ is preferable to communicating
by telephone. For an \email\ message has the immediacy of a telephone call
without any of the hassle of playing ``telephone tag,'' talking to voice
mail, or patiently explaining your quest to a secretary. Many of us find
that we send more {\it e}-mail than we do letters,\index{telephone tag}
and we use \email\ more often than we use the telephone.

Because the use of \email\ has become so prevalent, we must all learn some
basic etiquette of the \email\ system. As with many other activities in
life, \email\ is something that we can benefit\index{electronic mail!etiquette for} from if we give it just a few moments of reflection.

\begin{itemize}
\item Be sure that your \email\ messages go out with a complete header.
This header should include a ``From'' line, indicating your identity
and \email\  address, and a ``To'' line, indicating the identity
and \email\  address of the person to whom the message is being sent.
It is not mandatory, but is highly desirable, for you to include
a ``Subject Line'' in the header.  Many busy people receive 50 or
more \email\  messages per day;  you do such people a great favor
by\index{electronic mail!header in}
helping them quickly to identify and sort their\ \email.

\null \quad \ It is convenient to use {\tt Eudora} or {\tt Outlook} to manage \email.
These popular software devices will sort your email messages and put
them in pre-determined directories on your hard drive.   They can
be very useful.  Of course these software utilities depend on the Subject
line in order to do their jobs.

\null \quad \ Some systems allow you to strip away all or part of the
header of your \email\  message.  I urge you not to do this as such an
action is unprofessional and rude.  Sending anonymous \email\  is no better
than sending anonymous hate mail.

\item  I often receive \email\  messages that say ({\it in toto})
 ``Yes, I agree with you completely'' or ``Right on'' or ``There you
go again!''  I love fan mail as well as the next person, but I often
cannot tell what such \email\  is about.  Do yourself and your
correspondent a favor and either {\bf (i)} include the \email\  message
to\index{electronic mail!completeness of}
 which you are responding in your reply or {\bf (ii)} at least
include a sentence or two indicating to what you are responding.

\item  Sign your \email\  message with your full name.  Signing off
with ``See you later, alligator'' or ``That's all, Folks'' is momentarily 
amusing, but it often forces your recipient to search the header
of the\index{electronic mail!signature to}
message to determine whose pear-shaped tones he/she is reading.
Such a search is sometimes frustrating, and irksome to boot.

\null \quad \ The best possible ``signature'' to an {\it e}-mail message is something
like this:
\begin{verbatim}
*******************************************************
* Steven G. Krantz  (314) 935-6712  FAX (314) 935-6839*
* Department of Mathematics, Campus Box 1146          *
* Washington University in St. Louis                  *
* St. Louis, Missouri 63130-4899  sk@math.wustl.edu   *
******************************************************* 
\end{verbatim}
Of course you do not want to type out this mess each and
every time you send an \email\  message.  The operating system {\tt UNIX}
makes it easy for you to avoid such tedium.  {\tt Gmail} and
other mail utilities also make it straightforward to formulate
and install such a message.

\null \quad \ I have arranged for all my university \email\ to be
forwarded to my {\tt gmail} (i.e., Google mail) account.  I am happy
to say that {\tt gmail} is a powerful system that will allow
you to create a signature and to customize your emails in a number
of useful and attractive ways.  I note that {\tt gmail} has a terrific
spam filter.  And {\tt gmail} allows you to send attachments of 
any size.  It is a great system.

\item  The good news about \email\  is that it is a lot like
conversation.  It is spontaneous, natural, and candid.  The bad news
about \email\ is that it is a lot like conversation---without the
give and take of an interlocutor.  Thus we are tempted to type away
madly,\index{electronic mail!and conversation}
at high speed, having no care for corrections or
proofreading.  This is a big mistake.

\null \quad \ Proofread each \email\ message before it goes out. If the message is
important then proofread it several times. Most \email\ editors are easy
to use. In the {\tt UNIX} setting you have a choice: the {\tt PINE}
editor\index{electronic mail!proofreading of} is self-explanatory, and
much like a word processor; in the {\tt ELM} environment you can customize
the editing environment, using {\tt emacs},\index{emacs@{\tt emacs}} or {\tt
VI}, or another editor of your choosing. In any event, learn to use the
editor on your system and {\it use it.\/} Correct misspellings (many an
\email\ editor is equipped with a spell-checker) and misstatements. Clean
up your English. Some \email\ messages that you send will have the
permanence of a hard copy written letter. Send something that will reflect
well on you.

\null \quad \ In fact, when I am writing something of great importance, I compose
it on my home computer---on the text editor with which I am most
familiar (see Section 6.3 for a discussion of text editors).  I do
this in part for psychological reasons.  When I compose on my home
computer,\index{electronic mail!and the personal computer}
I do not worry about the system hanging or going down; I do
not worry about taking a break and being thrown off the system; and I
am using a writing environment with which I am thoroughly conversant.
I can use my spell-checker, my CD-ROM dictionary and thesaurus, and
other familiar resources to put the document in precisely the form
that I wish.  I also can sleep on the matter before I send the
document.

\null \quad The next morning, I bring the document to work on a flash drive, upload it
to the system (see Section 6.8), and then pull it into an {\it e}-mail
message using operating system commands. This methodology is a valuable
tool.

\item  Implicit in the preceding discussion is a major liability of
\email.  Too easily can you write something in haste in the \email\
environment\index{electronic mail!haste in}
and then just send it off---it only requires a key
stroke or two!---and then it is gone.  You cannot retrieve it.

\null \quad I once had a rather significant fight with another mathematician. He wrote
me a letter taking me to task for something that I had done. Fortunately,
this event occurred in the days before \email. I wrote a hasty and heated
response (in hard copy, for that was all that we had at the time) telling
this person that he was misguided and mean-spirited, and dropped it in the
department's outgoing mail tray. An hour or two later, I pulled the letter
from the mail (I had been stewing about it all the while), and penned a
milder version of the heated letter. This revision process repeated itself
throughout the day. By the end of the day, I had put in the mail a letter
of apology, acknowledging my error and thanking my correspondent for
calling it to my attention. I have always been happy for this outcome.
With \email\ the story would have ended differently, and badly.

\item Try to keep your \email\ messages brief. Of course I realize there
are times when you are circulating a report or writing a
detailed\index{electronic mail!brevity in} formal analysis of some
situation; in such circumstances, it may be appropriate to go on at some
length. But, most of the time, when writing \email, you are sending a
memo. Thus make it quick. Often, on the computer, we tend to do things
just because we can. Writing an \email\ message is a lot like talking, but
without the reality check of having someone interrupt you from time to
time. Thus you must show some good sense: say what you have to say, say it
cogently and completely and {\it concisely}, and then cease.

\item  You can easily forward any \email\  message that you receive
to anyone that you like.  I am astonished at the extent to 
which\index{electronic mail!forwarding of} this power is misused. 
When you receive a hard copy letter of recommendation in the
mail---for a tenure case, say---you probably do not make 50
photocopies of the letter and send them off to 50 different
mathematicians.  First, such an action would be rude; second, it
could have legal repercussions.  For a written letter, the sender
owns the contents and the recipient owns the piece of paper and that
particular {\it form\/} of its contents (that is the law).  A similar
legal protocol has been proposed for \email, although at this writing
the legislation has not been approved.  What I am discussing here is
not so much the law as common sense and common courtesy.

\null \quad People forward \email\  all over the place, with hardly a thought for
the consequences.  The courteous thing to do is to ask the author
before you forward anything.  Many people send me \email\ messages
that say ``Please delete this message after you have read it'' (the
implicit message here is ``Don't forward this to anyone!'').  I am
always punctilious about adhering to such a request, and I hope that
others are similarly considerate of my requests for discretion.

\item Electronic mail is not as secure as other forms of communication. Any
superuser on your system can eavesdrop on your\index{electronic mail!security of} 
\email, and computer bandits can break into the system and
perform all sorts of nasty deeds. Thus you need to exercise some restraint
with respect to what you say over \email. Many of us use hard copy letters
and the telephone for the most delicate matters.

\item In the early days of \email, a user had to be careful of line length:
lines longer than 80 characters were often truncated by\index{electronic
mail!line length in} either the sending or the receiving editor. Given
that some characters could be added in transit, it was best to keep lines
to 72 characters in length. Most editors and mail spoolers now can handle
longer lines, but careful users still keep lines no longer than 72
characters.

\null \quad Some mail spoolers and {\it e}-mail editors insert line breaks into {\tt
ASCII} files that they receive. (The likelihood of this inconvenience
increases if your lines are long.) Thus a perfectly good \TeX\ command
like \verb@\smallskip@ could be transmogrified to \verb@\smal@ at the end
of one line and \verb@lskip@ at the start of the next line. If you are
lucky, you will catch this glitch with a spell-checker. Of course you can
bullet proof your file by {\tt UUENCODE}-ing it or {\tt zip}-ing it before sending it.

\null \quad If you send a file to a friend with lines that are longer than 80
characters, and if he/she endeavors to print it out cold, then the lines
are likely to be chopped off in the hard copy. The industrious high-tech
recipient will reformat each paragraph before printing---using
\verb@<Esc>-q@ in {\tt emacs} or an analogous command on other systems.
Other recipients will miss a lot of information.

\null \quad Also avoid beginning any line with ``From'' or ``from,'' 
as this word is
proprietary\index{electronic mail!use of `from' in}
to \email\ (and will result in unwanted characters being
added to your document during the \email\ transmission process).
For example, in order to protect the special use of ``from,'' \email\
will replace it with ``\verb@>@from'' when it occurs at the beginning
of a line.
\end{itemize}

Electronic mail, or \email, is a marvelous tool.  It has affected the
mathematical infrastructure, and has altered the way that many of us
collaborate and communicate.  If each of us would exercise just a little
\email\ etiquette, then the annoyances attendant to \email\ would
be minimized.

%%%%%%%%%%%%%%%%%%%%%%%%%%%%%%%%%%%%%%%%%%%%%%%%%%%%%%%%%%

%% Chapter 5
\chapter{Books}

\begin{quote}
\footnotesize \sl Some books are to be tasted, others to be swallowed,
and some few to be chewed and digested.
\smallskip \hfill \break
\null \mbox{ \ \ } \hfill \rm Francis Bacon \break
\null \mbox{ \ \ } \hfill \rm {\it Essays\/} [1625], Of Studies
\end{quote}

\begin{quote}
\footnotesize \sl No man but a blockhead ever wrote except for money.
\smallskip \hfill \break
\null \mbox{ \ \ } \hfill \rm Samuel Johnson \break
\null \mbox{ \ \ } \hfill \rm quoted in Boswell's {\it Life of Samuel Johnson}
\end{quote}
\vspace*{.03in}

\begin{quote}
\footnotesize \sl I never think at all when I write \hfill \break
nobody can do two things at the same time \hfill \break
and do them both well.
\smallskip \hfill \break
\null \mbox{ \ \ } \hfill \rm Don Marquis 
\end{quote}

\begin{quote}
\footnotesize \sl A writer and nothing else is a man alone
in a room with the English language, trying to get
human feelings right.
\smallskip \hfill \break
\null \mbox{ \ \ } \hfill \rm John K.~Hutchens 
\end{quote}

\begin{quote}
\footnotesize \sl The writer who loses his self-doubt, who
gives way as he grows old to a sudden euphoria, to
prolixity, should stop writing immediately:  the
time has come for him to lay aside his pen.
\smallskip \hfill \break
\null \mbox{ \ \ } \hfill \rm Colette 
\end{quote}

\begin{quote}
\footnotesize \sl 
You can't polish cow chips.
\smallskip \hfill \break
\null \mbox{ \ \ } \hfill \rm paraphrased from Lyndon Johnson 
\end{quote}

%% Section 5.1
\markboth{CHAPTER 5.  BOOKS}{5.1.  WHAT CONSTITUTES A GOOD BOOK?}
\section{What Constitutes a Good Book?}
\markboth{CHAPTER 5.  BOOKS}{5.1.  WHAT CONSTITUTES A GOOD BOOK?}

Mathematics books are written all the time. Go to the library and pull one
at random off the shelf. Looks pristine, does it not? Or perhaps only the
first fifty pages show signs of reading. Many\index{book!math} an author
lavishes all his/her enthusiasm and creativity and energy on the first
part of his/her book; he/she then runs out of steam for the remainder.
Unfortunately, it is the reader who suffers the consequences.

Writing a good book requires more effort than many authors are willing to
give to the task. Writing a good {\it mathematics\/} book requires special
insights and skills. In my view, the hard work is worth it.
When\index{book!effort to write} you write a good mathematics paper, it
is only read by a small group of people. But write a good book and a lot
of people will see it. The book is a way of planting your flag, of putting
your stamp on the subject, of sharing with the world the fruits of your
hard labor.

My advice is not to consider writing a book until you have tenure and are
established somewhere. The task is just too time consuming, and is often
not construed as a positive contribution toward the tenure decision. Put
differently, and a bit simplistically, the view of the world is that an
Assistant Professor should be writing research papers and becoming
established in the research community. Once you have done that, and
achieved tenure status, then you have the leisure to consider other
pursuits.

Now let us consider what makes for a good book. First, and foremost, you
must have something to say. If you are only repeating, or paraphrasing,
what has been said before then you are contributing 
nothing\index{book!components of a good} to the subject. 
Second, you must have a plan for
saying it. The best method for writing a book is to immerse yourself
thoroughly in the subject. The book itself becomes your ``world'' for a
couple of years. A place to begin is to write a detailed outline of the
book. Begin by writing chapter headings. Then fill in some section
headings. After a while, the juices begin to flow and you will find
that\index{book!outline of} you cannot write fast enough to keep up with
the outline developing in your head.

Once the book outline is written, it should be emblazoned on your frontal
lobes. Carry it with you (in your head) all day long. I find, when
writing, that I am constantly jotting down thoughts or topics or phrases
that occur to me throughout the day. These can arise in conversation, or
in lectures, or while daydreaming. If you are thoroughly involved with the
project, then they come up.

Once you have a detailed plan of what you are going to do (and you are not
bound to this plan, for it will evolve as your work unfolds), then you
should begin to write. Write a chapter at a time. Completely immerse
yourself in each chapter. If, while writing Chapter 3, a thought occurs to
you about Chapter 6, then make a note.\index{book!method for writing} You
can, especially in the computer environment, jump from one chapter to
another. But the process can become confusing. Safest is to make a
note---in a notebook perhaps, or in a computer file that you can pull up
instantly. Then, when you begin work on Chapter 6, you have all your notes
to work from.

Remember, as you write, that you are taking material that you have
thoroughly digested and internalized and are presenting it to your
readers---many of whom are tyros. Thus you must perform a reverse
evolution to put yourself in the shoes of the learner. This may be hard to
do at first, but it is a worthwhile exercise: it helps you to see as a
whole how the subject is built and what questions it answers. It helps you
to understand motivation and foundations.

Keep in mind that organization is a powerful tool. I have seen too
many\index{book!organization of} math books that state lemmas
parenthetically. Here is an example:

\begin{quote}
We thus see that every pseudo-melange is a hyper-melange.  (We use
here the fact that every pseudo-melange is complete.  {\it Proof:}
Let ${\cal M}$ be a pseudo-melange.  Calculate its first
Sununu cohomology group, etc.)  \qquad \bad
\end{quote}

Here the author is writing a love letter to himself. If you write
such\index{book!ordering material} an epistle, then few will read it and
fewer still will derive anything from it. Especially when writing with a
computer, you can always add a lemma---wherever it is needed---and add
suitable connecting material as well. Do not succumb to the temptation to
skip this part of the writing regimen. Most of the process of developing a
book consists of attending to details like making sure that all your
lemmas and definitions are in place before you need them. You {\it must\/}
attend to these matters.

To recast what I have been trying to say in the last few paragraphs,
the first blush of writing can be lots of fun.  You organize a
subject in your head, or on paper.  In a flurry of enthusiasm, you
write a draft on paper.  You see the subject begin to shape up
as you, and only you, see it.  You begin to take possession of this
circle of ideas.  The process is exciting and stimulating.

But then the moment of truth arrives.  If you want to turn this
random sequence of meditations into a publishable book, one that
people will {\it read,\/} then some hard work lies ahead.  You must
go through the MS line by line, detail by detail, attending
to context, syntax, logic, motivation, and many other details
as well.  You will proofread the same passages over and over again.
Frequently, you will have to swallow your pride and rewrite an entire
section, or reorganize an entire chapter.  The revision process
is hard, tedious work and not for the faint of heart.

You must put yourself in the shoes of the first year graduate student, or
whoever represents the ground floor of those who might read your book.
Where\index{book!level of} will such a reader get hung up, and why? What
can you, as the author, do to help this person along?

Finding an original way to develop the proof of the latest theorem in
your subject is always a pleasure.  Reorganizing that material in a
new way, for your six or eight close buddies in the field, is
rewarding.  Much less stimulating is writing a chapter of motivation
and background material. But, thinking in terms of the longevity and
impact of your book, you must learn to admit that both of these tasks
are of paramount importance.  The latter is not going to have people
buying you drinks at the next conference, but it will help your book
to have an impact on the infrastructure of your subject.

To summarize, what makes for the writing of a good book is hard work
and unstinting attention to detail.  Frequently the work required is
tedious, and you will ask yourself why you cannot assign it to a
secretary or a graduate student.  The answer is that you are
producing {\it your book,\/} and it is for the ages, and you want it to
come out right.  

%% Section 5.2

\section{How to Plan a Book}

The business of planning a book has been touched on in the previous
section.  Here we flesh it out a bit.

A common way to develop\index{book!planning of} a mathematics book is
first to teach a course in the subject area. Indeed, teach it several
times. Develop detailed notes for the course. Polish them as you go. Get
your students and colleagues to read them, annotate them, criticize them.
Become a good observer: note\index{book!based on lecture notes} which
parts of the book make sense to your audience and which require additional
explanation from you. Use these notes and observations as a take-off point
for the book.

Mathematicians appear to be a shy, introspective 
lot.\footnote{An introverted mathematician is one who looks
at his shoes when he talks to you.  An extroverted mathematician
is one who looks at {\it your\/} shoes when he talks to you.}  It seems
to exhibit too much\index{introvert vs.\ extrovert} 
hubris for a mathematician to say ``Now I shall
write a book on thus and such.''  More often than not, the mathematician
sneaks into the task; and a good way to do this is to develop lecture
notes.

This lecture notes approach has several advantages over writing
the book cold.  First, you have the opportunity to classroom test
the material, to see in real time how students react to it, and to
modify it according to what you learn from the experience.   Second,
when you teach a course you are completely involved in the material,
and it is natural to develop it and revise it as you go.  Third, you
can show your lecture notes to colleagues---without much fear of
embarrassment because, after all, they are only lecture notes---and
learn from their comments and criticisms.  Fourth, if the material
does not seem to be developing expeditiously, you can abandon the
project without losing face.  After all, these were only lecture
notes.

It also helps to have a collaborator.  Imagine going to a colleague
at a conference or other group activity and saying ``You know, there
 ought to be a book on {\it badeboop badebeep}.''  If the colleague
indicates\index{book!written in collaboration}
assent, then you can begin to describe what material 
ought to be in the book.  Before long, you are swapping ideas, building
each other's enthusiasm.  Soon enough, you are writing a book
together.  Your collaborator is a reality check, and reassures you
that you have not set for yourself a long-term fool's errand (for
example, it would certainly be the pits to find out after two years
of hard work that your book topic ``Generalized  
Theory of Fluxions and Fluents'' was no longer a matter of
current interest).

Of course writing something as big as a book with a collaborator has
its down side too.  There will be periods when you are raring to go
and he/she is busy getting a divorce, or learning to chant ``Na myoho
rengae kyo,'' or moving into a yurt.  Or conversely.  Taking 
on a book collaborator is like adding a member to your
family.  And the family could become dysfunctional.

I want to leave you with one important thought about planning a book.
Try to have the entire vision of the book in
place before you launch full steam into the project.  Such planning
enables you to keep your sense of perspective, to know how much has
been accomplished and how much remains to be done.  It also helps to
prevent you from wandering off onto detours, or from developing
specious lines of investigation.  I have written books where I have
just started writing and let the course of events dictate where my
thoughts would lead me.  Sometimes this worked well; more often it did
not.  After writing many books, I can say with some confidence that
the planned approach is far superior.

%% Section 5.3

\section{The Importance of the Preface}

I have already indicated in Section 3.4 why the Preface to any
project is an important feature.  For something as grandiose as
a book,\index{book!preface for} the Preface is paramount.  
Writing the Preface is part of the planning
process, and it acts as your touchstone as you develop the project.

Indeed, while I am writing a book I often take a
break and spend some time staring at my Preface and my Table of Contents
(TOC).  It may well be that, at an advanced stage of the
writing, I no longer agree in detail with what the Preface and TOC
say.\index{book!Table of Contents for}  
But when I wrote the Preface\index{book!TOC for}
and TOC my thoughts were organized and galvanized and I
knew exactly what I was trying to accomplish.  Studying the Preface
and TOC is a way of reorienting myself.

And remember that your reviewers and your readers, if they are smart, will
study your Preface and TOC in detail. These two essential front matter
items will give them a preview of what they are about to read, and how to
go about reading it. Just as you write the introduction to a research
paper with the referee\index{book!front matter for} in mind, endeavoring
to answer or at least minimize all his observations and objections, so you
write your Preface and TOC with a view that you are deflecting all the
reader's {\it But\/}s.

Your Preface should not spare any detail.  You have obviously thought
about why existing books do not address or fill the need that your
book fills.  Spell this out in the Preface.  You have thought about
why your book has just the right level of detail and the right
prerequisites.  Say\index{book!Preface for} 
this in the Preface.  You have thought about
why your point of view is just the right one, and the points of
view in other books are either outdated or misguided.  Say so
(diplomatically) in the Preface.

Even if you were to write your Preface, polish it to perfection,
and then put it in the paper shredder, it would have been an important
and worthwhile exercise to write it.  Writing the Preface is your
(formal) way of working out exactly what you wish to accomplish
with your book.  

%% Section 5.4

\section{The Table of Contents}

In some sense, there is no way that you can know what will be in your book
before you have written it. But you certainly will know the milestones,
and the big ideas. In writing a novel, it may be 
possible\index{book!Preface for} to begin 
with ``It was the best of times, it was the worst of
times \dots'' and then let the ideas flow; however, technical writing
demands more deliberation. Somehow, writing ``Let $\epsilon > 0$'' does
not set one sailing into a disquisition on analysis. Mathematics is just
too technical and too complex; you must plan ahead.

Writing the TOC is part of the early process of developing your book. It
may hurt at first, and it may not feel like fun. But you will launch into
writing Chapter 1 more easily if you know in advance where you are headed;
conversely, if you do not know\index{book!Table of Contents for} where
you are headed, then how can you possibly begin? Treat the writing of the
TOC like working out on your NordicTrack:\reg just do it.

Make the TOC as detailed as you can.  The more thoroughly that you can
map out each chapter and each section, the more robust your
confidence will become.  That is, it will be much clearer that you
can and will write this book.  Always remember as you supply details
that you are not wedded to this particular form of the TOC.  You can,
and no doubt will, change it later.

If you find yourself unable to write the TOC, then maybe God is 
trying to tell you something.  Maybe you were not cut out to write
this book or, worse, maybe you have nothing to say.  Writing the TOC
is an acid test.  You will have to write it eventually.  What makes
you think that you will be able to write it {\it after\/} having
written all the chapters if you cannot write it before?  Does this
make any sense?  Write it now.

I may note that, when you are writing in \TeX, the trickiest feature
is formatting.  In particular, you may have trouble typesetting a Table
of Contents.  No worries.  \LaTeX\ will do it for you.  Suppose that
your source code \TeX\ file is {\tt myfile.tex}.  Simply enter the line
\begin{quote}
\verb@\tableofcontents@
\end{quote}
right after the \verb@\begin{document}@ line of your \TeX file and,
when you compile, \LaTeX\ will produce {\tt myfile.toc}.  That is
your Table of Contents.

%% Section 5.5

\markboth{CHAPTER 5.  BOOKS}{5.5.  TECHNICAL ASPECTS}
\section{Technical Aspects:  The Bibliography, the Index,           
   Appendices, etc.}
\markboth{CHAPTER 5.  BOOKS}{5.5.  TECHNICAL ASPECTS}

If you write your book using \LaTeX, or using the macros included
with the book [SK], then you have a number of powerful tools at
your disposal for completing some of the dreary tasks essential
to producing a good book.

In the old days, when an author created the index for a book, he/she
proceeded as follows.  (For effect, let me paint the whole
dreary picture from soup to nuts.)  First, the author sent his
manuscript\index{book!production of} 
into the publisher.  For a time, he/she would hear
nothing while the copy editor was working his/her voodoo on the MS.  Then
the publisher sent the author the copy-edited manuscript.  This gave
him/her the opportunity to reply to the editor's comments and
suggestions.  For example, the editor might have changed all 
the author's {\it
that\/}s to {\it which\/}s or vice versa.  The 
copy editor might have said
``You cannot call $G(x,y)$ `the Green's function' because that is
ungrammatical.''  Or ``you cannot refer to `Riemannian metrics' in
Chapter 10 because, when Riemann's name came up in earlier chapters,
it was not in adjectival form.'' (Both of these have happened to me;
in the penultimate example, I was advised to call $G(x,y)$ ``the
function of Mr.\ Green.'') 
In any event, the author slugged his/her way through the manuscript and
made his/her peace with the copy editor, sometimes via a shouting match
over the telephone.

At the next stage the author received ``galley proofs.''  These were
printouts of the typeset manuscript, but not broken for pages.
Galley proofs were often printed on paper that was 14 inches long or
more.\index{book!galley proofs for}  The author was supposed to read
the galleys with painstaking care, paying full attention to all
details.  The main purpose of this proofreading was to weed out any
errors---mathematical or linguistic or formatting or some other---that
were introduced by the typesetter.  In particular, one would check at
this stage that all the displayed mathematical formulas were set
correctly.

At the next, and final, stage the
author was sent ``page proofs.'' Now the author was receiving his
manuscript broken up into pages, and appearing more or less as it
would in the final book.  Space had been made for figures, and the
pages had running heads and actual page numbers.   At this propitious
moment, the author was (at least in theory) no longer checking for
mathematical, English, or typesetting errors.  In the best of all
possible worlds, at this stage a check was being made that the page
breaks did not alter the sense of the text, nor did 
they\index{book!page proofs for} result in figures being misplaced.

And it was at the page proof stage that the author made up the index.
First, he/she went through the page proofs and wrote each word to appear
in the index on a separate $3 \times 5$ card, together with the
correct page reference\index{book!making an index for} (which was
only {\it just now\/} available, since the author was working for the
first time with page proofs).  Then the author alphabetized all the
$3 \times 5$ cards.  Finally, the author typed up a draft of the
index.

In the modern, computer-driven environment for producing a book, the
production process is considerably streamlined.  If the manuscript is
submitted to the publisher in some form of \TeX, then usually
the\index{book!production in the computer age}
copy-edited manuscript stage and also the galley proof stage
are skipped.  The author works with page proofs only, and that is his/her
last ``pass'' over the manuscript.  The entire business of writing
words and page numbers on index cards, alphabetizing them, and then
typing up an index script is gone.  Here is the new methodology:

Imagine, for example, that you are using \LaTeX.  You can go through
your {\tt ASCII} source file and tag words. ({\it  You can do this at
any\index{book!production using \LaTeX} 
stage of your writing---indeed, you may do it rather naturally
``on the fly'' while you are creating the book.\/})  
For instance, suppose that
somewhere in your source file the word ``compact'' occurs, it is the
first occurrence, and you want that word to be in the index with that
particular page reference.  Then you put the code
\verb@\index{compact}@ immediately adjacent (with no intervening
space) to the occurrence in the text of the word ``compact'';  thus
\verb@\index{compact}@ now appears in your \TeX\ source file.  [This
additional \TeX\ code does not change the printed \TeX\ output.]  
There are modifications to the \verb@\index@ command to specify
subentries in the index, and also to allow you to index items
that are (strictly speaking) not words (for instance, you may
wish to have \verb@\begin{document}@ appear in the index if you are
writing a book about \TeX).

You place the command \verb@\makeindex@ in your \TeX\ source code file
right after the \verb@\begin{document}@ command.  Then, when you compile
the file \verb@myfile.tex@, a new file \verb@myfile.idx@ will be produced.
This is a raw form of your index, in which the entries appear in the order
in which they appear in the book---{\it not} alphabetized and not
with subentries in place.  But there will be a command in your \TeX\ system
that processes the file \verb@myfile.idx@ and produces yet another
file \verb@myfild.ind@.  {\it That} file is the final form of your
index that you can incorporate into the source code file for your book.

Just to repeat:  The indexing commands cause all the words that have been marked for
the index to appear in a single file, called myfile.idx (assuming
that the original file was myfile.tex), together with the appropriate
page references---{\it after\/} you have compiled the source file.  You
can then use the {\tt UNIX} command \verb@makeindex@ to alphabetize the
file MYFILE.IDX and to remove redundancies.  The procedure is
documented in the \LaTeX\ book [Lam], or in the file MAKEINDEX.TEX. 
(Alternatively, you can use operating system commands to alphabetize
the\index{book!making index using \LaTeX}
 index, and then do a little editing by hand to eliminate
repetitions and redundancies. The entire process usually takes just a
few hours.)   The disc that is included with the book [SK] also
includes macros that will assist in the making of an index.
The lovely book [MGBCR] has a detailed discussion of the \LaTeX\ commands
for making an index.

The reference [SG, pp.\ 76-96] treats all the technical aspects of
compiling a good\index{index, compiling a good} 
index.  The book [Lam] has a nice discussion of the
notion that you should index by {\it concept,\/} not by word.  The
former method allows the reader to find what he/she is looking for
quickly; the latter adds---unnecessarily---to the reader's labor.  A
good, and thorough, index adds immeasurably to the usefulness of a
book.  My claim is particularly true if your book is one to which a
typical reader will refer frequently and repeatedly---for example if
your book is meant to be a standard treatment of a mathematical
field.  Many otherwise fine mathematics books are flawed by lack of
an adequate index (or, for that matter, lack of an adequate
bibliography).\index{book!bibliography for}  

There are professional indexers who can produce a workable\index{indexers,
professional} index for any book. But nobody knows your book better than
you, the author. {\it You\/} should create the index. Given that modern
software makes the creation of an index relatively painless, there really
is no excuse for not creating one yourself.
 
Similar comments may be made about the Bibliography---this
procedure has already been discussed in detail in Section 2.6.  The
book [SK] tells you how to write \TeX\ macros to compile a glossary,
a table of notation, or any similar compendium.  The process is
rather technical, and I shall not describe it here.

I conclude with a few words about Appendices. You will sometimes come to a
point in your book where you feel that there is a calculation or a set of
lemmas that you know, deep down, must be included in the book; but it will
be painful to\index{book!appendices in} write them, and they will
interrupt the flow of your ideas. The solution then is to include them in
an Appendix. Just say in the text that, in order not to interrupt the
train of thought, you include details in Appendix III. Then you state the
result that you need and move on. This practice is smart exposition and
smart mathematics as well. It is also a way of managing your own psyche:
when you are attempting to tame technical material in the context of your
book proper, then it becomes a burden; if instead you isolate the same
material in an Appendix, then you loosen your fetters and the task becomes
much easier.

An Appendix also could include background results from undergraduate
mathematics, alternative approaches to certain parts of the material,
or just ancillary results that are important but too technical to
include in the text proper.  Appendices are a simple but important
writing device.  Learn to use them effectively.

%% Section 5.6

\markboth{CHAPTER 5.  BOOKS}{5.6.  HOW TO MANAGE YOUR TIME}
\section{How to Manage Your Time \hfill \break When Writing a Book}
\markboth{CHAPTER 5.  BOOKS}{5.6.  HOW TO MANAGE YOUR TIME}

Many a mathematics book is started with a bang, two-thirds of it is
written, the writer becomes bogged down in a 
struggle with\index{book!time management while writing} a 
piece of the exposition, or the
development of a particular theorem, and the book is never completed. I
cannot tell you how often this happens; perhaps more frequently than the
happy conclusion of the book sailing to fruition.  I imagine
that the same hangup can occur for the novelist, or for
the historical writer.

I would be naive, indeed silly, to suggest that those who cannot complete
their books are just insufficiently organized. Or that such people have
not read and digested my advice. Anyone can develop writer's block, or can
arrive at a point where the ideas being developed just do not work out, or
can just lose heart. We as mathematicians, however, are accustomed to this
dilemma. Most of the time, when we write a paper, things do not work out
as we anticipated. The hypotheses need to be adjusted, the conclusions
weakened,\index{book!writer's block while writing} the definitions
redeveloped. If you are going to write a book then you will have to apply
the same talents in the large. But you also need to think ahead to where
the difficulties will lie and how you will deal with them. One of the
advantages of doing mathematics is that nothing lies hidden. We can think
and plan the entire project through, if only we choose to do so.

People in twelve-step programs, with chemical dependencies, are
taught to live one day at a time.  Such people are taught to
concentrate on the ``now.''  If you are writing a book then, on the
one hand, you cannot afford this sort of shortsightedness.  You must
plan ahead, and have the entire project clearly in view.  If you kid
yourself about how Chapter 8 is going to work out then, when you get
to Chapter 8, you are going to pay.  By analogy, if you write a paper
in such a fashion that you shovel all the difficult ideas into Lemma
3 then, when it comes time to write and prove Lemma 3, you must face the
music.  You cannot fool Mother Nature.

But, having said this, and having (I hope) convinced you of the value of
planning, let me now put forth the advantages of tunnel vision. Once you
have done the detailed planning, and you are convinced that the book is
going to work, then develop an extremely narrow focus. Pick a section and
write it. You need not write the sections in logical\index{book!tunnel
vision when writing} order (though there is some sense to that). But, once
you have picked a section to work on, then focus on that one small task,
that one small section, and do it. If some worry about another section, or
another chapter, crops up then make a note of it and then press ahead with
the writing of your chosen section. Bouncing around from section to
section, and chapter to chapter---chasing corrections around a
never-ending vortex---is a sure path to disillusionment, depression, and
ultimate failure. You can always set up scenarios for defeat. Your
book-writing project can turn into a black hole, both for your time and
for your psychic energy. Writing a book is a huge task; nobody will blame
you if you give up, or abandon the effort. But with some careful planning,
with an incremental program for progress, and with some stamina, you can
make it to the end.

Paul Halmos\index{Halmos, Paul} [Ste] advocates, and describes in detail,
the ``spiral method'' for writing a book (or a paper, for that matter).
The idea is this: first you write Chapter 1, and then move on to Chapter
2.\index{book!spiral method for writing} After you have written Chapter 2,
you realize that Chapter 1 must be rewritten. You perform that rewrite,
re-examine Chapter 2, and then you move on to Chapter 3, after which you
realize that Chapters 1 and 2 must be rewritten. And so forth. If you are
writing by hand, with a pen on paper, then the spiral method takes place
in discrete steps as indicated. If, instead, you write with a computer
then the spiral method can take place in a more organic fashion: as you
are writing Chapter 3, and realizing that Chapter 1 needs modification,
you pull up Chapter 1 in another window and begin to make changes while
you are thinking about them. If those changes in turn necessitate a
massage of Chapter 2, then you pull it up in a third window. The advantage
of doing things in discrete steps, as described by Halmos,\index{Halmos,
Paul} is that you always know where you are and what you are doing; the
disadvantage of the organic approach is that you can become lost in a
vortex---caroming around among several chapters. The technique must be
used with care.

It can only improve your work to review Chapters $1$ through $(n-1)$ after
you have written Chapter $n$. On the other hand, if you do use the
out-of-the-box spiral method, as described and\index{Halmos, Paul}
recommended by Halmos, then one upshot will be that Chapter 1 of your book
will receive more attention than any other part, Chapter 2 will receive
the second greatest dose of attention, and so forth (for the proof, use
induction). As a result, your book {\it could\/} appear to the reader to
become looser and looser as it proceeds. Perhaps this is an acceptable
outcome, for only the die-hards will get to the end anyway. But when you
adopt a method for its good points, also be aware of its side effects.

Certainly choose a method that works for you---organic, inorganic, spiral,
or some other---and be sure to use it. \ \ If there is any time when it is
appropriate to be organized, methodical, indeed compulsive, that time is
when you are writing your book.

No matter what method you adopt for reviewing and modifying your work, keep
this in mind: only wimps revise their manuscripts; great authors throw
their work in the trash and rewrite. Such advice causes many to say ``That
is why I could never write a book;\index{book!revising the manuscript} it
is sufficient agony just to write a short paper.'' Rewriting is not so
difficult; in many ways it is easier than figuring out where to insert
words or to substitute passages. Treat your first try as just getting the
words out, for examination and consideration. Once the thoughts are lined
up in your head, then the first draft has served its purpose; you may as
well discard it (and {\it don't peek!\/}). The next go is your opportunity
to shape and craft the ideas so that they sing. The next round after that
allows you to polish the ideas so that they are compelling and forceful.
The final step allows you to buff them to a high sheen.

Use the advice of the last paragraph along with a dose of common
sense.  After you have struggled for a month to write down the 
proof of a difficult proposition, you are not going to throw it
in the trash and start again.  My advice here, as throughout this
book,\index{book!revision vs.\ writing anew} applies selectively.

Once you have arrived at (what appears to be) the end of the task of
writing your book, you still are not finished.  There remains a lot
of detail work.  You must prepare a good bibliography (Sections 2.6,
5.5).  You must prepare a good, detailed, index (the computer can
help a lot here---see Section 5.5).  If appropriate, you should
prepare a Table of Notation.  You\index{book!glossary for}
might consider building a Glossary. None\index{book!Table of Notation for}
of these tasks is a great deal of fun. But they will increase the value of
your book immeasurably. They can make the difference between an advanced
tract accessible to just a few specialists, or a book that opens up a
field.

%% Section 5.7

\markboth{CHAPTER 5.  BOOKS}{5.7.  WHAT TO DO WITH THE BOOK}
\section{What to Do with the Book Once It Is Written}
\markboth{CHAPTER 5.  BOOKS}{5.7.  WHAT TO DO WITH THE BOOK}

You have written your {\it magnum opus}, slaved over it for two or
more years, shown it to colleagues, received the praise of student
and mentor alike.  The manuscript is now polished to perfection.
There is no room for improvement.  Now what do you do with it?

The rules for submitting a book manuscript to a publisher are different
from those for\index{book!submitting to a publisher} 
submitting a research paper to a journal.  The hard and
fast rule for the latter is that you can only submit a research paper
to one journal at a time.  Most research journals tell you up front that,
by submitting a paper, you are representing that it has not been
submitted elsewhere.

Not so for books. You can submit a book manuscript simultaneously to
several different publishers. These days there are just a few mathematics
publishers---especially for advanced books. Get a feel for the different
publishers by looking at their book lists. You will see what quality of
books and authors they publish, and in what subject areas. Some
publishers, such as the AMS, CRC Press, Springer, 
Birkh\"{a}user, and the American Mathematical Society, have several book
series in mathematics. Familiarize yourself with all of them so that you
can make an informed choice. Talk to experienced authors to obtain the
sort of information that cannot be had from advertising copy.

If you want to jump-start the publication process, then you can begin long
before your book is completed. For example, if you are looking for a
typing grant or an advance, then you may wish to begin negotiations with
publishers after you have written just two or three chapters. Submit them,
along with a Preface or Prospectus\footnote{Like a Preface, the Prospectus
will describe what the book is about and why you have written it. Unlike a
Preface, the Prospectus\index{book!what to send the publisher} will
describe the audience, the competing texts, the types of courses that
could use the book, and the types of schools and departments that might
adopt the book.} (the marketing version of a Preface) and a TOC. And of
course include a brief cover letter saying who you are, what book you are
writing, and exactly what materials you are remitting.

Always send a manuscript to a publisher by either registered or certified
mail---return receipt requested. There are both practical and religious
reasons for doing so. First, it requires some effort and expense to
prepare a manuscript, plus the figures, plus the discs, for submission to
a publisher. You want to protect your investment of time and money; so
special mail services and even insurance are definitely in order. Less
obvious is an artifact of the way that publishing houses work: items that
arrive by regular post tend to get thrown into a pile; items that arrive
by registered or certified or express mail are given special treatment.
Stop and think about how many manuscripts, or how many pieces of mail, a
big publishing house will receive in any given business day. Now you will
understand why you should take pains to ensure that your manuscript
receives the particular attention that it deserves.

These days it is perfectly acceptable to send your book materials to a publisher
as a {\tt *.pdf} file in an \email\ attachment.  Describe in detail,
in the text of your \email, just what you are remitting---how many chapters,
what is the subject of the book, what books it should be compared to.   You
can even suggest some reviewers.  You should throw in a few sentences
about just who you are, what your background is, and why you are the
right author for this book.  It would not be out of place to attach
your Curriculum Vitae to this same \email.

In order to be able to negotiate intelligently with a publisher, be
sure\index{book!information for the publisher} 
to have the following information about your book under control:
\begin{enumerate}
\item Subject matter and working title
\item Level (graduate, undergraduate, professional, etc.)
\item Classes in which the book could be used
\item Existing books with which your book competes
\item Working length
\item Expected date of completion
\end{enumerate}

The publisher needs to know a subject area and working title for
in-house and developmental purposes.  The guys in the suits refer,
among themselves, to the ``Krantz project on fractals.''  So they
need a working title.  They need to know a working length and a
sketch of the potential market so that they can price out the
project.  They need to know an approximate due date so that they can
deal with scheduling (a non-trivial matter at a publishing house).

I am the consulting editor for a book series.  One of my earliest
authors completed his/her book two years late, with a book twice the
length originally projected; also the book was on a different subject
than that contracted, and with a different title.  And the
author wanted it to be published in two volumes!  I cannot tell
you how much trouble I had persuading the publisher to go ahead
with the project.  When you are dealing with a publishing house you
are dealing with business people.  You must endeavor to conform
to their view of the world.

If the publisher is interested in your project, then he/she will probably
solicit reviews. Some publishers will ask you to suggest reviewers
for\index{book!reviewers hired by publisher} your project. Most will not.
Expect the reviewing process to take three or four months. Expect to see
two to four reviews of your work.

One of the most difficult, and valuable, lessons that I have learned as an
author is to read reviews. By this I mean to read them
intensely\index{book!how to read reviews of} and dispassionately and to
learn what I can from them. Forget reacting to the criticisms. Forget
justifying yourself. Forget answering the reviewers' comments. Forget
melting down into an emotional puddle of goo. The point is this: even if
you cannot understand what the reviewer is thinking, what he/she describes is
nevertheless what he/she saw when reading the manuscript. The review
describes the impression that the manuscript made on him/her.  The main
question you should be asking yourself as you read the reviews is ``What
can I learn from these reviews?'' ``How can I use these comments to
improve my book?'' There is generally something of value in even the most
negative of reviews.

Usually the publisher has established an initial interest in
your project by looking at your Prospectus and TOC, and by
agreeing to undertake the cost of reviewing (unlike a referee
for a paper, a book reviewer is usually paid a modest
honorarium). If the consensus of the reviews is favorable,
then the publisher will most likely decide to publish your
book. He/she will then ask you to take the reviews under
advisement, and only that. The editor may want to discuss them
with you, and may even want your detailed reaction to them.
But few, if any, publishers will hold you accountable for each
comment made by each reviewer.\footnote{Note that these
remarks do not apply to the writing of a textbook at the lower
division\index{book!how to respond to reviews of} level, for
the so-called ``College Market.'' Such a project is more of a
team effort: you and the reviewers write the book together, in
a sort of Byzantine tug-of-war procedure. The process is best
learned by consenting adults in private, and I shall say
nothing more about it here.}

Remember this! And what I am about to say applies to research papers and to
books and to anything else that you submit for review: the reviewer is not
responsible for the accuracy and correctness of your work. There is only
one person who bears the ultimate responsibility,\index{book!accuracy in}
and that is you. Many reviewers will do a light reading, or an overview,
or will read the manuscript piecemeal, according to what interests them.
If the reviewers give you a ``pass,'' then that is good. But this ``pass''
is not a benediction, nor even a suggestion that everything you have written
is correct. You must check every word, and you yourself must certify every
word.

In any event, the period immediately following the review process is your
chance to take a couple of months and polish your manuscript yet again.
(You will also have the opportunity to make small changes later on in the
page proofs. But the post-review period is your last chance for
substantial rewriting.) {\it Treat this as a gift.\/} It would be
embarrassing to publish your book blind---with no reviews---and then to
have your friends point out all your errors and omissions, or (worse) that
your point of view is all wrong. The reviewing\index{book!revision of
manuscript} process, though not perfect, is a chance to collect some
feedback without losing face and without any repercussions.

After you have polished your MS to your satisfaction, and presumably 
shown it to some friends and students and colleagues, then you submit
the final, polished draft to the publisher.  Many publishers will
want this manuscript to be double or triple spaced, so that the 
various copy editors and typesetters will have room for their 
markings and queries.  The \TeX\ command \verb@\openup@$k$ \verb@\jot@,
where $k$ is a positive integer, will increase the between-line spacing
in your \TeX\ output by an amount proportional to $k$.

Nowadays almost all of the book-publishing process is conducted electronically,
and mostly over the Internet.  You submit your book as a {\tt *.pdf} file.
The publisher sends the {\tt *.pdf} file to the reviewers as an \email\ attachment.
Each reviewer sends in his/her report as an \email.  The publisher removes
any identifying lines from the reviews and passes them on to you (again by
\email).  You make the appropriate edits to your \TeX\ source file, declare
the book to be finished, and send both your {\tt *.tex} and {\tt *.pdf}
files to the publisher as \email\ attachments.  The copy editor works
on your {\tt *.pdf} file and marks edits on that file using ``electronic sticky
notes.''\footnote{Electronic sticky notes are a software utility that allows
you to paste little notes to any page of a {\tt *.pdf file}.  You simply
place your pointer where you want the note to be, right click on the mouse,
select ``electronic sticky note'' from the dropdown, and you get a little
yellow box in which to write your comment.  The little yellow boxes become part
of the file.  You will read the copy editor's comments in his/her electronic
sticky notes and you will respond with your own electronic sticky notes.}

Now here is one of the great myths that exists at large in the mathematical
community. People think that, in 2017, you send a flash drive or a CD-ROM, with \TeX\
code\index{book!mythical production of} on it, to the publisher. The
publisher puts the device in one end of a big machine and a box of books
comes out the other. Technologically this phenomenon is actually possible.
But a top-notch publishing house has a much more exacting procedure.

Here, instead, is what a good publishing house does with your manuscript
and disc. First, an editor decides whether your book is ready to go into
production. He/she may show your ``final manuscript'' to a member of
his/her editorial board, or\index{book!going into production} he/she may
make the decision on his/her own. But this hurdle must be jumped. Once the
book goes into production, some copy editing will be done. The actual
amount will vary from publishing house to publishing house. During the
copy editing process, your spelling, grammar, syntax, consistency of
style, and other nonmathematical aspects of your writing will be checked.
Depending on the density of corrections at this stage, you may or may not
be contacted. You may have to submit another manuscript.

One point that needs to be recorded is this.  With the advent of \TeX\ and the
Internet, more of the burden is placed on your shoulders.  When a copy
editor sends you the edits for your book, it is not enough for you
to say, ``OK, these edits are fine by me.''  You have to actually
go into your \TeX\ source code file and make the edits yourself (or at
least make the edits that you agree with).  When you are finished, then
you compile the source code file, produce a {\tt *.pdf} file, and send
that back into the publisher.

If you have never before written a book, then you may be surprised
at the many details that a copy editor will attend to when
handling your book.  Here are some of these:
\begin{itemize}
\item All displayed equations should be formatted in the same way.  
\item Left
and right page bottoms should align (this last task is something at
which Plain \TeX\ does not excel; \LaTeX\ handles the issue with the
\verb@\flushbottom@ command). 
\item No page should begin with a single line
that ends a paragraph (such an item is called an ``widow'').
\item No page should end with a single line that begins
a paragraph (these stragglers are called ``orphans'').  
\item Figures must be positioned properly, and rendered at the right size.
\item Running heads must be checked.
\item Page breaks must be checked.
\item Blank pages at the ends of chapters (when the last page
of text in the chapter is odd-numbered) must be completely blank.
\end{itemize}

If your project were typeset the old-fashioned way, with movable 
type---say that it has 400 pages---then the typesetting job would
cost \$15,000--\$20,000.  If instead you produce a \TeX\
file to a level of reasonable competence, then the adjustments
that\index{book!cost of producing} I 
described in the last paragraph will cost \$5,000 to \$7,000. So
\TeX\ {\it does\/} save money in the publishing process.  

After you have approved the page proofs, then that is the end of your role
in the publishing process (but see the {\it caveat\/} below about the
dreaded Marketing Questionnaire).  But the title page and
back cover copy are produced separately.  I suggest that you {\it insist\/} on
seeing the title page before the book goes to press. It happens---not
often---that an author's name is misspelled or an affiliation is rendered
incorrectly. Such an eventuality is embarrassing for everyone. It is best
to avert it.  And the back cover copy is a prominent advertisement for
your book; you want to be sure that it says the right things.

Even though your role is at an end, let me say a few words about what
happens next.  As is mentioned elsewhere in this book, when a \TeX\
file, consisting of {\tt ASCII} code, is compiled then the result is a {\tt *.dvi}
file.  Typically, this ``Device Independent File'' is then
translated, using software and without human intervention, to a {\tt
PostScript}\reg\ file. Why {\tt PostScript}\reg?  Many high-resolution
printers read {\tt PostScript}.\reg 
Once\index{book!production of} 
the files for the book have been translated into {\tt
PostScript}\reg, then the book is printed out at high resolution on
RC (resin coated) paper.  The result is a reproduction copy (or {\it
repro copy})\index{book!repro copy for} 
of your book printed on glossy, nonabsorbing paper, at
extremely\index{repro copy} 
high resolution.  All the smallest subscripts and superscripts
will be sharp and clear, even under magnification.

The repro copy of the book is then ``shot.''  Here, to be ``shot''
means to be photographed.  The pages of the book are photographed
onto film, in the fashion familiar to anyone who takes snapshots.  
But it is not printed onto photographic paper (what would
be\index{book!shooting of} the point of that?---it is {\it already\/}
on paper).  Instead, the negative is then exposed or ``burned'' into
chemically treated plates.  These plates are the masters from which
your book is printed.  (This process is becoming ever more
streamlined.  Today at the AMS, the ``repro copy'' step is skipped
altogether; the production department goes directly from the
electronic file to the negative.)

Once the printing, or lithographic, plates are prepared, then the rest of
the printing process---printing, cutting, and binding---is quite
automated. Good books are printed sixteen pages to a sheet, and then
folded and cut. This procedure results in the ``signatures'' that you can
see in the binding of any high-quality book (not a cheap paperback). (In
the old days the publisher did not cut the
signatures;\index{book!signatures in} a serious reader owned a book knife,
and did the cutting himself/herself.)

Interestingly, the physical cost of producing a book---that is the
printing, the binding, the cost of the paper---is well under \$5 per
volume; at least this is true if the print run is reasonably
large.\index{book!cost of producing} The difference in cost between
producing a paperback volume and a hardback volume is about \$2, depending
on the quality of the papers used. So why do math books cost so much?

The pricing question for books is all a matter of marketing. To be fair,
the publishing house has overhead. You remember the \$5,000 to \$7,000 for
the services of a \TeX nician? That is a cost that anyone can understand.
Then the salaries of the editor, the publisher, the company president, the
people in the production\index{book!marketing of} department, the costs of
marketing, the physical plant, and so forth must come out of money earned
from the sale of books. Most people, indeed most authors, are not
cognizant of the cost of warehousing books in a serviceable manner (so
that the books are readily accessible when an order comes in). Warehousing
is a fixed cost that adds noticeably to the expense of each and every book
that we buy. Taxes on inventory are {\it very} high. These last costs are
called ``overhead'' or ``plant costs,'' and play much the same role as the
overhead for an NSF Grant. Most publishing houses figure the cost of
producing a book by taking the up front, identifiable costs---technical
typesetting, any advance to the author, print costs (often the printing is
done by an outside firm), copy editing, composition, shooting---and then
adding a fixed percentage (from 30\% to 50\%) to cover the overhead that
was described above.

Then the editor does a simple arithmetic problem. He/she must make a
credible, conservative estimate as to how many copies your book will sell
in the first couple of years. Fifty years ago this was easy, since many
libraries had standing orders for all the major book series. (For example,
in the late 1960s, a company like Springer-Verlag or John Wiley could {\it
depend\/} on library sales of 1000-1200 copies for each book!) With
inflation, cutbacks, and other stringencies, libraries now pick and choose
each volume. Thus the editor must make an evaluation based on {\bf (i)}
whether the book is in a hot area, like dynamical systems or wavelets,
{\bf (ii)} whether people in disciplines outside mathematics (engineers,
for example) will buy it, {\bf (iii)} whether students will buy it,
{\bf (iv)}  whether the author has name recognition, and
{\bf (v)} whether the book can be used in any standard classes. Other
factors that figure in are {\bf (a)} Is this the first book in an
important field? {\bf (b)} Is there stiff competition from
well-established books? {\bf (c)} How much effort is the marketing
department willing to put into promoting the book? (You may suppose that
the marketing department will promote any book that the editorial
department sends in. On a {\it pro forma\/} level they will. But there is
a delicate dynamic between these two publishing house groups, and a
constant push and pull. A good editor takes pains to generate enthusiasm
among the marketing people for particular books.) Having evaluated these
factors, the editor writes a proposal for how many volumes of your book
the house can expect to sell within a reasonable length of time (a couple
of years). Then he/she figures in the company's standard profit
expectation. This gives rise to the wholesale price of the book.

As an example, suppose that you write a book on a fairly specialized
area of partial differential equations.  After an analysis of the
foregoing kind, the editor determines that the book is sure to sell
500 copies in the first two years.  The up front costs are \$15,000.
Add 50\% for overhead and that makes \$22,500.  Add 20\% for the
company's standard profit margin and that brings the total to
\$27,000.   The wholesale price of the book must, after sales of 500
copies, bring in that much money.  (If a given editor has several
books that fail to meet this simple criterion, then he/she is out of a
job.)  Now do the arithmetic.  You will find that the wholesale price
of this book must be \$54 per volume.  Thus a
bookstore will probably sell it for at least \$70 to \$80.  Now do you
understand why mathematics books cost what they do?

Incidentally, if the difference in cost between producing a hardcover copy
of a given book and a paperback copy of that same book is about \$2, then
what accounts for the large difference in cost between hardcover and
paperback books? The answer, apart from marketing voodoo, is that the
costs of producing the book tend to be covered by the sale of the
hardcover version. Thus the publisher has considerable latitude in pricing
the paperback edition. John Grisham novels\index{book!cost of hardcover
vs.\ paperback} stay in hardcover format for more than one year before the
paperback edition is released; usually, the hardcover edition sells
millions of copies. The production costs, and the huge advance that
Grisham\index{Grisham, John} garners for each of his books, are well
covered by the hardcover sales. Thus the publisher is ready to make real
money when the paperback edition is released. He/she can be imaginative
both in pricing and in production values---if the physical cost of
producing a volume is \$5-\$10, then he/she can price it for as little as
\$7-\$15 and expect to sell a great many copies. (Interestingly, the
entire notion of mass market paperbacks was invented in the early 1950s by
Mickey Spillane and his publishers---Dutton and Signet. By 1955,
Spillane\index{Spillane, Mickey} had written three of the five best
selling books in history---and he had only written three books! By
contrast, Margaret Mitchell's blockbuster {\it Gone with the Wind} [Mit]
sold fewer than a million copies in its first two years---all in
hardcover, of course. James Gleick's {\it Chaos} [Gle] has sold about the
same.)

Back to math books. In the preceding discussion there was an important
omission. How does the editor make the market determination that I
described three paragraphs ago? He/she can always consult his/her
editorial board and his/her trusted advisors. But let me reassure you that
he/she will\index{book!Marketing Questionnaire for} certainly study your
Prospectus and Preface, and he/she will pay close attention to your {\it
Marketing Questionnaire}.
 
The latter item bears some discussion.  Whenever you write a book for
a commercial publishing house, and often for a professional society
or a university press, you will be sent a Marketing Questionnaire
to complete.  I hate to complete these things, and you will too.  But
you must do it.  I have heard authors say ``I'll just phone the
editor and talk to him about this stuff.''  Sorry; that just will not
do.  You must complete the questionnaire, and carefully.

What is this mysterious object? First, the questionnaire is long---often 10
pages or more. Second, it asks a lot of embarrassing questions: What is
your hometown newspaper? Which professional societies might be interested
in your book? What are the ten strongest features of your book? What is
the competition? Why is your book better? In which classes can your book
be used? What is typical enrollment in those classes? How often are they
taught?

As mathematicians, we are simply not comfortable fielding questions such as
these. We do not think in these terms. But, if you have been attending to
the message of this section, then you can see how an editor can use this
information to help price out the book. So why can you not just go over
this stuff on the phone with the editor? One reason is that the editor
needs this information {\it in writing}---for the record, and to show that
he/she is working from information that {\it you\/} provided, and for
future reference. The other is that the questionnaire will be passed along
to the marketing department for the development of advertising copy and
marketing strategies for your book. Like it or not, the Marketing
Questionnaire is important. Take an hour and fill it out carefully.

When I was developing my first book, and negotiating with my publisher, I
asked the editor what I would be peeved about three years down the line.
He told me that I would be unhappy about the size\index{book!annoyances
with the publisher off} of the print run, and I would be unhappy with the
advertising. Then he explained to me how the world works. First, think
about the sales figures that I described above. And think about the fact
that a business must pay a substantial inventory tax for stock on hand.
Extra books sitting around are a liability. And today (with new printing
technology) small print runs are not so terribly expensive as they were
even ten years ago. Even print-on-demand is feasible in many cases.
So if the publisher thinks that your book will sell
500 copies in the first couple of years, then the first print run is
likely to be only 750. When that stock starts to run low, another 750 can
be generated easily. The money saved per unit with a print run of 1500 (as
opposed to 750) is relatively small, and is sharply offset by storage
costs and inventory tax.

And now a word about advertising. There is nothing that an author likes
better than to open the {\it Notices of the American Mathematical
Society}\index{book!advertising of}\index{book, promotion of} or the {\it
Mathematical Intelligencer\/} or the {\it American Mathematical Monthly\/}
and to see an ad for his/her book. Of course a full page ad is best (and
almost never seen), but a half page, or quarter page, or even an ad shared
with eleven other books, is just great. Typically, you will see such an ad
just once for your book. After that, your name and the title of your book
will appear in the company's catalogue. Of course there will be
advertising material OnLine. For a textbook there could be an entire Web
site devoted entirely to one book. Many publishers rely on ``card
decks''---stacks of $3'' \times 5''$ cards, each with a plug for a single
book---that are mailed in a block to mathematicians. Usually the potential
buyer can mail in a card, without money, and receive a copy of a
particular book for a 30 day examination period.

In the spirit of doing first things last, let me now say a few words about
book contracts. When a publishing house is interested in publishing your
book, then it will send you a contract. Typically,\index{book!contract
for} you will be offered a royalty rate of 10\% to 15\%. You will be given
a deadline, and this deadline is definitely negotiable.
Err\index{book!royalty for} on the conservative side (more time, rather
than less), so that you have a fighting chance of finishing the book on
time. If you do finish on schedule, then the publisher will take a shine
to your project. If you do not, and the project is six months late, then
most publishers will be forgiving; but, technically, a late project is no
longer under contract!

A rough page length will be specified in the contract, and a working title
given. Sometimes you will be offered an advance against royalties, or a
typing grant. Sometimes you will be asked to certify
that\index{book!contents of contract for} you will submit your manuscript
in some form of \TeX. Then there will be a lot of legal gobbledygook, most
of which seems to be slanted in favor of the publisher. For the most part,
it is. The publisher wants to be able to pull the plug on a project
whenever and wherever it deems such an action suitable. Honorable
publishers do not like to exercise this option, but they want to have the
option available.

I can tell you that many authors---especially first-time authors---are
quite uncomfortable with standard book contracts. This uneasiness stems,
for the most part, from lack of familiarity. The details of the contract
{\it can\/} be negotiated, and you should discuss with your editor any
passages or provisions that you do not like. If the\index{book!negotiating
the contract for} publisher wants {\it you\/} to render the artwork in
final form, and you cannot or will not do it, then negotiate. If you do
not like the deadline, then negotiate. If the number of gratis copies of
the work offered to the author is not adequate, then negotiate some more.
Usually such negotiations are fairly pleasant. You will find the editor
eager to cooperate---as long as your demands are within reason.

You may find it attractive to join the {\it Text and Academic Authors
Association} (TAA).\footnote{{\sl Text and Academic Authors Association},
P.\ O.\ Box 367, Fountain City, Wisconsin 54629.  The Web address
is \verb@www.taaonline.net@.}  This
organization was formed to defend the rights of authors, and will help
you in dealing with publishers.\index{Text and Academic Authors
Association}\index{TAA} It also has a rather informative newsletter.  And
membership gives you access to a number of useful discounts, so that
your dues are almost a wash.

I have dealt with many publishers.  Most of them are very good to their
authors (as well they should be) and most employ knowledgeable and 
competent editors.  However, forewarned is forearmed.  It is helpful
to be familiar with the publication process before you launch into it.

%% Section 2.8

%\markboth{CHAPTER 2.  TOPICS SPECIFIC TO MATHEMATICS}{2.7.  WHAT TO DO ONCE
%    THE PAPER IS WRITTEN}
\section{Royalties}
\markboth{CHAPTER 5.  TOPICS SPECIFIC TO MATHEMATICS}{5.8. ROYALTIES}

It makes sense that the author of a book will want to be compensated
for his/her efforts.  In other words, the author expects some royalties.
Of course a math book is not going to sell like a Tom Clancy novel.
But one can make a nontrivial amount of money from a math book.  As
an instance, calculus author Jim Stewart built a \$26 million dollar
home in Toronto with his royalties.

These days the royalty rate for an undergraduate text ranges
from 10\% to 15\%. It could be considerably more for a
well-established author with a popular book. For a graduate
text or monograph the royalty could be less. Here is the
passage from a recent contract for an upper division math text
going into its fourth edition (this is in fact a real analysis
text): 
\vspace*{.12in}

6. ROYALTIES
\begin{enumerate}
\item[{\bf (a)}]  The Publisher agrees to pay the Author (or someone
designated by the Author), and the Author shall accept as
payment in full for writing and delivering the Manuscript,
Illustrations, and index, for the performance of all of the
obligations of the Author hereunder, and for all the rights
granted to the Publisher pursuant to this Agreement, the
following amounts:
\begin{enumerate}
\item[{\bf (i)}] For copies in print or eBook format sold by the Publisher
in the United States of America, twelve percent (12\%) on the
first 750 copies and fifteen percent (15\%) thereafter of the
Publisher's net receipts (as defined in Paragraph 6(d) below).

\item[{\bf (ii)}] On translations, licensing sales, electronic database
sales, excerpts, abridgments, deep discount sales (sales at a
discount of fifty percent (50\%) or greater of the Publisher's
established list price of the Work), and on all sales of
copies of the Work outside the United States of America, the
Publisher shall pay royalties at one-half ($1/2$) the rate set
forth in Paragraph 6(a)(i) above in respect of the Publisher's
net receipts. In the event the Work is included in an
electronic database with other works, or is otherwise
exploited in combination with other works, royalties shall be
apportioned by Publisher in its sole discretion, exercised in
good faith.
\end{enumerate}
\item[{\bf (b)}] In the event the Publisher exercises any of the rights of
the Publisher pursuant to Paragraph 5 above and a royalty is
not specifically provided for, the royalty which shall be
payable to the Author shall be one-half ($1/2$) of the rate set
forth in Paragraph 6(a)(i) above in respect of the Publisher's
net receipts.

\item[{\bf (c)}] Notwithstanding the above, no royalty will be paid on
copies of the Work furnished gratis for review, advertising,
promotion, bonus, sample, or like purposes, or on copies of
the Work sold at less than Publisher's cost, or on any copies
returned to Publisher for any reason, or on copies of the Work
sold to the Author. Free use of the rights granted herein may
be made by the Publisher to promote the sale of copies of the
Work and the rights therein. The Publisher may create a
reasonable reserve for returns when calculating royalties.

\item[{\bf (d)}]  For purposes of this Agreement, the Publisher's ``net
receipts'' from sales shall mean monies received by the
Publisher from such sales less adjustments for discounts,
credits, and returns. Royalties will not be paid on prepaid
transportation, postage, insurance, and taxes. The Publisher's
"net receipts" from licensing or assignment shall mean monies
received by the Publisher less any specified costs of such
licensing or assignment.

\item[{\bf (e)}] All payments made under the terms of this Agreement will be
subject to Federal income tax withholding, as required by the
United States Internal Revenue Code.

\item[{\bf (f)}] All royalties and other income accruing to the Author under
this Agreement shall be credited to an account maintained on
the records of the Publisher (the ``Royalty Account''), which
Royalty Account will be charged for all amounts paid or
payable to Author, including any advance payments, and for all
amounts Author is charged, or obligated to pay, pursuant to
this Agreement.
\end{enumerate}

You can see that the publisher is careful to cover all possible
scenarios, and that the contract is written so that no
misunderstanding is possible.  The publisher is also very
explicit about royalty rates for {\it e}-books, for electronic
databases, and other high-tech versions of the book.

It is possible to negotiate the royalty rate with the
publisher. I once retained a publishing attorney to negotiate
publishing contracts for me. He got me some terrific royalty
rates, but afterwards the publishers were quite annoyed with
me for having indulged in this artifice.
\vfill
\eject

\hbox{ \ \ \ }

\thispagestyle{empty}

\newpage

%%%%%%%%%%%%%%%%%%%%%%%%%%%%%%%%%%%%%%%%%%%%%%%%%%%%%%%%%%

%% Chapter 6
\chapter{Writing with a Computer}

\begin{quote}
\footnotesize \sl Computers are useless.  They can only
give you answers.
\smallskip \hfill \break
\null \mbox{ \ \ } \hfill \rm Pablo Picasso 
\end{quote}

\begin{quote}
\footnotesize \sl If he wrote it he could get rid of it.  He had 
gotten rid of many things by writing them.
\smallskip \hfill \break
\null \mbox{ \ \ } \hfill \rm Ernest Hemingway \break
\null \mbox{ \ \ } \hfill \rm {\it Winner take Nothing\/} [1933].  Fathers and Sons 
\end{quote}

\begin{quote}
\footnotesize \sl Easy reading is damned hard writing.
\smallskip \hfill \break
\null \mbox{ \ \ } \hfill \rm Nathaniel Hawthorne  
\end{quote}

\begin{quote}
\footnotesize \sl In a very real sense, the writer writes
in order to teach himself, to understand himself,
to satisfy himself; the publishing of his ideas,
though it brings gratification, is a curious
anticlimax.
\smallskip \hfill \break
\null \mbox{ \ \ } \hfill \rm Alfred Kazin \break
\end{quote}

\begin{quote}
\footnotesize \sl On seeing a new piece of technology: \hfill \break

\noindent A science major says ``Why does it work?''  \hfill \break
An engineering major says ``How does it work?''  \hfill \break
An accounting major says ``How much does it cost?''   \hfill \break
A liberal arts major says ``Do you want fries with that?''  
\smallskip \hfill \break
\null \mbox{ \ \ } \hfill \rm Anon. 
\end{quote}

\begin{quote}

\footnotesize \sl [With reference to Germany]  One could almost
believe that in this people there is a peculiar sense of life as a
mathematical problem which is known to have no solution.
\smallskip \hfill \break
\null \mbox{ \ \ } \hfill \rm Isak Dinesen
\end{quote}

%% Section 6.1

\markboth{CHAPTER 6.   THE MODERN WRITING ENVIRONMENT}{6.1.  WRITING ON A COMPUTER}
\section{Writing on a Computer}
\markboth{CHAPTER 6.   THE MODERN WRITING ENVIRONMENT}{6.1.  WRITING ON A COMPUTER}

Today most every mathematician writes on a computer.  You may
find it cathartic to generate your early drafts writing by hand
with a pen.  That's fine.  All the best writers create their
work in that fashion.  But, in the end, your paper or book will
be rendered on a computer.  That is just the way that it is.

Clearly, when you are writing on a piece of paper with a pen
or pencil, then you can easily and naturally jump from one part
of the page to another.  You can, in a comfortable and intuitive
fashion, jot marginal notes and make insertions.  You can put
diacritical marks and editorial marks where appropriate.  You
can scan the current page, flip ahead or back to other pages,
sit under a tree with your entire MS clutched in your fist,
put Post-it\reg\ notes in propitious locations, tape addenda to 
pages, and so forth.

Now the fact is that almost all the ``old-fashioned'' devices described in
the last\index{writing!old-fashioned methodology} paragraph have analogues in the computer setting. And the
computer has capabilities that the traditional milieu lacks: magnificent
search facilities, unbeatable cut and paste features, the power to open
several different windows that either contain several different documents
or several different parts of the same document, and many others as well.
With a computer, you can have your text open in one window, a dictionary
open in another window, and the Internet open in a third window. What
could be better?

But you still must use the tools that work for you. If you have been
writing with a pen on paper for many years, then you may be disinclined to
change. At a prominent university on the east coast there was an eminent
and\index{writing environment!windows} 
prolific mathematician, who had access to any writing facilities that
one might wish, and who wrote by candlelight with a quill. \ \ That was
his choice, and it certainly worked for him.\footnote{He also
ate dinner every night wearing a tuxedo.} I also know
people who take lecture notes, directly in \TeX, on a notebook computer.
This I cannot imagine, but it works for them.\index{notebook computer}

In this section I want to say a few words about writing on the computer,
and what I find advantageous about it.  I am addicted to writing
on the computer.  It makes me more productive and efficient, and it saves
me a lot of time and aggravation.  

I will reserve comments about specific writing systems, like \TeX\ and {\tt
Word}, for a later section.

When writing on a computer, you can type as fast as you wish, never fearing
for spelling or other errors. When you become acclimated to the medium,
you can create text at least as fast as you would have with a pen
(assuming that you know how to type), and the text will always be legible.
You can make corrections, insertions, deletions, and move blocks of text
with blissful ease. You can print out beautiful paper copy of your work
(paper copy is called {\it hard copy}), and you can store your work on
your hard disc or hard drive (also known as the {\it fixed
disc}).\footnote{Today many computer systems have multiple external drives
(hard, flash, and other). I store all my work on external drives---never
on the fixed disc.} You never need worry about misplacing all or part of
your manuscript, since finding files on your hard drive is trivial. Even
if, weeks or months later, the only thing you can remember about your
document is a word or phrase in it, you will be able to find the file
containing the document in seconds.

To illustrate this last point, I often find myself printing out
another copy of a paper or chapter that I am working on, rather than
trying to find where I put my last paper copy.  I can find my file on
my hard disc and print it out in just a moment; the old approach, more
traditional and agonizing, of searching through my study for my hard
copy could take hours.  And remember this point:  any tool that
prevents your writing moods from being interrupted or jarred is a
valuable one.  My computer has eliminated, for me, the need to search
my office for the paper copy that I want to work on. It saves me
hours of time, and it saves me considerable irritation. Cherish
those tools that make your\index{computer!advantages of writing with}
life easier, and learn to use them well.

When working on a computer, you can easily keep every single version
of a document you are writing.  Suppose, for example that you are writing
an article about diet fads among troglodytes.  The first version of your
article could be called {\tt TROG.001}.  After you modify it,
the second version could be called {\tt TROG.002}.  The third would
then be called {\tt TROG.003}.  And so forth.  All these would be neatly
stored, and accessible, on your hard disc.   Compare with the situation,
in a paper office, in which you had thirty-two 
different versions of a document.
How would you store them all?   How would you keep track of and
differentiate among them?  How would you access them?
Note that a computer also\index{writing!file management during}
assigns a time and date stamp to each
file you process.  Thus, when you do a directory reading,
you would see something like this:
$$
\begin{array}{lllrr}
\hbox{\tt TROG} & \tt 001 & \qquad  \tt 2357 & \qquad \hbox{\tt 9-21-94} & \qquad \hbox{\tt 11:15pm}  \\ 
\hbox{\tt TROG} & \tt 002 & \qquad  \tt 3309 & \qquad \hbox{\tt 9-22-94} & \qquad \hbox{\tt 2:31pm} \\   
\hbox{\tt TROG} & \tt 003 & \qquad  \tt 3944 & \qquad \hbox{\tt 9-24-94} & \qquad \hbox{\tt 10:42pm} \\   
\hbox{\tt TROG} & \tt 004 & \qquad  \tt 4511 & \qquad \hbox{\tt 9-29-94} & \qquad \hbox{\tt 9:11am} \\   
\hbox{\tt TROG} & \tt 005 & \qquad  \tt 3173 & \qquad \hbox{\tt 10-2-94} & \qquad \hbox{\tt 2:04am} \\ 
\end{array}
$$
\noindent We see here five versions of the paper.  The third column
shows the number of bytes in each version.  The fourth shows the
date on which the editing of that version was completed.  The
last column shows the exact time of completion.

Note that the versions grew in size until, in the wee hours of
October 2, the author decided to discard more than 1300 bytes of the
document; this resulted in version 005.  Is it not reassuring to know
that all the old versions are available, just in case the author
decides to resuscitate one of his/her old turns of phrase?  Whether or
not you are in the habit of examining old drafts of your work, you
will find it psychologically helpful to have all the old versions. 
When the work is complete you can, if you wish, discard all the
drafts but the final one.  But the fact is that mass storage space is
so cheap and plentiful these days that every draft of every one of
your works, even if you are Stephen King and Tom Clancy rolled into
one, will only take up a small fraction of your hard disc.   
	
If you do your writing with a first-rate text editor, as I do, then
you\index{writing!using a text editor during} have powerful tools at your
disposal (see [SK] for a discussion of text editors). You\index{text
editor} can open several files simultaneously, have several different
portions of the same file open at the same time, and have a bibliographic
resource file open; with an environment like {\tt Windows},\reg you also
can have an optical drive or OnLine thesaurus and dictionary open and be connected to the
Internet---and you can jump from one setting to\index{windows@{\tt Windows}} the
other effortlessly. Given that any trip to the dictionary could take ten
minutes the old-fashioned way, and ten seconds the electronic way, think
of how much time you will save over a period of several years. Again---and
here is {\it the\/} most important point---by using technology you
circumvent the danger of your thought processes and your creative juices
being interrupted.

Even though I am addicted to writing with a computer, hard copy plays an
important role in my writing process. For, after I have written a draft, I
print it out, lounge in my most comfortable chair, and proofread and edit.
There exist methods of proofreading and editing directly on
the\index{writing!important of hard copy of} computer---I shall not go
into them here. But, because of my age and my training, I find that there
is nothing like a paper copy and a red pen to stimulate critical thinking.
You will have to decide for yourself what works for you.

There is a down side to writing with a computer; you can work your way past
this one, but you had best know about it in advance. When\index{computer!disadvantage of writing with} you create a document on a computer
system---especially if you use a sophisticated computer typesetting system
like \TeX---then the printed copy looks like a finished product. This
makes it even more difficult than usual for you, the author, to see the
flaws that are present. Even with handwritten copy you will have
difficulty seeing that certain paragraphs must go and others must be
rearranged or rewritten. But, when the MS is typeset, the product looks
etched in stone. One cannot imagine how it could be any more perfect.
Believe me, it can always stand improvement. You will have to retrain
yourself to read your typeset work critically.

(For the flip side of the last paragraph, consider this.  I was
recently asked, by an important publishing house, to evaluate a
manuscript for a textbook that they were considering developing.  The
manuscript was {\it handwritten}.  This flies in the face of all that
is\index{book!handwritten manuscript for}
holy; a manuscript that is going to a publisher should always be typed or
word processed. In any event, I took what I was given and wrote my report.
But this was a difficult process for me. I had to keep telling myself that
this was {\it not\/} a rough set of notes, that it was a polished
manuscript---even though it was handwritten and {\it looked\/} like a
rough set of notes. Play this paragraph off against the last one for
a\index{writing!form over substance in} lesson about form over substance.)

And now a coda on backups. If you use a computer for your work, then
develop the habit of doing regular backups. The ``by the book''
method\index{computer!backups} for doing backups is the ``modified Tower
of Hanoi'' protocol. This gives you access to any configuration that your
hard disc has had for the past several weeks. Not all of us are up to that
level of rigor. But do {\it something}. At least once per week, back up
all your critical files to an external hard drive, a flash drive, or some
other mass storage device. Losing your {\tt C:} drive is
analogous to having your house burn down. It is an experience that you can
well do without. Regular backups are a nearly perfect insurance against
such a calamity.

These days it is a good yoga to back up your files to the cloud.   On the one
hand, the cloud is merely\index{cloud!backup to} somebody else's hard drive.  No different from your
hard drive.  On the other hand, the cloud is very well maintained and is 
accessible from anywhere---with a notebook computer, a tablet, or a smart
cell phone.  And it is extra security for your data.

%% Section 6.2

\section{Word Processors}

I have already indicated in Section 6.1 the advantage
that working on a computer has to offer.  Next I shall specialize
down to word processors and what they do.  (I do this in part so
that, when you read Section 6.5 about \TeX, you will appreciate
the differences.)
 
A word\index{word!processor} processor is a piece of software; 
you use this device for entering text on the computer
screen, and for saving the text on a storage device (usually
a disc).  You engage in this process by striking keys on a keyboard---very
similar to typing.  The word\index{word!processor, uses of}
processor performs many useful functions for you:
\begin{enumerate}
\item  When you get to the end of a line, the word processor jumps
to a new line---you do not have to listen for a bell, or keep one
eye on the text, as you did in the days of typewriters.
\item  The word processor allows you to insert or delete text,
or to move blocks of text from one part of the document to another,
with ease and convenience.  You can create a new document (such as
a letter) by making a few changes to an existing document.
\item  The word processor right justifies (evens up the right margin)
of your document.  This process
results in a more polished look.  
\item  The word processor can check your spelling.
\item  The word processor communicates with your printer, and ensures
that the document is printed out just as it appears on the screen
(this is what we call {\tt WYSIWYG}, or ``What you see is what you get'').
\item  The word processor enables you, if you wish, to incorporate
graphics into your document.
\item  The word processor allows you to\index{wysiwyg@{\tt WYSIWYG}} perform ``global search
and replace'' functions.  For example, if you are writing a paper
about mappings, and you decide to change the name of your mapping 
from $F^*$ to $G_k$, then this can be done {\it throughout the paper\/}
with a few keystrokes.
\item  The word processor allows you to 
select from among several different
fonts:  roman, boldface, italic, typewriter-like, and so forth.
\end{enumerate}							    

These days, most professional people prepare their documents on a
word processor.  Using a word processor saves time, money, and
manpower.  From the point of view of a mathematician, a word
processor is not entirely satisfactory.  The primary reason is that a
word processor will not typeset mathematics in an acceptable
fashion.  A typical word processor can display {\it some\/}
mathematics, but not in a form similar to what you would see in a
high quality book.  The word processor cannot treat complicated
mathematical expressions:  a commutative diagram, the quotient of a
matrix by an integral, or a matrix with entries that are themselves
matrices. Even for simple mathematical expressions, such as a
character with both a superscript and a subscript, the output from a
word processor is\index{word!processor, deficiencies of} nowhere
near the quality that one would see in a typeset book.  There are patches
you can buy---for Microsoft {\tt Word}, for instance---that enable
some mathematical formulas.  Among these are {\tt MathType} and
{\tt MathML}.  But they are of nowhere near the
quality that \TeX\ outputs.  Certainly not the quality of a finished
book.  When it comes to delicate matters of kerning and other spacing
and formatting issues, word processors are limited in their abilities.
And this is in the nature of things, just because
a word processor is {\tt WYSIWYG}.

Outside of mathematics---in the {\it text\/}---word processors fall
short in that they do not {\it kern\/} the letters in words;
many word processors use monospaced\index{kerning}\index{monospace type}
fonts, just like a typewriter. 
This fact means that the word processor does not perform the delicate
spacing between letters---spacing that {\it depends\/} on which two
letters are adjacent---that is standard in the typesetting process. 
The word processor does not offer the variety of fonts, in the
necessary range of sizes, that is ordinarily used in typesetting.
The word processor does not have sufficient power to adjust horizontal
and vertical spacing on the page---both essential for the demands
of quality page composition.   

Put in other terms, a word processor is constrained by the fact that it is
{\tt WYSIWYG} (``what you see is what you get'').\index{wysiwyg@{\tt WYSIWYG}}
Even a high quality screen is no more than 150 pixels (dots) per inch, while high
quality printing is 2400 dots or more per inch. Since a word processor
prints {\it exactly} what appears on the screen, it can format with no
more precision than what can be displayed on the screen. \TeX, by
contrast, is a markup language. It gives typesetting and formatting
commands. It can position each character on the page within an accuracy of
$10^{-6}$ inches.

Because a word processor is {\tt WYSIWYG},
any file produced by a word processor will contain hidden formatting
commands.  One side effect of this simple fact is that if you
cut out a piece of text from a word processor file and move it
to another part of the file then it may not
format properly.  As an instance, suppose that you have a displayed
quotation (such a display usually has text with wider margins and
space above and below).  Snip that out using standard commands
for your word processor and drop it in elsewhere; it will not
format correctly and you will waste a lot of time fixing it up.
Because \TeX\ is a markup language, it does not suffer this
formatting malady.  One of the beautiful features of \TeX\ is
that you can cut and move a fantastically complicated display
and it will not change one iota.

Finally, no word processor is universal.  There are too many word
processing systems.  They are all compatible to a degree, but not in
the way that they treat mathematics.  Thus, again, if one is doing
mathematics\index{word!processing, non-portability of}
using a word processor on the Internet then one will be
hindered.

%% Section 6.3

\section{Using a Text Editor}

Text editors are, primarily, for the use of programmers.  A
programmer wants an environment for entering computer code; the code
will later\index{text editor}
be {\it compiled\/} by a JAVA compiler, a PHP compiler,
a C$++$ compiler, or some other compiler.  Thus a text editor should
not perform value-added features to the code that has been entered: 
there should be no formatting commands, no instructions for the
printer, no hidden bytes, or any other secondary data.  A file
created with a text editor should comprise only the original {\tt ASCII}
code, together with space and line break commands.

In today's world there are rather more sophisticated tools,
like the open source environment {\tt Eclipse}.  In the\index{eclipse@{\tt Eclipse}}
words of {\tt Eclipse}
\begin{quote}
Eclipse provides integrated development environments (IDEs) 
and platforms for nearly every language and
architecture. We are famous for our Java IDE, C/C++, JavaScript and PHP
IDEs built on extensible platforms for creating desktop, Web and cloud
IDEs. These platforms deliver the most extensive collection of add-on
tools available for software developers.
\end{quote}
For the purposes of doing \TeX, a standard text editor
is sufficient for most people.  The text editor that
comes with PC-\TeX\ is customized for \TeX\ users.  And
a text editor like {\tt Crisp} can be customized by the user
for \TeX\ or for any other application.

A document printed directly from a file created with a text editor
would look just like what you see on your computer
screen---typewriter-like font and all---with a ragged right
edge and with old-fashioned monospacing.  Such a document might be
acceptable for an in-house memo, but is not formatted in a manner that
would be suitable for public use.  Thus why would a mathematician
want to use a text editor?

Today, \TeX\ is the document creation utility of choice for
mathematicians (see Section 6.5).  Apart from its flexibility and the
extremely\index{tex@\TeX!vs.\ word processors}
high\index{tex@\TeX!and text editors} 
quality of its output, \TeX\ is also infinitely portable
and it is the one system that you can depend on most (and soon all)
mathematicians knowing.  If you want to work with a mathematician in
Germany, using the Internet, and if you were to say to your
collaborator ``let's use {\tt OpenOffice}\reg,'' then you would be
laughed right off the stage.  The only choice is \TeX\ (or one of its
variants, such as \LaTeX\ or AMS-\TeX).  And the point is this: \TeX\
is a high level computing language (and also a {\it markup\/}
language---see Section 6.5).  You create a \TeX\ document using a
text editor. 

Many a \TeX\ system comes bundled with its own text editor.  Usually,
such a bundled editor has many useful features that make it
particularly easy to create \TeX\ documents.  If you are a PC user,
however, then you are accustomed to selecting your own software.  {\tt Windows}
comes\index{notepad@{\tt NotePad}} with a perfectly serviceable text editor called {\tt NotePad};
many popular word processors, such as {\tt OpenOffice}\reg\
and {\tt Word}\reg, have a ``text editor'' mode.  But much more
sophisticated text editors may be purchased commercially.  One of the
best is {\tt Crisp}\reg\ (which is a version of the {\tt UNIX} 
editor {\tt emacs}\index{emacs@{\tt emacs}}\index{crisp@{\tt Crisp}} 
that has been adapted for the PC)  
A good text editor can be customized for specific
applications, allows you to open several documents and several
windows\index{text editor!customization of}
 at once, has sophisticated search and cut-and-paste
operations, and will serve you as a useful tool.

%% Section 6.4

\markright{6.4.  SPELL-CHECKERS, ETC.}
\section{Spell-Checkers, Grammar Checkers, and the Like}
\markright{6.4.  SPELL-CHECKERS, ETC.}

The great thing about a document created on a computer is that the
document is stored on your hard disc as a computer file.
Thus\index{spell-checkers} your document has become a sequence of
bytes.  In most cases, your document in electronic form will consist
primarily of {\tt ASCII} code---{\tt ASCII} is the
international code for describing the characters that appear on your
computer keyboard and your computer screen.  

A computer file, consisting of a sequence of bytes, is grist for
your computer's mill.  The file is data ready to be manipulated.  Apart
from sending the data to a screen or to a printer, what else can
your\index{computer!files, processing of}
computer do with it?  Here are some options:

\begin{itemize}
\item It can check the words in the file for spelling.
\item It can check for repeated words, misused words, omitted words.
\item It can check grammar and syntax.
\item It can check style.
\end{itemize}

At this writing, spell-checkers are highly sophisticated tools. A good
spell-checker can zip through a 10,000 word document of ordinary text in a
minute or two. It will flag a word that it does not recognize, suggest
alternatives, and ask you what you want to do about it. It will catch many
standard typographical errors, such as typing ``naet'' for ``neat,'' or
such as typing the word ``the'' at the end of line $n$ and also at the
beginning of line $(n+1)$. Of course it will\index{spell-checkers} also
flag most proper names, archaic spellings, and many foreign words and
mathematical terms. As you use your spell-checker, you can augment its
vocabulary (which is performed semiautomatically, so requires little
labor), hence your spell-checker becomes more and more accustomed to {\it
your particular writing}. Given that a spell-checker requires very little
effort to learn and use, and that it can only add to the precision of your
document (it suggests changes, and makes only those that you approve),
\hbox{you would be foolish} not to use a spell-checker. \ {\it However:} 
Never allow the spell-checker to lull you into a false sense of
security. To wit, the ultimate responsibility for correct spelling lies
with you (see below for more on the limitations of spell-checkers).

If you use a garden variety spell-checker on a \TeX\ document, then
you will be most unhappy.  The spell-checker will flag every
\TeX\ command (words beginning with \verb@ \ @) and every math
formula (set off by \verb@ $ @ signs).  You will find processing
even a short \TeX\ document to be an agony.  Fortunately, the spell-checker 
{\tt MicroSpell}\reg\ has a ``\TeX\ mode''; in that
mode, {\tt MicroSpell}\reg\ knows to ignore \TeX\ commands and\index{microspell@{\tt MicroSpell}}
math formulas.

Do not use a spell-checker foolishly.  If you intend to write the
word ``unclear'' and instead write ``ucnlear'' (a common
transposition error), then the spell-checker will certainly tell you,
and that is useful information.  But if you intend ``unclear'' and
instead write ``nuclear,'' then the spell-checker will forge blithely
ahead---because ``nuclear'' is a {\it word\/}, and a spell-checker will
only flag non-words.  If you mean to say that someone
is ``weird'' and instead you say he/she is ``wired,'' then your message
may still trickle through; but your spell-checker will not help you
to get it right.  The lesson is clear (rather than unclear):  if your
document passes the spell-checker, then you know that certain
rudimentary errors are not present; however certain other, more
sophisticated, errors could be present.  Will you have to catch them
yourself, with old-fashioned proofreading?  Read on.

An interesting operating system tool, 
available also in {\tt UNIX} and other computer
environments, is the ``word counter.''  Run an {\tt ASCII} file through the
word\index{word!counters} 
counter and this utility reports {\bf (i)} how many words
there are in your document, {\bf (ii)} what is the most frequently
used word and how often the word is used, {\bf (iii)} what is the second
most frequently used word and how often is it used, etc.  This device
can easily be construed as an example of using a computer to do
something just because it can be done.  But we all have personal
foibles, and the word counter can help to detect them.
I, for example, tend to overuse the word ``really''; I had
to make special passes through this manuscript to weed out many
occurrences of that word.  A more informed opinion about which words
you overuse can be made if you use word-counting software.  If you
use the word ``really'' more often than you use the word ``the,''
then you may be in trouble; however, if you use it less often than you
use the word ``flagellate,'' then a different conclusion is in order.

I read of a professional author being stymied for a period of a year
as a result of using word counting and other software.  He ran one of
his famous stories through the software, and it pointed out certain
words and phrases that he overused.  Thereafter, whenever the author
was writing and began to use one of the pegged words or phrases, he
panicked.  The matter became worse and worse, and he eventually
developed a writer's block.  It took him considerable therapy, not to
mention stress and hard work, to defeat this block.

One of the big events in the world of finance in the last many years
is the invention of the Black/Scholes option pricing scheme.  This
very sophisticated technology uses stochastic integrals---{\it very} advanced
mathematics.  It\index{Black/Scholes option pricing scheme} 
won the Nobel Prize for Myron Scholes (Fischer Black
had died).  Naturally Scholes's school, Stanford University, wanted
to make a big deal out of their distinguished faculty member's encomium,
so an article was written for the in-house newsletter.  Unfortunately,
some foolish editor ran a spell-checker on the article and ended
up changing every occurrence of ``Myron Scholes'' to ``moron schools.''
And that is the way that the article appeared.  Woe is us.

We conclude with another anecdote, courtesy of G. B. Folland.  One
of Folland's publishers used a spell-checker that recognized the
word ``homomorphism'' but not the word ``homeomorphism.'' The result?
The copy editor changed every instance of the latter to the former.
The original manuscript contained several dozen of each.  Now do
you see how a spell-checker can get a person into trouble? 

%% Section 6.5

\section{What Is \TeX\ and Why Should You Use It?}

\TeX, created by Donald Knuth in the early 1980s, is an electronic
typesetting\index{tex@\TeX!what is} system.  Designed by a mathematician,
specifically for the creation of mathematical documents, it also
is a versatile tool in other typesetting tasks.  In fact \TeX\
is used in many law offices, and is also used to typeset {\it TV Guide}.
The reference [Kn] tells something of the philosophy behind
the creation of \TeX.

What makes \TeX\ such a powerful tool?  First, \TeX\ is almost
infinitely portable.  A \TeX\ document created with a Macintosh
computer in Hong Kong can be sent over {\it e}-mail to a PC user
in Sheboygan who in turn can send it on to the user of a Cray I
in Bielefeld.  During this process, there are never any problems
with compiling, printing, or viewing.

The book [SK] already contains this author's efforts at describing
the inner workings of \TeX\ and how to learn them.  I shall not
repeat that material here.  Instead, I shall say just a few words
about how \TeX\ is used.  

\TeX\ is {\it not\/} a word processor.  Instead \TeX\ 
is what is called a ``markup
language.'' ``Markup language'' means 
that, in your \TeX\ document (created with a\index{tex@\TeX}
text editor---not a word processor), you enter commands that tell
\TeX\ what you wish to have\index{tex@\TeX!vs.\ word processors} 
appear on each page, and in what
position.  \TeX\ allows you to position each character on the page to
within $10^{-6}$ inch accuracy.

If you think about all the material that appears on a typeset page,
then what is described in the last paragraph sounds arduous---like
it is simply too much trouble.  Fortunately, \TeX\ performs most
typesetting tasks automatically.

If you are typesetting ordinary prose, then you simply type the words
on the screen, with spaces between consecutive words.  With \TeX, you
can leave any amount of space between successive words in your source
code; you can also put any number of words on each line of code. 
\TeX\ will choose the correct spacing, and the correct number of
words for each line, when it compiles the document.  You indicate a
new paragraph by leaving a blank line. There is almost nothing more
to say about typesetting text:  \TeX\ spaces letters correctly,
it chooses the correct amount of
space to put between words, how to put space between paragraphs, and
so forth.  It makes each line come out flush right, and ensures that
each page contains the correct number of lines---not too many and not
too few.\index{tex@\TeX!spacing}\index{tex@\TeX!kerning}

For mathematics, there are English-language-like commands that
tell \TeX\ just what you want.  I will\index{tex@\TeX!sample code}
present just one example:  The code

\begin{verbatim}
Now it is time to do some mathematics---a task for 
which, given that we have spent many years at the 
university, we are eminently well prepared.  Our work
is inspired by the identity $X(1 + X) = X + X^2$.

Let us consider the equation
$$
\int_X^{X^2 - X} {{\alpha^3 
    + 17{{\alpha} \over {\alpha-2}}} \over 
     {{\alpha-5} \over {\alpha + 1}}} \, d\alpha
      = \hbox{det} \, \left ( \matrix 
     {X^2 & 3X & X \cr
   {{X^3 - 4} \over {X + 1}} & \sin X & \log X \cr
   {{X} \over {X+1}} & \hbox{erf}\, X & \sqrt{X} \cr
       }  \right ) 
$$
which has been a matter of great interest in recent years.
\end{verbatim}

\noindent would typeset as
\medskip \hfill \break

Now it is time to do some mathematics---a task for which, given that
we have spent many years at the university, we are eminently well
prepared.  Our work is inspired by the identity $X(1 + X) = X + X^2$.

Let us consider the equation
$$
\int_X^{X^2 - X} {{\alpha^3 
    + 17{{\alpha} \over {\alpha-2}}} \over {{\alpha-5} \over {\alpha + 1}}} \, d\alpha
        = \hbox{det} \, \left ( \matrix 
       {X^2 & 3X & X \cr
        {{X^3 - 4} \over {X + 1}} & \sin X & \log X \cr
         {{X} \over {X+1}} & \hbox{erf}\, X & \sqrt{X} \cr
       }  \right ) 
$$
which has been a matter of great interest in recent years.

I should stress that you should {\it not} use a word processor to
create your \TeX\ source file.  A word processor has a very large
number of hidden binary commands for formatting and visualizing.
These hidden commands will confuse the \TeX\ compiler, and give
you output that is quite different from what you want.  See Section
6.3 on text editors to find out what a text editor is and how to
get one.
   
Even though you may not know \TeX, you should have little difficulty
seeing the correspondence between the code that is entered and
the resulting output.  (Note that the single dollar signs signify
material to be typeset in ``in-text'' math mode; the double dollar signs
tell \TeX\  first to enter, and then to exit, ``displayed'' math mode.)

After you have created an {\tt ASCII} file with your text editor,
call it {\tt myfile.tex}, then you compile it with the command
{\tt tex myfile} (variants are {\tt latex myfile} and {\tt amstex myfile}).
This\index{tex@\TeX!files, compiling}
creates the ``device independent file,'' called {\tt myfile.dvi}.
The {\tt dvi} file can be ported to a printer, to a screen, or translated
to {\tt Postscript}\reg.\index{postscript@{\tt Postscript}}

As you can see from the preceding example, \TeX\ does a magnificent
job of typesetting mathematics.  Usually no human intervention
is required in order to obtain the quality and precision that
you desire.  

One interesting feature of \TeX\ is that you cannot expect to see on
the screen exactly what you will obtain in your printed output.  For
even a high quality screen has resolution about 150 or so
pixels\index{tex@\TeX!previewing}
per inch. Today, printers have a resolution of 2400 or more
dots per inch. 
The {\tt Preview}\reg\ programs\index{preview@{\tt Preview}}
that come with \TeX\
allow you to view your document to the extent of seeing where the
various elements appear on the page---sufficient for doing elementary
editing.  But, to view the final output accurately, you must print a
hard copy.

It is fairly straightforward to import a graphic into
a \TeX\ document.  An example of the relevant command is\index{tex@\TeX!graphics in}

\begin{verbatim}
\begin{figure}
\centering				     
\includegraphics[height=2.25in, width=2.75in]
    {k:/pubbooks/harmonic/figures/fig1-1.eps}
\caption{A ``pop'' or ``click''.}
\end{figure}
\end{verbatim}
You can see that we use the \verb@\begin{figure}@ command to create a ``float.''
This means that the figure does not have a fixed position but is floated around
for a best fit.  The \verb@\includegraphics@ command specifies height and
width for the figure, and also calls in the specific graphics file.  There
is also a caption command.

\TeX\ was originally designed with the notion of maximum power and
flexibility in mind; Knuth planned that each discipline would develop its
own style files to tailor \TeX\ to its own uses. The variant\index{tex@\TeX!flexibility of} \LaTeX, created by Leslie Lamport, endeavors to serve all
end users. More specialized style\index{tex@\TeX!vs.\ \LaTeX} files are
available from the American Mathematical Society (to give just one
example); these enable the AMS-\TeX\ user to typeset a paper in the style
of any of the AMS primary journals.

There is a whole new world of document-preparation tools available
today.  As a semi-neanderthal, I would be more than sympathetic if
you do not want to dive into all the graphic and typesetting and
electronic features that I have described here.  In fact these tools
are best learned in gradual stages.  The learning curve for \TeX\
alone is rather steep, although the book [SK] makes strides toward
jump-starting the learning process.  My recommendation is to begin by
learning some form of \TeX.  \LaTeX\ is a particularly popular form
of \TeX, and one favored by publishers (because it is more structured
and steers the author toward more standard formatting styles than
does Plain \TeX).  The reference [SK] creates an
accessible bridge between Plain \TeX\ (the most flexible \TeX\ tool)
and \LaTeX\ (the least flexible \TeX\ tool).  Most mathematics
departments have the hardware, the software, and the expertise to
make it easy for you to learn \TeX.  This software is one of today's
standard mathematical tools. You are shooting yourself in the foot
not to learn it.  

For graphics, you may find that Adobe {\tt Illustrator} or Corel {\tt
DRAW!} or (on a {\tt UNIX} system) {\tt xfig} is a useful utility. Any of
these devices will output graphics in {\tt *.pdf} format or {\tt *.bmp}
format or {\tt *.eps} format or dozens of other popular graphics formats.
\index{adobeillustrator@Adobe {\tt Illustrator}}\index{coreldraw@Corel {\tt DRAW!}}

There are various front ends available to make \TeX\ more user friendly.
For myself, \LaTeX\ is sufficiently friendly.  I can work with it
comfortably.

There are a number of variants of \TeX. First, one should note that there
are two obvious places to get your own copy of \TeX.
One\index{tex@\TeX!variants of} is the company PC-\TeX, located in San
Francisco. The other is the freeware version of \TeX\ called Mik\TeX. Both
are very good.

As previously noted, Donald Knuth did not market \TeX.  In fact he
{\it gave} it to the American Mathematical Society.  The AMS conscripted
Michael Spivak to creat AMS-\TeX, which is the AMS's version of \TeX.
There is also LAMS\TeX.  This is the AMS version of \LaTeX.   The fact
remains that most mathematicians, and most publishers, use \LaTeX.  It
has become the {\it lingua franca} of the \TeX\ world.

Why does the AMS need its own version of \TeX?  Well, from the AMS point
of view, Knuth's version of \TeX\ has certain liabilities.  For one thing,
Knuth does not like the blackboard bold font, so he did not include it
in \TeX.  These are characters\index{amstex@AMS\TeX} l
ike $\RR, \CC, \ZZ$.   Most mathematicians
are wedded to this font.  The AMS added other accoutrements to \TeX\ to
make it more useful to mathematicians.  Similar comments can be made
about LAMS\TeX.\index{blackboard bold fonts}

Michael Spivak invested a great deal of time and money and effort
in developing a new set of fonts (the MathTime Professional fonts) 
that\index{Spivak, M.}\index{Spivak, M.!fonts} are an alternative to the Computer Modern
font\index{fonts!MathTime} of Knuth.  Spivak's observations included that
\begin{itemize}
\item Some of the Greek letters in Computer Modern were hard to
distinguish from others.
\item The large parentheses in Computer Modern are not rounded as
they should be.
\item The root signs in Computer Modern are not formulated properly.
\end{itemize}
There are a number of other technical ways in which the MathTime fonts
are an improvement over the Computer Modern.  Spivak's fonts may be
purchased on the Internet.  There is a trial version of the fonts
that is free.

The world of \TeX\ has become a way of life for many people.  The 
\TeX\ Users Group (or TUG) is an organization dedicated to promoting
and developing \TeX.  The TUG newsletter is a fascinating read for
those interested in \TeX.\index{tex@\TeX\ Users Group}

We conclude this section with an amusing story about \TeX.  Mathematician
Pete Casazza of the University of Missouri wrote a book several years
ago.  Naturally he wrote it in \TeX.  At some point in the publication
process, his publisher sent him a package containing edited page proofs
and a disk with the \TeX\ file on it.  Now the University of Missouri
is in Columbia, Missouri.  But the U. S. Postal Service sent the
package to the {\it country} of Colombia.  Well, the authorities
in Colombia opened the package and were most curious as to what
was on the disk.  They examined the \TeX\ file and saw all the
dollar signs (recall that dollar signs are used to format math
formulas).  They rapidly concluded that these were the books for
an illegal drug cartel.  After a few months of hard study, 
they were unable to determine just what the file was telling them.
So they enlisted the help of the FBI.  The FBI was able to discover
that\index{Casazza, P.} this was in fact a \TeX\ file containing nothing but mathematics.
So, after a delay of a good many months, Pete Casazza finally got
his package and was able to proceed with the publication of his book.

%% Section 6.6

\section{Graphics}

As indicated elsewhere, the most common
method for including graphics in a book is still to create them {\it
separately}, each on its own page.  The drawings could be created by
hand, with pen and ink.  Or they could be produced with Corel {\tt
DRAW!}\reg, or {\tt MacDraw}\reg, or Adobe {\tt Illustrator}\reg, 
or {\tt xfig}, or
any\index{graphics created with a computer}
number of other packages.  To repeat, each figure should be on a
separate piece of high quality drawing paper (available from any
store that carries art supplies) and drawn in dense black ink.  Use a
proper drawing pen---not a ball point, or a rolling writer, or a
pencil.  Best is to draw the figures (considerably) larger than
they will actually appear in the book, in thick dark strokes.  When
they are photographically reduced to fit, then the pen strokes come
out sharper, denser, and darker.

Each figure should be labeled clearly:  a typical label might
be
\begin{quote}
\bf Chapter 3 \ \ Section 2 \ \ Figure 5
\end{quote}

\noindent Correspondingly,\index{figure!labels}
 somewhere in Section 2 of Chapter 3 there
should be a space set aside for this figure, and it should be labeled
``Figure 5.''  (I am assuming here, for simplicity, that you are
producing your document in some version of \TeX; if not, then forget
about leaving a space in the document, but {\it do\/} put a label in
the margin.)  And be sure that the text contains a specific reference
to each figure by name (label); do not leave it to the reader to determine what
figure goes with which set of ideas.  (The same remark applies, of
course, to tables.)  It helps, though it is not mandatory, to
give each figure a title and a caption.

Drawing good illustrations for your work is an art.  
A good figure is not too busy, does
not have extraneous information or extraneous penstrokes, and
displays\index{figures!creating informative}
its message prominently and clearly.  The books [Tuf1] and
[Tuf2] by Edward Tufte will give you a number of useful pointers on
how to develop powerful graphics for your work.

Of course we all know\index{Tufte, E.} that there are copious electronic tools for
creating artwork in your manuscript.  Just as an example, many
versions of \TeX\ have simple commands, such as \verb@\psfig@,
that\index{figures!software for creating}
allow you to import an encapsulated {\tt Postscript}\reg\ file
into your document.  In one\index{postscript@{\tt Postscript}} common 
scenario, a \verb@\special@
command\index{specialcommand@\special command} 
insets raw printer commands into the file that will
communicate with your printer.  The result is that your {\tt
Postscript}\reg\ figure\index{postscript@{\tt Postscript}} appears right on the printed page ({\it
provided} that you have a {\tt Postscript}\reg\ printer or know how
to use {\tt Ghostscript} to make {\tt Postscript}\reg\ talk to a 
non-{\tt Postscript}\reg\ printer).  Some\index{postscript@{\tt Postscript}} 
versions of \TeX---such as
{\tt Personal \TeX}\reg---understand several other graphics languages
as well.  For example, the Hewlett-Packard language {\tt PCL} is a graphics
protocol designed for use with certain HP printers.  And many
graphics programs give you a choice of several different graphics
output languages; these could include {\tt ps}, {\tt eps}, 
{\tt bmp}, or {\tt wmf} graphics
images. The documentation for your \TeX\ software (for instance {\tt
Personal \TeX}\reg) will explain precisely which graphics languages
it can handle and how it does so.

And now a caveat about {\tt Mathematica}\reg, {\tt
Maple}\reg, and the like. These, too, are small miracles.  If you
need\index{mathematica@{\tt Mathematica}}\index{maple@{\tt Maple}}
 to draw a hyperboloid of one sheet, or the graph of $z = \log\bigl
(|\sin (x^2 + y^3)|\bigr )$, then there is nothing to beat {\tt
Mathematica}\reg.  I recommend that you use it. {\tt
Mathematica}\reg\ will output your figure in encapsulated {\tt
Postscript}\reg, for\index{postscript@{\tt Postscript}} storage on your hard disc, and in principle this
file can be imported into your document.

A final note:  ask {\tt Mathematica}\reg\ to graph a horrendously
complicated function of two variables, and it will do so in an
instant.  Such tasks are what {\tt Mathematica}\reg\ is designed to 
perform.  And it
will provide the labels on the axes automatically.  But endeavor to
draw a rectangle or triangle, and to label the vertices in your own
fashion, and it may take you an hour.  Conversely, I can hand draw the
rectangle or triangle and provide the labels in five minutes.  But it
could take me hours to graph the function.  Instead you
should draw the triangle or rectangle using Corel {\tt DRAW!}
or Adobe {\tt Illustrator}.  {\it Use the proper tools
in the proper context.}

%% Section 6.7

\markright{6.7.  THE INTERNET AND {\tt HYPERTEXT}}
\section{The Internet and {\tt hypertext}}
\markright{6.7.  THE INTERNET AND {\tt HYPERTEXT}}

Just a word about {\tt hypertext}, and about electronic publishing in
general.  The spirit of electronic publishing is to bypass the
traditional hard copy of published materials, and instead make the
materials available on the Internet and the World Wide Web.  Readers
would\index{internet}\index{hypertext@{\tt hypertext}} 
be identified and would pay either by buying a password or by
paying the publisher to make materials available to a {\it
particular} CPU with a particular identification number (the IP
address---given\index{IP address} by four octets of code).

A part of this new electronic publishing environment is {\tt hypertext}.
With {\tt hypertext}, certain words or phrases in the electronic document
appear in an accented form---often in a different color or
underlined.  If the reader ``clicks'' on the accented word, then he/she
is ``jumped'' to a cognate item.  For instance, if you are reading a
book on the function theory of several complex variables, you come
across the word ``pseudoconvex,'' and you cannot recall what it
means, then---instead of madly flipping through the book trying to
find the definition (that is the old way)---you click on the word and
are jumped either to the passage that contains the definition, or
perhaps to a lexicon, or perhaps to a menu that offers you several
options. Alternatively, you could click on a reference to another
book or paper and you would be jumped to the reference---to the {\it
actual text of the reference}---no matter where in the world the
source is.  Or you could click on the name of a mathematician who is
mentioned in the text and you would be jumped to his/her home page, or to
his/her publication list.   Yet another scenario is that you could click
on an icon or a button and bring up animated graphics.

Clearly {\tt hypertext} is an amazing device, and the possibilities
that it offers are vast and amazing.  In the coming years, as the
world decides what role electronic publishing will have in our lives,
how to charge for it, how to market it, how to archive electronic
documents, and so forth, we will see more and more electronic books
and documents.  For now, matters are in a developmental stage.
 
There now exist many electronic journals.  An electronic journal is
one in which all transactions---submission, remanding to a referee,
referee's report, editorial decision, and publication---are 
executed over the Internet.  No hard copies of the journal are produced,
nor archived.

One interesting features of all-electronic journals is that they
make the notion of an ``issue'' or a ``volume'' pass\'{e}.  Papers can
be published sequentially rather than in batches.

Several of these new electronic journals
are ``startup'' journals, run ``for love'' by an individual from
his/her office computer.  Others are\index{journal!electronic}
institutionalized, but are still free.  Still others require paid
subscription.

The notion of {\it Open Access} has begun to play a major role
in journal publishing---not just in mathematics but across the
sciences.  There are now over 10,000 open access journals.
Open access journals have been heavily promoted by NIH (the National
Institutes of Health).  In fact any scientist who receives NIH funding
is {\it required} to publish his/her\index{Open Access journals} results in an Open Access journal.

Just what is Open Access (OA)?  The basic principle of an OA journal
is that anyone can read the articles in the journal without paying
for a subscription.  There are typically no subscriptions to an OA
journal.  This is fine for a journal run by a group of volunteers on
a personal computer.  But if a commercial publisher chooses to publish
an OA journal then someone has to pay the publisher's expenses, and the
publisher needs to make a profit.  So a fee is levied against the author,
and the fee is generally nontrivial.  Often the fee is in the thousands
of dollars.

Proponents of OA hope that the model for journals will shift:  from the
university paying many thousands of dollars per year in subscriptions to
instead the university paying many thousands of dollars per year in author
fees.  Some universities have in fact bought in to this new model.  Most
have not.  In some cases the author will have a grant that will pay the fee.
In most cases (especially in mathematics) this is not so.   Some (see [Ewi]) have
argued that OA is turning scholarly journals into vanity presses.  Others have
observed that dishonest publishers will accept a great many substandard papers
just to generate a strong cash flow.

The main point here is that the Internet has opened up a world of new
possibilities for scholarly publishing.  OA is just one of these, and it
is something that we all need to learn to live with.   The article [Ewi]
presents a cogent analysis of the OA movement.

There are different flavors of Open Access, and a good place to read about
the details is [Sub]. OA delivered by journals is called {\it gold}. OA
delivered by repositories is called {\it green}. {\it Libre} OA involves
removing permissions.  The Open Access Newsletter, created by Peter Suber,
will give you the chapter and verse on OA from the point of view
of its partisans.  Its URL is
\begin{verbatim}
http://legacy.earlham.edu/~peters/fos/newsletter/archive.htm
\end{verbatim}

At the beginning of 1995, the AMS (American Mathematical Society)
eliminated the ``Research Announcements'' section from the {\it
Bulletin} of the AMS and is instead
publishing research announcements---with essentially the same
editorial\index{ams@AMS Research Announcements} policies and publishing standards---in electronic form. The
AMS is also making archival/disaster-backup copies of all the startup
electronic journals, strictly as a service to the mathematical
community and at no charge.  As of this writing, the AMS has
initiated several brand new subject-area electronic journals for
which subscribers will pay a modest fee. 

Most electronic journals are run with the same editorial procedures as
for\index{electronic journals!advantages of} a paper journal (these procedures are described in Section 2.7).
The primary differences are two:  
\begin{enumerate} 
\item  With an electronic journal, page limitations are not 
important; therefore an electronic journal can publish longer 
articles and more of them.
\item  The avowed goal of most electronic journals is to generate no
paper---none whatsoever.  Therefore papers are submitted by {\it
e}-mail and forwarded to the editor and to the referee by {\it
e}-mail.  The referee's report is sent back to the editor,
and then to the author, by {\it e}-mail;
and any revisions are resubmitted by {\it e}-mail.  The paper is
published electronically, either on a bulletin board or on a server.
There is never a hard copy of anything. 
The published paper is posted in an output language, such as {\tt
*.dvi} (a Device Independent file generated by \TeX) or {\tt *.eps}
(an encapsulated {\tt Postscript}\reg\ file ) or {\tt *.pdf} 
(an Acrobat\reg\ file).\index{postscript@{\tt Postscript}} 
{\it There are no reprints}, and there are no hard copies of the
journal.  In some cases, CD-ROM versions of the electronic journal
are available for archival purposes.  Usually the
end user can download individual papers for (compiling and) hard copy
printing for personal use.  
\end{enumerate}
                   
Some hard copy journals are now simultaneously
publishing an electronic version.  One interesting
innovation is that some traditional journals make any mathematical
paper available electronically, for a modest charge, as soon as the
paper has been accepted.

The are several advantages of electronic-only journals:  {\bf (1)}  the
journals take no shelf space (a fact that is of immense importance
to librarians), {\bf (2)}  the journals cannot be lost or stolen,
{\bf (3)}  an unlimited number of readers can access the articles at
the same time, {\bf (4)}  (in many cases) individuals can print
out their own hard copies of any given article.

The world of electronic publishing is just opening up, and 
promises new frontiers of publishing activity and also of legal
complications.\index{electronic publishing}  
As an example, the copyright law issues connected
with electronic publishing are immense [Oke].  

Some authors are making entire books available at no charge on the
World Wide Web.  Commercial publishers are also exploring the
publication of electronic forms of books.  In fact some publishers
will propose to an author that a home page be set up for his/her book,
and\index{electronic books!available at no charge} that not simply the book but also a variety of ancillaries appear
on the Web site.  These ancillaries could include relevant papers, a
bibliographic database, exercise books, lecture manuals---you name
it.   In some scenarios, a publisher may develop a version of a book
to which readers may contribute interactively.  
%% Section 6.9

\def\email{{\it e}-mail}

\markright{6.9.  COLLABORATION BY E-MAIL}
\section{Collaboration by {\it e}-Mail; Uploading and Downloading}
\markright{6.9.  COLLABORATION BY E-MAIL}

Writing a collaborative mathematical work is a source of great
pleasure.  It is especially fun when you use {\it e}-mail as a tool. 
Entire chapters can be zapped around the world in an instant.  You
get immediate feedback on your ideas.  In many (but not all) ways, 
collaborating by {\it e}-mail
is like\index{electronic mail!collaboration over} having your
partner in the office next door.\index{email@{\it e}-mail}

Many of us do our work on the computer system at school.  The school
system is probably a network---most likely a {\tt UNIX} system.  If
this describes your working environment, then Internet collaborative
work is a breeze.  You work on your document---using an editor like
{\tt vi} or {\tt emacs}.  When you are ready to share it (call the
paper {\tt ourpaper.tex}) with your collaborator, you will send
it as an \email\ attachment.  Usually it is best to send both
a {\tt *.tex} file and a {\tt *.pdf} file.  Of course you will still
retain the master copy on your system's hard disc.  

If you wish, you can type comments at the beginning of, or in
the middle of, the \TeX\ document. If you precede each line of
the comment material with a \verb@ % @ symbol, then \TeX\
will\index{comments in a \TeX\ document} ignore those lines.
\index{tex@\TeX, comments in a document}

Perhaps you do your work at home, on a PC or a Macintosh, and then
bring the files to work for further processing and {\it e}-mailing.
Thus you copy the files to a flash drive and must ``upload'' the files
from the flash drive to the computer system at school.  You will have to
consult your local guru for the details of this uploading process.
But\index{uploading files}
 note this:  there are differences in the ways that files are
formatted on a microcomputer as compared with a mini or mainframe
computer that uses the {\tt UNIX} operating system.  

The public domain operating system {\tt LINUX} for the personal
computer is a popular choice these days.  {\tt LINUX} is a version
of {\tt UNIX} that is designed for PC-type computers (there is also
a version of {\tt UNIX} designed for Macintosh computers).  In
{\tt LINUX} you can\index{unix@{\tt UNIX}} open either a {\tt DOS}
window or a Macintosh window, and the different operating systems can
talk to each other.  These days, all Macintosh computers run on
LINUX.  One benefit of LINUX is that hackers tend to ignore LINUX machines.
They attach Windows machines because there are many more of them.

\section{Mathematical Collaboration in Today's World}

In the old days---say 100 years ago---mathematical collaboration
was relatively rare.  There were only several hundred research mathematicians
in the world, and each of these sat alone in his or her office
and applied himself/herself to the proving of theorems.  Hardy
and Littlewood, who together wrote more than 100 papers, were certainly
the exception.	 In those days, the only two modes of discourse
were snail mail or meeting face-to-face.  There was the telephone,
but\index{Littlewood, J.\ E.}\index{Hardy, G. H.} 
long-distance calls were considered to be prohibitively expensive.

In today's world more papers are written collaboratively than not.
And there are so many devices to enable this collaboration.  Certainly
collaboration by \email\ is quick and convenient.  Drafts of papers,
written in \TeX\ or \LaTeX\ can easily be sent as \email\ attachments.
One can use {\tt Skype} to have multi-hour rap sessions with a collaborator
on the other side of the world---with no cost to anyone.  

{\tt FaceTime} is an Apple product that works on {\tt iPhone}s,
tablets, and Macintosh computers.  It allows you to speak
to\index{facetime@{\tt FaceTime}} a friend or collaborator and see him/her at the same time.
So, in principle, both of you could be writing on a blackboard
or white board and each could see what the other is doing.
Collaborating with {\tt FaceTime} is almost the same as being
in the same room together.

Of course attending conferences, workshops, and research institutes is
a terrific way to hook up with people who have interests similar
to your own.  After you have established a working relationship with
some of these people, then you can go back to your home institutions
and communicate by one of the methods described above.

There are also many opportunities for students.  The Mathematical
Sciences Research Institute in Berkeley sponsors summer workshops
for graduate students each year.  The National Science Foundation
sponsors the Research Experience for Undergraduates program.  There
are also summer internships.

\section{If You are Not a Native English Speaker}

On the one hand, it is rapidly becoming the case that English
is the default language for writing in mathematics.  If you
want your ideas to be read around the world, you write in English.
French and German are fine, but they do not have the universality
of English.

But many of us are not native English speakers.  Our English may be
serviceable, but it is not perfect.  We may need some help to get
our prose up to snuff for publication.  What are the options?

First, many publishers can offer help with English.  They have staff
professionals who are skilled at helping non-native speakers sharpen
their prose.  They know enough about other languages that they
know which bugs to look for and how to fix them.  Do be sure
to ask your publisher whether they can provide such assistance if
you need it.\index{non-native speakers}

Second, there are professional private writing coaches.  Of course they
are for hire, so you will have to find the resources to remunerate them.
But the expense will probably be worth it.  And you may, in the process,
gain\index{writing!coaches} a valuable ally for future writing projects.

It is also possible that you can get one of your students---either an
undergraduate or a graduate student---who is a native English speaker
to help you with your writing.  The student will likely be thrilled
to be asked to assist with a scholarly task.  And you will enjoy
working with a student on a worthwhile project.

Finally, you may have a generous and friendly colleague who is willing
to give some time to helping you with your writing.  After all, this
is in part why we have colleagues.  And you will feel quite comfortable
working side-by-side with a colleague of your own age and with
similar training.  

As time goes on, if you work at it, your English will become better
and better.  Some of the best writers-in-English that I know are
not native English speakers.  A good example is Benoit Mandelbrot, who
was an extraordinarily gifted writer.  Indeed a good deal of his
success can be attributed to his writing ability.  Another example
is Elias M. Stein at Princeton.  His analysis books have been enormously
influential, and are widely read with great pleasure.
\vfill
\eject

\hbox{ \ \ \ }

\thispagestyle{empty}

\newpage

%% Chapter 7
\chapter{The World of High-Tech Publishing}

\begin{quote}
\footnotesize \sl Never spend more than a year on anything.
\smallskip \hfill \break
\null \mbox{ \ \ } \hfill \rm Jeff Ullman 
\end{quote}

\begin{quote}
\footnotesize \sl The commonest thing is delightful if only one hides it.
\smallskip \hfill \break
\null \mbox{ \ \ } \hfill \rm Oscar Wilde
\end{quote}

\begin{quote}
\footnotesize \sl Not of the letter, but of the spirit:  for the
letter killeth, but the spirit giveth life.
\smallskip \hfill \break
\null \mbox{ \ \ } \hfill \rm The Holy Bible, the New Testament \break
\null \mbox{ \ \ } \hfill \rm The Second Epistle of Paul the Apostle \break
\null \mbox{ \ \ } \hfill \rm to the Corinthians.  Chapter 3, Verse 6.
\end{quote}

\begin{quote}
\footnotesize \sl The road to hell is paved with works-in-progress.
\smallskip \hfill \break
\null \mbox{ \ \ } \hfill \rm  Philip Roth
\end{quote}

\begin{quote}
\footnotesize \sl The road to hell is paved with adverbs.
\smallskip \hfill \break
\null \mbox{ \ \ } \hfill \rm Stephen King
\end{quote}

\begin{quote}
\footnotesize \sl Who wants to become a writer? And why? 
Because it's the answer to everything. It's the streaming 
reason for living. To note, to pin down, to build up, 
to create, to be astonished at nothing, 
to cherish the oddities, to let nothing go down the drain, 
to make something, to make a great flower out of 
life, even if it's a cactus.
\smallskip \hfill \break
\null \mbox{ \ \ } \hfill \rm Enid Bagnold
\end{quote}

\begin{quote}
\footnotesize \sl To gain your own voice, you have to forget about having it heard.
\smallskip \hfill \break
\null \mbox{ \ \ } \hfill \rm Allen Ginsberg, WD
\end{quote}

\begin{quote}
\footnotesize \sl Cheat your landlord if you can and must, but do not 
try to shortchange the Muse. It cannot be done. 
You can't fake quality any more than you can fake a good meal.
\smallskip \hfill \break
\null \mbox{ \ \ } \hfill \rm William S. Burroughs
\end{quote}

\begin{quote}
\footnotesize \sl All readers come to fiction as willing accomplices to your lies. 
Such is the basic goodwill contract made the moment we pick up 
a work of fiction.
\smallskip \hfill \break
\null \mbox{ \ \ } \hfill \rm Steve Almond, WD
\end{quote}

The advent of the computer, and particularly of the Internet,
has completely changed the face of modern mathematical publishing.
There are many new artifacts and features of this world.  And many
new forces at play.  In this chapter we attempt to describe the
key new components of our publishing life.

In the Middle Ages and the Renaissance, books were expensive and rare
and really only available to the priveleged few.  Today many
books are available for free on the Web.  Many encyclopedias
(especially, but not exclusively, {\tt Wikipedia}) are also freely
available.  Many {\tt Google} tools are free.  This truly
is the information age.\index{google@{\tt Google}}\index{wikipedia@{\tt Wikipedia}}

\section{Preprint Servers}

In the old days, when you wrote a math paper, it was typed up by a
manuscript typist using an IBM Selectric typewriter. This typewriter was
special because it had ``element balls'' with special characters such as
math symbols, letters from the Greek alphabet, and special braces and
brackets. It was still necessary to render some symbols by hand with an
inkpen, but the Selectric did most of the work.\index{IBM Selectric
typewriter}\index{photocopy machine}

You would have many copies of this work reproduced on the 
photocopy machine, and you would mail these (with snail mail)
to your colleagues all over the world.  This is how a mathematician
would establish his/her reputation and make his/her mark on
the profession.  You could not afford to wait for the paper
to be published; this could cause a delay of a few years, and
your likelihood of getting scooped was nontrivial.  You
had to get the word out right away.  This is how it was done.\footnote{Of
course you could also give seminars and speak at conferences.  This was
an important part of the profession.  It was also quite common to send
out postcards announcing results.  But electronic media
were not at all available fifty years ago.}

The trouble with the system just described is that it meant that
well established people at the top universities heard all the
new developments first.  More obscure mathematicians who were
not well connected were generally left out of the loop.  They
could go to conferences and get some hints about new developments.
But they did not get their information in a timely fashion.

Now things have changed. There are a good many preprint servers
that serve as repositories for new mathematics. What is a
preprint server? It is a Web site where you can post your new
paper.\index{preprint server} Now we must understand clearly that a paper posted on a
preprint server is not refereed or vetted in any way. It is
just posted for all the world to see. And, indeed, absolutely
anyone can view or download or print the papers posted on a preprint server.
And virtually anyone can post on a preprint server (although
some, like {\tt arXiv}, have an entry level for submission).
It is an observed fact that {\tt arXiv} is the most popular and
prevalent preprint server for mathematics.  More will be said
about\index{arxiv@{\tt arXiv}} this tool in what follows.

There are a number of specialized preprint servers for
particular research areas of mathematics. But the most
prominent and widely used preprint server is {\tt arXiv}.
Developed by Paul Ginsparg in 1991, {\tt arXiv} started as
a\index{preprint server!specialized} physics preprint server. But now it handles mathematics,
computer science, statistics, quantitative finance, and quantitative biology as well. As many as
10,000 papers, in the six indicated fields, per month are posted on
{\tt arXiv}. It would be foolish to assert, as many people do,
that ``all math papers are now posted on {\tt arXiv}.'' What is
more accurate is to say that the number of papers posted on
{\tt arXiv} is approaching 30\%. And it is growing. But there
are plenty of older mathematicians who do not give a hoot
about {\tt arXiv}. And there are a number of other
mathematicians who prefer to post their work on specialized
preprint servers that are dedicated to particular areas of
mathematics.  And still others who just cannot be bothered.

The Web site \verb@http://www.arxiv.org@ offers statistics on
the use of {\tt arXiv}.  According to the latest data, there
are 1,213,827 articles now posted on {\tt arXiv}.  If we estimate
that 400,000 of these are in mathematics (and that is a generous
estimate), and we note that about 100,000 math articles are produced per year,
and finally we note that {\tt arXiv} is 25 years old, then it is easy
to see that {\tt arXiv} has\index{arxiv@{\tt arXiv}} not yet taken over.  But it could.

Plenty of mathematicians are tired of dealing with obstinate referees
and arrogant editors.  They feel that, having posted their work on {\tt arXiv},
they have published it and that is all that they owe to the mathematics
community.  One could argue the point.  One could claim that the traditional
refereeing and publishing process guarantees the robustness and longevity
of our work.  That displaying mathematics as an undifferentiated melange
of non-reviewed work is neither productive nor useful.  But these ideas
are still very much in the air and still very much being debated.

The great thing about {\tt arXiv} is that it is very easy to type
in an author's name and get a listing of all his/her most recent
papers.  And it is equally easy to download any of them.  You can
also tell {\tt arXiv} which areas of mathematics you are interested in
and it will send you an \email\ notification each day of what
new papers have been posted.\index{arxiv@{\tt arXiv}, search on}

The server {\tt arXiv} has become so well established that it is now possible,
with many journals, to submit a paper by just providing a pointer
to your {\tt arXiv} posting.  Most professional journals are fairly
free and easy about {\tt arXiv}.  They will not insist that you
take down your {\tt arXiv} posting as soon as your paper
is accepted by the journal.   Book publisher are different, and they
often {\it will} ask you to take your book down from {\tt arXiv} once
it is officially published.\index{arxiv@{\tt arXiv}, removing posting from}

An interesting feature of {\tt arXiv} is that it only accepts
\TeX\ submissions.  That is, when you upload your paper, it
must be in raw \TeX\ or \LaTeX\ form.  Not {\tt *.pdf}, not {\tt *.docx}.
Only \TeX.  And {\tt arXiv} {\it compiles} your paper right on
the spot.  If it succeeds, then you can proceed with the submission
process.  If it fails then you are dead in the water.  This is
another motivation for you to learn to make your \TeX\ files
self-contained.	 In fact the {\tt arXiv} Web site is quite
explicit in stating that it prefers \LaTeX 2e.  And, since
that is the driving form of \TeX\ behind Mik\TeX, that preference
makes some sense.\index{arxiv@{\tt arXiv}, use of \TeX\ in}

If your paper has separate graphics files, then you may find it
tricky to get {\tt arXiv} to compile and accept your paper.  But it
can be done.  I have done it.\index{arxiv@{\tt arXiv}, use of graphics in}

As noted elsewhere in this book, once you write something then it is
immediately copyrighted to you. That is still true when you put a paper on
{\tt arXiv}. When you next submit the paper to a journal, it is quite
standard for the journal to ask you to sign a copyright transfer
agreement. When you sign it, then the copyright moves to the
publisher.\footnote{Some mathematicians prefer to retain the copyright to
themselves.\index{copyright!to you} This is because, for instance, it may happen years later that
someone wants to put together a volume of historically influential papers
in a certain subject area. If your paper is chosen for this volume, and if it
is copyrighted to some other publisher, then nasty negotiations may ensue.
And nasty fees. If you feel strongly about retaining the copyright to your
work, you may have to negotiate with your journal publisher.}

Just for the record, here is a fairly friendly journal publisher's policy
toward posting papers on the {\tt arXiv}:
\begin{quote}
The ASL hereby grants to the Author the non-exclusive right to
reproduce the Article, to create derivative works based upon
the Article, and to distribute and display the Article and any
such derivative work by any means and in any media, provided
the provisions of clause (3) below are met. The Author may
sub-license any publisher or other third party to exercise
those rights.
\end{quote}
And here is a slightly less friendly policy, which still allows the \hfill \break
author to post on {\tt arXiv}:
\begin{quotation}
I understand that I retain or am hereby granted (without the
need to obtain further permission) rights to use certain
versions of the Article for certain scholarly purposes, as
described and defined below (Retained Rights"), and that no
rights in patents, trademarks or other intellectual property
rights are transferred to the journal.

The\index{arxiv@{\tt arXiv}, publisher policy towards} 
Retained Rights include the right to use the Pre-print or
Accepted Authors Manuscript for Personal Use, Internal
Institutional Use and for Scholarly Posting; and the Published
Journal Article for Personal Use and Internal Institutional
Use.
\end{quotation}

Now the truth is that {\tt arXiv}, in its raw form, is rather
stodgy and difficult to use.  Fortunately for us, Greg Kuperberg
has\index{front@{\tt Front} for the {\tt arXiv}}\index{Kuperberg, G.} 
created a front end for {\tt arXiv} called {\tt Front}.  The
URL for {\tt Front} is 
\begin{verbatim}
http://front.math.ucdavis.edu/
\end{verbatim}
{\tt Front} is very user-friendly and easy to use.  I recommend it.

A nasty problem that we all have to deal with in the modern world is this.
Once I write a paper, there are soon many versions of it floating around.
There could be a dozen versions on my school computer, another dozen
versions on my home computer, a version on {\tt arXiv}, a version
on my Home Page, and so forth.  Which is the definitive version of the
paper?  There is no clear and easy answer to this question.  You may
want to give the matter some thought and establish a definitive policy for
yourself.\index{paper, multiple copies of}

\section{\tt MathSciNet}

We have mentioned {\tt MathSciNet} at several earlier junctions
in the book.\index{mathscinet@{\tt MathSciNet}}  Here we treat the topic more discursively.

In 1869 Felix M\"{u}ller and Carl Ohrtmann created the periodical
{\it Jahrbuch \"{u}ber die Fortschritte der Mathematik}. Its purpose was
to archive the mathematical literature. The {\it Jahrbuch} was published
by Walter de Gruyter in one volume per year until 1943. A total of 68
volumes, containing records of 200,000 publications, appeared in the {\it
Jahrbuch}. One of the wonderful things about the mathematical literature is
that it never goes out of date. The articles in {\it Jahrbuch} are still
of value. Therefore, more recently, Bernd Wegner, Keith Dennis, and Elmar
Mitter have created an OnLine version of the {\it Jahrbuch} called ERAM
(Electronic Research Archive for Mathematics).  Because of events described
below, most living mathematicians have never seen (nor perhaps even heard
of) the {\it Jahrbuch}.

In 1939 Otto Neugebauer, who had in 1931 created {\it Zentralblatt
f\"{u}r Mathematik} in Germany, fled from the Nazis and moved to
Brown\index{Neugebauer, O.} University and created {\it Mathematical Reviews}.  The purpose
of both {\it Zentralblatt} and {\it Math Reviews} was to archive the 
mathematical literature.\index{mathematical reviews@{\it Mathematical Reviews}}  Each of these journals publishes a brief
review\index{zentralblattfurmathematik@{\it Zentralblatt f\"{u}r Mathematik}} 
(about a paragraph or so) about most of the papers published
in most of the math journals around the world.  In actuality, journals
are classified by type:  for some journals, all articles in all issues
are reviewed; for other journals, some articles in each issue are reviewed;
in still other journals, no articles are reviewed.

For the American Mathematical Society, {\it Math Reviews} is a big enterprise.
Situated in an old brewery in Ann Arbor, Michigan, at least 75 people are employed
in the production of {\it Math Reviews}.  A great deal of care is put into
sorting out names (so that all the different John Smith's are distinguished),
sorting out articles with similar titles, and getting all the bibliographic
information correct.   And a huge amount of effort is devoted to
requisitioning and classifying and typesetting and organizing the
reviews of the individual papers and books.

Unlike most other abstracting databases, {\tt MathSciNet} takes care
to\index{mathreviews@{\it Math Reviews}} identify authors properly.  Its author search allows the
user to find publications associated with a given author
record, even if multiple authors have exactly the same name.
Mathematical Reviews personnel will sometimes even contact
authors directly to ensure that the database has correctly attributed
their papers. On the other hand, the general search menu uses
string matching in all fields, including the author field.
This functioning is needed for the database to access some old
reviews (before 1940) which have not yet been completely
integrated and thus cannot be found by searching for the
author first.

In 1980, {\it Math Reviews} was converted to an OnLine database, and this
eventually evolved into {\tt MathSciNet} in the 1990s. It is safe to say that {\tt
MathSciNet} has revolutionized the mathematics profession. Now virtually
any mathematician can, from virtually any location, look up papers and
books in the mathematical sciences, assemble bibliographies and reading
lists,\index{mathscinet@{\tt MathSciNet}} and become acquainted with the literature. It is now relatively
straightforward to assemble bibliographies and reference lists.

All sorts of interesting searches can be done in {\tt
MathSciNet}. You can search for author(s), title, journal,
{\tt MathSciNet} ID number, and many other choices. {\tt
MathSciNet} provides a Bib\TeX\ entry with each review. It is
linked to {\tt MathJax}.\footnote{{\tt MathJax} is a
cross-browser {\tt JavaScript} library that displays
mathematical notation in web browsers, using {\tt MathML},
\LaTeX, and {\tt ASCIIMathML} markup. MathJax is released as
open-source software under the Apache License.}\index{mathscinet@{\tt MathSciNet}\ search}

You can also, on {\tt MathSciNet}, calculate
your collaboration distance to another mathematician.  So, for instance, if you wrote
a paper with Riemann who wrote a paper with Gauss, then your ``Gauss number''
is 2, and {\tt MathSciNet} can calculate that for you.\index{collaboration distance}

If you look up {\tt Steven Krantz} on {\tt MathSciNet}, then it
will tell you {\bf (i)}  how many papers and books Krantz has
written, {\bf (ii)} how many citations there have been of his
work, and {\bf (iii)} who his principal collaborators are.  There
is a wealth of information on {\tt MathSciNet} and this is a tool
that is certainly worth mastering.\index{collaborators, number of}\index{citations, number of}

\section{Mathematical Blogs and Related Ideas}

A blog (a truncation of the expression {\tt weblog}) is a
discussion or informational website published on the World
Wide Web consisting of discrete, often informal diary-style
text entries (these are usually called ``posts''). Posts are typically displayed in
reverse\index{blogs}\index{chat rooms}\index{wikis@{\tt wiki}s} chronological order, so that the most recent post
appears first, at the top of the web page. Until 2009, blogs
were usually the work of a single individual,
occasionally of a small group, and often covered a single
subject or topic.  Since 2010, ``multi-author blogs'' 
have developed, with posts written by large numbers of authors
and sometimes professionally edited. The rise of {\tt Twitter} and other
``microblogging'' systems helps integrate multiple author blogs and single-author
blogs into the news media.  ``Blog'' can also be used as a verb,
meaning to maintain or add content to a blog. The emergence
and growth of blogs in the late 1990s coincided with the
advent of Web publishing tools that facilitated the posting of
content by non-technical users who did not have much
experience with {\tt HTML} or computer programming. Previously, a
knowledge of such technologies as {\tt HTML} and File Transfer
Protocol had been required to publish content on the Web, and
as such, early Web users tended to be hackers and computer
enthusiasts. Since 2010, the majority of blogs are interactive Web 2.0
websites,\index{blogs vs.\ social networking sites} allowing visitors to leave online comments and even
message each other via GUI widgets on the blogs, and it is
this interactivity that distinguishes them from other static
websites.  In that sense, blogging can be seen as a form of
social networking service (although, in the next section, we
are careful to distinguish blogs from social networking
sites). Indeed, bloggers do not only
produce content to post on their blogs, but also build social
relations with their readers and other bloggers.  

By definition, a {\it blog} is a discussion or informational
Web site published on the Web and consisting of discrete,
often informal, text entries (called {\it posts}). A {\it chat
room}, by contrast, can have many contributors. A {\tt Wiki}
allows most anyone to {\it edit} the material being posted.

In mathematics, blogs have become very popular. See [Bae] for
an enthusiastic discussion of math blogs. John Baez's math
blog can be accessed at \verb@johncarlosbaez.wordpress.com@.

A math blog can be about a particular mathematical topic, or
even a single mathematical research problem.  Participating in
a math blog is very much like participating in a coffeeroom
discussion, with many participants from all over the world.
There are many a story of important problems being solved by
the participants in a math blog or chat room.

Of course many new complexities can arise from such an event.
If 200 people contribute to the solution of a problem, then
who writes it up? Whose name goes on the paper? Who decides
where to submit it? If you submit to an Open Access journal
and there is an (often nontrivial) author fee, then who pays it?
\index{collaboration!complications with mass}

Some of the most famous and popular math blogs and chat
rooms and {\tt Wiki}s were created
by Fields Medalists---notably Timothy Gowers and Terence Tao.
Gowers has also created {\tt polymath} and {\tt MathOverFlow}. These are Web
sites specifically designed to bring together large groups of people
to\index{mathoverflow@{\tt MathOverFlow}} work on specific mathematical research problems.  \hfill \break
{\tt MathOverFlow} has\index{polymath@{\tt polymath}} had a number of notable successes.

It is particularly easy to participate in {\tt MathOverFlow}.
Go to \hfill \break
\verb@mathoverflow.net@ and you will be immediately introduced
to current problems under discussion.  Further down the page
you are asked to contribute your own comments.  And so now
you are a member of the gang!

It is equally easy to become involved in {\tt polymath}.  Just
go to \verb@https://polymathprojects.org/@ and you are off
and running.  Also check out \verb@\michaelnielsen.org@ for
insights into {\tt polymath}.

One advantage of a blog is that it does not have to meet the
usual scholarly publication standards, and does not have to
fit into the purview of any of the standard scholarly journals.  
It does not have to undergo any refereeing or vetting.  You
do not have to deal with tiresome referees or officious editors.
You can create a blog about
any topic that you think to be of interest, or that you think
others\index{blog!topics for} will find to be of interest.  This could include
\begin{enumerate}
\item[{\bf (a)}]  your experience as a project {\tt NeXT} fellow,
\item[{\bf (b)}]  how to teach a certain unusual topic,
\item[{\bf (c)}]  how to handle tricky situations with students,
\item[{\bf (d)}]  a discussion of an interesting paper that you read recently,
\item[{\bf (e)}]  how to study for an exam,
\item[{\bf (f)}]  how to prepare for a job interview,
\item[{\bf (g)}]  how to come up with examples on your own,
\item[{\bf (h)}]  an attack on a specific research problem,
\item[{\bf (i)}]  partial results on a particular research problem.
\end{enumerate}

There are a number of Web sites that will help you to create
a blog of your own.  Among these are {\tt siteblog}, {\tt SiteBuilder},
{\tt website}, {\tt HostGater}, and {\tt iPage}.  These will
help you to implement the sort of interactivity that makes a blog
effective.

The term {\it chat room}, or chatroom, is primarily used to describe any form of
synchronous conferencing, occasionally even asynchronous conferencing. The
term can thus mean any technology ranging from real-time online chat and
online\index{chat room} interaction with strangers (e.g., online forums) to fully immersive
graphical social environments. The primary use of a chat room is to share
information via text with a group of other users. Generally speaking, the
ability to converse with multiple people in the same conversation
differentiates chat rooms from instant messaging programs, which are more
typically designed for one-to-one communication. The users in a particular
chat room are generally connected via a shared interest or other similar
connection, and chat rooms exist to cater to a wide range of subjects.
New technology has enabled the use of file sharing and webcam to be
included in some programs. This would be considered a chat room.

Another type of Web page that is becoming very popular is the
{\tt wiki} page. A {\tt wiki} is a Web site
that\index{wiki@{\tt wiki}\ page} provides collaborative modification of its content and
structure directly from the Web browser. In a typical {\tt wiki},
text is written using a simplified markup language (known as
``wiki markup''), and often edited with the help of a rich-text
editor. A wiki is run using wiki software, otherwise known
as a wiki engine. There are dozens of different wiki engines
in use, both standalone and part of other software, such as
bug tracking systems.  The Web site
\begin{verbatim}
https://en.wikipedia.org/wiki/Wikipedia:How_to_create_a_page
\end{verbatim}
will tell you how to create a {\tt wiki} page.

As you may know, {\tt Wikipedia} is a {\it very large} OnLine
encyclopedia---one with\index{wikipedia@{\tt Wikipedia}} millions of entries. And it is written
by the readership. The {\tt Wikipedia} organization does a
fairly careful job of monitoring the articles, and they are
generally of good quality. It is something of an encomium to
have a {\tt Wikipedia} article written about oneself.

You may actually want to consider writing an article for {\tt
Wikipedia}. If, for instance, you work in several complex
variables, and you are interested in domains of finite type or
boundary orbit accumulation points (both subjects of current
intense interest), then you may be disappointed to find that
{\tt Wikipedia} has no article on either of these topics. So
you may like to write one. The URL in the preceding paragraph
but one will tell you how to do so. The {\tt Wikipedia}
article will {\it not} identify you as the author, and it {\it
will} allow others to correct and augment your words. It is a
collaborative process, and it can be fun.

There is also {\tt WikiMath}, which is a {\tt wiki} page designed
specifically for mathematics and mathematical questions.  It has
been described as follows:\index{wikimath@{\tt WikiMath}}
\begin{quote}
This wiki collects tasks and topics from mathematics, including their solutions. This
is for everyone who, by himself/herself, feels a need for
mathematical help.  You can look here for your task. Either you find a
solution directly, or you can hope that maybe the next person
interested in WikiMath will discuss your project on a new page. 
\end{quote}

Blogs and chat rooms and {\tt wiki}s and related utilities have expanded our
ability to communicate with people all over the world.  They
have augmented the already exploding activity of mathematical
collaboration.  They are a significant new part of life.

\section{{\tt FaceBook}, {\tt Twitter}, {\tt Instagram}, and the Like}

Social media have become a significant factor in the modern
world.\index{social media}\index{facebook@{\tt FaceBook}} They foster social relationships, they have played a
signficant\index{twitter@{\tt Twitter}} role in helping people meet mates, and also in
helping people to re-gain contact with others whom they
haven't seen or communicated with in many years. {\tt
FaceBook}, {\tt Twitter}, and other utilities play a major
role in advertising. Newly elected President of the United
States Donald Trump frequently uses {\tt Twitter} to
promulgate his ideas and opinions.\index{Trump, D.}

Social media are computermediated technologies that allow the creating and
sharing of information, ideas, career interests, and other forms of
expression via virtual communities and networks. The variety of standalone
and builtin social media services currently available 
challenges a succinct definition.

However,\index{social media!common features} there are some common features:
\begin{enumerate}
\item[{\bf 1.}]  Social media are interactive Web 2.0 Internet based
applications.
\item[{\bf 2.}]  User-generated
content, such as text posts or comments, digital photos
or videos, and data generated through all online interactions are the
lifeblood of social media.
\item[{\bf 3.}]  Users create service-specific
profiles for the {\tt Web}site or app that are
designed and maintained by the social media organization.
\item[{\bf 4.}]  Social media facilitate the development of online social networks by
connecting a user's profile with those of other individuals and/or
groups.
\end{enumerate}

Social media use Web-based and mobile technologies on smartphones and tablet
computers to create highly interactive platforms through which
individuals, communities, and organizations can share, cocreate, discuss,
and modify user-generated content or premade content posted OnLine. They
introduce substantial and pervasive changes to communication between
businesses, organizations, communities, and individuals. Social media
change the way businesses, organizations, communities, and individuals interact.
Social media change the way individuals and large organizations
communicate. These changes are the focus of the emerging field of
technoself studies. In America, a survey reported that 84\% of
adolescents in have a Facebook account.

What role can social media play in mathematics?  Mathematician Michael Jury
was perhaps the first person ever to use social media to solve a research
problem.\index{social media!and mathematics}  How could he have done this?

{\tt FaceBook}, for instance, makes it very easy to communicate with a group
of friends---even a very broad group.  So Mike sent a message to a large
group of mathematical friends telling of a place where he was stuck in
a problem that he was working on.  Within the same day he had an answer
to\index{Jury, M.} his question.  He solved the problem and wrote a nice paper.

Certainly it would be possible to send a {\tt Tweet} (a {\tt Twitter} message)
announcing that you have proved a nice new theorem.   This is not currently
the\index{tweet@{\tt Tweet}} default way to announce a new result.  Instead people post their work
on {\tt arXiv}, they submit a research announcement to {\it Research Announcements
of the AMS}, or they give a talk at a conference.  Or they might send around
a mass email.  But the world is changing around us and social media
may eventually play a more prominent role in mathematician's lives.

What are some of the things that social media could do for mathematics?
Here are some partial answers:
\begin{itemize}
\item One could, in principle, use {\tt FaceBook} as a tool to promote
collaboration.  I frankly do not know anyone who does this.  Chat rooms
would perhaps be more appropriate.  Most people that I know collaborate using
\email---just sending drafts back and forth as attachments.

\item A social media utility could be used to augment a mathematics class.  One could
post anecdotes about famous mathematician, bits of mathematical trivia,
interesting math facts, challenge problems, and the like.  One caution
is that it is dangerous territory to make only some of your students
your ``Friends.''  Avoid that temptation.  

\item Social media could be used to remind students about upcoming tests and
other class events (review sessions, films, special presentations, and
the like).  

\item Social media could be used quite effectively to announce and promote
upcoming talks, upcoming conferences, new workshops, and other mathematical
events.  

\item You could easily post emendations and errata to your lectures on
a socail media site.

\item You could, if you like, {\tt Tweet} a link and make a {\tt FaceBook}
post to your latest paper on {\tt arXiv}.
\end{itemize}
	     
Social\index{social media!vs. blogs} media are often confused with blogs and other Internet utilities.
The distinguishing feature of genuine social media are these:
\begin{enumerate}
\item {\bf User accounts:}  If a site allows visitors to create their own accounts that they can log into, then that's a good sign there's
going to be social interaction. You can't really share information or interact with others online without doing it through a
user account.

\item {\bf Profile pages:}  Since social media is all about communication, a profile page is often necessary to represent an individual. It
often includes information about the individual user, like a profile photo, bio, website, feed of recent posts,
recommendations, recent activity, and more.

\item {\bf Friends, followers, groups, hashtags and so on:} Individuals use
their accounts to connect with other users. They can also use them to
subscribe to certain forms of information.

\item {\bf News feeds:} When users connect with other users on social
media, they're basically saying, ``I want to get information from these
people.'' That information is updated for them in real time via their news
feed.

\item {\bf Personalization:} Social media sites usually give users the
flexibility to configure their user settings, customize their profiles to
look a specific way, organize their friends or followers, manage the
information they see in their news feeds, and even give feedback on what
they do or don't want to see.

\item {\bf Notifications:} Any site or app that notifies users about
specific information is definitely playing the social media game. Users
have total control over these notifications and can choose to receive the
types of notifications that they want.

\item {\bf Information updating, saving, or posting:} If a site or an app
allows you to post absolutely anything, with or without a user account,
then it is social. It could be a simple text-based message, a photo upload,
a {\tt YouTube} video, a link to an article, or anything else.

\item {\bf Like buttons and comment sections:} Two of the most common ways
we interact on social media are via buttons that represent a ``like'' plus
comment sections where we can share our thoughts.

\item {\bf Review, rating, or voting systems:} Besides liking and
commenting, lots of social media sites and apps rely on the collective
effort of the community to review, rate, and vote on information that they
know about or have used. Think of your favorite shopping sites or movie
review sites that use this social media feature. 
\end{enumerate}
								  
There are philosophically related items on the Web that do not exactly fit
the moniker ``social media.'' An example is {\tt reddit}. {\tt Reddit} is
a social news aggregation, web content rating, and discussion website.
{\tt Reddit}'s registered community members can submit content, such as
text posts or direct links. Registered users can then vote submissions up
or\index{reddit@{\tt Reddit}} down to organize the posts and determine their position on the site's
pages. The submissions with the most positive votes appear on the front
page or the top of a category. Content entries are organized by areas of
interest called ``sub{\tt reddit}s.'' The sub{\tt reddit} topics include
news, science, gaming, movies, music, books, fitness, food, and
image-sharing, among many others.

There are a good many people who spend several hours everyday recording
all the details of their lives on social media.  And this is accompananied
by copious photographs and other graphics.  They use social media to
promulgate their political opinions, their sexual opinions, and their
social opinions.  It is too soon to tell what role social media may
play in the mathematical sciences.  But the potential is there.

\section{Print-on-Demand Books}

In the old days---fifty years ago let us say---the production of a book
was very formulaic.  Photo engravers used sulphuric acid to
create copper plates which were used as the printing plates in
a high-speed printing press.  And a print run had to be at
least 1000 units in order to be cost effective.  Printing
just\index{print-on-demand books} a few copies of a book was virtually infeasible.  No more.

Nowadays print-on-demand is both feasible and cost effective.  Because
of electronic media, there is no longer a notion of a book
``going out of print.''  The electronic data for a book---stored on
a hard drive or a flash drive or a tape---is always there.  And backed
up on remote devices as well.  There are numerous companies---including
Book1One, Xulon Publishers, and altagraphics---which can produce print-on-demand
books for you.\index{Xulon Publishers}\index{Book1One}\index{altagraphics}

The Espresso Book Machine (EBM) is a print-on-demand machine that 
prints,\index{Espresso Book Machine} collates, covers, and binds a single book in a few minutes.
The EBM is small enough to fit in a retail book store or small library 
room, and as such it is targeted at retail and library markets. 
The EBM can potentially allow readers to obtain any book title, 
even books that are out of print. The machine takes a {\tt *.pdf} 
file for input and prints, binds, and trims the reader's selection as 
a paperback book.

Jason Epstein gave a series of lectures in 1999 about 
his\index{Epstein, J.} experiences in publishing. Epstein mentioned in 
his speech that a future was possible in which customers 
would be able to print an out-of-stock title on the spot, 
if a book-printing machine could be made that would fit in a 
store. He founded 3BillionBooks with Michael Smolens, a Long 
Island\index{Smolens, M.} entrepreneur in Russia, and Thor Sigvaldason, 
a consultant at Price Waterhouse Coopers. At the time, Jeff 
Marsh, a St Louis engineer\index{Marsh, J.} and inventor, had already 
constructed a prototype book printer that could both photocopy and 
bind. Marsh was working on this project for Harvey Ross, 
who held a patent for such a machine. Peter Zelchenko, a 
Chicago-based technologist and a partner of Ross in a related 
patent effort, worked with Marsh to prove the concept and 
also helped bring Marsh and other players together with 
several venture interests.

Ultimately Epstein, together with Dane Neller, former President
and CEO of Dean and Deluca, licensed Marsh's invention and
founded On Demand Books. The first Espresso\index{On Demand Books} 
Book Machine was
installed and demonstrated on June 21, 2007 at the New York
Public Library's Science, Industry and Business Library. For a
month, the public was allowed to test the machine by printing
free copies of public domain titles provided by the Open
Content Alliance, a non-profit organization with a database of
over 200,000 titles.

The direct-to-consumer model supported by the Espresso Book Machine 
eliminates the need for shipping, warehousing, returns and 
pulping of unsold books; it allows simultaneous global availability 
of millions of new and backlist titles.

Unfortunately the Espresso Book Machine costs about \$150,000.  Or you
can lease it for \$5,000 per month.  So this is out of reach
for most people.  But a number of bookstores and libraries have one,
and let their customers use it for a modest fee.  It is possible,
at least in principle, to produce a hard copy of a {\tt Google} digitized
book (which is, by definition, open access) for about \$8.   

Now it is conceivable to have your working seminar at University X put
together a book gathering together the thoughts that you have been 
developing for the past few years and have several copies printed
up for use by the group (and for the students as well).  

Remember that, as soon as you write something, then it is automatically
copyrighted to you.  So you need not worry about protecting your
book once it is printed.

Today {\tt Amazon} is the world's largest book seller---and
seller of everything else as well. {\tt Amazon} has
revolutionized the book business in many ways. {\tt CreateSpace} 
is an artifact of {\tt Amazon} that allows you to
create your\index{createspace@{\tt CreateSpace}}\index{amazon@{\tt Amazon}} 
own electronic book to be posted on {\tt Amazon}.  And
{\tt Amazon} can produce the book in hard copy as well.

It has been six years since Amazon acquired {\tt CreateSpace}, an
on-demand publishing platform, and almost four years since
they announced the free online setup for self-publishing.
While four years seems like a long time in our fast-paced
world, self-publishing still has not reached the mass audience.
Even the biggest social media gurus still take the traditional
route, only choosing to self-publish when they've been
rejected by mainstream publishing houses.

The truth is that print-on-demand publishing is the fastest, most
profitable and easiest way to get your written thoughts out
there. Today, self-published books are even distributed to
traditional outlets like Barnes \& Noble and academic
libraries.\index{self-publishing}\index{print-on-demand} 

Of course self-publishing means you do not get the marketing
resources that come with a traditional publishing deal, but in
our world of social media, that can be easily fixed. So if
self publishing is so easy, why do we not see more authors
using it? Most people are simply not aware of the low barrier
to entry.   This could be the wave of the future.
\vfill
\eject

\hbox{ \ \ \ }

\thispagestyle{empty}

\newpage

%%%%%%%%%%%%%%%%%%%%%%%%%%%%%%%%%%%%%%%%%%%%%%%%%%%%%%%%%%%%%%%

%% Chapter 8
\chapter{Closing Thoughts}

\begin{quote}
\footnotesize \sl Of all those arts in which the wise excel, \hfill \break
Nature's chief masterpiece is writing well.
\smallskip \hfill \break
\null \mbox{ \ \ } \hfill \rm John Sheffield, Duke of Buckingham and Normanby \break
\null \mbox{ \ \ } \hfill \rm Essay on Poetry [1682]
\end{quote}

\begin{quote}
\footnotesize \sl Great prizes are reader interest and understanding; all else is
secondary.  Graceful prose, imagery, wit, even orthography and
grammar are only means to more important ends.  This observation 
makes writing and reading more of a colloquy and less a lonely or
isolating business.
\smallskip \hfill \break
\null \mbox{ \ \ } \hfill \rm from the dust jacket of {\it Mathematical Writing\/} [KnLR]
\end{quote}
\vspace*{.03in}

\begin{quote}
\footnotesize \sl
England has forty-two religions and two sauces.
\smallskip \hfill \break
\null \mbox{ \ \ } \hfill \rm Voltaire     
\end{quote} 

\begin{quote}
\footnotesize \sl A writer needs three things, experience, 
observation, and imagination, any two of which, at times any
one of which, can supply the lack of the others.
\smallskip \hfill \break
\null \mbox{ \ \ } \hfill \rm William Faulkner 
\end{quote}

\begin{quote}
\footnotesize \sl 
Isaac Newton invented his theory of gravity when he was 21.
I'm 32, and I just found out that Garfield and Heathcliffe are
two different cats.
\smallskip \hfill \break
\null \mbox{ \ \ } \hfill \rm Anon. 
\end{quote}

\begin{quote}
\footnotesize \sl Anything that helps communication is good.  
Anything that hurts is bad.
\smallskip \hfill \break
\null \mbox{ \ \ } \hfill \rm Paul Halmos\index{Halmos, Paul}
\end{quote}

%% Section 7.1

\markboth{CHAPTER 7.  CLOSING THOUGHTS}{7.1.  WHY IS WRITING IMPORTANT?}
\section{Why Is Writing Important?}
\markboth{CHAPTER 7.  CLOSING THOUGHTS}{7.1.  WHY IS WRITING IMPORTANT?}

The case for writing, indeed for good writing, has been made
throughout this book. Writing is our tool for communicating our
ideas, and for leaving a legacy for future generations.  One of the
marvels of genuinely outstanding writing is its longevity.  In many
ways, the\index{writing!the case for good}
writings of Herodotus, of Descartes, of Plato, or of
Faulkner are as vibrant and important today as when they were first
penned.  

Writing at the very highest level is often painstaking and tedious.
A good author can spend an entire day agonizing over a word, a comma,
or a phrase.  He/she will revise the work mercilessly.  For the working
mathematician, I am not recommending this sort of writing.  Do it if
you like; but this level of precision and artistry
is not what our profession either demands or needs. 
In fact the sort of clear, cogent, precise writing that I am
promoting here requires little more effort than lousy writing
requires.  Like the ability to scuba dive, the ability to write well
is in truth a matter of becoming conversant with the basic principles and
then practicing.  Once you become comfortable with
the process, then writing becomes less of a chore and more of a
pleasurable pastime.  It allows you to view your written work as an
accomplishment to be proud of, rather than another agony that you
have\index{writing!the pleasure of} slogged your way through.

We all grow up speaking English (or some other native language).
After a while, we convince ourselves that we are able to express our
thoughts verbally---regardless of our technical facility with
grammar, usage, and syntax.  As we grow older, a corollary of such
reasoning is that we all think that we know how to write.  A result of this
process is that it is more difficult to teach people to write (and,
in turn, for them to learn to write) than it is to teach people
calculus. When a student has his/her calculus work marked incorrect,
then\index{writing!learning the art of} 
he/she is inclined to say ``Apparently I don't know how to do this
kind of problem.  I'd better get some help.''   But when a student
has his/her writing marked up and criticized, then he/she is liable to go to
the instructor and say ``Well, just what is it that you want?''

Learning to write well is a yoga; it is a manner of being trained in
self-criticism and self-instruction. Fields Medalist Enrico Bombieri has
observed to me that his artistic activities, particularly his painting,
have helped him to see things more clearly, and in greater detail. Just
so, learning to write well will\index{writing!as yoga} sharpen your
thoughts, develop your skills at ratiocination, and help you to
communicate more effectively.

Developing an ability to write effectively will give you an appreciation of
the writing, and the thinking, of others. And you will learn from their
writing---both what to do and what not to do. It will add a new dimension
to your life. I hope that it is a happy one.
\vfill
\eject

\hbox{ \ \ \ }

\thispagestyle{empty}

\vfill
\eject

\hbox{ \ \ }

\thispagestyle{empty}

\newpage

\begin{theindex}

  \item abbreviations, 60
  \item abstract
    \subitem how to write, 76
    \subitem notation in, 77
    \subitem of paper, 65
    \subitem references in, 77
    \subitem simplicity of, 77
  \item accuracy, 10
  \item acronym, 60
    \subitem use of in bibliographic \hfill \break
             references, 79
  \item adjectives, 49
  \item adverbs, 49
  \item affiliation
    \subitem of author, 65
    \subitem on work, 19
  \item Agnew, Spiro, 13
  \item `All', `Any', `Each', `Every', 35
  \item alliteration, 13, 14
  \item `alternate' vs.\ `alternative', 50
  \item AMS Standard Cover Sheet, 153
  \item assessing your work, 85
  \item audience, 2, 3
  \item author
    \subitem disputes, 93, 94
    \subitem mutual respect of, 94
    \subitem name, 64

  \indexspace

  \item bibliographic
    \subitem reference, unacceptable, 81
    \subitem references, 22
    \subitem references in \LaTeX, 80
    \subitem style, for a book, 79
    \subitem style, for a paper, 79
  \item bibliography, 22
    \subitem how to write, 77
    \subitem in \TeX, 78
    \subitem style of, 78
  \item big words vs.\ small words, 31
  \item book
    \subitem accuracy in, 182
    \subitem advertising of, 189
    \subitem annoyances with the publisher off, 189
    \subitem appendices in, 175
    \subitem based on lecture notes, 169
    \subitem bibliography for, 175
    \subitem components of a good, 166
    \subitem contents of contract for, 190
    \subitem contract for, 190
    \subitem cost of hardcover vs.\ paperback, 187
    \subitem cost of producing, 185, 186
    \subitem effort to write, 166
    \subitem front matter for, 171
    \subitem galley proofs for, 173
    \subitem glossary for, 179
    \subitem going into production, 183
    \subitem handwritten manuscript for, 199
    \subitem how to read reviews of, 181
    \subitem how to respond to reviews of, 182
    \subitem information for the publisher, 180
    \subitem level of, 168
    \subitem making an index for, 173
    \subitem making index using \LaTeX, 174
    \subitem marketing of, 186
    \subitem Marketing Questionnaire for, 188
    \subitem math, 166
    \subitem method for writing, 167
    \subitem mythical production of, 183
    \subitem negotiating the contract for, 190
    \subitem ordering material, 167
    \subitem organization of, 167
    \subitem outline of, 166
    \subitem page proofs for, 173
    \subitem planning of, 169
    \subitem Preface for, 171
    \subitem preface for, 170
    \subitem production in the computer age, 173
    \subitem production of, 172, 185
    \subitem production using \LaTeX, 174
    \subitem repro copy for, 185
    \subitem review, check list for, 130
    \subitem review, how to write a, 129
    \subitem review, negative, 131
    \subitem review, positive, 131
    \subitem review, what should be in a, 129
    \subitem reviewers hired by \hfill \break
             publisher, 181
    \subitem revising the manuscript, 178
    \subitem revision of manuscript, 182
    \subitem revision vs.\ writing anew, 178
    \subitem royalty for, 190
    \subitem shooting of, 185
    \subitem signatures in, 186
    \subitem spiral method for writing, 177
    \subitem submitting to a publisher, 179
    \subitem Table of Contents for, 170, 171
    \subitem Table of Notation for, 179
    \subitem time management while writing, 176
    \subitem TOC for, 170
    \subitem tunnel vision when writing, 177
    \subitem what to send the publisher, 180
    \subitem writer's block while writing, 176
    \subitem writing of a, 108
    \subitem written in collaboration, 169
  \item book, promotion of, 189
  \item brevity, 36
  \item British English, 32

  \indexspace

  \item cf., 37
  \item choice of words, 32
  \item citation
    \subitem styles, 82
  \item citations
    \subitem on-the-fly, 82
    \subitem Princeton style, 82
  \item claim
    \subitem use of, 72
  \item `clearly', 43
  \item clutter on desk, 22
  \item collaboration
    \subitem as adventure, 94
    \subitem etiquette, 94
  \item collaborative work, 93
  \item colored pens, 21
  \item commas, overuse of, 29
  \item comments in a \TeX\ document, 218
  \item communicating author, 89
  \item communication skills, 157
  \item `compare' and `contrast', 51
  \item computer
    \subitem advantages of writing with, 197
    \subitem backups, 199
    \subitem disadvantage of writing with, 199
    \subitem files, processing of, 204
  \item concentration, 23
  \item confidence, 2
  \item contractions, 38
  \item contractions, use of, 33
  \item cover letter, sample of, 156
  \item cover letter, what to include in a, 156

  \indexspace

  \item damning with faint praise, 115
  \item date on paper, 65
  \item date on work, 19
  \item deferred proof, 72
  \item definition
    \subitem how to state, 74
    \subitem statement of, 74
    \subitem use of ``if and only if'' in, 76
  \item definitions
    \subitem how many to include, 75
    \subitem how many to supply, 74
    \subitem importance of, 74
    \subitem placement of, 68
    \subitem which not to include, 75
  \item `denote', 38
  \item desk supplies, 21
  \item detail, 12
  \item `different from' vs.\ `different than', 51
  \item displayed math vs.\ in-text math, 30
  \item `due to' vs.\  `because' vs.\ `through', 52
  \item Dyson, Freeman, 110

  \indexspace

  \item e.g., 37
  \item editor, dealing with, 90
  \item electronic mail
    \subitem and conversation, 161
    \subitem and the personal computer, 161
    \subitem brevity in, 162
    \subitem collaboration over, 218
    \subitem collaboration via, 158
    \subitem completeness of, 160
    \subitem etiquette for, 159
    \subitem forwarding of, 163
    \subitem haste in, 162
    \subitem header in, 159
    \subitem line length in, 163
    \subitem proofreading of, 161
    \subitem security of, 163
    \subitem signature to, 160
    \subitem significance of, 158
    \subitem use of `from' in, 164
    \subitem uses of, 158
  \item electronic publishing, 217
  \item {\tt emacs}, 161, 204
  \item English ability of candidate, 126
  \item English, role of in mathematics, 30
  \item exposition
    \subitem how to write, 101
    \subitem importance of, 100, 109
    \subitem level of difficulty in, 101
    \subitem use of examples in, 102
    \subitem what is, 100

  \indexspace

  \item faddish prose, 32
  \item `farther' vs.\ `further', 52
  \item figure labels, 213
  \item figures
    \subitem creating informative, 213
    \subitem software for creating, 213
  \item Flaubert, Gustave, 14
  \item flippancy, 32
  \item foreign words and phrases, 31

  \indexspace

  \item galley proofs, return of, 92
  \item geometric measure theory, 69
  \item Gettysburg address, 18
  \item given, 29
  \item Goethe, Johann, 14
  \item good
    \subitem sense, 52
    \subitem taste, 52
  \item grammar checkers, 204
  \item grant
    \subitem credibility of your proposal for, 150
    \subitem dissemination in, 151
    \subitem page limitation for application, 151
    \subitem program officer of, 151
    \subitem proofreading of the proposal, 152
    \subitem prospectus for, 149
    \subitem self-evaluation in, 151
  \item grants
    \subitem availability of, 149
  \item graphics created with a computer, 212
  \item Grisham, John, 188

  \indexspace

  \item Halmos, Paul, 10, 16, 28, 60, 76, 177, 178, 243
  \item Hardy, G.\ H., 159
  \item Harvard system for bibliographic references, 79
  \item hence, 15
  \item homily, 9
  \item `hopefully', 52
  \item Hugo, Victor, 23
  \item {\tt hypertext}, 214

  \indexspace

  \item i.e., 37
  \item if and only if, 27
  \item if-then, 26
  \item iff, 28
  \item `if' vs.\ `whether', 41
  \item index, compiling a good, 175
  \item indexers, professional, 175
  \item `infer' and `imply', 41
  \item infinitives, splitting of, 53
  \item Instructions
    \subitem for Submission, 88
    \subitem to Authors, 77, 88
  \item `in terms of', 54
  \item internet, 214
  \item introvert vs.\ extrovert, 169
  \item intuitionistic ethics, 14
  \item IP address, 214
  \item `its' and `it's', 41
  \item `I' vs.\ `we' vs.\ `one', 33

  \indexspace

  \item jargon, 60
  \item job
    \subitem applying for a, 152
    \subitem cover letter for application, 152
    \subitem interview for, 155
    \subitem offer of a, 155
    \subitem short list for, 154
    \subitem talk given in candidacy for, 155
    \subitem teaching credentials for, 154
    \subitem what to include in the application for a, 154
  \item journal
    \subitem backlog of, 88
    \subitem electronic, 215
    \subitem ranking of, 87
    \subitem selecting one for your paper, 87, 88

  \indexspace

  \item kerning, 202
  \item key words, 65

  \indexspace

  \item language as a weapon, 8
  \item \LaTeX, 78
  \item `lay' and `lie', 42
  \item lemma
    \subitem use of, 72
  \item `less' and `fewer', 42
  \item Lincoln, Abraham, 18
  \item linear ordering of written discourse, 7
  \item Littlewood, J.\ E., 159
    \subitem precepts, 133
  \item love letter, 2, 3
  \item love letter to your self, 66

  \indexspace

  \item manipulative language, 9
  \item {\tt Maple}, 214
  \item math paper, components of, 64
  \item Math Reviews number, 78
  \item {\tt Mathematica}, 214
  \item Mathematical Reviews, 65
  \item MathSciNet, 78, 82
  \item {\it modus ponendo ponens}, 26
  \item MR subject classification numbers, 65

  \indexspace

  \item n.b., 37
  \item name on work, 18
  \item `need only', 55
  \item notation, 15
    \subitem abuse of, 24
    \subitem at the beginning of a sentence, 24
    \subitem choosing, 76
    \subitem consistency of, 25
    \subitem good, 24
    \subitem importance of, 76
    \subitem overuse of, 16, 25, 26
    \subitem planning, 25
    \subitem unnecessary, 15
    \subitem use of, 24
  \item notebook computer, 196
  \item numbering
    \subitem pages, 18
    \subitem schemes, 19
    \subitem systems, 20
  \item numbers
    \subitem and numerals, 42
    \subitem use of in bibliographic references, 79

  \indexspace

  \item obscure expression, 5
  \item `obviously', 43
  \item old-fashioned prose, 32
  \item Olivier, Laurence, 11
  \item one side of paper, writing on, 20
  \item opinion piece
    \subitem arguments in, 106
    \subitem conclusion of, 106
    \subitem having something to say in an, 105
    \subitem how to write an, 104
    \subitem research for, 105
    \subitem summation in, 106
    \subitem thesis of an, 106
  \item organization, 10
  \item outline, importance of, 101
  \item overused words, 43

  \indexspace

  \item paper
    \subitem acknowledgement of receipt, 89
    \subitem contents of, 66
    \subitem cover letter for, 89
    \subitem dealing with referee's report, 90
    \subitem first paragraph of, 67
    \subitem galley proofs for, 92
    \subitem how to submit, 85
    \subitem how to write, 83
    \subitem introduction to, 67
    \subitem order of author names, 93
    \subitem ordering material in, 66
    \subitem organization of, 68
    \subitem publishable, 83
    \subitem referee of, 67
    \subitem where to submit, 83, 85, 86
  \item parallel structure, 55
  \item participial phrases, 55
  \item passive voice, 17, 18
  \item pen vs.\ pencil, 21
  \item plural forms of foreign nouns, 44
  \item position paper, research for, 105
  \item possessives, 44
  \item postal address of author, 65
  \item {\tt Postscript}, 209, 213, 214, 217
  \item precise use of language, 6
  \item precision and custom, 45
  \item preface
    \subitem as writing compass, 108
    \subitem importance of the, 107
    \subitem purpose of the, 107
    \subitem questions that it answers, 107
    \subitem what is a, 107
  \item prepositions, ending a sentence with, 56
  \item {\tt Preview}, 209
  \item primary sources, use of, 77
  \item priority disputes, 93
  \item progress in mathematics, 85
  \item proof
    \subitem by contradiction, 73
    \subitem organizing, 68, 72
  \item proofreading, 10
  \item prose vs.\  mathematics, 30
  \item pseudo-proofs, 102

  \indexspace

  \item q.v., 37
  \item Quine, W.\ V.\ O., 27
  \item quotations, 56

  \indexspace

  \item recommendation letter
    \subitem binary comparisons in, 115, 117, 118
    \subitem brevity in, 129
    \subitem declining to write, 125
    \subitem describing scholarly work in, 115
    \subitem enthusiasm in, 121
    \subitem errors in, 121
    \subitem for a non-research school, 120
    \subitem for your own Ph.D.\ student, 125
    \subitem inflation in, 123
    \subitem introductory paragraph, 114
    \subitem mistakes in, 118
    \subitem not writing, 113
    \subitem professional, 112
    \subitem reasons not to write, 113
    \subitem scientific work described in, 116
    \subitem significance of a negative, 121
    \subitem specificity in, 115, 124
    \subitem summation of, 114
    \subitem things to say in, 118
    \subitem truth in, 126
  \item recommendation letter for a student, 124
  \item recommendation, letter of, 112
  \item redundancy, 57, 69
  \item referee
    \subitem dealing with, 90
    \subitem errors by, 90
    \subitem thanking, 91
  \item referee's report
    \subitem how to write a, 133
    \subitem key points in, 134
    \subitem keying to a given journal, 134
    \subitem level of detail in a, 133
    \subitem what to include in a, 133
  \item references, which to cite, 78
  \item repetition, 69
  \item repetitive sounds, 13
  \item repro copy, 185
  \item Riemann hypothesis, 85, 139, 150
    \subitem proof of, xiv, 105
  \item rules of grammar, 34
    \subitem flexibility of, 59
    \subitem strictness of, 59
  \item run-on sentences, 14

  \indexspace

  \item Saturday Night Live, 73
  \item say something, 2
  \item sense, 10
  \item sense of audience, 108
  \item serial comma, 42
  \item Shakespeare, William, 14
  \item ``shall'' and ``will'', 58
  \item short paragraphs vs.\ long \hfill \break
        paragraphs, 31
  \item short sentences vs.\ long sentences, 26
  \item simple sentence
    \subitem structures, 31
    \subitem vs.\ complex sentence, 31
  \item singular constructions vs.\ plural, 29
  \item so, 15
  \item sound, 10
    \subitem and sense, 11, 14, 34
  \item spacing on page, 21
  \item \verb@\special@ command, 213
  \item spell-checkers, 204, 205
  \item spelling, 10
  \item Spillane, Mickey, 188
  \item status of your statements, 72
  \item stopping places for reader, 30
  \item subject
    \subitem and object, 48
    \subitem and verb, agreement of, 45
    \subitem classification numbers, 65
  \item `suffices to', 55
  \item survey
    \subitem bibliography for, 102
    \subitem conclusion of, 102
    \subitem giving credit in a, 103
    \subitem how to write, 102
    \subitem imprecision in, 104
    \subitem purpose of, 104
    \subitem writing of, 102

  \indexspace

  \item TAA, 191
  \item Table of Contents
    \subitem importance of, 108
    \subitem role in the writing process, 108
  \item talk
    \subitem apocryphal audience for a, 135
    \subitem blackboard use in, 142
    \subitem checklist for, 137
    \subitem conclusion of a, 136
    \subitem dreaming in, 139
    \subitem entry points, 140
    \subitem exit points, 140
    \subitem flexibility of a, 134
    \subitem focus of, 137
    \subitem giving credit in, 141
    \subitem how to give a, 134
    \subitem inflexibility of a, 134
    \subitem informality in, 139
    \subitem ingredients of a, 135
    \subitem mantra for, 141
    \subitem organization of, 136
    \subitem overhead slides in, 143
    \subitem preparation of, 140, 141
    \subitem proofs in, 138
    \subitem specificity in, 138
    \subitem structuring of, 138
    \subitem time management in, 135, 140, 144
    \subitem time segments in, 136
    \subitem title of, 138
    \subitem types of, 136
    \subitem use of examples in, 135
  \item teaching philosophy, statement of, 153
  \item technology and the bibliography, 82
  \item telephone tag, 159
  \item terminology as organizational tool, 75
  \item terseness, 12
  \item \TeX, 207
    \subitem and text editors, 203
    \subitem files, compiling, 209
    \subitem flexibility of, 210
    \subitem previewing, 209
    \subitem sample code, 208
    \subitem vs.\ \LaTeX, 210
    \subitem vs.\ word processors, 203, 207
    \subitem what is, 207
  \item Text and Academic Authors Association, 191
  \item text editor, 198, 203
    \subitem customization of, 204
  \item thanks to granting agencies, 65
  \item ``that'' and ``which'', 58
  \item ``then prove that'', 27
  \item theorem
    \subitem how to prove, 71
    \subitem how to state, 69
    \subitem statement of, 68
    \subitem stating in one sentence, 71
  \item therefore, 15
  \item `this' and `that', 46
  \item title, importance of, 64
  \item `to be', 50
  \item TOC, 108
  \item tone, 32
  \item `trivially', 43

  \indexspace

  \item {\tt UNIX}, 218
  \item uploading files, 218

  \indexspace

  \item Vietnam war, 13
  \item Vita
    \subitem honesty in, 147
    \subitem importance of the, 145
    \subitem parts of the, 146
    \subitem references in, 153
    \subitem sample of, 148
    \subitem things not to include in a, 157
    \subitem what not to include in a, 149

  \indexspace

  \item walk in the woods, Halmos style, 16
  \item Wermer, John, 135
  \item `we', use of, 33
  \item when to stop writing, 8
  \item `where', 47
  \item `who' and `whom', 48
  \item {\tt Windows}, 82
  \item {\tt Windows95}, 198
  \item word
    \subitem counters, 206
    \subitem order, 28
    \subitem processing, non-portability of, 202
    \subitem processor, 200
    \subitem processor, deficiencies of, 201
    \subitem processor, uses of, 200
    \subitem treated as an object, 27
  \item writer's block, 7
  \item writing
    \subitem and thought, 4
    \subitem as yoga, 245
    \subitem file management during, 198
    \subitem form over substance in, 199
    \subitem important of hard copy of, 199
    \subitem large, 21
    \subitem learning the art of, 244
    \subitem mathematics vs.\ writing English, 23
    \subitem the case for good, 244
    \subitem the pleasure of, 244
    \subitem using a text editor during, 198
  \item {\tt WYSIWYG}, 202

\end{theindex}

\markboth{INDEX}{INDEX}

\addcontentsline{toc}{chapter}{\thechapter Index}

\end{document}